\documentclass[journal]{IEEEtran}
\ifCLASSINFOpdf
\else
\fi
\RequirePackage{amsmath,amsthm,amsfonts,amssymb}
\RequirePackage[numbers]{natbib}
\RequirePackage{xr-hyper}  
\RequirePackage[colorlinks,citecolor=blue,urlcolor=blue]{hyperref} 
\RequirePackage{graphicx} 

\usepackage{adjustbox}
\usepackage{booktabs}
\usepackage{enumitem}
\usepackage{makecell}
\numberwithin{equation}{section}
\theoremstyle{plain}
\newtheorem{theorem}{Theorem}
\numberwithin{theorem}{section}
\newtheorem{lemma}[theorem]{Lemma}
\newtheorem{proposition}[theorem]{Proposition}
\newtheorem{corollary}[theorem]{Corollary}
\theoremstyle{remark}
\newtheorem{remark}{Remark}
\numberwithin{remark}{section}
\newtheorem{example}{Example} 
\numberwithin{example}{section}
\newtheorem{definition}{Definition}
\numberwithin{definition}{section}
\usepackage{xspace}
\usepackage{bbm}
\usepackage{mathtools}
\usepackage{dsfont}
\usepackage{lineno}
\usepackage[mathscr]{euscript}

\newcommand{\suchthat}{\mathchoice{\colon}{\colon}{:\mspace{1mu}}{:}}

\def\Var{\mathop{\rm Var}\nolimits}%

\newcommand{\Lc}{\mathcal{L}}
\newcommand{\Ps}{\mathscr{P}}
\newcommand{\Qs}{\mathscr{Q}}
\newcommand{\Xc}{\mathcal{X}}

\newcommand{\ft}{{\tilde{f}}}
\newcommand{\ph}{{\hat{p}}}
\newcommand{\pt}{{\tilde{p}}}
\newcommand{\qh}{{\hat{q}}}

\newcommand{\rb}{{\mathbf r}}
\newcommand{\Sb}{{\mathbf S}}
\newcommand{\xb}{{\mathbf x}}
\newcommand{\Xb}{{\mathbf X}}
\newcommand{\yb}{{\mathbf y}}
\newcommand{\Yb}{{\mathbf Y}}
\newcommand{\zb}{{\mathbf z}}

\def\a{\alpha}
\def\b{\beta}
\def\d{\delta}
\def\eps{\epsilon}
\def\ups{\upsilon}
\def\th{\theta}
\def\Th{\Theta}
\def\Rho{\mathrm{P}}

\def\vth{\vartheta}

\DeclareMathOperator{\Leb}{\lambda_{\mathrm{Leb}}}
\DeclareMathOperator\E{\mathsf{E}}
\let\P\relax
\DeclareMathOperator{\P}{\Ps}
\newcommand{\Pt}{\tilde{\P}}
\DeclareMathOperator{\Q}{\Qs}

\let\digamma\relax
\DeclareMathOperator\digamma{\Psi}

\DeclareMathOperator{\Beta}{\mathsf{B}}
\DeclareMathOperator{\Binom}{\mathsf{Bin}}
\DeclareMathOperator{\Normal}{\mathsf{N}}
\DeclareMathOperator{\Poisson}{\mathsf{Poi}}

\DeclareMathOperator{\GammaDist}{\mathsf{G}}

\DeclareMathOperator{\Unif}{\mathsf{Unif}}

\newcommand\eg{e.g.,\xspace}
\newcommand\ie{i.e.,\xspace}
\def\textiid{i.i.d.\@\xspace}
\newcommand\iid{\ifmmode\text{ i.i.d. } \else \textiid \fi}
\newcommand{\Real}{\mathbb{R}}
\newcommand{\Natural}{\mathbb{N}}
\newcommand{\Integer}{\mathbb{Z}}
\newcommand{\ones}{\mathds{1}}

\newcommand{\half}{\frac{1}{2}}

\DeclarePairedDelimiterX{\infdivx}[2]{(}{)}{%
  #1\;\delimsize\|\;#2%
}
\newcommand{\D}{D\infdivx*}

\DeclarePairedDelimiterX{\norm}[1]{\lVert}{\rVert}{#1}
\DeclarePairedDelimiterX{\abs}[1]{\lvert}{\rvert}{#1}

\newcommand*\diff{\mathop{}\!\mathrm{d}}

\newcommand{\defeq}{\mathrel{\mathop{:}}=}

\newcommand{\supp}{\textnormal{supp}}

\newcommand{\at}{{\tilde{a}}}
\newcommand{\bt}{{\tilde{b}}}

\newcommand{\ut}{{\tilde{u}}}
\newcommand{\vt}{{\tilde{v}}}

\newcommand{\qtil}{\tilde{q}}

\newcommand{\Bb}{\mathbb{B}}
\renewcommand{\E}{\mathbb{E}}
\renewcommand{\Sb}{\mathbb{S}}
\newcommand{\That}{\hat{T}}
\newcommand{\Ft}{\widebar{F}}
\newcommand{\Ut}{\widebar{U}}
\newcommand{\Tt}{\widebar{T}}

\newcommand{\taut}{\tilde{\tau}}
\newcommand{\nut}{\tilde{\nu}}

\newcommand{\psit}{\widebar{\psi}}
\newcommand{\psitt}{\widebarbar{\psi}}

\newcommand{\phibar}{\widebar{\phi}}

\newcommand{\Renyi}{R\'enyi}
\newcommand{\JSD}{D_{\mathsf{JS}}}

\newcommand{\bigO}{O}
\newcommand{\bigOtilde}{\tilde{\bigO}}

\DeclarePairedDelimiter\ceil{\lceil}{\rceil}
\DeclarePairedDelimiter\floor{\lfloor}{\rfloor}

\DeclareMathOperator{\polyln}{\mathrm{poly}\ln}
\DeclareMathOperator{\rvol}{\varrho}
\def\maximal{M}
\def\minimal{m}
\def\maxf{W}
\def\minf{w}
\def\expvolfunc{R}

\DeclareFontFamily{OT1}{pzc}{}
\DeclareFontShape{OT1}{pzc}{m}{it}{<-> s * [1.100] pzcmi7t}{}
\DeclareMathAlphabet{\mathpzc}{OT1}{pzc}{m}{it}

\newcommand{\numberthis}{\addtocounter{equation}{1}\tag{\theequation}}

\def\textin{\text{in}}
\def\textout{\text{out}}

\newcommand{\Cab}{C}
\newcommand{\Cbe}{c}

\newcommand{\revision}[1]{\textcolor{black}{#1}}
\newcommand{\secondrevision}[1]{\textcolor{black}{#1}}

\usepackage{cuted}
\makeatletter
\let\save@mathaccent\mathaccent
\newcommand*\if@single[3]{%
  \setbox0\hbox{${\mathaccent"0362{#1}}^H$}%
  \setbox2\hbox{${\mathaccent"0362{\kern0pt#1}}^H$}%
  \ifdim\ht0=\ht2 #3\else #2\fi
  }
\newcommand*\rel@kern[1]{\kern#1\dimexpr\macc@kerna}
\newcommand*\wideaccent[2]{\@ifnextchar^{{\wide@accent{#1}{#2}{0}}}{\wide@accent{#1}{#2}{1}}}
\newcommand*\wide@accent[3]{\if@single{#2}{\wide@accent@{#1}{#2}{#3}{1}}{\wide@accent@{#1}{#2}{#3}{2}}}
\newcommand*\wide@accent@[4]{%
  \begingroup
  \def\mathaccent##1##2{%
    \let\mathaccent\save@mathaccent
    \if#42 \let\macc@nucleus\first@char \fi
    \setbox\z@\hbox{$\macc@style{\macc@nucleus}_{}$}%
    \setbox\tw@\hbox{$\macc@style{\macc@nucleus}{}_{}$}%
    \dimen@\wd\tw@
    \advance\dimen@-\wd\z@
    \divide\dimen@ 3
    \@tempdima\wd\tw@
    \advance\@tempdima-\scriptspace
    \divide\@tempdima 10
    \advance\dimen@-\@tempdima
    \ifdim\dimen@>\z@ \dimen@0pt\fi
    \rel@kern{0.6}\kern-\dimen@
    \if#41
      #1{\rel@kern{-0.6}\kern\dimen@\macc@nucleus\rel@kern{0.4}\kern\dimen@}%
      \advance\dimen@0.4\dimexpr\macc@kerna
      \let\final@kern#3%
      \ifdim\dimen@<\z@ \let\final@kern1\fi
      \if\final@kern1 \kern-\dimen@\fi
    \else
      #1{\rel@kern{-0.6}\kern\dimen@#2}%
    \fi
  }%
  \macc@depth\@ne
  \let\math@bgroup\@empty \let\math@egroup\macc@set@skewchar
  \mathsurround\z@ \frozen@everymath{\mathgroup\macc@group\relax}%
  \macc@set@skewchar\relax
  \let\mathaccentV\macc@nested@a
  \if#41
    \macc@nested@a\relax111{#2}%
  \else
    \def\gobble@till@marker##1\endmarker{}%
    \futurelet\first@char\gobble@till@marker#2\endmarker
    \ifcat\noexpand\first@char A\else
      \def\first@char{}%
    \fi
    \macc@nested@a\relax111{\first@char}%
  \fi
  \endgroup
}
\makeatother

\newcommand\doubleoverline[1]{\overline{\overline{#1}}}

\newcommand\widebar{\wideaccent\overline}
\newcommand\widebarbar{\wideaccent\doubleoverline}

\newcommand*\patchAmsMathEnvironmentForLineno[1]{%
   \expandafter\let\csname old#1\expandafter\endcsname\csname #1\endcsname
   \expandafter\let\csname oldend#1\expandafter\endcsname\csname end#1\endcsname
   \renewenvironment{#1}%
      {\linenomath\csname old#1\endcsname}%
      {\csname oldend#1\endcsname\endlinenomath}}%
\newcommand*\patchBothAmsMathEnvironmentsForLineno[1]{%
   \patchAmsMathEnvironmentForLineno{#1}%
   \patchAmsMathEnvironmentForLineno{#1*}}%
\AtBeginDocument{%
\patchBothAmsMathEnvironmentsForLineno{equation}%
\patchBothAmsMathEnvironmentsForLineno{align}%
\patchBothAmsMathEnvironmentsForLineno{flalign}%
\patchBothAmsMathEnvironmentsForLineno{alignat}%
\patchBothAmsMathEnvironmentsForLineno{gather}%
\patchBothAmsMathEnvironmentsForLineno{multline}%
}

\begin{document}
%
\title{Nearest Neighbor Density Functional Estimation\\
from Inverse Laplace Transform}

\author{J.~Jon~Ryu\IEEEauthorrefmark{1},~\IEEEmembership{Student Member,~IEEE,}
        Shouvik~Ganguly\IEEEauthorrefmark{1},~\IEEEmembership{Member,~IEEE,}
        Young-Han~Kim,~\IEEEmembership{Fellow,~IEEE,}
        Yung-Kyun~Noh,~\IEEEmembership{Member,~IEEE}
        and~Daniel~D.~Lee~\IEEEmembership{Fellow,~IEEE}
\thanks{Manuscript received July 31, 2020; revised January 14, 2021; accepted January 26, 2022. 
This work was supported in part by the National Science Foundation under Grant CCF-1911238.
Y.-K. Noh was partly supported by NRF/MSIT (No. 2017R1E1A1A03070945, 2021M3E5D2A01019545), IITP/MSIT Artifcial Intelligence Graduate School Program for Hanyang University (2020-0-01373).
}
\thanks{\IEEEauthorrefmark{1}J.~J.~Ryu and S.~Ganguly contributed equally to this work.}
\thanks{J.~J.~Ryu and Y.-H.~Kim are with the Department of Electrical and Computer Engineering, University of California, San Diego, La Jolla, CA 92093 USA (e-mail: \href{mailto:jongharyu@ucsd.edu}{jongharyu@ucsd.edu}; \href{mailto:yhk@ucsd.edu}{yhk@ucsd.edu}).}%
\thanks{S.~Ganguly was with the Department of Electrical and Computer Engineering, University of California, San Diego, La Jolla, CA 92093 USA. He is now affiliated at XCOM Labs, San Diego, CA 92121 USA
(e-mail: \href{mailto:sganguly@xcom-labs.com}{sganguly@xcom-labs.com}).}%
\thanks{Y.-H.~Kim is with the Department of Electrical and Computer Engineering, University of California at San Diego, La Jolla, CA 92093 USA and Gauss Labs Inc, Seoul, South Korea (e-mail: \href{mailto:yhk@ucsd.edu}{yhk@ucsd.edu}).}%
\thanks{Y.-K.~Noh is with Department of Computer Science, Hanyang University, Seoul 04763, Republic of Korea
and School of Computational Sciences, Korea Institute for Advanced Study, Seoul 02455, Republic of Korea
(e-mail: \href{mailto:nohyung@hanyang.ac.kr}{nohyung@hanyang.ac.kr}).}
\thanks{D.~D.~Lee is with Cornell Tech, New York, NY 10044 USA and Global AI Center for Samsung Research (e-mail: \href{mailto:ddl46@cornell.edu}{ddl46@cornell.edu}).}
}

\maketitle

\begin{abstract}
A new approach to $L_2$-consistent estimation of a general density functional using $k$-nearest neighbor distances is proposed, where the functional under consideration is in the form of the expectation of some function $f$ of the densities at each point. The estimator is designed to be asymptotically unbiased, using the convergence of the normalized volume of a $k$-nearest neighbor ball to a Gamma distribution in the large-sample limit, and naturally involves the inverse Laplace transform of a scaled version of the function $f.$
Some instantiations of the proposed estimator recover existing $k$-nearest neighbor based estimators of Shannon and R\'enyi entropies and Kullback--Leibler and R\'enyi divergences, and discover new consistent estimators for many other functionals such as logarithmic entropies and divergences. 
The $L_2$-consistency of the proposed estimator is established for a broad class of densities for general functionals,
and the convergence rate in mean squared error is established as a function of the sample size for smooth, bounded densities.
\end{abstract}


%
\IEEEpeerreviewmaketitle

\section{Introduction}
\IEEEPARstart{T}{his} paper studies the problem of estimating \revision{an entropy functional} of the form
\begin{equation*}
T_f(p)\defeq \E_{\Xb\sim p}[f(p(\Xb))]
= \int f(p(\xb))p(\xb)\diff\xb,
\end{equation*}
where $f\suchthat \Real_+\to\Real$ is a given function and $p$ is a probability density over $\Real^d$.
Table~\ref{table:estimator_functions_single} lists examples of $f$ and the corresponding functional $T_f$. The goal is to estimate $T_f(p)$ based on independent and identically distributed (\iid) samples $\Xb_{1:m}=(\Xb_1,\ldots,\Xb_m)$ from $p$ by forming an estimator $\That_f^m(\Xb_{1:m})$ that converges to $T_f(p)$ in $L_2$ as the sample size $m$ grows to infinity, that is,
\[
\lim_{m\to\infty}\E\bigl[\bigl(\That_f^{m}(\Xb_{1:m})-T_f(p)\bigr)^2\bigr] = 0.
\]

More generally, let $f\suchthat \Real_+\times \Real_+\to\Real$ and consider a \revision{divergence} functional
\begin{equation*}
T_f(p,q)\defeq\E_{\Xb\sim p}[f(p(\Xb),q(\Xb))]
=\int f(p(\xb),q(\xb))p(\xb)\diff\xb
\end{equation*}
of a pair of probability densities $p$ and $q$ over $\Real^d.$ 
Table~\ref{table:estimator_functions_two} lists examples of $f$ and the corresponding $T_f$.
In this case, the main problem is to construct an estimator $\That_f^{m,n}(\Xb_{1:m},\Yb_{1:n})$ based on \iid samples $\Xb_{1:m}$ from $p$ and $\Yb_{1:n}$ from $q,$ independent of each other, such that
\begin{equation*}
\lim_{m,n\to\infty}\E\bigl[\bigl(\That_f^{m,n}(\Xb_{1:m},\Yb_{1:n})-T_f(p,q)\bigr)^2\bigr] = 0.
\end{equation*}

Consistent estimation of such quantities, such as Shannon's differential entropy ($f=\ln(1/p)$), (exponentiated) \Renyi{} $\a$-entropies ($f=p^{\a-1}$), Kullback--Leibler (KL) divergence ($f=\ln(p/q)$), Hellinger distance ($f = \sqrt{q/p}$), (exponentiated) \Renyi{} $\a$-divergences ($f=p^{\a-1}q^{-\a}$), and Jensen--Shannon divergence (see Table~\ref{table:estimator_functions_two}), is a problem of considerable practical interest, having wide-ranging applications in parameter estimation~\cite{Weidemann--Stear1969,Wolsztynski--Thierry--Pronzato2005}, goodness-of-fit testing~\cite{Girardin--Lequesne2017,Crzcgorzewski--Wirczorkowski1999,Goria--Leonenko--Mergel--Inverardi2005}, quantization~\cite{Marano--Matta--Willett2007}, independent component analysis~\cite{Kraskov--Stogabauer--Grassberger2004, Learned-Miller--Fisher-III2003, Boukouvalas--Mowakeaa--Fu--Adali2016}, texture classification~\cite{Hero--Ma--Michel--Gorman2002, Susan--Hanmandlu2013}, design of experiments~\cite{Liepe--Filippi--Michal--Stumpf2013, Lewi--Butera--Paninski2007}, pattern recognition~\cite{Hero--Michel1999, Neemuchwala--Hero--Carson2005, Lajevardi--Hussain2009,Shan--Gong--McOwan2005}, clustering and feature selection~\cite{Lajevardi--Hussain2009, Aghagolzadeh--Soltanian-Zadeh--Araabi--Aghagolzadeh2007, Peng--Long--Ding2005, Sotoca--Pla2010}, and statistical inference~\cite{Giet--Lubrano2008}.  
In addition, divergence estimates can be used as  measures of distance between two distributions and thus can generalize distance-based algorithms for metric spaces to the space of probability distributions; see, for example, \cite{Oliva--Poczos--Schneider2013,Henderson--Gallagher--Eliassi-Rad2015} and the references therein.

\begin{table*}
\caption{Examples of functionals of one density and their estimator functions $\phi_{k}(u)$.
A reference is given whenever an estimator already exists in the literature.
The last column presents a pair of exponents $(a_k,b_k)$ of the polynomial envelope of the estimator function $\phi_{k}(u)$. 
The constant $\eps$, if any, can be chosen as an arbitrarily small positive number.
For the first three examples, $k>-a_k$ is required to guarantee the existence of the corresponding inverse Laplace transform.
Here, $\digamma(\a)$ denotes the digamma function~\citep{Korn--Korn2000}; see also Example~\ref{ex:differential_entropy}.
}
\label{table:estimator_functions_single}
\centering
\resizebox{\textwidth}{!}{
\begin{tabular}{l l l l}
 \toprule
 Name
    & \makecell[l]{
            $T_f(p)=\E_p[f(p)]$
            }
    & \makecell[l]{$\displaystyle \phi_{k}(u)=\frac{\Gamma(k)}{u^{k-1}}\Lc^{-1}\Bigl\{\frac{f(p)}{p^k}\Bigr\}(u)$} & \makecell[l]{$(a_k, b_k)$}\\
    \midrule
    \makecell[l]{
    Differential entropy\\
    \cite{Kozachenko--Leonenko1987,Singh--Misra--Hnizdo--Fedorowicz--Demchuk2003,Goria--Leonenko--Mergel--Inverardi2005}\\
    (Examples~\ref{ex:differential_entropy}, \ref{ex:consistency_differential_entropy}, \ref{ex:estimating_differential_entropy}, \ref{ex:estimating_differential_entropy2}, \ref{ex:estimating_differential_entropy_varying_k})}
        & $\displaystyle
        \E\Bigl[\ln\frac{1}{p}\Bigr]$
        & $\displaystyle \ln u-\digamma(k)$
        & \makecell[l]{
            $(-\eps,\eps)$\\
        }\\
    \midrule[0.01pt]
    \makecell[l]{
        $\a$-entropy~\cite{Leonenko--Pronzato--Savani2008}\\
        ($\a\ge 0$)\\
        (Examples~\ref{ex:renyi_entropy}, \ref{ex:consistency_renyi_entropy}, \ref{ex:estimating_renyi_entropy}, \ref{ex:estimating_renyi_entropy2}, \ref{ex:estimating_renyi_entropy_varying_k})
        }
        & $\displaystyle
        \E[p^{\a-1}]$
        & \makecell[l]{$\displaystyle \frac{\Gamma(k)}{\Gamma(k-\a+1)} \Bigl(\frac{1}{u}\Bigr)^{\a-1}$
        } & \makecell[l]{$(1-\a,1-\a)$} \\
    \midrule
    \makecell[l]{
    Logarithmic $\a$-entropy\\
    ($\a> 0$)\\
    (Example~\ref{ex:generalized_alpha_entropy})
    }
        & $\displaystyle
        \E\Bigl[
        p^{\a-1}\ln\frac{1}{p}
        \Bigr]$
        & \makecell[l]{$\displaystyle
        \frac{\Gamma(k)}{\Gamma(k-\a+1)}
        u^{-\a+1}
        (\ln u-\digamma(k-\a+1))$
        } & \makecell[l]{
        $(1-\a-\eps,1-\a+\eps)$
        }\\
    \midrule[0.01pt]
    \makecell[l]{
    Exponential $(\a,\b)$-entropy\\
    ($\a>0,\b\ge 0$)\\
    (Example~\ref{ex:exponential_alphabeta_entropy})}
        & $\displaystyle \E[p^{\a-1}e^{-\b p}]$
        & \makecell[l]{$\displaystyle \frac{\Gamma(k)}{\Gamma(k-\a+1)}\frac{(u-\b)^{k-\a}}{u^{k-1}}\ones_{[\b,\infty)}(u)$
        } & \makecell[l]{
        $(0,1-\a)$ for $k\ge \a$
        }\\
\bottomrule
\end{tabular}
}
\end{table*}

One of the most basic and prominent nonparametric approaches is the $k$-nearest neighbor ($k$-NN) based method, which is appealing since
\revision{its hyperparameter tuning is relatively simple}
and is computationally efficient, especially when $k$ is held fixed, independent of the sample sizes $m$ and $n$.
\revision{In this paper, we propose a new, universal design principle of a $L_2$-consistent $k$-NN based estimator for a wide class of the density functionals $T_f(p)$ and $T_f(p,q)$ based on the inverse Laplace transform, which generalizes many existing estimators which have been developed and analyzed separately. 
Based on the proposed mathematical framework, we establish the consistency and the rate of convergence in MSE of the density functional estimator under fairly general regularity conditions, by extending and simplifying the existing analyses of the KL estimator by \citet{Bulinski--Dimitrov2019a,Bulinski--Dimitrov2019b} and \citet{Gao--Oh--Viswanath2018tit}.}

\subsection{The proposed single-density functional estimators}

Suppose that a metric $\rho\suchthat \Real^d\times\Real^d\to\Real_+$ is associated with the $d$-dimensional space $\Real^d$.
Given samples $\Xb_{1:m}$ and a point $\xb\in\Real^d$, we denote the $k$-NN distance of $\xb$ from the samples by $r_{km}(\xb)\defeq r_{k}(\xb|\Xb_{1:m})$ for $k\le m$.
\revision{Here, $r_k(\xb|A)$ denotes the $k$-NN distance of $\xb$ from a set $A\subseteq \Real^d$, where the distance tie is broken arbitrarily.}
The key statistic in this paper is a normalized volume 
\begin{equation}
\label{eq:def_ukm}
U_{km}(\xb)\defeq U_k(\xb|\Xb_{1:m})\defeq m\Leb(\Bb(\xb,r_{k}(\xb|\Xb_{1:m})))
\end{equation}
of the $k$-NN ball centered at $\xb$ with respect to $\Xb_{1:m}$.
Here and henceforth, $\Leb$ denotes the Lebesgue measure over $\Real^d$, $\Bb(\xb,r)\defeq \{\yb\in\Real^d\suchthat \rho(\xb,\yb)< r\}$ denotes the open ball of radius $r>0$ centered at $\xb\in\Real^d$, and $\widebar{\Bb}(\xb,r)$ denotes the closure of $\Bb(\xb,r)$.
When the $k$-NN distance $r_k$ is evaluated at one of the samples $\xb=\Xb_i$ $(1\le i\le m)$, we define it as $r_k(\Xb_i|\Xb_{1:i-1}\Xb_{i+1:m})$ to exclude the trivial zero distance.
Consequently, we use the convention
\begin{align*}
U_{km}(\Xb_i)\defeq 
(m-1) \Leb(\Bb(\xb,r_k(\xb|\Xb_{1:i-1}\Xb_{i+1:m}))).
\end{align*}
Note that under this convention, we have 
\begin{align}\label{eq:convention}
U_{km}(\Xb_m)=U_{k,m-1}(\Xb_m).    
\end{align}

\begin{table*}
\caption{Examples of functionals of two densities and their estimator functions $\phi_{kl}(u,v)$. The absolute continuity $\P\ll\Q$ is assumed implicitly unless stated otherwise.
A reference is given whenever an estimator already exists in the literature.
The last column presents pairs of exponents $(a_{kl},b_{kl})$ and $(\at_{kl},\bt_{kl})$ of the polynomial envelopes of the estimator function $\phi_{kl}(u,v)$ in $u$ and $v$, respectively.
The constant $\eps$, if any, can be chosen as an arbitrarily small positive number.
For each case, $k>-a_{kl}$ and $l>-\at_{kl}$ is required to guarantee the existence of the corresponding inverse Laplace transform.}
\label{table:estimator_functions_two}
\centering
\begin{adjustbox}{center}
\resizebox{\textwidth}{!}{
\begin{tabular}{ l l l l}
\hline
 Name
    & \makecell[l]{
            $T_f(p,q)=\E_p[f(p,q)]$
            }
    & \makecell[l]{
    $\displaystyle \phi_{kl}(u,v)=\frac{\Gamma(k)\Gamma(l)}{u^{k-1}v^{l-1}}\Lc^{-1}\Bigl\{\frac{f(p,q)}{p^kq^l}\Bigr\}(u,v)$} 
    & \makecell[l]{$(a_{kl}, b_{kl})$;\\$(\at_{kl},\bt_{kl})$}\\
 \midrule
    \makecell[l]{
    	KL divergence~\cite{Wang--Kulkrani--Verdu2009}\\
    	(Examples~\ref{ex:kl_divergence}, \ref{ex:consistency_kl_divergence}, \ref{ex:estimating_kl_divergence}, \ref{ex:estimating_kl_divergence_varying_kl},
    	\ref{supp:ex:kl})
        }
        & $\displaystyle
        \E\Bigl[\ln\frac{p}{q}\Bigr]$
        & $\displaystyle \ln \frac{v}{u}+\digamma(k)-\digamma(l)$ 
        & \makecell[l]{
        $(-\eps,\eps)$;\\$(-\eps,\eps)$
        }\\
    \midrule[0.01pt]
    \makecell[l]{
        $\a$-divergence~\cite{Poczos--Schneider2011}\\
        ($\a>0$)\\
        (Examples~\ref{ex:polynomial_functional}, \ref{ex:consistency_renyi_divergence}, \ref{ex:estimating_renyi_divergence}, \ref{ex:estimating_renyi_divergence_varying_kl}, \ref{supp:ex:poly})
	    }
        & $\displaystyle 
        	\E\Bigl[
            \Bigl(\frac{p}{q}\Bigr)^{\a-1}
            \Bigr]
            $
        & \makecell[l]{$\displaystyle \frac{\Gamma(k)\Gamma(l)}{\Gamma(k-\a+1)\Gamma(l+\a-1)}\Bigl(\frac{v}{u}\Bigr)^{\a-1}$
        }
        & \makecell[l]{
            $(1-\a,1-\a)$;\\$(\a-1,\a-1)$
        } \\
    \midrule
    \makecell[l]{
  	    Logarithmic $\a$-divergence\\
  	    ($\a>0$)\\
  	    (Examples~\ref{ex:logarithmic_alpha_divergence}, \ref{supp:ex:log_alpha_div})
  	}
    & $\displaystyle
    \E\Bigl[
    \Bigl(\frac{p}{q}\Bigr)^{\a-1}\ln\frac{p}{q}
    \Bigr]$
    & \makecell[l]{$\displaystyle
    \frac{\Gamma(k)\Gamma(l)}{\Gamma(k-\a+1)\Gamma(l+\a-1)}
    \Bigl(\frac{v}{u}\Bigr)^{\a-1}\times$\\
    $\qquad\qquad\displaystyle\bigl(\ln\frac{v}{u}+\digamma(k-\a+1)-\digamma(l+\a-1)\bigr)$
    }
    & \makecell[l]{
        $(1-\a-\eps,1-\a+\eps)$;\\$(\a-1-\eps,\a-1+\eps)$
    }\\
    \midrule[0.01pt]
    \makecell[l]{
        Le Cam distance\\
        (Examples~\ref{ex:asymp_nn_classfication}, \ref{ex:estimating_asymp_nn_classfication}, \ref{supp:ex:NNclassification})
        }
		& $\displaystyle \E\Bigl[\frac{(p-q)^2}{2p(p+q)}\Bigr]$
        & 
        \newcommand{\smallbinom}[2]{\biggl(\genfrac{}{}{0pt}{}{#1}{#2}\biggr)}
        \makecell[l]{
        $\displaystyle
        2\smallbinom{k+l-2}{k-1}^{-1}\Bigl\{\sum_{j=0}^{l-1}\smallbinom{k+l-2}{k-1+j}\Bigl(-\frac{u}{v}\Bigr)^{j}-$
        \\
        \qquad\qquad\quad
        $\displaystyle
        \Bigl(-\frac{u}{v}\Bigr)^{l-1}\Bigl(1-\frac{v}{u}\Bigr)^{k+l-2}\ones_{[v,\infty)}(u)\Bigr\}$
        }
        & \makecell[l]{
            $(-k+1,l-1)$;\\$(-l+1,k-1)$
        }\\
    \midrule[0.01pt]
        \makecell[l]{
        Entropy difference\\ ($\Q\ll\P$)\\
        (Example~\ref{supp:ex:ent_diff})
    }
    & $\displaystyle
    \E\Bigl[
        \ln \frac{1}{p}-\frac{q}{p}\ln \frac{1}{q}
    \Bigr]$
    & \makecell[l]{$\displaystyle\frac{(l-1)}{k}\frac{u}{v}(\digamma(l-1)-\ln v)-(\digamma(k)-\ln u)$
    } 
    & \makecell[l]{
        $(-\eps,1)$;\\$(-1-\eps,-1+\eps)$
    } \\
    \midrule[0.01pt]
    \makecell[l]{
        Reverse KL divergence\\ ($\Q\ll\P$)\\
        (Example~\ref{supp:ex:rev_kl})
    }
        & $\displaystyle
            \E\Bigl[\frac{q}{p}\ln \frac{q}{p}\Bigr]$
        & \makecell[l]{$\displaystyle \frac{l-1}{k}\frac{u}{v}\Bigl(\ln \frac{u}{v}+\digamma(l-1)-\digamma(k+1)\Bigr)$
        }  
        & \makecell[l]{
            $(1-\eps,1+\eps)$;\\
            $(-1-\eps,-1+\eps)$
            } \\
    \midrule[0.01pt]
    \makecell[l]{
        Jensen--Shannon divergence\\ ($\Q\ll\P$)\\
        (Examples~\ref{ex:JSD}, \ref{ex:estimating_JSD}, \ref{supp:ex:JSD})
    }
        & 
        $\displaystyle
        \E\Bigl[\half\ln \frac{2p}{p+q}
		+\frac{q}{2p}
         \ln\frac{2q}{p+q}
        \Bigr]$
        & See Example~\ref{ex:JSD}.
        & \makecell[l]{
            $(-k+1, l-1)$;\\
            $(-l+1, k-1)$
            }\\
 \bottomrule
\end{tabular}
}
\end{adjustbox}
\end{table*}

Let $\GammaDist(\a,\b)$ denote the Gamma distribution with shape parameter $\a>0$ and rate parameter $\b>0$, whose density is
\begin{align*}
    \frac{\b^\a}{\Gamma(\a)}u^{\a-1}e^{-\b u},\quad u\ge 0.
\end{align*}
Here $\Gamma(\a)\defeq\int_{0}^{\infty }x^{\a-1}e^{-x}\diff x$ denotes the Gamma function.
The following fact on the asymptotic distribution of $U_{km}(\xb)$ is well known~\cite{Singh--Misra--Hnizdo--Fedorowicz--Demchuk2003,Goria--Leonenko--Mergel--Inverardi2005,Leonenko--Pronzato--Savani2008}. The proof is presented in Appendix~\ref{supp:sec:technical_lemmas:knn_stats} for completeness.

\begin{proposition}\label{prop:knndist}
Suppose that $k\ge 1$ is a fixed integer, and let $\Xb_{1:m}$ be \iid samples drawn from $p$ on $\Real^d$.
Then, for almost every $\xb$,
$U_{km}(\xb)$ converges to a $\GammaDist(k,p(\xb))$ random variable in distribution as $m$ goes to infinity.
\end{proposition}

This general convergence result is the cornerstone of the design of our estimator. 
To be more specific, 
for functionals of one density $p$,
consider an estimator of the form 
\begin{equation}
\label{eq:estimator_single}
\That_f^{(k)}(\Xb_{1:m}) = \frac{1}{m}\sum_{i=1}^m\phi_k(U_{km}(\Xb_i))
\end{equation}
that depends on the samples only through the $k$-NN distance evaluated at each of them. 
As a necessary condition for the $L_2$-consistency of this estimator, the function $\phi_k$
should be chosen such that
\[
\lim_{m\to\infty} \E[\That_f^{(k)}] = T_f(p),
\]
that is, the estimator is asymptotically unbiased.
On the one hand, since $\Xb_{1:m}$ are identically distributed,
we have, from \eqref{eq:convention} and \eqref{eq:estimator_single}, that $\E[\That_f^{(k)}] = \E[\phi_k(U_{k,m-1}(\Xb_m))]$, and thus the desired asymptotic unbiasedness for a \emph{fixed} $k$ can be expressed equivalently as
\begin{equation}\label{eq:asymp_unbiasedness}
\lim_{m\to\infty}\E[\phi_k(U_{k,m-1}(\Xb_m))]
= T_f(p) 
= \int p(\xb)f(p(\xb))\diff\xb.
\end{equation}
On the other hand, from Proposition~\ref{prop:knndist}, we expect that under certain regularity conditions,
\begin{align*}
\lim_{m\to\infty} \E[\phi_k(U_{k,m-1}(\Xb_m))]
&= \E[\phi_k(U_{k\infty}(\Xb))]
\numberthis
\label{eq:convergence_in_expectation}\\
&= \int p(\xb) \E[\phi_k(U_{k\infty}(\xb))]\diff \xb,
\end{align*}
where $U_{k\infty}(\xb)$ is a $\GammaDist(k,p(\xb))$ random variable,
independent of
$\Xb\sim p$ for every $\xb$.
We choose $\phi_k(u)$ so as to equate the integrands in~\eqref{eq:asymp_unbiasedness} and \eqref{eq:convergence_in_expectation}, \ie
for every $p > 0$,
if $U \sim \GammaDist(k,p)$,
then
\begin{align*}
f(p)
&= \E[\phi_k(U)] \\
&= \int_0^\infty\phi_k(u)\frac{p^k}{\Gamma(k)}u^{k-1}e^{-up}\diff u \\
&= \frac{p^k}{\Gamma(k)}\Lc\{u^{k-1}\phi_k(u)\}(p),\numberthis\label{eq:desired_relation}
\end{align*}
where $\Lc\{\cdot\}$ represents the \emph{one-sided Laplace transform} (see, \eg \cite[Ch.~29]{Korn--Korn2000}), defined as
\begin{equation*}
\Lc\{g(u)\}(p)\defeq \int_0^\infty g(\ut)e^{-p\ut}\diff\ut.
\end{equation*}
Rearranging the terms in \eqref{eq:desired_relation}, we obtain the key equation of this paper via inverse Laplace transform
\begin{align}\label{eq:estimator_function_single}
\phi_k(u)=\frac{\Gamma(k)}{u^{k-1}}\Lc^{-1}\Bigl\{\frac{f(p)}{p^k}\Bigr\}(u),
\end{align}
which we refer to as the \emph{estimator function} $\phi_k$ for $f$ with parameter $k$.
In general, inverse Laplace transform $\Lc^{-1}\{\cdot\}(\cdot)$ can be obtained by the \emph{Bromwich integral}, which is the contour integral
\begin{align*}
\Lc^{-1}\{f(p)\}(u) = \frac{1}{2\pi i}\lim_{T\to\infty}\int_{\gamma-iT}^{\gamma+iT} e   ^{pu}f(p)\diff p,
\end{align*}
where $\gamma$ is chosen so that all singularities of $f(p)$ lie to the left
of the vertical line $\mathrm{Re}(p)=\gamma$ in the complex plane and that $f(p)$ is bounded on the line
(see, \eg \cite[Ch.~2]{Cohen2007}).
For most cases of our interest (see Tables~\ref{table:estimator_functions_single} and \ref{table:estimator_functions_two}), however, inverse Laplace transforms can be computed using known transforms of elementary functions~\cite{Korn--Korn2000}, along with several properties of Laplace transform, such as linearity, time-scaling, and convolution. 
The reader is referred to Table~\ref{supp:table:inverse_laplace} in Appendix~\ref{supp:sec:exmps} for a list of elementary Laplace transforms.
Note, for example, that by the linearity of the inverse Laplace transform, if $\phi_k$ is the estimator function for $f$, then the estimator function for $af+b$ is $a\phi_k+b$ for any $a,b\in\Real$.
Concrete examples of estimator functions for different choices of $f$ are presented in Table~\ref{table:estimator_functions_single}.
See Appendix~\ref{supp:sec:exmps} for detailed derivation of these examples.

The main contributions of
this paper, for single-density functionals, are as follows:
By establishing the asymptotic unbiasedness
condition in~\eqref{eq:asymp_unbiasedness} and~\eqref{eq:convergence_in_expectation} of
the proposed estimator~\eqref{eq:estimator_single},
the necessity of which was first observed in the Ph.D. thesis of one of the authors~\cite[Ch.~5]{Noh2011},
and by establishing that
the variance of the estimator
also vanishes asymptotically,
we show that the proposed estimator is $L_2$-consistent under mild regularity conditions on densities.
\revision{The general statement (Corollary~\ref{cor:consistency_single}) capture the hardness of estimating a given functional based on $k$-NN statistics as a polynomial tail behavior of its corresponding inverse Laplace transform.}
For smooth, bounded densities, we also establish the polynomial convergence rate in mean-squared error (MSE) by carefully bounding nonasymptotic error terms. 
\revision{Informally, under certain regularity conditions, we establish that
\[
\E[(\hat{T}_f^{(k)}-T_f(p))^2] = \bigOtilde(m^{-\lambda(\sigma_p,a,k)}) + \bigO(m^{-1/2}),
\]
where $\sigma_p$ is the order of smoothness of the underlying distribution $p$, $a$ quantifies how much the functional $T_f$ is affected by \emph{high} densities (see~\eqref{eq:polynomial_tail}), and $\lambda(\sigma,a,k)$ is the bias rate exponent defined in \eqref{eq:lambda_classP_w_sigma}; see Section~\ref{sec:single:main_results} and Corollary~\ref{cor:single_mse} for details. }
\secondrevision{For example, when the densities are sufficiently smooth, \ie $\sigma_p\ge 1$, the rate exponent becomes $\lambda\approx 1/d$ for $k$ sufficiently large, implying the approximate MSE rate of $\tilde{O}(m^{-1/\max\{d,2\})})$.}

\subsection{The proposed double-density functional estimators}
For functionals of two densities, we naturally extend the same idea to the Laplace transform in two dimensional spaces. For $g:\Real_+^2\to\Real,$ we use $(u,v)$ and $(p,q)$ to denote ``time domain'' and ``frequency domain'' variables, respectively, and define
\begin{align*}
\Lc\{g(u,v)\}(p,q)\defeq \int_0^{\infty}\int_0^{\infty} g(\ut,\vt)e^{-p\ut}e^{-q\vt}\diff\ut\diff\vt.
\end{align*}
Note we keep dummy variables such as $u$ and $v$ in $\Lc\{g(u,v)\}(p,q)$ explicit, so as to avoid any confusion on which function is being transformed.
We define 
the \emph{estimator function} $\phi_{kl}$ for $f$ with parameters $(k,l)$, computed through the two-dimensional inverse Laplace transform, as
\begin{align}
\label{eq:estimator_function_double}
\phi_{kl}(u,v)
&=\frac{\Gamma(k)\Gamma(l)}{u^{k-1}v^{l-1}}
    \Lc^{-1}\Bigl\{\frac{f(p,q)}{p^kq^l}\Bigr\}(u,v).
\end{align}
When $T_f(p,q)$ is in the form of divergence, i.e., $f(p,q)$ is a function of $p/q$,  
the corresponding estimator function $\phi_{kl}(u,v)$ is also a function of $u/v$; see Proposition~\ref{supp:prop:gandh} in Appendix~\ref{supp:sec:exmps}.
Concrete examples of estimator functions for different choices of $f$ are presented in Table~\ref{table:estimator_functions_two}.
See Appendix~\ref{supp:sec:exmps} for detailed derivations of these examples.
Given two sets of samples $\Xb_{1:m}$ from $p$ and $\Yb_{1:n}$ from $q$, we further define
\begin{equation*}
V_{ln}(\xb)\defeq V_{l} (\xb|\Yb_{1:n}) \defeq n\Leb(\Bb(\xb,r_l(\xb|\Yb_{1:n}))).
\end{equation*}
We then propose a $(k,l)$-NN estimator of the form
\begin{equation}\label{eq:estimator_double}
\That_f(\Xb_{1:m},\Yb_{1:n}) 
= \frac{1}{m}\sum_{i=1}^m\phi_{kl}(U_{km}(\Xb_i), V_{ln}(\Xb_i)).
\end{equation}
As in the single-density case, we establish
the $L_2$-consistency
and  MSE convergence rate
of our estimator~\eqref{eq:estimator_double} under respective regularity conditions.

Throughout the paper, we assume the Euclidean distance, \ie $\rho(\xb,\yb)=\|\xb-\yb\|$, but the results will continue to hold for the $p$-norm ($p\ge 1$) with minor modifications; see Section~\ref{sec:conclusion} for related remarks.

\textbf{Notation.}
We use $\rvol_d(v)\defeq (v/\ups_d)^{1/d}$ to denote the radius of a $d$-dimensional ball of a volume $v$ and $\ups_d(r)\defeq \rvol_d^{-1}(r)=\Leb(\Bb(0,r))$ to denote the volume of ball of radius $r$.
We further use $\ups_d\defeq\ups_d(1)=2^{d}\Gamma(1+\frac{1}{2})^d\Gamma(1+\frac{d}{2})^{-1}$ to denote the volume of the unit ball $\Bb(0,1)$.
We denote the density of a random variable $U$ as $\rho_U(u)$.
We use the calligraphic letters $\P$ and $\Q$ to denote the probability measures corresponding to the density $p$ and $q$, respectively, and denote the support of a density $p$ as
\[
\supp(p)\defeq 
\{\xb\in\Real^d\suchthat \P(\Bb(\xb,r))>0,~\forall r>0\}.
\]
We use $\P\ll \Q$ to denote
the absolute continuity of $\P$ with respect to $\Q$.
For nonnegative functions $A(x)$ and $B(x)$ of $x\in\Xc$, we write $A(x)\lesssim_{\a} B(x)$ if there exists $C(\a)>0,$ depending only on some parameter $\a,$ such that $A(x)\le C(\a)B(x)$ for all $x\in\Xc.$
We use the standard Bachmann--Landau notation $\bigO$ and $\Theta$ (see, \eg \cite{Cormen--Leiserson--Rivest--Stein2009}) throughout the paper, and write $f(n) = \bigOtilde(g(n))$ to represent the polylogarithmic order $f(n) = \bigO(g(n)(\ln g(n))^k)$ for some $k\in\Real.$ 
We use the shorthand notation $a\wedge b=\min\{a,b\}$ and $a\vee b=\max\{a,b\}$.
Finally, $\ones_A$ stands for the indicator function of a set $A$.

\paragraph*{Organization} The rest of the paper is organized as follows. 
Section~\ref{sec:related_work} discusses the relevant literature and positions our contributions in that context.
We analyze the proposed estimator for functionals 
of one density (cf.~\eqref{eq:estimator_single} and \eqref{eq:estimator_function_single}) in Section~\ref{sec:single}
and of two densities (cf.~\eqref{eq:estimator_double} and \eqref{eq:estimator_function_double}) in
Section~\ref{sec:double}.
We discuss the convergence rate of the estimators with adaptive choices of $k$ and $l$ in Section~\ref{sec:adaptive_choice_k}.
We present in Section~\ref{sec:exp} numerical results to demonstrate the proposed estimator for a few synthetic examples.
Section~\ref{sec:conclusion} concludes the paper.

\section{Related work}\label{sec:related_work}
One of the most straightforward estimators of the density functional $T_f(p)=\E_{\Xb\sim p}[f(p(\Xb))]$ is the ``plug-in'' estimator that first forms a density estimate $\xb\mapsto \ph(\xb)$ from the samples $\Xb_{1:m}$, such as
the standard $k$-NN density estimate
\begin{align}\label{eq:knn_density}
    \ph_{km}(\xb)=\ph(\xb)=\frac{k/m}{\Leb(\Bb(\xb,r_{km}(\xb)))},
\end{align}
then plugs it in as
\begin{align}\label{eq:plugin}
    \tilde{T}_f(\ph) = \frac{1}{m}\sum_{i=1}^m f(\ph(\Xb_i)).
\end{align}
Building on the consistency of the $k$-NN density estimate $\ph_{km}$ when $k$ increases sublinearly with $m$ \cite{Loftsgaarden--Quesenberry1965,Biau--Devroye2015}, one can establish the consistency and finite-sample analysis of the plug-in estimator when $k\to\infty$~\cite{Sricharan--Raich--Hero2012,Sricharan--Wei--Hero2013,Moon--Hero2014a,Moon--Hero2014b}.
For estimating the double-density functional $T_f(p,q)=\E_{\Xb\sim p}[f(p(\Xb),q(\Xb))]$,  \citet{Berrett--Samworth2019} recently proposed a weighted version of the plug-in $(k,l)$-NN estimators of the form 
\begin{align}\label{eq:plugin_double}
    \tilde{T}_f(\ph,\qh) = \frac{1}{m}\sum_{i=1}^m f(\ph(\Xb_i),\qh(\Xb_i)),
\end{align}
with the $k$-NN density estimate $\ph_{km}$ and the $l$-NN density estimate $\qh_{kn}$ based on the samples $\Xb_{1:m}$ from $p$ and $\Yb_{1:n}$ from $q$, respectively.
They proved its efficiency by establishing a tight local asymptotic minimax lower bound and established a corresponding central limit theorem, given that $k$ and $l$ of the weighted-averaged plug-in estimators grow to infinity.

For a \emph{fixed} $k$,
however,
an appropriate ``bias correction'' is necessary for the plug-in estimator in~\eqref{eq:plugin} to be asymptotically unbiased, since the fixed-$k$-NN density estimator in~\eqref{eq:knn_density} is not consistent for a finite $k$.
A fixed-$k$ plug-in estimator with bias correction was first studied by \citet{Kozachenko--Leonenko1987}, who applied $1$-NN distances to estimate differential entropies of densities on $\Real^d$ based on an idea of \citet{Dobrushin1958}, and established the $L_2$-consistency of their estimator.
Subsequently, 
\citet{Singh--Misra--Hnizdo--Fedorowicz--Demchuk2003} and \citet{Goria--Leonenko--Mergel--Inverardi2005} 
generalized 
the 1-NN Kozachenko--Leonenko estimator to $k \ge 1$ as
\begin{align}
\label{eq:kozachenko_leonenko_estimator}
\That_{\textsf{KL}}^{(k)}(\Xb_{1:m}) 
&= \tilde{T}_f(\ph_{km})+ \ln k -\digamma(k) \\
&= \frac{1}{m}\sum_{i=1}^m 
    \ln \frac{1}{\ph_{km}(\Xb_i)} 
    + \ln k -\digamma(k),
\nonumber
\end{align}
where 
$\digamma(x)\defeq\Gamma'(x)/\Gamma(x)$ denotes the digamma function~\citep{Korn--Korn2000}.
As the canonical fixed-$k$ density functional estimator, the Kozachenko--Leonenko estimator has been investigated extensively in the literature. 
Beyond the $L_2$-consistency, \citet{Tsybakov--van-der-Meulen1996} first established $\sqrt{m}$-consistency, \ie the
$L_2$-convergence rate of $\bigO(m^{-1})$, of a truncated version of the $1$-NN Kozachenko--Leonenko estimator in $\Real$, which
was extended by
\citet{Gao--Oh--Viswanath2018tit} to $k \ge 1$ and $d \ge 1$.
Some recent developments include a central limit theorem~\cite{Delattre--Fournier2017}, results on large-$k$ behavior~\cite{Berrett--Samworth--Yuan2019}, and minimax optimality~\cite{Han--Jiao--Weissman--Wu2017,Jiao--Gao--Han2018}.

Along the same line, 
$L_2$-consistent fixed-$k$ or fixed-$(k,l)$ plug-in estimators with proper additive or multiplicative bias
correction were proposed%
\footnote{%
As pointed out in
\cite{Pal--Poczos--Szepesvari2010}, 
there are slight errors in the original analyses in
\cite{Kozachenko--Leonenko1987,Goria--Leonenko--Mergel--Inverardi2005,Wang--Kulkrani--Verdu2009,Leonenko--Pronzato--Savani2008}
when invoking asymptotic theory to establish
$L_2$-consistency. 
Correct proofs were given later in
\cite{Bulinski--Dimitrov2019a,Bulinski--Dimitrov2019b,Leonenko--Pronzato2010}.}
for KL divergence~(\citet{Wang--Kulkrani--Verdu2009}),
\Renyi{} entropies~(\citet{Leonenko--Pronzato--Savani2008}), \Renyi{} divergences~(\citet{Poczos--Schneider2011}), and several other divergences of a specific polynomial form~(\citet{Poczos--Xiong--Sutherland--Schneider2012}).
These plug-in estimators can be
expressed in general as
\begin{align}\label{eq:plugin_affine_bias_correction}
\tilde{T}_f^{\mathsf{aff}}(\ph)
&=a_k\tilde{T}_f(\ph)+b_k,
\intertext{or}
\tilde{T}_f^{\mathsf{aff}}(\ph,\qh)
&=a_{kl}\tilde{T}_f(\ph,\qh)+b_{kl},
\label{eq:plugin_affine_bias_correction2}
\end{align}
where $\ph$ is 
the fixed-$k$-NN density estimator
from $\Xb_{1:m}$ in~\eqref{eq:knn_density}, 
$\qh$ is the fixed-$l$-NN density
estimate similarly obtained from $\Yb_{1:n}$,
and
$(a_k, b_k)$ and $(a_{kl}, b_{kl})$ determine
functional-specific bias correction, respectively. 
Many density functionals beyond the special examples mentioned earlier,
however, do not allow such affine bias correction. For example,
a plug-in estimator for the
logarithmic $\a$-entropy in Table~\ref{table:estimator_functions_single}
cannot be made unbiased, even asymptotically,
by any affine bias correction.

A more general approach to 
correcting bias of the fixed-$k$
plug-in estimator was
proposed 
by 
\citet{Singh--Poczos2016}
as
\begin{align}\label{eq:singh_poczos}
\tilde{T}_{b \circ f}(\ph)=\frac{1}{m}\sum_{i=1}^m b_{km}(f(\ph_{km}(\Xb_i))),
\end{align}
which obviously subsumes affine bias correction.
This estimator was shown to be $L_2$-consistent for a fixed $k$ with definite convergence rate if there exists
a bias-correcting
function
$b_{km}$ that satisfies
\begin{align}\label{eq:singh_poczos_requirement}
\E[b_{km}(f(\ph_{km}(\xb)))]
=\E[f(\widebar{p}_{km}(\xb))]
\end{align}
for every $m$ and any
underlying density $p$, and for $\P$-a.e.\@ $\xb$, where \begin{equation*}
\widebar{p}_{km}(\xb)
=\frac{\P(\Bb(\xb,r_{km}(\xb)))}{\Leb(\Bb(\xb,r_{km}(\xb)))}
\end{equation*}
is the average density 
over the $k$-NN ball $\Bb(\xb,r_{km}(\xb))$.
Despite the general form of this estimator,
however, 
the existence of $b_{km}$ satisfying
the stringent condition of equality
in \eqref{eq:singh_poczos_requirement}
for every $m$
could be established
only for differential entropy
(and only for KL divergence in case of
functionals of two densities).

In contrast to the existing literature, our estimator
\begin{align}\label{eq:estimator_k_phi}
\That_f^{(k)}(\Xb_{1:m}) 
&= \frac{1}{m}\sum_{i=1}^m\phi_k(U_{km}(\Xb_i))\\
&= \frac{1}{m}\sum_{i=1}^m \phi_k\Bigl(\frac{k}{\ph_{km}(\Xb_i)}\Bigr)\nonumber
\end{align}
bypasses the whole bias correction issue of the plug-in approach by specifying the estimator function $\phi_k$ directly via the inverse Laplace transform~\eqref{eq:estimator_function_single}. 
Here, we identified that $U_{km}(\xb)=k/\ph_{km}(\xb)$ by the respective definitions in \eqref{eq:def_ukm} and \eqref{eq:knn_density}.
Our approach naturally 
unifies all existing estimators of the form~\eqref{eq:plugin_affine_bias_correction}
or~\eqref{eq:plugin_affine_bias_correction2},
and finds new estimators for logarithmic entropies and divergences that cannot be obtained even in the most general bias-corrected form~\eqref{eq:singh_poczos} of the traditional plug-in estimator~\eqref{eq:plugin}.
For example, our estimator for the logarithmic $\a$-entropy ($f(p)=p^{\a-1}\ln(1/p)$) is characterized by the estimator function
\begin{align}\label{eq:exmp_transformed_function}
\phi_k(u)
&=\phi_k\Bigl(\frac{k}{p}\Bigr) \\
&=\frac{\Gamma(k)}{\Gamma(k-\a+1)}k^{-\a+1}p^{\a-1}\Bigl(\ln\frac{k}{p} -\Psi(k-\a+1) \Bigr),
\nonumber
\end{align}
which cannot be expressed as a function $b_{km}(f(p))$ for some $b_{km}$.

\revision{We comment on} how 
analysis techniques of the proposed estimators are
related to those in the literature.
Through the design of our estimator functions~\eqref{eq:estimator_function_single} and \eqref{eq:estimator_function_double} via inverse Laplace transform, we can naturally
extend and simplify existing analyses for differential entropy and KL divergence by \citet{Bulinski--Dimitrov2019a,Bulinski--Dimitrov2019b},
and 
establish the asymptotic unbiasedness of our estimators~\eqref{eq:estimator_single} and \eqref{eq:estimator_double} for a general functional. 
By adapting the nonasymptotic analysis for differential entropy in \citet{Gao--Oh--Viswanath2018tit}, we can also establish the bias convergence rate of the estimator for a general functional, but without truncation.
For variance analysis, we deviate from the aforementioned work \cite{Bulinski--Dimitrov2019a,Bulinski--Dimitrov2019b,Gao--Oh--Viswanath2018tit} for simplicity and deploy a technique for the Euclidean space used by \citet{Singh--Poczos2016}; see also \cite[Ch.~7]{Biau--Devroye2015}.
Note, however, that 
the established variance results of our estimator continue to hold under the $p$-norm; see Remark~\ref{rem:variance_analysis}.
Our consistency analysis
(unbiasedness and vanishing variance)
strengthens and simplifies many existing ones
including those for
\Renyi{} entropies~\cite{Leonenko--Pronzato--Savani2008}, \Renyi{} divergences~\cite{Poczos--Schneider2011}, and divergences of polynomial form~\cite{Poczos--Xiong--Sutherland--Schneider2012}.
The convergence rates for the functionals
in Tables~\ref{table:estimator_functions_single} and \ref{table:estimator_functions_two} are established in this paper for the first time, except the Kozachenko--Leonenko estimator~\cite{Tsybakov--van-der-Meulen1996,Gao--Oh--Viswanath2018tit,Singh--Poczos2016,Jiao--Gao--Han2018} and the KL divergence estimator~\cite{Singh--Poczos2016}.

\revision{In a different direction of investigation, kernel density estimator (KDE)-based approaches have been widely studied in the literature for estimation of smooth density functionals, which also include many of the examples presented in Sections~IV and VI as special cases. 
\citet{Birge--Massart1995} established a minimax optimal rate $O(m^{-\frac{8\sigma}{d+4\sigma}}+m^{-1})$ on convergence rates in MSE of estimators of certain integral functionals involving the density and its derivatives under H\"older smoothness of order $\sigma$ (Definition~\ref{def:Holder}) on the density and demonstrated that the parametric rate $O(1/m)$ is achievable if the density is sufficiently smooth, say, $\sigma\ge d/4$. 
For estimating polynomial divergence functionals, \citet{Krishnamurthy--Kandasamy--Poczos--Wasserman2014} proposed plug-in estimators corrected through estimating higher-order terms in the von Mises expansions, which may require computationally demanding numerical integration, and established a minimax lower bound $\Omega(m^{-\frac{8\sigma}{4\sigma+d}}+m^{-1})$ under H\"older smoothness of order $\sigma>0$.
\citet{Kandasamy--Krishnamurthy--Poczos--Wasserman--Robins2015} generalized this approach to more general functionals and mutual information and established similar rates. 
In another line of work, extending the boundary-corrected plug-in estimator for mutual information of \cite{Liu--Lafferty--Wasserman2012}, \citet{Singh--Poczos2014icml,Singh--Poczos2014nips} established the MSE rate $O(m^{-\frac{2\sigma}{\sigma+d}}+m^{-1})$ for a kernel-based plug-in estimator of a class of density functionals under certain regularity conditions; we remark that this approach commonly requires a prior knowledge on the support.}

\revision{Convergence of $k$-NN distance-based estimators of density functionals can be improved by using the so-called ``ensemble method'', where a convex combination of estimators with different $k$ values is used. 
\citet{Moon--Sricharan--Hero2017} studied the ensemble method for estimation of the mutual information between two continuous random variables, and demonstrated that under certain broad regularity conditions on the density, the optimal convex combination, which can be computed by solving a convex optimization problem, yields the \emph{parametric} MSE rate $O(1/m)$ provided that the density is sufficiently smooth.
In a similar spirit, \citet{Moon--Sricharan--Greenewald--Hero2018}, \citet{Noshad--Moon--Sekeh--Hero2017}, and \citet{Wisler--Moon--Berisha2018} obtained the MSE rate $O(1/m)$ for estimating the KL divergence, $f$-divergences, and a wider class of density functionals including $f$-divergences, respectively, using the ensemble method. 
Analyzing the ensemble version of the proposed estimators is beyond the scope of this paper.}

\revision{We finally remark that \citet{Nguyen--Wainwright--Jordan2010} studied the estimation of $f$-divergences through minimization of empirical risk, by formulating the problem as a convex program. They established convergence rates when the likelihood ratio between the two distributions belongs to a reproducing kernel Hilbert space. It seems, however, quite nontrivial to compare these assumptions with those on smoothness used in the present work.}

\section{Functionals of one density}\label{sec:single}
Recall 
that we define the estimator function $\phi_{k}\suchthat \Real_+\to \Real$ for a given $f\suchthat \Real_+\to\Real$, with parameter $k\in\Natural$ as
\begin{align*}
\phi_k(u)
=\frac{\Gamma(k)}{u^{k-1}}\Lc^{-1}\Bigl\{\frac{f(p)}{p^k}\Bigr\}(u),
\tag{\ref{eq:estimator_function_single}}
\end{align*}
whenever the inverse Laplace transform exists, and then define the estimator as
\begin{align*}
\That_f^{(k)}(\Xb_{1:m}) = \frac{1}{m}\sum_{i=1}^m\phi_k(U_{km}(\Xb_i)).
\tag{\ref{eq:estimator_single}}
\end{align*}

\begin{remark}\label{rem:adaptive_choice_k}
One can check that, for all the examples in Table~\ref{table:estimator_functions_single}, 
\begin{align}\label{eq:tilde_f_in_k_limit}
\lim_{k\to\infty}\phi_k\Bigl(\frac{k}{p}\Bigr)=f(p)
\end{align} 
for each $p>0$. 
In light of \eqref{eq:estimator_k_phi}, this observation heuristically indicates that our estimator becomes closer to the plug-in estimator~\eqref{eq:plugin} as we use larger, fixed $k$. This observation is consistent with our intuition that we do not need any bias correction for the plug-in estimator with very large $k$, since the plugged-in $k$-NN density estimate~\eqref{eq:knn_density} becomes consistent as $k\to\infty$ in the sample limit~\cite{Loftsgaarden--Quesenberry1965}.
\end{remark}

To analyze the proposed estimator for general functionals $T_f(p)$ in a unified manner, we abstract  polynomial tail behaviors 
of each estimator function $\phi_k(u)$
as $u\downarrow 0$ and $u\uparrow\infty$ by a pair of constants $(a_k,b_k)\in\Real^2$ such that $\abs{\phi_k(u)}\lesssim \psi_{a_k,b_k}(u)$, where we define a piecewise polynomial function $\psi_{a,b}\suchthat \Real_+\to\Real$ for $a,b\in \Real$ as
\begin{align}
\psi_{a,b}(u)\defeq
\begin{cases}
u^{a} &\text{if }0<u\le 1,\\
u^{b} &\text{if }u > 1.
\end{cases}\label{eq:polynomial_tail}
\end{align}

\revision{Note that as $a$ gets larger and $b$ gets smaller, the piecewise polynomial function $\psi_{a,b}(u)$ decays faster as $u\downarrow 0$ and as $u\uparrow\infty$, respectively. 
Therefore, $a$ and $b$ quantify the amount of contribution of low and high density values to the estimator function $\phi_k(u)$, respectively.
Consistent with the observation that such extreme density values typically make the density functional estimation problem harder, we will establish stronger statements for functionals with larger $a$ and smaller $b$.
Below we present the estimator functions for a few representative functionals.}

\begin{example}[Differential entropy~\cite{Kozachenko--Leonenko1987}]
\label{ex:differential_entropy}
For $f(p)=\ln(1/p)$ and any $k\ge 1$, we can compute, as detailed in Example~\ref{supp:ex:kl} in Appendix~\ref{supp:sec:exmps},
\[
\phi_k(u)=\ln u-\Psi(k).
\] 
Note that we can write $\digamma(k)=H_{k-1}-\gamma$ for $k\in \Natural$, where $H_k=\sum_{i=1}^k(1/i)$ denotes the $k$-th harmonic number and $\gamma\defeq \lim_{k\to\infty} (H_k-\ln k)$ denotes the Euler--Mascheroni constant~\citep{Korn--Korn2000}.
As a bound on the estimator function $\phi_k(u)$, we consider
\[
\abs{\phi_k(u)}
\lesssim |\ln u| + 1
\lesssim \psi_{-\eps,\eps}(u)
\] 
for any arbitrarily small $\eps>0$ throughout the paper.
A finer analysis without relying on the polynomial bound $\psi_{-\eps,\eps}(u)$ may lead to a marginal improvement in the resulting performance guarantee~\cite{Gao--Oh--Viswanath2018tit,Bulinski--Dimitrov2019a,Bulinski--Dimitrov2019b}, but we do not pursue that in this paper.
\end{example}
\begin{example}[$\a$-entropy~\cite{Leonenko--Pronzato--Savani2008}]\label{ex:renyi_entropy}
For $f(p)=p^{\a-1}$ ($\a\ge 0$), we refer to the density functional $T_f(p)=\int p^{\a}(\xb)\diff\xb$ as the \emph{$\a$-entropy}. In the literature, this functional appears in \citet{Renyi1961} entropy $h_\a(p)=(\ln T_f(p))/(1-\a)$ and \citet{Harvda--Charvat1967} or \citet{Tsallis1988} entropy $\tilde{h}_\a(p)=(1-T_f(p))/(\a-1)$.
For any $k\in\Natural$ such that $k>\a-1$, we can compute, as verified in Example~\ref{supp:ex:poly} in  Appendix~\ref{supp:sec:exmps},
\[
\phi_k(u)=\frac{\Gamma(k)}{\Gamma(k-\a+1)}\Bigl(\frac{1}{u}\Bigr)^{\a-1},
\]
which allows the tight polynomial bound
\[
\abs{\phi_k(u)}\lesssim\psi_{1-\a,1-\a}(u).\]
\end{example}

\begin{example}[Logarithmic $\a$-entropy]
\label{ex:generalized_alpha_entropy}
For $f(p)=p^{\a-1}\ln(1/p)$ ($\a>0$), we refer to the density functional $T_f(p)=\int p^{\a}(\xb)\ln(1/p(\xb))\diff\xb$ as the \emph{logarithmic $\a$-entropy}.
For any $k\in\Natural$ such that $k>\a-1$, we can compute, as verified in Example~\ref{supp:ex:log_alpha_div} in Appendix~\ref{supp:sec:exmps},
\[
\phi_k(u)=\frac{\Gamma(k)}{\Gamma(k-\a+1)}u^{-\a+1}(\ln u -\Psi(k-\a+1)),
\]
and we consider
\[
|\phi_k(u)|\lesssim u^{-a+1}(|\ln u|+1) \lesssim \psi_{1-\a-\eps,1-\a+\eps}
\]
for any arbitrarily small $\eps>0$ as its polynomial bound.
\end{example}

\begin{example}[Exponential $(\a,\b)$-entropy]
\label{ex:exponential_alphabeta_entropy}
For $f(p)=p^{\a-1}e^{-\b p}$ ($\a>0$, $\beta\ge0$), we refer to the density functional $T_f(p)=\int p^{\a}(\xb)e^{-\b p(\xb)}\diff\xb$ as the \emph{exponential $(\a,\b)$-entropy}.
For any $k\in\Natural$ such that $k>\a-1$, we can compute 
\[
\phi_k(u)=\frac{\Gamma(k)}{\Gamma(k-\a+1)}\frac{(u-\b)^{k-\a}}{u^{k-1}}\ones_{[\b,\infty)}(u)
\]
using time shifting property of Laplace transform from the estimator function expression of the $\a$-entropy.
The estimator function $\phi_k$ can be bounded as
\[
|\phi_k(u)|\lesssim \psi_{0,1-\a}(u)
\]
for $k\ge \a$ and cannot be bounded by a piecewise polynomial function if $k<\a$.
\end{example}

In our subsequent analysis,
regularity conditions for the consistency and convergence rate of the proposed estimator
depend on $k$ and $f$ via 
the lower tail exponent $a$ and the upper tail exponent $b$. 
By \eqref{eq:estimator_k_phi}, extreme values of $\ph_{km}$ are amplified more via $\phi_k$ as $a$ decreases and and $b$ increases.
Hence, intuitively, when $a$ is large and $b$ is small,
the regularity conditions are milder and the estimator converges faster.

\subsection{Consistency}\label{sec:single:consistency}
Focusing solely on the
asymptotic behavior of our
estimator,
we can establish the $L_2$-consistency for general functionals
under mild assumptions on densities.
To state the results rigorously, we first define certain technical conditions.
For future use in Section~\ref{sec:double:consistency} for functionals of two densities, we state the conditions in terms of two densities $p$ and $\pt$ such that $\P\ll \tilde{\P}$.
Later, we identify $\pt$ as the density $p$ for samples $\Xb_{1:m}$ or the density $q$ for samples $\Yb_{1:n}$.

\revision{For the sake of easy analysis of density functional estimators, the standard simplifying assumptions are global upper- and lower-boundedness on the underlying density $p$, \ie there exist $c>0$ and $C>0$ such that $c\le p(\xb)\le C$ for any $\xb\in\supp(p)$; note that the boundedness of the support follows from the lower boundedness of the density.
In what follows, to establish the asymptotic consistency of the proposed estimators for a larger class of densities, we will consider weaker conditions than the boundedness assumptions, similar to those in~\cite{Bulinski--Dimitrov2019a,Bulinski--Dimitrov2019b}.}

For each $r>0$, we define the local maximal operator $\maximal_r$ on $\Real^d$ for a density $p$ by
\begin{align*}
\maximal_r{p}(\xb) &\defeq \sup_{r'\in (0,r]} \frac{\P(\Bb(\xb,r'))}{\Leb(\Bb(\xb,r'))}.
\end{align*}
Similarly, for each $r>0$, we define the local minimal operator $\minimal_r$ on $\Real^d$ for a density $p$ by
\begin{align*}
\minimal_r{p}(\xb) &\defeq \inf_{r'\in (0,r]} \frac{\P(\Bb(\xb,r'))}{\Leb(\Bb(\xb,r'))}.
\end{align*}
For each $r>0$, $\xb\mapsto \maximal_r{p}(\xb)$ and $\xb\mapsto \minimal_r{p}(\xb)$ are lower- and upper-semicontinuous, respectively, and so are Borel measurable~\cite{Bulinski--Dimitrov2019a,Bulinski--Dimitrov2019b}.
\revision{In particular, $\maximal_rp(\xb)$ and $\minimal_rp(\xb)$ are pointwise upper and lower bounds, respectively, on the density $p.$}

Given a \revision{non-decreasing} function $\xi\suchthat\Real_+\to\Real_+$, for densities $p$ and $\pt$, we define the functionals 
\begin{align*}
\maxf(p,\pt; \vth,r)&\defeq \int p(\xb)(\maximal_r{\pt}(\xb))^\vth\diff\xb,\\
\minf(p,\pt;\xi,\vth,r)&\defeq \int p(\xb)\xi((\minimal_r{\pt}(\xb))^{-\vth})\diff\xb,
\end{align*}
and
\begin{align*}
\expvolfunc(p,\pt; \xi,\vth,r)&\defeq \iint_{\rho(\xb,\yb)>r} p(\xb)\pt(\yb)\xi(\ups^\vth(\rho(\xb,\yb)))\diff\xb\diff\yb
\end{align*}
for each $\vth>0$ and $r>0$.
\revision{Here we define these quantities with possibly different densities $p$ and $\pt$ for the future use with double-density functionals; for single-density functionals, the readers can simply assume $p=\pt$.
In place of the upper- and lower- boundedness assumptions on the density $\pt$, we will impose the finiteness of the expected values $\maxf(p,\pt;\vartheta,r)$ and $\minf(p,\pt;\xi,\vartheta,r)$, respectively.
Further, $\expvolfunc(p,\pt; \xi,\vth,r)$ roughly quantifies how fast $p$ and $\pt$ decay to zero in there tails.
Observe that $\expvolfunc(p,\pt; \xi,\vth,r)\to 0$ as $r\to\infty$. 
Intuitively, as the tails of $p$ and $\pt$ decay faster, the speed of convergence of $\expvolfunc(p,\pt; \xi,\vth,r)$ will be faster.
In particular, if both $p$ and $\pt$ have bounded support, then $\expvolfunc(p,\pt; \xi,\vth,r)=0$ for $r$ sufficiently large. 
Note further that $\maxf$, $\minf$, and $\expvolfunc$ become larger as $\vartheta$ increases. }

Given $k\in\Natural$ and $(a,b)\in\Real^2$, 
\revision{consider the following conditions.}

\begin{enumerate}[label=\textbf{(U$_{p\pt}$; $k,a$)},leftmargin=*]
\item \label{cond:low}
Either $a\ge 0$, or if $a<0$, 
then there exists $r>0$ such that $\maxf(p,\pt;k,r)<\infty$.
\end{enumerate}
%
%
\begin{enumerate}[label=\textbf{(L$_{p\pt}$; $\xi,b$)},leftmargin=*]
\setcounter{enumi}{1}
\item \label{cond:high}
Either $b\le 0$, or if $b>0$, then there exists $r>0$ such that $\minf(p,\pt;\xi,b,r) < \infty$
and
\begin{equation*}
\limsup_{m\to\infty} \xi(m^b)\expvolfunc\bigl(p,\pt; \xi,b,\varrho\bigl(\frac{\kappa_m}{m}\bigr)\bigr) <\infty
\numberthis\label{eq:Rcondition}
\end{equation*}
for some $\kappa_m$ such that $\kappa_m/m\to\infty$ and $(\ln\kappa_m)/m\to 0$ as $m\to\infty$.
\end{enumerate}
\revision{Recall that the polynomial tail exponents $a$ and $b$ of the the $k$-NN estimator function~\eqref{eq:estimator_k_phi} of a given density functional quantify the amount of contribution of high and low density values to the estimator, respectively.
Hence, $a$ is coupled with  $\maxf$ that captures the upper boundedness of the density, while $b$ is pertinent to $\minf$ and $\expvolfunc$ that quantify the lower boundedness.
We note that as $a$ gets larger, $k$ gets smaller, and $b$ gets smaller, conditions~\ref{cond:high} and \ref{cond:low} become weaker, thus encompassing a larger class of densities.}

Let $\Xi$ be the class of non-decreasing functions $\xi\suchthat\Real_+\to\Real_+$ such that $\xi(t)/t\to\infty$ as $t\to\infty$, that $\xi(t_1t_2)\le \xi(t_1)\xi(t_2)$ for any $x,y>t_0$ for some $t_0\in\Real_+$, and that $\omega(\xi)\defeq \inf\{\eta>1\suchthat \xi(t)/t^\eta\to0 \text{ as }t\to\infty\}<\infty$.
For example, $\xi_1(t)=(t\ln t)\vee 0 \in\Xi$ with $t_0=e$ and $\omega(\xi_1)=1$, and $\xi_2(t)=t^\alpha\in\Xi$ for $\alpha>1$ with $t_0=0$ and $\omega(\xi_2)=\alpha$.

We are now ready to state the $L_2$-consistency results.
We show separately that the bias and variance converge to zero under certain regularity conditions.
Note that all estimator functions presented in Table~\ref{table:estimator_functions_single} are continuous.
Throughout, we consider a fixed $(a,b) \in \Real^2$ for a target functional $T_f(\cdot)$ that satisfies 
$\abs{\phi_k(u)}\lesssim \psi_{a,b}(u)$, provided
that the estimator function $\phi_k(u)$ exists for $k > -a$.

\begin{theorem}[Vanishing bias]\label{thm:vanishing_bias:single}
For a target functional $T_f(\cdot)$, 
if the estimator function $\phi_k$ is continuous and the underlying density $p$ satisfies
\hyperref[{cond:low}]{\textbf{(U$_{pp}$; $k,a$)}} and 
\hyperref[{cond:high}]{\textbf{(L$_{pp}$; $\xi,b$)}}
with some function $\xi\in\Xi$,
then the estimator~\eqref{eq:estimator_single} with fixed $k>-\omega(\xi) a$ is asymptotically unbiased.
\end{theorem}

\begin{theorem}[Vanishing variance]
\label{thm:vanishing_variance:single}
For a target functional $T_f(\cdot)$, if the underlying density $p$ satisfies \hyperref[{cond:low}]{\textbf{(U$_{pp}$; $k,a$)}} and
\hyperref[{cond:high}]{\textbf{(L$_{pp}$; $\xi,b$)}} with $\xi(t)=t^2$, 
the variance of the estimator~\eqref{eq:estimator_single} with fixed $k>-2a$ converges to zero as $m\to\infty$.
\end{theorem}

Combining Theorems~\ref{thm:vanishing_bias:single} and \ref{thm:vanishing_variance:single}, the $L_2$-consistency readily follows as a corollary.
\begin{corollary}[Consistency]
\label{cor:consistency_single}
For a target functional $T_f(\cdot)$,
if the estimator function $\phi_k$ is continuous and the underlying density $p$ satisfies
\hyperref[{cond:low}]{\textbf{(U$_{pp}$; $k,a$)}} and \hyperref[{cond:high}]{\textbf{(L$_{pp}$; $\xi,b$)}} with $\xi(t)=t^2$, 
then the estimator~\eqref{eq:estimator_single} with fixed $k>-2a$ is $L_2$-consistent.
\end{corollary}

In the following examples, we illustrate how Corollary~\ref{cor:consistency_single} can be instantiated for a few representative functionals.

\begin{example}[Differential entropy; Example~\ref{ex:differential_entropy} contd.]\label{ex:consistency_differential_entropy}
Recall that for any $k\in\Natural$, $|\phi_k(u)|\lesssim \psi_{-\eps,\eps}(u)$ for arbitrarily small $\eps>0$. 
By Corollary~\ref{cor:consistency_single}, the estimator~\eqref{eq:estimator_single} is $L_2$-consistent if the underlying density $p$ satisfies that \hyperref[{cond:low}]{\textbf{(U$_{pp}$; $k,-\eps$)}} and \hyperref[{cond:high}]{\textbf{(L$_{pp}$; $\xi,\eps$)}} with $\xi(t)=t^2$ for some $\eps>0$.
We note that the condition~\eqref{eq:Rcondition} in \hyperref[{cond:high}]{\textbf{(L$_{pp}$; $\xi,\eps$)}} can be relaxed to a milder condition in which there exist some $\d,R>0$ such that
\begin{align*}
\iint_{\rho(\xb,\yb)>R} p(\xb)p(\yb)|\ln \ups(\rho(\xb,\yb))|^{\d}\diff\xb\diff\yb
<\infty
\end{align*}
by performing a similar analysis based on the upper bound $|\phi_k(u)|\lesssim |\ln u| + 1$, \ie without invoking the polynomial bound $\psi_{-\eps,\eps}(u)$ for an arbitrarily small $\eps>0$. 
This recovers a similar result reported in \cite{Bulinski--Dimitrov2019b}.
\end{example}

\begin{example}[$\a$-entropy; Example~\ref{ex:renyi_entropy} contd.]
\label{ex:consistency_renyi_entropy}
Recall that for any $k\in\Natural$, $|\phi_k(u)|\lesssim \psi_{1-\a,1-\a}(u)$. 
For $\a>1$, since $b=1-\a<0$, the estimator with fixed $k>2(\alpha-1)$ is $L_2$-consistent if $p$ satisfies \hyperref[{cond:low}]{\textbf{(U$_{pp}$; $k,a$)}},
which slightly generalizes
the upper-boundedness condition and the requirement $k>2\a-1$ assumed in
\citet{Leonenko--Pronzato--Savani2008}.
For $\a<1$, since $a=1-\a>0$, the estimator with fixed $k\ge 1$ is $L_2$-consistent if $p$ satisfies \hyperref[{cond:high}]{\textbf{(L$_{pp}$; $\xi,b$)}} with $\xi(t)=t^2$, for examples, if $p$ is bounded away from zero and supported over a hyperrectangle.
We remark that \citet{Leonenko--Pronzato2010} reported the $L_2$-consistency of the estimator for densities
satisfying alternate conditions when $\a<1$.
\end{example}

\subsubsection{Proof of Theorem~\ref{thm:vanishing_bias:single} (vanishing bias)}
\label{sec:proof:thm:bias:single}
If the estimator function $\phi_k$ is continuous, by the continuous mapping theorem and Proposition~\ref{prop:knndist}, we have the convergence of the statistic $\phi_k(U_{k,m-1}(\Xb_m))$ to $\phi_k(U_{k\infty}(\Xb))$ in distribution as $m\to \infty$, where $U_{k\infty}(\xb)$ is a $\GammaDist(k,p(\xb))$ random variable, independent of $\Xb\sim p$ for $\P$-a.e.\@ $\xb$.
Hence, if the sequence of random variables $(\phi_k(U_{k,m-1}(\Xb_m)))_{m\ge 1}$ is uniformly integrable, we readily establish the asymptotic unbiasedness:
\begin{align*}
\lim_{m\to\infty} \E[\That_f^{(k)}(\Xb_{1:m})]
&=\lim_{m\to\infty} \E[\ph_k(U_{k,m-1}(\Xb_m))]\\
&=\E[\phi_k(U_{k\infty}(\Xb))] 
=T_f(p).
\end{align*}
To show the uniform integrability of $(\phi_k(U_{k,m-1}(\Xb_m)))_{m\ge 1}$, we invoke the following lemma.
\begin{lemma}[De la Vall\'{e}e Poussin theorem~{\cite[Theorem~1.3.4]{Borkar1995}}]
\label{lem:de_la_vallee}
A collection of random variables $(X_i)_{i\in I}$ is uniformly integrable if and only if there exists a non-decreasing function $\xi\suchthat \Real_+\to\Real_+$
such that $\sup_{i\in I} \E[\xi(\abs{X_i})] < \infty$ and $\xi(t)/t\to\infty$ as $t\to\infty$.
\end{lemma}

Observe that we have
\begin{align*}
&\E[\xi(\abs{\phi_k(U_{k,m-1}(\Xb_m))})]\\
&=\int p(\xb)\E[\xi(\abs{\phi_k(U_{k,m-1}(\xb))})]\diff\xb\\
&\lesssim\int p(\xb)\E[\xi(\psi_{a,b}(U_{k,m-1}(\xb)))]\diff\xb\\
&=\int p(\xb)\int_0^{\infty}\xi(\psi_{a,b}(u))\diff F_{km}(u|\xb)\diff\xb.
\end{align*}
Since $\xi\in\Xi$, we have $-\int_0^1 u^k\diff\xi(u^{a\wedge 0})<\infty$ for $k>-\omega(\xi)a$ and $\int_0^\infty e^{-t}\xi(t^{b\vee 0})\diff t<\infty$, and thus we can apply Lemma~\ref{supp:lem:generic:boundedness} in Appendix~\ref{supp:sec:technical_lemmas:generic}, which yields
\begin{align*}
\limsup_{m\to\infty}
\E[\xi(\abs{\phi_k(U_{k,m-1}(\Xb_m))})]
<\infty.
\end{align*}
This ensures the uniform integrability of $(\phi_k(U_{k,m-1}(\Xb_m)))_{m\ge 1}$ by the de la Vall\'{e}e Poussin theorem (Lemma~\ref{lem:de_la_vallee}), and thus concludes the proof.
\qed

\subsubsection{Proof of Theorem~\ref{thm:vanishing_variance:single} (vanishing variance)}
\label{sec:proof:thm:variance:single}
By Lemma~\ref{supp:lem:generic_variance:single} for the Euclidean space $(\Real^d,\|\cdot\|)$, we have
\begin{align*}
\Var(\That_f^{(k)})
\le\frac{2(1+k\gamma_d)}{m}
    \{& (2k+1)\E[\phi_{k}^2(U_{k,m-1}(\Xb_m))]\\
    &+2k\E[\phi_{k}^2(U_{k+1,m-1}(\Xb_m))]
\},
\end{align*}
\revision{where $\gamma_d$ is a constant which only depends on $d$; see Lemma~\ref{supp:lem:generic_variance:single}.}
Since $\xi(t)=t^2$ and $k>-2a$ imply that $-\int_0^1 u^k\diff\xi(u^{a\wedge 0})<\infty$ and $\int_0^\infty e^{-t}\xi(t^{b\vee 0})\diff t<\infty$, we can apply Lemma~\ref{supp:lem:generic:boundedness}, which ensures for $k'\in\{k,k+1\}$ that 
\[
\limsup_{m\to\infty}\E[\phi_k^2(U_{k',m-1}(\Xb_m))]
<\infty.
\]
It establishes $\Var(\That_f^{(k)})=\bigO(m^{-1})$ for $m$ sufficiently large.
\qed

\begin{remark}
\label{rem:variance_analysis}
The variance analysis relies on the Efron--Stein inequality (Lemma~\ref{supp:lem:efron_stein}) and a covering lemma (Lemma~\ref{supp:lem:knn_cone_covering}) that only applies to the Euclidean space; see Appendix~\ref{supp:sec:tech_lem_var}.
An idea for the generic variance bound (Lemma~\ref{supp:lem:generic_variance:single}) first appeared in \citet{Singh--Poczos2016} as a generalization of a technique for analyzing the 1-NN Kozachenko--Leonenko estimator by \citet[Ch.~7]{Biau--Devroye2015}, and has been employed in the literature to bound the variance of $k$-NN based estimators; see, \eg \citet{Moon--Sricharan--Hero2017}.
We note that one can attain the same rate (up to polylogarithmic factors) under the $p$-norm, by instead adapting the analysis in \citet{Gao--Oh--Viswanath2018tit}. 
As it demands a rather involved argument to bound a covariance term, however, we present a simpler approach in this paper.
\end{remark}

\subsection{Convergence rates for smooth, bounded densities}
\label{sec:single:main_results}
So far, we have established the $L_2$-consistency of the proposed estimator for general functionals under mild assumptions on densities.
Under rather stronger assumptions such as smoothness and boundedness,
we can actually establish the convergence rate of the proposed estimator in MSE.
Specifically, we consider certain regularity conditions adapted from \cite{Gao--Oh--Viswanath2018tit}.

First, we assume that
\begin{enumerate}[label=\textbf{(U$_p$)},leftmargin=*]
\item\label{cond:bounded_above}
there exists $0<\Cab_p<\infty$ such that $p(\xb)\le \Cab_p$ almost everywhere (a.e.).
\end{enumerate}
Further, we impose a few conditions related to lower-boundedness of the density, that is,
\begin{enumerate}[label=\textbf{(L\arabic*$_p$)},leftmargin=*]
\setcounter{enumi}{0}
\item\label{cond:bounded_below}
there exists $\Cbe_p>0$ such that $p(\xb)\ge \Cbe_p$ for $\xb\in\supp(p)$,
\item\label{cond:bounded_support}
the support of $p$ is bounded, and
\item\label{cond:regular_support}
there exists $r>0$ such that
\begin{align}
\eta_p\defeq \inf_{\xb\in\supp(p)}\inf_{r'\in(0,r]} \frac{\Leb(\Bb(\xb,r')\cap \supp(p))}{\Leb(\Bb(\xb,r'))}>0.
\nonumber
\end{align}
\end{enumerate}
The last condition~\ref{cond:regular_support} is called the $(\eta_p,r)$-regularity of $\supp(\mu)$ in the literature~\cite{Audibert--Tsybakov2007}.

\begin{remark}\label{rem:boundedness_implies_consistency}
The upper-boundedness condition~\hyperref[{cond:bounded_above}]{\textbf{(U$_{p}$)}} implies the condition~\hyperref[{cond:low}]{\textbf{(U$_{pp};k,a$)}}, since $\maximal_r{p}(\xb)\le \Cab_{p}<\infty$ for every $\xb\in\Real^d$ and any $r>0$. 
Also, the conditions \hyperref[{cond:bounded_below}]{\textbf{(L1$_{p}$)}}, \hyperref[{cond:bounded_support}]{\textbf{(L2$_{p}$)}}, and \hyperref[{cond:regular_support}]{\textbf{(L3$_{p}$)}} on lower-boundedness of $p$ imply the condition~\hyperref[{cond:high}]{\textbf{(L$_{pp};\xi,b$)}} for any nonnegative function $\xi$, since for $b>0$ we have
\begin{align*}
\minf(p,p;\xi,b,r)
&=\int p(\xb)\xi((\minimal_r{p}(\xb))^{-b})\diff\xb \\
&\le \int p(\xb)\xi((\eta_{p}\Cbe_{p})^{-b})\diff\xb = \xi((\eta_{p}\Cbe_{p})^{-b})<\infty
\end{align*}
for some $r>0$ by \hyperref[{cond:bounded_below}]{\textbf{(L1$_{p}$)}} and \hyperref[{cond:regular_support}]{\textbf{(L3$_{p}$)}}, and $\expvolfunc(p,p;\xi,b,\varrho(\kappa_m/m)))=0$ for $m$ sufficiently large by the boundedness of the support of $p$ from \hyperref[{cond:bounded_support}]{\textbf{(L2$_{p}$)}}.
\end{remark}
We recall the following notion of H\"{o}lder continuity for smoothness of the density $p$,
which is assumed commonly in nonparametric statistics; see, \eg \cite{Birge--Massart1995,Krishnamurthy--Kandasamy--Poczos--Wasserman2014,Singh--Poczos2016,Han--Jiao--Weissman--Wu2017,Jiao--Gao--Han2018}.

\begin{definition}\label{def:Holder}
For $\sigma>0$, a function $g\suchthat\Real^d\to \Real$ is said to be \emph{$\sigma$-H\"older continuous
over an open subset $\Omega\subseteq \Real^d$} if
$g$ is continuously differentiable over $\Omega$ up to order $\kappa\defeq \ceil{\sigma}-1$ and
\begin{align}\label{eq:Holder_constant}
L(g;\Omega)\defeq
\sup_{\substack{\rb\in\Integer_+^d\\ 
				\abs{\rb}=\kappa}}
\sup_{\substack{\yb,\zb\in\Omega\\ 
                \yb\neq \zb}}
	\frac{\abs{\partial^{\rb}g(\yb)-\partial^{\rb}g(\zb)}}{\norm{\yb-\zb}^{\b}}
    <\infty,
\end{align}
where $\b\defeq \sigma-\kappa$. 
Here we use a multi-index notation (see, \eg \cite[Ch.~8]{Folland2013}), 
that is, $\abs{\rb}\defeq r_1+\cdots+r_d$ for $\rb\in \Integer_+^d$ and $\partial^\rb g(\xb)\defeq \partial^\kappa g(\xb)/(\partial x_1^{r_1}\cdots\partial x_d^{r_d}).$ 
\end{definition}

Since the density is not smooth on the boundary of the support due to the lower-boundedness condition~\ref{cond:bounded_below},
we assume a smoothness
condition on the underlying density only over the interior of its support
and 
impose a separate regularity condition on the boundary:
\begin{enumerate}[label=\textbf{(S$_p$)},leftmargin=*]
\item\label{cond:smoothness_int}
The density $p$ is $\sigma_p$-H\"older continuous over the interior of $\supp(p)$ for $\sigma_p\in(0,2]$, and
\end{enumerate}
\begin{enumerate}[label=\textbf{(B$_p$)},leftmargin=*]
\item\label{cond:boundary}
the boundary of $\supp(p)$ has finite $(d-1)$-dimensional Hausdorff measure~\cite{Folland2013}.
\end{enumerate}

Truncated versions of well-known distributions such as exponential, Gaussian, and Cauchy distributions, as well as distributions with bounded support, such as uniform distribution and beta distributions with parameters $\a,\b\ge 1,$ satisfy these conditions with $\sigma_p=2$, and the truncated Laplace distribution satisfies the conditions with $\sigma_p=1$; see Appendix~\ref{supp:sec:examples_densities} for details on these examples. 
For densities of \emph{unbounded} support,
we provide a separate
treatment using a variant of our estimator; see Section~\ref{sec:convergence_rate:truncated:unbdd_support}.

Equipped with these regularity conditions,
we upper bound the MSE
of our estimator by considering its bias and variance separately.

\begin{theorem}[Bias rate]\label{thm:bias_class1_fixed_k_single}
For a target functional $T_f(\cdot)$, if the underlying density $p$ satisfies the conditions \ref{cond:bounded_above}, \ref{cond:bounded_below}, \ref{cond:bounded_support}, \ref{cond:regular_support}, \ref{cond:smoothness_int}, and \ref{cond:boundary},
then the estimator~\eqref{eq:estimator_single} with fixed $k>-a$ satisfies
\begin{align}\label{eq:bias_bound_generic_expression}
\bigl|\E[\That_f^{(k)}]-T_f(p)\bigr|
=\bigOtilde
(m^{-\lambda(\sigma_p,a,k)})
\end{align}
as $m\to\infty$, where 
\begin{align}
\lambda(\sigma,a,k)
&=\begin{cases}
\frac{1}{d}(\sigma\wedge 1)(\frac{k+a}{k-1})
	& \text{if }a\le -\frac{\sigma}{d}-1,\\
\frac{1}{d}(\sigma\wedge \frac{k+a}{k-1})
	&\text{if }-\frac{\sigma}{d}-1 < a \le -1,\\
\frac{1}{d}(\sigma\wedge 1)
    &\text{if }a > -1.
\end{cases} 
\label{eq:lambda_classP_w_sigma}
\end{align}
\end{theorem}

\begin{remark}
Since $k > -a$ is required to apply Theorem~\ref{thm:bias_class1_fixed_k_single},
when $a \le -1$
(for example, the $2$-entropy),
our estimator is well-defined 
and $\lambda$ in
\eqref{eq:lambda_classP_w_sigma}
is positive only for $k > 1$.
Conversely, our bias bound holds for
1-NN estimators of any functional $T_f(p)$ with estimator function $\phi_1(u)$ of lower tail exponent $a > 1$, the examples of which include differential entropy, the $\a$-entropy with $\a<2$, the logarithmic $\a$-entropy with $\a<2$, and exponential $(\a,\b)$-entropy with $\a\le 1$ in Table~\ref{table:estimator_functions_single}.
\end{remark}

\begin{remark}
\label{rem:rate_increasing_ak}
\secondrevision{The rate exponent $\lambda$ increases as the lower-tail-polynomial exponent $a$ increases, or equivalently, the estimator function $\phi_k(u)$ converges to 0 faster as $u\downarrow 0$.
If $a$ is independent of $k$, the rate exponent $\lambda$ becomes larger with larger $k$. 
In Section~\ref{sec:adaptive_choice_k}, we show that a properly growing $k$ in sample size can guarantee the largest rate exponent in~\eqref{eq:lambda_classP_w_sigma}. 
Note, however, that if $a$ decreases as $k$ increases, which is the case for some exceptional cases (Examples~\ref{ex:estimating_asymp_nn_classfication}~and~\ref{ex:estimating_JSD}), the rate exponent could become slower with larger $k$. 
This is in contrast to the large-$k$ requirement for \emph{plug-in} estimators, to guarantee the underlying $k$-NN density estimate to be consistent. 
We remind that our estimator is designed to be asymptotically unbiased for every fixed $k$, without appealing to the consistency of the $k$-NN density estimator, and it thus does not contradict the behavior of plug-in estimators. }
\end{remark}

\begin{remark}
\secondrevision{The upper tail exponent $b$ appears only in the exponent of polylogarithmic factors $\bigO(\polyln(m))$ in the rate, and thus is hidden by $\bigOtilde$ in \eqref{eq:lambda_classP_w_sigma}. 
At a finer scale, the rate increases as $b$ decreases; see the proof of Theorem~\ref{supp:lem:generic_outer_bias_case1} and Lemmas~\ref{supp:lem:generic_inner_bias:single_case1} and \ref{supp:lem:generic_outer_bias_case1} in Appendix~\ref{supp:sec:proof:thm:bias_class1_fixed_k_single}.}
\end{remark}

The variance of the estimator can be
bounded without the smoothness conditions.
\begin{theorem}[Variance rate]\label{thm:variance_rate:single}
For a target functional $T_f(\cdot)$, 
if the underlying density $p$ satisfies 
\ref{cond:bounded_above}, \ref{cond:bounded_below}, 
\ref{cond:bounded_support},
and
\ref{cond:regular_support},
then the estimator~\eqref{eq:estimator_single} with fixed $k>-2a$ satisfies
\begin{align}\label{eq:variance:single}
\Var(\That_f^{(k)})
=\bigO(m^{-1}). 
\end{align}
\end{theorem}

Combining Theorem~\ref{thm:bias_class1_fixed_k_single} on bias and Theorem~\ref{thm:variance_rate:single} on variance, we can obtain the convergence rate in MSE and establish the $L_2$-consistency  of the estimator.

\begin{corollary}[Convergence rate]\label{cor:single_mse}
Under the same assumptions in Theorem~\ref{thm:bias_class1_fixed_k_single}, then the estimator~\eqref{eq:estimator_single} with fixed $k>-2a$ satisfies
\begin{align}
\E\bigl[\bigl(\That_f^{(k)}-T_f(p)\bigr)^2\bigr]
=\bigOtilde(m^{-2\lambda(\sigma_p,a,k)}+m^{-1}).
\label{eq:mse_rate_classP}
\end{align}
\end{corollary}

\begin{remark}
For $d\ge 2,$
the bias bound always dominates the variance bound so that the MSE is bounded as $\bigOtilde(m^{-2\lambda})$.
For $d=1$, the variance bound may dominate the bias bound, depending on $\sigma_p$ and $a$.
\end{remark}

\begin{remark}
\label{rem:smoothness}
\revision{We note that the bias rate of the proposed estimators under H\"older smoothness of order $\sigma>0$ is at most $O(m^{-(\sigma\wedge 1)/d})$; it may be improved to $O(m^{-(\sigma\wedge 2)/d})$ if the \emph{boundary bias} is ignored, as remarked in \cite{Gao--Oh--Viswanath2018tit}, but it still suffers the curse of dimensionality. 
As pointed out in \citet{Jiao--Gao--Han2018}, it is an inherent problem with any \emph{positive}-kernel-based estimator that a higher smoothness $\sigma>2$ cannot be exploited in density functional estimation~\cite[Chapter~1]{Tsybakov2009}.
In particular, the key component in our analysis is Lemma~\ref{supp:lem:GOVLemma4_smoothing} from \cite{Jiao--Gao--Han2018}, which cannot be improved for $\sigma>2$.
See \cite{Han--Jiao--Weissman--Wu2017} for an extensive deliberation on this issue and see \cite{Delattre--Fournier2017,Sricharan--Raich--Hero2012,Sricharan--Wei--Hero2013,Moon--Sricharan--Hero2017,Moon--Sricharan--Greenewald--Hero2018,Noshad--Moon--Sekeh--Hero2017,Wisler--Moon--Berisha2018} for a solution based on the jackknife idea for some density functionals.
Providing a remedy to the limitation of the proposed estimators is left as an open problem.}
\end{remark}

\begin{remark}
\label{rem:minimax}
\revision{An estimator of a given density functional is said to be \emph{minimax} optimal if its MSE for the worst-case density is no larger than that of any other estimator. 
In general, the established convergence rates in MSE, including the rates for divergence functional estimators in Corollaries~\ref{cor:double_mse}, are not minimax optimal~\cite{Singh--Poczos2014icml,Singh--Poczos2014nips,Krishnamurthy--Kandasamy--Poczos--Wasserman2014,Kandasamy--Krishnamurthy--Poczos--Wasserman--Robins2015} due to the suboptimal bias rates; see, \eg Example~\ref{ex:estimating_differential_entropy}.
Since our main focus is on providing unified consistency and convergent rate analyses of the proposed generic estimators, we leave proving minimax optimality under proper regularity conditions with or without modifications of the proposed estimators as important future directions.
For the special case of differential entropy, we note that \citet{Jiao--Gao--Han2018} established an asymptotic minimax optimality of the Kozachenko--Leonenko estimator~\cite{Jiao--Gao--Han2018} for for smooth densities of order $\sigma\in(0,2]$ \emph{over a torus} (no boundary condition), matching the lower bound of \cite{Han--Jiao--Weissman--Wu2017} up to a polylogarithmic factor.}
\end{remark}

\begin{example}[Differential entropy; Example~\ref{ex:differential_entropy} contd.]\label{ex:estimating_differential_entropy}
Recall from Example~\ref{ex:differential_entropy} that $|\phi_k(u)|\lesssim \psi_{-\eps,\eps}(u)$ for any arbitrarily small $\eps>0$.
Suppose that the underlying density $p$
satisfies the conditions~
\ref{cond:bounded_above}, \ref{cond:bounded_below}, \ref{cond:bounded_support}, \ref{cond:regular_support}, \ref{cond:smoothness_int}, and \ref{cond:boundary}, in Theorem~\ref{thm:bias_class1_fixed_k_single}
with some $\sigma_p \in (0,2]$. 
Then we have the bias exponent
$\lambda = {\sigma_p}/{d}$ as in the third
case of~\eqref{eq:lambda_classP_w_sigma}
and 
the variance exponent of $1$ from~\eqref{eq:variance:single}.
Consequently, by Corollary~\ref{cor:single_mse} the 
MSE of our estimator
is bounded as
$\bigOtilde(m^{-2(\sigma_p\wedge 1)/d}+m^{-1})$.
This result recovers the same MSE rate of a truncated Kozachenko--Leonenko estimator in \cite{Gao--Oh--Viswanath2018tit} for $\sigma_p=2$.
\revision{We remark that \citet{Gao--Oh--Viswanath2018tit} reported a lower bound $\Omega(m^{-\frac{16}{d+8}}+m^{-1})$ for estimating differential entropy under $\sigma=2$ and hence,
the convergence rate is not minimax optimal.}
\end{example}

\begin{example}[$\a$-entropy; Example~\ref{ex:renyi_entropy} contd.]\label{ex:estimating_renyi_entropy}
Recall from Example~\ref{ex:renyi_entropy} that $\abs{\phi_k(u)}\lesssim\psi_{1-\a,1-\a}(u)$ for any $k\in\Natural$ such that $k>\a-1$.
Hence, for densities satisfying the conditions \ref{cond:bounded_above}, \ref{cond:bounded_below}, \ref{cond:bounded_support}, \ref{cond:regular_support}, \ref{cond:smoothness_int}, and \ref{cond:boundary}, the MSE of our estimator~\eqref{eq:estimator_single} with fixed $k>2(\a-1)$
is bounded as \eqref{eq:mse_rate_classP} with the bias rate exponent
\begin{align}\label{eq:bias_rate_exp_renyi}
\lambda(\sigma_p,a,k)
=\begin{cases}
\frac{1}{d}(\sigma_p\wedge 1) &\text{if }\a<2,\\
\frac{1}{d}(\sigma_p\wedge\frac{k+1-\a}{k-1}) &\text{if }2\le \a<2+\frac{\sigma_p}{d},\\
\frac{1}{d}(\sigma_p\wedge 1)(\frac{k+1-\a}{k-1}) &\text{if }\a \ge 2+\frac{\sigma_p}{d}.
\end{cases}
\end{align}
Note that similar convergence rates can be established for the logarithmic $\a$-entropy and the exponential $(\a,\b)$-entropy.
\end{example}

\subsubsection{Proof of Theorem~\ref{thm:bias_class1_fixed_k_single} (bias rate)}\label{sec:proof:thm:bias_rate:single}
First note that $U_{km}(\Xb_1),\ldots,U_{km}(\Xb_m)$ are identically distributed, and
$U_{km}(\Xb_m)=U_{k,m-1}(\Xb_m)$ by definition; see~\eqref{eq:convention}.
Hence, we can write
\begin{align*}
\E[\That_f^{(k)}
]
&=\E[\phi_k(U_{k,m-1}(\Xb_m))]\\
&=\int
	\E[\phi_k(U_{k,m-1}(\Xb_m))|\Xb_m=\xb]p(\xb)
    \diff \xb\\
&=\int
	\E[\phi_k(U_{k,m-1}(\xb))]p(\xb)
    \diff \xb,\numberthis\label{eq:expectation_estimator_single}
\end{align*}
where the last equality holds since $\Xb_m$ and $\Xb_{1:m-1}$ are independent.
Recall from Proposition~\ref{prop:knndist} that $U_{km}(\xb)$ converges to a $\GammaDist(k,p(\xb))$ random variable $U_{k\infty}(\xb)$ for $\P$-a.e.\@ $\xb$.
Thus, by the construction~\eqref{eq:desired_relation} of the estimator function $\phi_k(u)$, we can express the density functional as
\begin{align*}
T_f(p)
&=\int f(p(\xb))p(\xb)\diff \xb
=\int
	\E[\phi_k(U_{k\infty}(\xb))]
    p(\xb)
    \diff \xb.
\end{align*}
Applying the triangle inequality, we first have
\begin{align*}
&\bigl|\E[\That_f^{(k)}]-T_f(p)\bigr|\\
&\le\int
	p(\xb)\abs*{
    \E[\phi_k(U_{k,m-1}(\xb))
    -\phi_k(U_{k\infty}(\xb))]
    }
    \diff\xb\\
&= \int p(\xb)\Bigl|\int_0^{\infty}
\phi_k(u)(\rho_{U_{k,m-1}(\xb)}(u)-\rho_{U_{k\infty}(\xb)}(u))\diff u \Bigr| \diff\xb.
\numberthis\label{eq:target_functional_single}
\end{align*}
For some real numbers $\tau_m$ and $\nu_m$ such that $0\le \tau_m \le 1 \le \nu_m <\infty$, which are to be determined later as functions of $k,a,d$, and $\sigma_p$,
we break the inner integral and apply the polynomial bound $|\phi_k(u)|\lesssim \psi_{a,b}(u)$ with the triangle inequality to obtain
\begin{align*}
\bigl|\E[\That_f^{(k)}]-T_f(p)\bigr|
&\lesssim 
I_{\textout,1}
+I_{\textin,1}
+I_{\textin,2}
+I_{\textout,2},
\numberthis\label{eq:bias_decomposition}
\end{align*}
where
\begin{align*}
I_{\textout,1}
&\defeq \E_{p}[I_{\textout,1}(\Xb)]\\
&=\E_{p}\Bigl[\int_0^{\tau_m} 
\psi_{a,b}(u)(\rho_{U_{k,m-1}(\Xb)}(u)+\rho_{U_{k\infty}(\Xb)}(u))
\diff u\Bigr],\\
I_{\textin,1}
&\defeq \E_{p}[ I_{\textin,1}(\Xb)] \\
&=\E_{p}\Bigl[ \int_{\tau_m}^1 
\psi_{a,b}(u)|\rho_{U_{k,m-1}(\Xb)}(u)
-\rho_{U_{k\infty}(\Xb)}(u)|
\diff u \Bigr],\\
I_{\textin,2}
&\defeq \E_{p}[ I_{\textin,2}(\Xb)] \\
&=\E_{p}\Bigl[\int_1^{\nu_m} 
\psi_{a,b}(u)|\rho_{U_{k,m-1}(\Xb)}(u)
-\rho_{U_{k\infty}(\Xb)}(u)|
\diff u \Bigr],
\end{align*}
and
\begin{align*}
I_{\textout,2}
&\defeq \E_{p}[ I_{\textout,2}(\Xb) ] \\
&=\E_{p}\Bigl[ \int_{\nu_m}^{\infty}
\psi_{a,b}(u)(\rho_{U_{k,m-1}(\Xb)}(u)+\rho_{U_{k\infty}(\Xb)}(u))
\diff u \Bigr].    
\end{align*}
The \emph{inner bias} terms $I_{\textin,1}$ and $I_{\textin,2}$ can be bounded by Lemma~\ref{supp:lem:generic_inner_bias:single_case1} under the conditions \ref{cond:bounded_above}, \ref{cond:smoothness_int}, and \ref{cond:boundary},
and the \emph{outer bias} terms $I_{\textout,1}$ and $I_{\textout,2}$ can be bounded by Lemma~\ref{supp:lem:generic_outer_bias_case1} under the conditions \ref{cond:bounded_above}, \ref{cond:bounded_below}, \ref{cond:bounded_support}, and \ref{cond:regular_support}.
After putting the bounds from Lemmas~\ref{supp:lem:generic_inner_bias:single_case1} and \ref{supp:lem:generic_outer_bias_case1} together, a proper choice of the break points $(\tau_m,\nu_m)$ concludes the proof; see Appendix~\ref{supp:sec:proof:thm:bias_class1_fixed_k_single} for the details.
\qed

\begin{remark}
The key step in this analysis is the decomposition in \eqref{eq:bias_decomposition}, which is based on the construction of the estimator~\eqref{eq:desired_relation} from its asymptotic unbiasedness.
Moreover, by considering only the polynomial tail behavior of each estimator function and using \eqref{eq:bias_decomposition}, our analysis can deal with a general functional in a simple, unified manner.
The rest of the bias analysis,
that is, bounding
the four bias terms, 
closely follows and naturally extends that of \cite{Gao--Oh--Viswanath2018tit} for a truncated version of the Kozachenko--Leonenko estimator of differential entropy. 
\end{remark}

\subsubsection{Proof of Theorem~\ref{thm:variance_rate:single} (variance rate)}
\label{sec:proof:thm:variance_rate:single}
Since the boundedness conditions \ref{cond:bounded_above}, \ref{cond:bounded_below}, 
\ref{cond:bounded_support},
and
\ref{cond:regular_support} imply \hyperref[{cond:low}]{\textbf{(U$_{pp}$; $k,a$)}} and
\hyperref[{cond:high}]{\textbf{(L$_{pp}$; $\xi,b$)}} (see Remark~\ref{rem:boundedness_implies_consistency}), the variance rate directly follows from the proof of Theorem~\ref{thm:vanishing_variance:single} in Section~\ref{sec:proof:thm:variance:single}.
\qed

\subsection{Convergence rates for smooth densities of unbounded support}
\label{sec:convergence_rate:truncated:unbdd_support}
Theorem~\ref{thm:bias_class1_fixed_k_single} establishes the bias rate of the proposed estimator for smooth, bounded densities that inherently assume nonsmooth boundary.
In this section, we establish convergence rate of a truncated version of the estimator for densities of unbounded support.

For functionals of one density, we define a truncated version of the estimator~\eqref{eq:estimator_single} as
\begin{align}\label{eq:estimator_single:truncated}
\Tt_f^{(k)}(\Xb_{1:m})
\defeq
\frac{1}{m}\sum_{i=1}^m
\phibar_{k}(U_{km}(\Xb_i);\tau_m,\nu_m),
\end{align}
where we define the truncated estimator function
\[
\phibar_k(u;\tau,\nu)\defeq \phi_k(u)\ones_{(\tau,\nu)}(u)
\]
and the \emph{lower and upper truncation points} $\tau_m,\nu_m\in \Real_+$ are hyperparameters such that $0\le \tau_m\le 1\le \nu_m<\infty$ that are to be determined based on the function $f$,
the dimension $d$, the number of nearest neighbors $k$, and/or the smoothness order of the underlying density $p$.

We assume the following condition on the tail behavior of the underlying density, which is
more general than
\ref{cond:bounded_below}:
\begin{enumerate}[label=\textbf{(L\arabic*$_p^{\prime}$)},leftmargin=*]
\setcounter{enumi}{0}
\item \label{cond:exp_integral_poly}
There exist $\th>0$ and $D_0>0$ such that $\int p(\xb)e^{-\b p(\xb)}\diff\xb \le D_0\b^{-\th}$ for all $\b>1$.
\end{enumerate}
This tail condition with $\th=1$ was originally considered by \citet{Tsybakov--van-der-Meulen1996} for their analysis in $\Real$.
As pointed out in \cite{Tsybakov--van-der-Meulen1996},  densities with strictly sub-exponential tails, such as Gaussian distributions, satisfy \ref{cond:exp_integral_poly} with $\th=1$.
It can also be shown that densities with polynomially decaying tails satisfy condition \ref{cond:exp_integral_poly} for some $0<\th<1$.

We additionally introduce
the following functional-dependent condition on
the behavior of the estimator
function for small density values:
\begin{enumerate}[label=\textbf{(L\arabic*$_p$)},leftmargin=*]
\setcounter{enumi}{3}
\item\label{cond:additional_regularity}
There exists $\d>0$ such that $\int p(\xb)(p(\xb))^{-(1+\d)b}\diff\xb <\infty$.
\end{enumerate}

Finally, as we consider densities with unbounded support, we assume that
\begin{enumerate}[label=\textbf{(S$_p^\prime$)},leftmargin=*]
\item \label{cond:smoothness_Rd}
the density $p$ is $\sigma_p$-H\"older continuous over $\Real^d$ \revision{for $\sigma_p\in(0,2]$},
\end{enumerate}
in place of \ref{cond:smoothness_int}.

Exclusively for the following proposition, we additionally assume that $\phi_k(u)$ satisfies $\abs{\phi_k(u)}\lesssim \psi_{a,b}(u)$, $\phi_{k}(u)$ is differentiable at any $u>0$, and $\abs{\phi_k'(u)}\lesssim \psi_{a-1,b-1}(u)$, which hold for all the examples in Table~\ref{table:estimator_functions_single}.
\begin{proposition}[Bias rate for smooth densities of unbounded support]\label{prop:bias_class2_fixed_k_single}
For a target functional $T_f(\cdot)$, if the underlying density $p$ satisfies the conditions \ref{cond:bounded_above},
\ref{cond:exp_integral_poly},
\ref{cond:additional_regularity}, and
\ref{cond:smoothness_Rd}, then the truncated estimator \eqref{eq:estimator_single:truncated} with $-a<k<-b+\th+1$ and truncation points
\begin{align}
\tau_m
&=\begin{cases}
\Th(m^{-\frac{\sigma_p}{d}\frac{1}{k-\frac{\sigma_p}{d}-1}}) & \text{if }a\le-\frac{\sigma_p}{d}-1,\\
\bigO(m^{-\frac{\sigma_p}{d}\frac{1}{k+a}}) & \text{o.w.}
\end{cases}
\label{eq:alpha_classQ}
\end{align}
and 
\begin{align}\label{eq:beta_classQ}
&\nu_m\\
&=\begin{cases}
\Th(m^{(\frac{\sigma_p}{d}\wedge 1)\frac{1}{\th-k-b+1}})& 
    \text{if }k\le -b-1, b\le-\frac{\sigma_p}{d}-1,\\
\Th(m^{\frac{\sigma_p}{d}\frac{1}{\th-k+\frac{\sigma_p}{d}+2}}\bigr)
    & \text{if }k\le -b-1, b>-\frac{\sigma_p}{d}-1,\\
\Th(m^{\frac{1}{\th+2}})
    & \text{if } k>-b-1, b\le-\frac{\sigma_p}{d}-1,\\
\Th(m^{(\frac{\sigma_p}{d}\wedge 1)\frac{1}{\th+2}})
    & \text{if } k>-b-1, b>-\frac{\sigma_p}{d}-1,
\end{cases}\nonumber
\end{align}
with $\nu_m=o(\sqrt{m})$ as $m\to\infty$
satisfies
\begin{align}
\bigl|\E\bigl[\Tt_f^{(k)}\bigr]-T_f(p)\bigr|
&=\bigO
\bigl(m^{-\lambda_\tau\wedge\lambda_\nu}\bigr),
\nonumber
\end{align}
where
\begin{align}
\lambda_\tau
=\begin{cases}
\frac{\sigma_p}{d}\frac{k+a}{k-\frac{\sigma_p}{d}-1} & \text{if }a\le-\frac{\sigma_p}{d}-1,\\
\frac{\sigma_p}{d} & \text{o.w.},
\end{cases}
\label{eq:lambda_tau_classQ}
\end{align}
and
\begin{align}
&\lambda_\nu\label{eq:lambda_beta_classQ}
\\
&= \begin{cases}
    \frac{\sigma_p}{d}\wedge 1
    & \text{if }k\le-b-1, b\le-\frac{\sigma_p}{d}-1,\\
    (\frac{\sigma_p}{d}(1-\frac{b+\frac{\sigma_p}{d}+1}{\th-k+\frac{\sigma_p}{d}+2}))\wedge 1
    & \text{if }k\le-b-1, b>-\frac{\sigma_p}{d}-1,\\
    \frac{\sigma_p}{d}\wedge(1-\frac{k+b+1}{\th+2})
    & \text{if }k>-b-1, b\le-\frac{\sigma_p}{d}-1,\\
    (\frac{\sigma_p}{d}\wedge 1) (1-\frac{k+b+1}{\th+2})
    & \text{if }k>-b-1, b>-\frac{\sigma_p}{d}-1.
\end{cases}\nonumber
\end{align}
\end{proposition}

We can establish the variance rate with truncation under only the upper-boundedness condition, without explicitly imposing the condition $k>-2a$ as required in Theorem~\ref{thm:variance_rate:single}.
\begin{proposition}[Variance rate of truncated estimator]\label{proposition:variance:single:truncated}
For a target functional $T_f(\cdot)$,
if the underlying density $p$ satisfies \ref{cond:bounded_above}, then
the estimator~\eqref{eq:estimator_single:truncated} with $k>-a$ satisfies
\begin{align}\label{eq:variance:single:truncated}
\Var\bigl(\Tt_f^{(k)}\bigr)
=\bigO
\Bigl(
\frac{k^2}{m}\bigl(k^{-k}\tau_m^{(k+2a)\wedge 0}+\nu_m^{2b\vee 0}\bigr)
\Bigr).
\end{align}
\end{proposition}

Combining Propositions~\ref{prop:bias_class2_fixed_k_single} and \ref{proposition:variance:single:truncated}, we can obtain a corresponding consistency result as in Corollary~\ref{cor:single_mse}, the formal statement of which is omitted.

At face value, 
Proposition~\ref{prop:bias_class2_fixed_k_single} enlarges considerably the class of densities under the purview of our analyses. 
On the flip side, however, it requires the underlying density to be smooth over the whole of $\Real^d$ and this rules out, for example, the uniform distribution, which is covered by \ref{cond:bounded_below}.
Thus, Proposition~\ref{prop:bias_class2_fixed_k_single}
and Theorem~\ref{thm:bias_class1_fixed_k_single} complement each other.

The stringent requirement $k<-b+\th+1$ in Proposition~\ref{prop:bias_class2_fixed_k_single} is due to a bias term $\bigO(\nu_m^{b+k-1-\th})$ that appears in the analysis; a smaller $k$, which is, of course, still larger than $-a$, gives a tighter bound on this term, whereas a larger $k$ is desired to reduce the bias due to the lower truncation. 
Proposition~\ref{prop:bias_class2_fixed_k_single} thus cannot guarantee the $L_2$-consistency of the estimator when $k$ grows as $m\to\infty$, as the condition $k<\th-b+1$ is violated. 

\begin{example}[Differential entropy; Example~\ref{ex:differential_entropy} contd.]\label{ex:estimating_differential_entropy2}
For estimating differential entropy, recall that $|\phi_k(u)|\lesssim \psi_{-\eps,\eps}(u)$ for arbitrarily small $\eps>0$.
Consider densities that satisfy the conditions \ref{cond:bounded_above},
\ref{cond:exp_integral_poly},
\ref{cond:additional_regularity}, and
\ref{cond:smoothness_Rd} for some $0<\th\le 1$.
Since Proposition~\ref{prop:bias_class2_fixed_k_single} requires $k<\th+1-\eps$, we need to choose $k=1$ to guarantee the $L_2$-consistency of our estimator.
We obtain a bias bound $\bigO(m^{-\frac{\th-\eps}{\th+2}(\frac{\sigma_p}{d}\wedge 1)})$, a variance bound $\bigO(m^{-(1-\d)})$ for arbitrarily small $\d>0$ from Proposition~\ref{proposition:variance:single:truncated}, and thus the MSE rate $\bigO(m^{-\frac{2(\th-\eps)}{\th+2}(\frac{\sigma_p}{d}\wedge 1)})$.
In particular, for one-dimensional densities with $\sigma_p\ge 1$ and $\th=1$, we obtain the MSE rate $\bigO(m^{-\frac{2(1-\eps)}{3}})$.
Note that this rate is slightly worse than $\bigO(m^{-1})$, as obtained by \citet[Section~2, pp.\@ 77--78]{Tsybakov--van-der-Meulen1996} under different regularity conditions with a faster growing upper truncation point $\nu_m=\Th(\sqrt{m})$.
\end{example}

\begin{example}[$\a$-entropy; Example~\ref{ex:renyi_entropy} contd.]
\label{ex:estimating_renyi_entropy2}
Consider estimating the $\a$-entropy ($\a\neq 1$) of densities that satisfy the conditions \ref{cond:bounded_above},
\ref{cond:exp_integral_poly},
\ref{cond:additional_regularity}, and
\ref{cond:smoothness_Rd} with some $\th>0$.
Since $|\phi_k(u)|\lesssim \psi_{1-\a,1-\a}(u)$, we need to use $k\in(\a-1,\a+\th)$ for our estimator to apply Proposition~\ref{prop:bias_class2_fixed_k_single}.
By setting the truncation points as
\begin{align*}
&(\tau_m,\nu_m)\\
&=\begin{cases}
(\bigO(m^{-\frac{\sigma_p}{d}\frac{1}{k-\a+1}}),\Theta(m^{(\frac{\sigma_p}{d}\wedge 1)\frac{1}{\th+2}})), &\text{if }\a<\frac{\sigma_p}{d}+2,\\
(\Th(m^{-\frac{\sigma_p}{d}\frac{1}{k-\frac{\sigma_p}{d}-1}}), \Th(m^{\frac{1}{\th+2}})), &\text{if }\a\ge\frac{\sigma_p}{d}+2,
\end{cases}
\end{align*}
our estimator achieves the bias rate $\bigO(m^{-(\lambda_\tau\wedge \lambda_\nu)})$,
where
\begin{align*}
(\lambda_\tau,\lambda_\nu)
&=\begin{cases}
(\frac{\sigma_p}{d}, (\frac{\sigma_p}{d}\wedge 1)\frac{\th+\a-k}{\th+2}) &\text{if }\a<\frac{\sigma_p}{d}+2,\\
(\frac{\sigma_p}{d}\frac{k-\a+1}{k-\frac{\sigma_p}{d}-1}, \frac{\sigma_p}{d}\wedge \frac{\th+\a-k}{\th+2}) &\text{if }\a\ge\frac{\sigma_p}{d}+2.
\end{cases}
\end{align*}
From Proposition~\ref{proposition:variance:single:truncated}, we can bound the variance of our estimator as $\bigO(m^{-\lambda_{\mathrm{v}}})$, where
\begin{align*}
\lambda_{\mathrm{v}}
&=\begin{cases}
1-(\frac{\sigma_p}{d}\wedge 1)\frac{2(1-\a)\vee 0}{\th+2} &\text{if }\a<\frac{\sigma_p}{d}+2,\\
1-\frac{\sigma_p}{d}\frac{(2\a-k-2)\vee 0}{k-\frac{\sigma_p}{d}-1} &\text{if }\a\ge\frac{\sigma_p}{d}+2,
\end{cases}
\end{align*}
and thus we establish the MSE rate $\bigO(m^{-2(\lambda_\tau\wedge\lambda_\nu)}+m^{-\lambda_{\mathrm{v}}})$.
\end{example}

\begin{remark}
We remark in passing on the consistency of the truncated estimator (without convergence rate analysis). 
With lower truncation point $\tau_m$ such that $\tau_m^{k+2a}=o(m)$, the conditions $k>-2a$ can be relaxed to $k>-a$ in Corollary~\ref{cor:consistency_single}. 
Moreover, a very mild upper truncation of speed $\nu_m=e^{o(m)}$ can relax the condition~\hyperref[{cond:high}]{\textbf{(L$_{pp}$)}} assumed in the consistency results to a milder one, \ie
\begin{enumerate}[label=\textbf{(L$_{p\pt}^\prime$; $\xi,b$)},leftmargin=*]
\item \label{cond:high:milder}
Either $b\le 0$, or if $b>0$, then there exists $r>0$ such that $\minf(p,\pt;\xi,b,r) < \infty$
\end{enumerate}
with $\pt=p$.
\end{remark}

\section{Functionals of two densities}
\label{sec:double}
We now consider estimating a functional $T_f(p,q)$ of two densities $p$ and $q$.
Henceforth, we assume that $\P\ll\Q$.
Recall that for fixed $k,l\in\Natural$ and a given $f\suchthat \Real_+^2\to\Real$,
we define the \emph{estimator function} $\phi_{kl}\suchthat \Real_+^2\to\Real$ of $f$ with parameters $k,l$ as
\begin{align*}
\phi_{kl}(u,v)
&=\frac{\Gamma(k)\Gamma(l)}{u^{k-1}v^{l-1}}
    \Lc^{-1}\Bigl\{\frac{f(p,q)}{p^kq^l}\Bigr\}(u,v),
\tag{\ref{eq:estimator_function_double}}
\end{align*}
whenever the inverse Laplace transform exists, and then define the estimator as
\[
\That_f^{(k)}(\Xb_{1:m},\Yb_{1:n}) = \frac{1}{m}\sum_{i=1}^m\phi_{kl}(U_{km}(\Xb_i),V_{ln}(\Yb_i)).
\tag{\ref{eq:estimator_double}}
\]
Here we define \[
V_{ln}(\xb) \defeq U_{l}(\xb|\Yb_{1:n}) = n\Leb(\Bb(\xb,r_l(\xb|\Yb_{1:n}))).\] 

\begin{remark}
Similar to the observation made in Remark~\ref{rem:adaptive_choice_k}, an analogous limiting behavior 
\begin{align*}
\lim_{k,l\to\infty} \phi_{kl}\Bigl(\frac{k}{p},\frac{l}{q}\Bigr)
=f(p,q)
\end{align*}
can be verified for all the examples in Table~\ref{table:estimator_functions_two} except Le Cam distance and Jensen--Shannon divergence.
\end{remark}

As for the single-density case, a polynomial tail behavior of the estimator function $\phi_{kl}(u,v)$ affects the convergence rate of each instantiated estimator.
We describe a tail behavior of $\phi_{kl}(u,v)$ by a quadruple $(a_{kl},b_{kl},\at_{kl},\bt_{kl}) \allowbreak\in\Real^4$ such that $\abs{\phi_{kl}(u,v)}\lesssim \psi_{a_{kl},b_{kl}}(u)\psi_{\at_{kl},\bt_{kl}}(v)$. 
This characterization allows us to handle the convergence of $U_{km}(\xb)$ and $V_{ln}(\xb)$ separately so that we can extend the analysis for the single-density case in a straightforward manner.
Note that for all the examples presented in Table~\ref{table:estimator_functions_two}, $(a_{kl},b_{kl},\at_{kl},\bt_{kl})$ can be found as constants independent of $k$ and $l$, except Le Cam distance and Jensen--Shannon divergence.
Also note that all the estimator functions $\phi_{kl}(u,v)$ presented in Table~\ref{table:estimator_functions_two} are continuous.

\begin{example}[KL divergence~\cite{Wang--Kulkrani--Verdu2009}]
\label{ex:kl_divergence}
For $f(p,q) = \ln(p/q),$ we can compute, as shown in Example~\ref{supp:ex:kl} in Appendix~\ref{supp:sec:exmps},
\[
\phi_{kl}(u,v) = \ln\frac{v}{u}+H_{k-1}-H_{l-1}.
\]
As a bound on the estimator function $\phi_{kl}(u,v)$, we consider 
\begin{align*}
|\phi_{kl}(u,v)|
&\lesssim 1+|\ln u|+|\ln v|\\
&\lesssim (1+|\ln u|)(1+|\ln v|)
\lesssim \psi_{-\eps,\eps}(u)\psi_{-\eps,\eps}(v)
\end{align*}
for any arbitrarily small $\eps>0$.
\end{example}
\begin{example}[Polynomial functional~\cite{Poczos--Schneider2011,Poczos--Xiong--Sutherland--Schneider2012}]
\label{ex:polynomial_functional}
For $f(p,q)=p^{\a-1}q^{\b}$ $(\a>0,\b>1-\a)$ and any $k,l\in\Natural$ such that $k>\a-1$ and $l>\b$, we can compute, as shown in Example~\ref{supp:ex:poly} in Appendix~\ref{supp:sec:exmps},
\begin{align*}
\phi_{kl}(u,v)=\frac{\Gamma(k)\Gamma(l)}{\Gamma(k-\a+1)\Gamma(l-\b)}u^{1-\a}v^{-\b},
\end{align*}
which allows the tight polynomial bound
\[
\abs{\phi_k(u)}\lesssim\psi_{1-\a,1-\a}(u)\psi_{-\b,-\b}(v).
\]
This class of polynomial functionals includes many important functionals. 
For the special instance of $\b=1-\a$, we refer to the density functional $T_f(p,q)=\int p^{\a}(\xb)q^{1-\a}(\xb)\diff\xb$ as the \emph{$\a$-divergence}, which appears in the literature in a few different forms; see, \eg \citet{Renyi1961} and \citet{Cichocki--Lee--Kim--Choi2008}.
\end{example}

\begin{example}[Logarithmic $\a$-divergence]
\label{ex:logarithmic_alpha_divergence}
For $f(p,q)=(p/q)^{\a-1}\ln(p/q)$ $(\a>0)$, we refer to the density functional $T_f(p,q)=\int p^{\a}(\xb)q^{1-\a}(\xb)\ln(p(\xb)/q(\xb))\diff\xb$ as the \emph{logarithmic $\a$-divergence}.
For any $k,l\in\Natural$ such that $k>\a-1$ and $l>1-\a$, we can compute, as shown in Example~\ref{supp:ex:log_alpha_div} in Appendix~\ref{supp:sec:exmps},
\begin{align*}
\phi_{kl}(u,v)
&=\frac{\Gamma(k)\Gamma(l)}{\Gamma(k-\a+1)\Gamma(l+\a-1)}\\
&\quad\times u^{-\a+1}\Bigl(\ln\frac{v}{u}+\digamma(k-\a+1)-\digamma(l+\a-1)\Bigr).
\end{align*}
As a bound on the estimator function $\phi_{kl}(u,v)$, we consider
\begin{align*}
|\phi_{kl}(u,v)|
&\lesssim u^{-\a+1}v^{\a-1}(1+|\ln u|+|\ln v|)\\
&\lesssim u^{-\a+1}v^{\a-1}(1+|\ln u|)(1+|\ln v|)\\
&\lesssim \psi_{1-\a-\eps,1-\a+\eps}(u)\psi_{\a-1-\eps,\a-1+\eps}(v)
\end{align*}
for any arbitrarily small $\eps>0$.
\end{example}

\begin{example}[Le Cam distance]\label{ex:asymp_nn_classfication}
\revision{For $f(p,q)=(p-q)^2/(2p(p+q))$, the corresponding divergence functional
\begin{align*}
D_{\mathsf{LC}}(p,q)
&=\half\int \frac{(p(\xb)-q(\xb))^2}{p(\xb)+q(\xb)}\diff\xb\\
&=1-\int\frac{2 p(\xb)q(\xb)}{p(\xb)+q(\xb)}\diff\xb
\end{align*}
is called Le Cam distance~\citep[p.~47]{LeCam2012} in the literature~\citep{Polyanskiy--Wu2019}.
We note in passing that this functional has a connection to the nearest neighborhood binary classification rule: 
it is well known that the asymptotic error of the nearest neighborhood binary classification for equiprobable classes is given as $\half (1-T_f(p,q))$~\cite{Cover--Hart1967}.}
For any $k,l\in\Natural$, we can compute, as shown in Example~\ref{supp:ex:NNclassification} in Appendix~\ref{supp:sec:exmps},
\begin{align*}
\phi_{kl}(u,v)&=
2\binom{k+l-2}{k-1}^{-1}\Bigl(-\frac{u}{v}\Bigr)^{l-1}
\times\\
&\qquad
\biggl\{\sum_{i=0}^{l-1}\binom{k+l-2}{i}\Bigl(-\frac{v}{u}\Bigr)^i \\ 
&\qquad\qquad - \Bigl(1-\frac{v}{u}\Bigr)^{k+l-2}\ones_{[v,\infty)}(u)\biggr\}-1.
\end{align*}
As a bound on the estimator function $\phi_{kl}(u,v)$, we have
\begin{align*}
|\phi_{kl}(u,v)|
&\lesssim \psi_{-k+1,l-1}(u)\psi_{-l+1,k-1}(v).
\end{align*}
\end{example}

\begin{example}[Jensen--Shannon divergence]\label{ex:JSD}
When $\Q\ll\P$, we can write Jensen--Shannon divergence as
\[
\JSD(p,q)
=\half\Bigl(D\Bigl(p~\Big\|~\frac{p+q}{2}\Bigr) + D\Bigl(q~\Big\|~\frac{p+q}{2}\Bigr)\Bigr)
=T_f(p,q)
\]
for 
\begin{equation*}
f(p,q)=\half\Bigl(\frac{q}{p}+1\Bigr)\ln\frac{2}{(q/p)+1}+\frac{q}{2p}\ln \frac{q}{p},
\end{equation*} 
where $\D{p}{q}$ denotes the KL divergence between $p$ and $q$.
For any $k\ge 1$ and $l\ge 2$, we can compute, as shown in Example~\ref{supp:ex:JSD} in Appendix~\ref{supp:sec:exmps},
\begin{align*}
\phi_{kl}(u,v)
&= \half\Bigl\{\ln 2
+\frac{l-1}{k}\frac{u}{v} \Bigl( \Psi(l-1)-\Psi(k+1)+\ln 2\frac{u}{v}\Bigr)
\\
&\qquad\quad
+B_{kl}(u,v) + \frac{l-1}{k}\frac{u}{v} B_{k+1,l-1}(u,v)\Bigr\},
\end{align*}
where $B_{kl}(u,v)$ is defined in \eqref{eq:jsd_estimation_function}.
\begin{table*}
\centering
\begin{minipage}{\textwidth}
\begin{align*}
B_{kl}(u,v)
&= \begin{cases}
\displaystyle
\binom{k+l-2}{k-1}^{-1}\sum_{j=1}^{l-1}\binom{k+l-2}{k-1+j}\frac{(-u/v)^{j}}{j} & 
\displaystyle
\text{if } \frac{u}{v}< 1,\\
\displaystyle
-\ln \frac{u}{v} + \binom{k+l-2}{k-1}^{-1}\Bigl\{
-\sum_{j=-k+1}^{-1}\binom{k+l-2}{k-1+j}\frac{(-u/v)^{j}}{j}
+\sum_{\substack{j=-k+1\\j\neq 0}}^{l-1}\binom{k+l-2}{k-1+j} \frac{(-1)^{j}}{j} \biggr\}
& \displaystyle\text{if } \frac{u}{v}\ge 1.
\end{cases}
\numberthis\label{eq:jsd_estimation_function}
\end{align*}
\end{minipage}
\medskip
\hrule
\end{table*}
As a polynomial bound, we have
\[
|\phi_{kl}(u,v)| \lesssim \psi_{-k+1,l-1}(u)\psi_{-l+1,k-1}(v).
\]
\end{example}


\subsection{Consistency}
\label{sec:double:consistency}
As in Section~\ref{sec:single:consistency}, we can establish the $L_2$-consistency of the estimator of functionals of two densities under mild regularity conditions.
Throughout,
we consider a fixed 
$(a,b,\at,\bt)\in\Real^4$
for a target functional $T_f(\cdot,\cdot)$ whose estimator function $\phi_{kl}$
satisfies 
$\abs{\phi_{kl}(u,v)}\lesssim \psi_{a,b}(u)\psi_{\at,\bt}(v)$,
provided 
that the estimator function $\phi_{kl}$ exists for $k>-a$ and $l>-\at$.

\begin{theorem}[Vanishing bias]\label{thm:vanishing_bias:double}
For a target functional $T_f(\cdot,\cdot)$, if the estimator function $\phi_{kl}(u,v)$ is continuous and the underlying densities $p$ and $q$ satisfy
\hyperref[{cond:low}]{\textbf{(U$_{pp}$; $k,a$)}},
\hyperref[{cond:high}]{\textbf{(L$_{pp}$; $\xi^2,b$)}},
\hyperref[{cond:low}]{\textbf{(U$_{pq}$; $l,\at$)}}, and
\hyperref[{cond:high}]{\textbf{(L$_{pq}$; $\xi^2,\bt$)}} for some function $\xi\in\Xi$,
then the estimator~\eqref{eq:estimator_double} with $k>-2\omega(\xi)a$ and $l>-2\omega(\xi)\at$ is asymptotically unbiased as $m,n\to\infty$.
\end{theorem}

\begin{theorem}[Vanishing variance]
\label{thm:vanishing_variance:double}
For a target functional $T_f(\cdot,\cdot)$, if the underlying densities $p$ and $q$ satisfy
\hyperref[{cond:low}]{\textbf{(U$_{pp}$; $k,a$)}},
\hyperref[{cond:high}]{\textbf{(L$_{pp}$; $\xi^2,b$)}},
\hyperref[{cond:low}]{\textbf{(U$_{pq}$; $l,\at$)}}, and
\hyperref[{cond:high}]{\textbf{(L$_{pq}$; $\xi^2,\bt$)}} with $\xi(t)=t^2$,
then then the variance of the estimator~\eqref{eq:estimator_double} with fixed $k>-4a$ and fixed $l>-4\at$ converges to zero as $m,n\to\infty$.
\end{theorem}

\begin{corollary}[Consistency]
\label{cor:consistency_double}
For a target functional $T_f(\cdot,\cdot)$, if the estimator function $\phi_{kl}(u,v)$ is continuous and the underlying densities $p$ and $q$ satisfy
\hyperref[{cond:low}]{\textbf{(U$_{pp}$; $k,a$)}},
\hyperref[{cond:high}]{\textbf{(L$_{pp}$; $\xi^2,b$)}},
\hyperref[{cond:low}]{\textbf{(U$_{pq}$; $l,\at$)}}, and
\hyperref[{cond:high}]{\textbf{(L$_{pq}$; $\xi^2,\bt$)}}
with $\xi(t)=t^2$,
then the estimator~\eqref{eq:estimator_double} with fixed $k>-4a$ and fixed $l>-4\at$ is $L_2$-consistent.
\end{corollary}

In the following examples, we illustrate how Corollary~\ref{cor:consistency_double} can be instantiated for a few representative functionals. 

\begin{example}[KL divergence; Example~\ref{ex:kl_divergence} contd.]
\label{ex:consistency_kl_divergence}
Recall that for estimating differential entropy,
$|\phi_{kl}(u,v)|\lesssim \psi_{-\eps,\eps}(u)\psi_{-\eps,\eps}(v)$ for arbitrarily small $\eps>0$ and for any $k,l\in\Natural$. 
By Corollary~\ref{cor:consistency_double}, the estimator~\eqref{eq:estimator_double} with fixed $k\ge 1$ and $l\ge 1$ is $L_2$-consistent if the underlying densities $p$ and $q$ satisfy 
\hyperref[{cond:low}]{\textbf{(U$_{pp}$; $k,-\eps$)}},
\hyperref[{cond:high}]{\textbf{(L$_{pp}$; $\xi^2,\eps$)}},
\hyperref[{cond:low}]{\textbf{(U$_{pq}$; $l,-\eps$)}}, and
\hyperref[{cond:high}]{\textbf{(L$_{pq}$; $\xi^2,\eps$)}}
with $\xi(t)=t^2$.
As discussed in Example~\ref{ex:consistency_differential_entropy}, a finer analysis recovers a similar consistency result established in \cite{Bulinski--Dimitrov2019b}.
\end{example}

The proofs of the main results (Theorems~\ref{thm:vanishing_bias:double}, \ref{thm:vanishing_variance:double}, \ref{thm:bias_class1_fixed_kl_double}, and \ref{thm:variance_rate:double}) in this section follow with minor extensions to those of the single-density case, and are deferred to Appendix~\ref{supp:sec:proof:main_results}.

\begin{example}[$\a$-divergence; Example~\ref{ex:polynomial_functional} contd.]
\label{ex:consistency_renyi_divergence}
Recall that for estimating the $\a$-divergence ($\a\neq1$), we have $|\phi_{kl}(u,v)|\lesssim \psi_{1-\a,1-\a}(u)\psi_{\a-1,\a-1}(v)$ for any $k,l\in\Natural$ such that $k>\a-1$ and $l>1-\a$. 
For $\a>1$, since $b=1-\a<0$ and $\at=\a-1>0$, the estimator with fixed $k>4(\a-1)$ and $l\ge 1$ is $L_2$-consistent if the underlying densities $p$ and $q$ satisfy that 
\hyperref[{cond:low}]{\textbf{(U$_{pp}$; $k,1-\a$)}} and
\hyperref[{cond:high}]{\textbf{(L$_{pq}$; $\xi^2,\a-1$)}}
with $\xi(t)=t^2$. 
For $\a<1$, since $a=1-\a>0$ and $\bt=\a-1<0$, the estimator with $k\ge 1$ and $l>4(1-\a)$ is $L_2$-consistent if the underlying densities $p$ and $q$ satisfy that 
\hyperref[{cond:high}]{\textbf{(L$_{pp}$; $\xi^2,b$)}} and
\hyperref[{cond:low}]{\textbf{(U$_{pq}$; $l,\at$)}} 
with $\xi(t)=t^2$.
This consistency result covers a strictly larger class of densities than an earlier result by \citet{Poczos--Schneider2011}, whereby the $L_2$-consistency of the estimator with $l=k$ is established under rather stronger assumptions such as boundedness and uniform continuity of densities. 
Moreover, Propositions~\ref{cor:consistency_double} and \ref{cor:consistency_double} strengthen the $L_2$-consistency result established in \citet{Poczos--Xiong--Sutherland--Schneider2012} for a polynomial functional (see Example~\ref{ex:polynomial_functional}), which subsumes $\a$-divergence.
\end{example}

\subsection{Convergence rates for smooth, bounded densities}
\label{sec:double:main_results}

\begin{theorem}[Bias rate]\label{thm:bias_class1_fixed_kl_double}
For a target functional $T_f(\cdot,\cdot)$, if the underlying density $p$ satisfies \ref{cond:bounded_above}, \ref{cond:bounded_below}, \ref{cond:bounded_support}, \ref{cond:regular_support}, \ref{cond:smoothness_int}, and \ref{cond:boundary}, and $q$ satisfies \hyperref[{cond:bounded_above}]{\textbf{(U$_q$)}}, \hyperref[{cond:bounded_support}]{{\textbf{(L1$_q$)}}}, \hyperref[{cond:bounded_below}]{{\textbf{(L2$_q$)}}}, \hyperref[{cond:regular_support}]{{\textbf{(L3$_q$)}}}, \hyperref[{cond:smoothness_int}]{{\textbf{(S$_q$)}}}, and \hyperref[{cond:boundary}]{{\textbf{(B$_q$)}}},
then the estimator~\eqref{eq:estimator_double} with fixed $k>-a$ and $l>-\at$
satisfies
\begin{align*}
\bigl|\E[\That^{(k,l)}_f]-T_f(p,q)\bigr|
  &=\bigOtilde\bigl(
    m^{-\lambda(\sigma_p,a,k)}
    +n^{-\lambda(\sigma_q,\at,l)}\bigr),
\end{align*}
as $m,n\to\infty$, 
where the rate exponent function $\lambda(\sigma,a,k)$ is as defined in \eqref{eq:lambda_classP_w_sigma}.
\end{theorem}

\begin{theorem}[Variance rate]
\label{thm:variance_rate:double}
For a target functional $T_f(\cdot,\cdot)$, if the underlying density $p$ satisfies \ref{cond:bounded_above}, \ref{cond:bounded_below}, \ref{cond:bounded_support}, and \ref{cond:regular_support}, and $q$ satisfies
\hyperref[{cond:bounded_above}]{\textbf{(U$_q$)}}, \hyperref[{cond:bounded_below}]{\textbf{(L1$_q$)}}, \hyperref[{cond:bounded_support}]{\textbf{(L2$_q$)}}, and \hyperref[{cond:regular_support}]{\textbf{(L3$_q$)}},
then the estimator~\eqref{eq:estimator_double} with fixed $k>-2a$ and fixed $l>-2\at$ satisfies
\begin{align}\label{eq:variance:double}
\Var\bigl(\That_f^{(k,l)}\bigr)
=\bigO(m^{-1}).
\end{align}
\end{theorem}

Combining Theorems~\ref{thm:bias_class1_fixed_kl_double} and Theorem~\ref{thm:variance_rate:double}, we obtain the convergence rate in MSE and conclude the $L_2$-consistency of the estimator.

\begin{corollary}[Convergence rate]\label{cor:double_mse}
Under the same assumptions in Theorem~\ref{thm:bias_class1_fixed_kl_double}, then the estimator~\eqref{eq:estimator_double} with fixed $k>-2a$ and fixed $l>-2\at$ satisfies
\begin{align}
\label{eq:mse_rate_classP_double}
&\E\bigl[\bigl(\That^{(k,l)}_f -T_f(p,q)\bigr)^2\bigr]\nonumber\\
&=\bigOtilde
\bigl(
m^{-2\lambda(\sigma_p,a,k)}
+n^{-2\lambda(\sigma_q,\at,l)}
+m^{-1}
\bigr)
\end{align}
and thus is $L_2$-consistent.
\end{corollary}

\begin{remark}
Similar to the single-density case, if $d\ge 2$, the bias bound dominates the variance bound. 
\end{remark}

\begin{example}[KL divergence; Example~\ref{ex:kl_divergence} contd.]\label{ex:estimating_kl_divergence}
For estimating KL divergence, recall that $|\phi_{kl}(u,v)|\lesssim \psi_{-\eps,\eps}(u)\psi_{-\eps,\eps}(v)$ for any arbitrarily small $\eps>0.$ 
It can be shown, using Theorems~\ref{thm:bias_class1_fixed_kl_double} and \ref{thm:variance_rate:double}, that for estimating the (forward) KL or reverse KL divergences between any two densities $p$ and $q$ such that $\P\ll\Q$, each of which is either the uniform distribution, or one of the truncated Gaussian, Cauchy, Laplace, or exponential distributions, we obtain a bias bound of $\bigOtilde(m^{-1/d})$ and a variance bound of $\bigO(m^{-1}),$ and therefore, the MSE rate of $\bigOtilde(m^{-2/d}+n^{-2/d}+m^{-1})$ as established in Corollary~\ref{cor:double_mse}. 
\end{example}

\begin{example}[$\a$-divergence; Example~\ref{ex:polynomial_functional} contd.]\label{ex:estimating_renyi_divergence}
For estimating the $\a$-divergence ($\a>0$), recall that $|\phi_{kl}(u,v)|\lesssim \psi_{1-\a,1-\a}(u)\psi_{\a-1,\a-1}(v)$ for any $k,l\in\Natural$ such that $k>\a-1$ and $l>1-\a$.
Hence, if $p$ satisfies \ref{cond:bounded_above}, \ref{cond:bounded_below}, \ref{cond:bounded_support}, \ref{cond:regular_support}, \ref{cond:smoothness_int}, and \ref{cond:boundary}, and $q$ satisfies \hyperref[{cond:bounded_above}]{\textbf{(U$_q$)}}, \hyperref[{cond:bounded_support}]{{\textbf{(L1$_q$)}}}, \hyperref[{cond:bounded_below}]{{\textbf{(L2$_q$)}}}, \hyperref[{cond:regular_support}]{{\textbf{(L3$_q$)}}}, \hyperref[{cond:smoothness_int}]{{\textbf{(S$_q$)}}}, and \hyperref[{cond:boundary}]{{\textbf{(B$_q$)}}}, then the MSE of the estimator~\eqref{eq:estimator_double} with $k>2(\a-1)$ and $l>2(1-\a)$ is bounded as \eqref{eq:mse_rate_classP_double} with the bias rate exponents
\begin{align*}
\lambda(\sigma_p,a,k)=
    \begin{cases}
    \frac{1}{d}(\sigma_p\wedge 1) &\text{if }\a<2,\\
    \frac{1}{d}(\sigma_p\wedge\frac{k+1-\a}{k-1}) &\text{if }2\le \a<2+\frac{\sigma_p}{d},\\
    \frac{1}{d}(\sigma_p\wedge 1)(\frac{k+1-\a}{k-1}) &\text{if }\a \ge 2+\frac{\sigma_p}{d}.
    \end{cases}
\end{align*}
and
\begin{align*}
\lambda(\sigma_q,\at,l)=
    \frac{1}{d}(\sigma_q\wedge 1).
\end{align*}
This result also holds for the logarithmic $\a$-divergence.
\end{example}

\subsection{\secondrevision{Le Cam distance and Jensen--Shannon divergence:} Performance guarantee with truncation}\label{sec:convergence_rate:truncated}
\secondrevision{The statements in the previous section do not apply to the estimators for Le Cam distance (Example~\ref{ex:asymp_nn_classfication}) and Jensen--Shannon divergence (Example~\ref{ex:JSD}). 
The difficulty arises from the fact that the estimator function $\phi_{kl}$ for these divergences have lower-polynomial-tail exponents $(a,\at)=(-k+1,-l+1)$ which become smaller with larger $k$ and $l$. Therefore, while the bias guarantees (Theorems~\ref{thm:vanishing_bias:double}~and~\ref{thm:bias_class1_fixed_kl_double}) are still applicable, we cannot control the variance of the estimator using Theorems~\ref{thm:vanishing_variance:double}~or~\ref{thm:variance_rate:double}, as $(a,\at)=(-k+1,-l+1)$ does not meet the requirements \{$k>-4a$, $l>-4\at$\} or \{$k>-2a$, $l>-2\at$\}.}

\secondrevision{To handle the variance of the estimator for these exceptional cases, we consider a truncated version of the estimator~\eqref{eq:estimator_double}.}
For functionals of two densities, we define the truncated estimator as 
\begin{align}\label{eq:estimator_double:truncated}
&\Tt_f^{(k,l)}(\Xb_{1:m},\Yb_{1:n})\nonumber\\
&\defeq
\frac{1}{m}\sum_{i=1}^m
\phibar_{kl}(U_{km}(\Xb_i),V_{ln}(\Xb_i);\tau_m,\nu_m,\taut_n,\nut_n),
\end{align}
where we define the truncated estimator function
\[
\phibar_{kl}(u,v;\tau,\nu,\taut,\nut)\defeq \phi_{kl}(u,v)\ones_{(\tau,\nu)}(u)\ones_{(\taut,\nut)}(v)
\] 
and the \emph{truncation points} $\tau_m,\nu_m,\taut_n,\nut_n$ are hyperparameters such that $0\le \tau_m\le 1\le \nu_m\le \infty$ and $0\le \taut_n\le 1\le \nut_n\le \infty$.
\secondrevision{As noted earlier, we do not require the upper-truncation points in contrast to Section~\ref{sec:convergence_rate:truncated:unbdd_support} and thus only consider a \emph{lower-truncated estimator} with $\nu_m=\infty$ and $\nut_n=\infty$ in this section.}

We can first establish the consistency of the lower-truncated estimator.
\begin{proposition}[Consistency]
\label{prop:consistency_double:truncated}
For a target functional $T_f(\cdot,\cdot)$, if the estimator function $\phi_{kl}(u,v)$ is continuous and the underlying densities $p$ and $q$ satisfy
\hyperref[{cond:low}]{\textbf{(U$_{pp}$; $k,a$)}},
\hyperref[{cond:high:milder}]{\textbf{(L$_{pp}^\prime$; $\xi^2,b$)}},
\hyperref[{cond:low}]{\textbf{(U$_{pq}$; $l,\at$)}}, and
\hyperref[{cond:high:milder}]{\textbf{(L$_{pq}^\prime$; $\xi^2,\bt$)}}
with $\xi(t)=t^2$,
then the lower-truncated estimator~\eqref{eq:estimator_double:truncated} with fixed $k>-a$ and $l>-\at$ and with \secondrevision{lower-truncation points such that $\tau_m^{(k+4a)\wedge 0}\taut_n^{(l+4\at)\wedge 0}=o(m)$ is $L_2$-consistent.}
\end{proposition}

We can also establish convergence rate of the truncated estimator~\ref{eq:estimator_double:truncated} for functionals of two densities.
Define a lower truncation point function as
\begin{align}
&\tau(m,\sigma,a,k)\label{eq:alpha_classP_w_sigma}\\
&= 
\begin{cases}
\Theta\bigl(m^{-\frac{\sigma\wedge 1}{d(k-1)}}\bigr)
	& \text{if }a\le -\frac{\sigma}{d}-1,\\
\Theta\bigl(m^{-\frac{1}{d(k-1)}}\bigr)
	&\text{if }-\frac{\sigma}{d}-1 < a \le -1,\\
O\bigl(m^{-\frac{1}{d(a+1)}}\bigr)
	&\text{if }a > -1.
\end{cases}\nonumber
\end{align}

\begin{proposition}[\secondrevision{Convergence rate}]\label{prop:double_mse:truncated}
For a target functional $T_f(\cdot,\cdot)$, if the underlying density $p$ satisfies the conditions \ref{cond:bounded_above}, \ref{cond:bounded_below}, \ref{cond:smoothness_int}, and \ref{cond:boundary}, and $q$ satisfies the conditions \hyperref[{cond:bounded_above}]{\textbf{(U$_q$)}}, 
\hyperref[{cond:bounded_below}]{{\textbf{(L1$_q$)}}}, 
\hyperref[{cond:smoothness_int}]{{\textbf{(S$_q$)}}}, and
\hyperref[{cond:boundary}]{{\textbf{(B$_q$)}}}, the truncated estimator~\eqref{eq:estimator_double:truncated} with fixed $k>-a$ and $l>-\at$ satisfies
\begin{align*}
&\E\bigl[\bigl(\Tt^{(k,l)}_f
    -T_f(p,q)\bigr)^2\bigr]\\
  &=\bigOtilde
  \bigl(
    m^{-2\lambda(\sigma_p,a,k)}
    +n^{-2\lambda(\sigma_q,\at,l)}
    +m^{-1}\tau_m^{(2a+k)\wedge 0}\taut_n^{(2\at+l)\wedge 0}
    \bigr),
\end{align*}
as $m,n\to\infty$, and thus is $L_2$-consistent.
\end{proposition}

\begin{example}[Le Cam distance; Example~\ref{ex:asymp_nn_classfication} contd.]
\label{ex:estimating_asymp_nn_classfication}
For estimating $T_f(p,q)$ with $f(p,q)=q/(p+q)$, recall that $|\phi_{kl}(u,v)|\lesssim \psi_{-k+1,l-1}(u)\psi_{-l+1,k-1}(v)$ for any $k\ge 1$ and $l\ge 1$. 
\secondrevision{For densities $p$ and $q$ satisfying conditions in Proposition~\ref{prop:consistency_double:truncated}, the lower-truncated estimator~\eqref{eq:estimator_double:truncated} for Le Cam distance is $L_2$-consistent. 
In particular, the estimator with $k=l=1$ is consistent even without lower truncation, since $\tau_m^{(k+4a)\wedge 0}\taut_n^{(l+4\at)\wedge 0}=\tau_m^0\taut_n^0=0$ with $\tau_m=\taut_n=0$ and $k=l=1$. 
If the underlying densities $p$ and $q$ satisfy the conditions in Proposition~\ref{prop:double_mse:truncated},
then the lower-truncated estimator with fixed $k\ge1$ and $l\ge1$ and truncation points
$\tau_m=\tau(m,\sigma_p,-k+1,k)$, and $\taut_n=\tau(n,\sigma_q,-l+1,l)$
satisfies}
\begin{align*}
&\E\bigl[\bigl(\That^{(k,l)}_f
    -T_f(p,q)\bigr)^2\bigr]\numberthis\label{eq:rate_nnce_truncated}\\
  &=\bigOtilde\bigl(
    m^{-2\lambda_k(\sigma_p)}
    +n^{-2\lambda_l(\sigma_q)}
    +m^{-1}\tau_m^{(-k+2)\wedge 0}\taut_n^{(-l+2)\wedge 0}\bigr),
\end{align*}
as $m,n\to\infty$, where $\lambda_p=\lambda_k(\sigma_p)$ and $\lambda_q=\lambda_l(\sigma_q)$, where
\[
\lambda_k(\sigma)\defeq \lambda(\sigma,-k+1,k)
=\begin{cases}
\frac{1}{d}(\sigma\wedge 1) & \text{if }k=1,\\
\frac{1}{d}(\sigma\wedge\frac{1}{k-1}) & \text{if }2\le k<2+\frac{\sigma}{d},\\
\frac{1}{d}\frac{\sigma\wedge 1}{k-1}, &\text{if }k>2+\frac{\sigma}{d}.
\end{cases}
\]
\secondrevision{Based on this rate-exponent expression and the additional factor of $\tau_m^{(2a+k)\wedge 0}\tilde{\tau}_n^{(2\at+l)\wedge 0}$ in the variance rate which only worsens the rate with larger $k$ and $l4$, one would expect that the convergence becomes only slower as $k$ and/or $l$ become large, and thus, the fastest rate achieved is
$\bigOtilde(m^{-\frac{2}{d}(\sigma_p\wedge 1)}+n^{-\frac{2}{d}(\sigma_q\wedge 1)}+m^{-1})$, when $k=1$ and $l=1$ with lower truncation points $\tau_m=0$ and $\taut_n=\Th(n^{-\frac{1}{d}})$.
This is in contrast with Remark~\ref{rem:rate_increasing_ak}, where we observed faster convergence with larger values of $k$ when $a$ does not decrease in $k$.
We note that the experiments with synthetic data in Section~\ref{sec:exp} show that the estimator performs well even for large values of $k$ and $l$, suggesting that the detrimental effect of the lower tail exponents might be removed with a tighter analysis.}
\end{example}

\begin{example}[Jensen--Shannon divergence; Example~\ref{ex:JSD} contd.]\label{ex:estimating_JSD}
For estimating Jensen--Shannon divergence, recall that $|\phi_{kl}(u,v)| \lesssim \psi_{-k+1,l-1}(u)\psi_{-l+1,k-1}(v)$ for any $k\ge 1$ and $l\ge 2$. 
\secondrevision{For densities $p$ and $q$ satisfying conditions in Proposition~\ref{prop:consistency_double:truncated}, the lower-truncated estimator~\eqref{eq:estimator_double:truncated} for Jensen--Shannon divergence is $L_2$-consistent. Also, we do not require the lower-truncation $\tau_m$ for $k=1$, by the same argument in the previous example.}
If the underlying densities $p$ and $q$ satisfy the conditions in Proposition~\ref{prop:double_mse:truncated} and \secondrevision{additionally} $\Q\ll\P$,
then the estimator~\eqref{eq:estimator_double} with fixed $k\ge1$ and $l\ge2$ and the same truncation points in Example~\ref{ex:estimating_asymp_nn_classfication}
satisfies \eqref{eq:rate_nnce_truncated}.
\secondrevision{The established rate seems to get only slower as $k$ and/or $l$ become large, and thus achieves its fastest rate
$\bigOtilde(m^{-\frac{2}{d}(\sigma_p\wedge 1)}+n^{-\frac{2}{d}(\sigma_q\wedge 1)}+m^{-1})$ when $k=1$ and $l=2$ with lower truncation points $\tau_m=0$ and $\taut_n=\Th(n^{-\frac{1}{d}})$.
Note, however, this conclusion might not hold in practice; see Example~\ref{ex:estimating_asymp_nn_classfication}.}
\end{example}

\section{Adaptive choices of \texorpdfstring{$k$}{k} and \texorpdfstring{$l$}{l}}
\label{sec:adaptive_choice_k}
In Section~\ref{sec:single}, we established the convergence rate of the proposed estimator~\eqref{eq:estimator_single} for fixed $k$.
Since $\E[\phi_k(U_{k\infty}(\xb))]=f(p(\xb))$ for each valid $k\in\Natural$ by design, we can choose any valid
$k$ without violating the asymptotic unbiasedness.
\secondrevision{In Remark~\ref{rem:rate_increasing_ak},} we observed that a larger \emph{fixed} $k$ in general leads to a larger rate exponent in~\eqref{eq:lambda_classP_w_sigma},
and thus, a faster convergence rate. 
This prompts the 
question of whether increasing $k \to \infty$ along with $m$ improves the convergence rate upon fixed $k$.
The following proposition answers this in the affirmative. 
The proof is deferred to Appendix~\ref{supp:sec:proofs:adaptive:single}.

\begin{proposition}[Convergence rate and $L_2$-consistency with increasing $k$]\label{prop:bias_class1_increasing_k_single}
For a target functional $T_f(\cdot)$, 
if the underlying density $p$ satisfies \ref{cond:bounded_above}, \ref{cond:bounded_below}, \ref{cond:bounded_support}, \ref{cond:regular_support}, \ref{cond:smoothness_int}, and \ref{cond:boundary},
then the estimator~\eqref{eq:estimator_single} with $k=\Th((\ln m)^{1.1})$
satisfies
\begin{align}\label{eq:bias_rate_increasing_k_single}
\bigl|\E\bigl[\That_f^{(k)}\bigr]-T_f(p)\bigr|
&=\bigOtilde
\bigl(m^{-\frac{\sigma_p\wedge 1}{d}}\bigr)
\end{align}
as $m\to\infty$.
Furthermore, the estimator~\eqref{eq:estimator_single} satisfies
\begin{align}\label{eq:mse_rate_increasing_k}
\E\bigl[\bigl(\That_f^{(k)}-T_f(p)\bigr)^2\bigr]
=\bigOtilde\bigl(m^{-\frac{2(\sigma_p\wedge 1)}{d}}+m^{-1}\bigr)
\end{align}
and thus is $L_2$-consistent.
\end{proposition}

\begin{remark}
As  expected heuristically, the bias rate exponent $(\sigma_p\wedge 1)/d$ in \eqref{eq:bias_rate_increasing_k_single} equals  the limit of the finite-$k$ rate exponent in \eqref{eq:lambda_classP_w_sigma} as $k\to\infty$.
\end{remark}

\begin{remark}
There is no consensus on the optimal choice of $k$ for functional estimation in the literature.
For example, \citet{Singh--Poczos2016} analyzed $k=\bigO(1)$, whereas \citet{Berrett--Samworth--Yuan2019} suggested \revision{$k=\bigO((\ln m)^5)$ for asymptotic efficiency of the estimator, a slightly faster choice than the previous theorem}, for differential entropy. 
\citet{Perez-Cruz2009} discussed some relevant empirical results on the choice of $k$. 
\end{remark}

\begin{remark}
\secondrevision{While our main focus in this paper is to establish consistency and convergence rates for the proposed estimators with fixed $k$ (and $l$), we point out that a tighter analysis on the dependence on $k$ may lead to a better asymptotic convergence rate. Note that the analysis of Kozachenko--Leonenko estimator by \citet{Berrett--Samworth--Yuan2019} allows polynomial growth of $k$ in the sample size. 
The loose dependence on $k$ in our analysis can be traced back to Lemma~\ref{supp:lem:GOVLemma2_pdf_gap_bound}, which quantifies the gap between densities of the normalized volume of $k$-NN ball $U_{km}(\xb)$ and its limiting Poisson random variable $U_{k\infty}(\xb)$.
To tighten the bound, one needs to sharpen Lemma~\ref{supp:lem:GOVLemma5_poisson} on the speed of convergence of a Poisson binomial random variable to a Poisson random variable.}
\end{remark}

\begin{example}[Differential entropy; Example~\ref{ex:estimating_differential_entropy} contd.]
\label{ex:estimating_differential_entropy_varying_k}
Applying Proposition~\ref{prop:bias_class1_increasing_k_single} on differential entropy with $k = \Th((\ln m)^{1.05})$, we obtain the MSE rate~\eqref{eq:mse_rate_increasing_k}.
This rate is the same as the fixed-$k$ case in Example~\ref{ex:estimating_differential_entropy}.
\end{example}

\begin{example}[$\a$-entropy; Example~\ref{ex:estimating_renyi_entropy} contd.]
\label{ex:estimating_renyi_entropy_varying_k}
Applying Proposition~\ref{prop:bias_class1_increasing_k_single} on $\a$-entropy with $k = \Th((\ln m)^{1.05})$, we obtain the bias rate exponent $(\sigma_p\wedge 1)/d$, which is greater than or equal to that in Example~\ref{ex:estimating_renyi_entropy} with $k$ fixed.
\end{example}

Similarly to the single-density case, we can establish the convergence rate when $k$ and $l$ vary polylogarithmically with $m$ and $n$, provided that $m$ and $n$ grow to infinity in the same speed, \ie $m\asymp n$.
The following proposition can be proved by extending the proof of Proposition~\ref{prop:bias_class1_increasing_k_single} to the double-density case as in the proofs of Theorems~\ref{thm:bias_class1_fixed_kl_double} and \ref{thm:variance_rate:double}, and thus is omitted.

\begin{proposition}[Convergence rate and $L_2$-consistency with increasing $k$ and $l$]\label{prop:bias_class1_increasing_kl_double}
For a target functional $T_f(\cdot,\cdot)$, if the underlying densities $p$ and $q$ satisfy the conditions \ref{cond:bounded_above}, \ref{cond:bounded_below}, \ref{cond:bounded_support}, \ref{cond:regular_support}, \ref{cond:smoothness_int}, \ref{cond:boundary},
\hyperref[{cond:bounded_above}]{\textbf{(U$_q$)}}, 
\hyperref[{cond:bounded_below}]{{\textbf{(L1$_q$)}}}, 
\hyperref[{cond:bounded_support}]{{\textbf{(L2$_q$)}}}, 
\hyperref[{cond:regular_support}]{{\textbf{(L3$_q$)}}}, 
\hyperref[{cond:smoothness_int}]{{\textbf{(S$_q$)}}}, and
\hyperref[{cond:boundary}]{{\textbf{(B$_q$)}}},
then the estimator~\eqref{eq:estimator_double} with $k=\Theta((\ln m)^{1.1})$ and $l=\Theta((\ln n)^{1.1})$ satisfies
\begin{align*}
\bigl|\E[\That^{(k,l)}_f]-T_f(p,q)\bigr|
  &=\bigOtilde(m^{-\frac{\sigma_p\wedge 1}{d}}+n^{-\frac{\sigma_q\wedge 1}{d}}),
\end{align*}
as $m,n\to\infty$ with $m\asymp n$.
Furthermore, the estimator~\eqref{eq:estimator_double} satisfies
\begin{align*}
&\E\bigl[(\That^{(k,l)}_f
    -T_f(p,q))^2\bigr]\\
&=\bigOtilde(m^{-\frac{2(\sigma_p\wedge 1)}{d}}+n^{-\frac{2(\sigma_q\wedge 1)}{d}}+m^{-1}),
\numberthis\label{eq:mse_rate_increasing_kl}
\end{align*}
and thus is $L_2$-consistent, provided that $m\asymp n$.
\end{proposition}

\begin{remark}
For $d\ge 2$, if $k$ and $l$ increase as in Proposition~\ref{prop:bias_class1_increasing_kl_double}, the bias bound always dominates the variance bound so that the MSE is bounded as $\bigO(m^{-1})$.
For $d=1$, the variance bound may dominate the bias bound depending on $\sigma_p,\sigma_q$, $d$, and/or the choices of $k$ and $l$.
\end{remark}

\begin{example}[KL divergence; Example~\ref{ex:estimating_kl_divergence} contd.]
\label{ex:estimating_kl_divergence_varying_kl}
Letting $k$ and $l$ increase as $k = \Theta((\ln m)^{1.05}))$ and $l = \Theta((\ln n)^{1.05})),$ we obtain the MSE rate~\eqref{eq:mse_rate_increasing_kl} for estimating KL divergence.
As a complementary asymptotic result, \citet{Wang--Kulkrani--Verdu2009} showed that the $(k,l)$-NN KL divergence estimator with $k=k_m$ and $l=l_n$ such that $k_m/m\to 0$ and $k_m/(\ln m)\to\infty$ as $m\to\infty$ and $l_n/n\to 0$ and $l_n/(\ln n)\to\infty$ as $n\to\infty$ converges to the true KL divergence almost surely for uniformly continuous densities bounded from below on their support.
\end{example}

\begin{example}[$\a$-divergence; Example~\ref{ex:estimating_renyi_divergence} contd.]
\label{ex:estimating_renyi_divergence_varying_kl}
Letting $k$ and $l$ increase as $k = \Theta((\ln m)^{1.05}))$ and $l = \Theta((\ln n)^{1.05})),$
the MSE of our estimator
is bounded as \eqref{eq:mse_rate_increasing_kl}.
\end{example}

\section{Numerical results}\label{sec:exp}
The performance of the proposed estimators~\eqref{eq:estimator_single} and \eqref{eq:estimator_double} for several density functionals were simulated over $500$ runs for sample sizes ranging from $100$ till $25600$.\footnote{The code is available at \url{https://github.com/jongharyu/knn-functional-estimation}.}
For each dimension $d$ from $1$ through $5$, we considered the uniform density $\Unif([0,1]^d)$, the Gaussian density $\Normal(0,I_d)$ restricted to $\|\xb\|\le 3$, and the Gaussian density $\Normal(0,I_d)$ as the density $p$.
For double-density functionals, we considered $\Unif([0,2]^d$),
$\Normal(0,4I_d)$ restricted to $\Bb(0,3)$, and $\Normal(0,4I_d)$ as the density $q$.\footnote{As an exception for the experiment with the Jensen--Shannon divergence estimator, instead of $\Unif([0,1]^d)$ and $\Unif([0,2]^d)$, we used piecewise constant densities $p$ and $q$ supported on $[0,1]^d$, which are defined as follows:
\[
p(\xb)=\begin{cases}
3/2 &\text{if }0\le x_1\le 1/2,\\
1/2 &\text{if }1/2<x\le 1,
\end{cases}
\quad\text{and}\quad
q(\xb)=\begin{cases}
1/2 &\text{if }0\le x_1\le 1/2,\\
3/2 &\text{if }1/2<x\le 1.
\end{cases}
\]}
Note that all the functionals considered in these simulations can be expressed in closed form up to incomplete gamma function, except the exponential entropies, Le Cam distance, and Jensen--Shannon divergences for Gaussian densities. We estimated the latter using Monte Carlo approximation.
Polynomial rates of convergence were observed for all cases, and in each case, the exponent was calculated by ordinary least-squares linear regression between the logarithms of the sample sizes and the MSE.
We considered $k\in\{1,2,3,4,5,10,15\}$ and,
for double-density functional estimators,
$l=k$  for simplicity.

\begin{figure*}[htp]
\centering
\begin{minipage}[b]{0.28\linewidth}
  \centerline{\includegraphics[width=\textwidth]{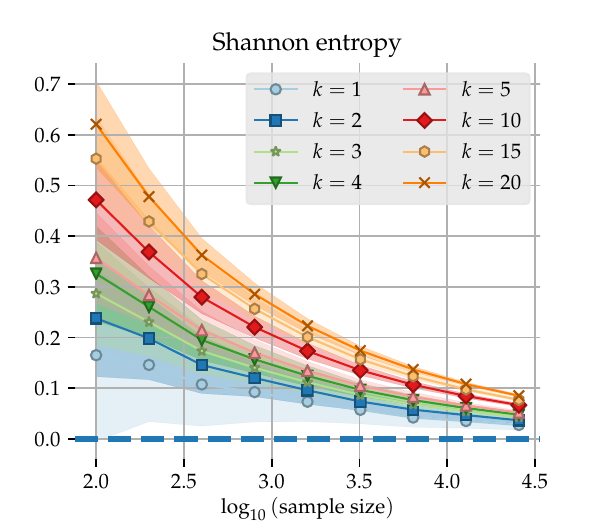}}
\end{minipage}
\qquad
\begin{minipage}[b]{0.28\linewidth}
  \centerline{\includegraphics[width=\textwidth]{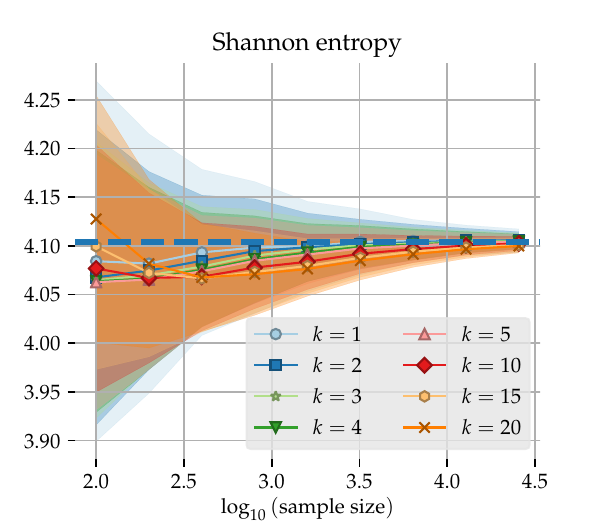}}
\end{minipage}
\qquad
\begin{minipage}[b]{0.28\linewidth}
  \centerline{\includegraphics[width=\textwidth]{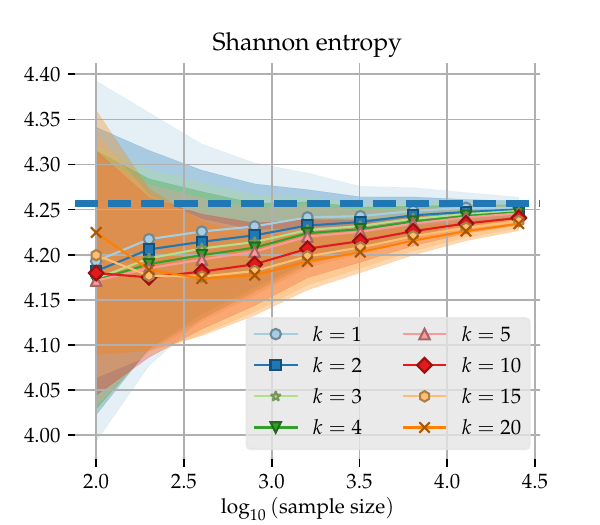}}
\end{minipage}\\
\begin{minipage}[b]{0.28\linewidth}
  \centerline{\includegraphics[width=\textwidth]{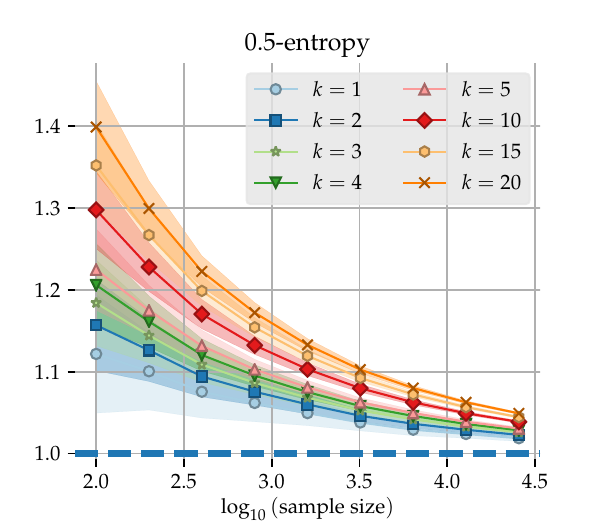}}
\end{minipage}
\qquad
\begin{minipage}[b]{0.28\linewidth}
  \centerline{\includegraphics[width=\textwidth]{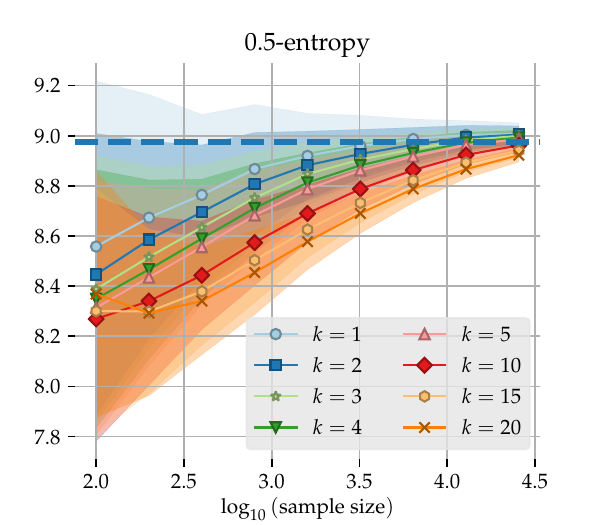}}
\end{minipage}
\qquad
\begin{minipage}[b]{0.28\linewidth}
  \centerline{\includegraphics[width=\textwidth]{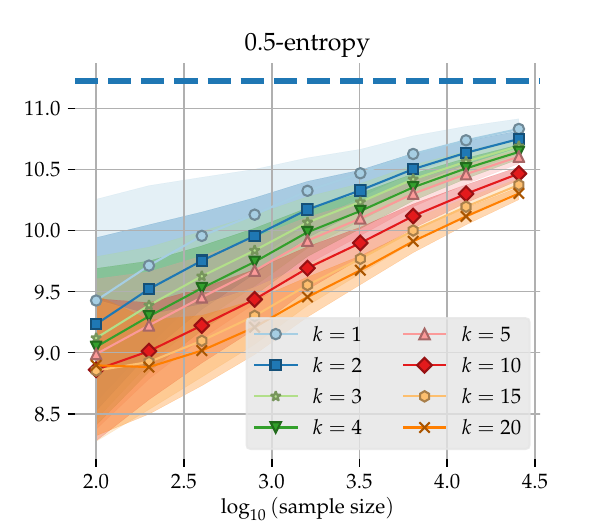}}
\end{minipage}
\begin{minipage}[b]{0.28\linewidth}
\centerline{\includegraphics[width=\textwidth]{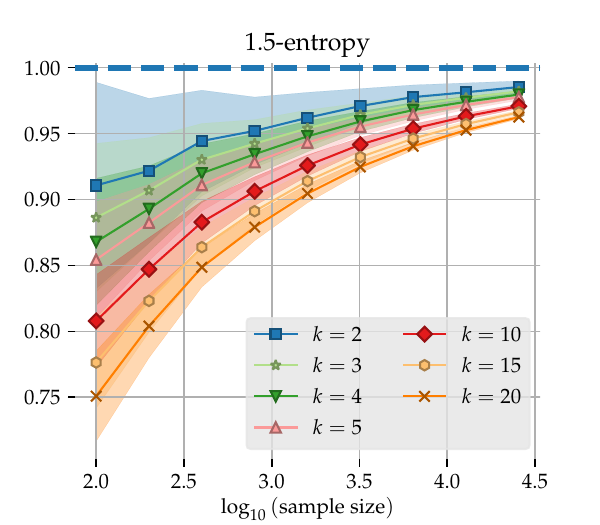}}
\end{minipage}
\qquad
\begin{minipage}[b]{0.28\linewidth}
\centerline{\includegraphics[width=\textwidth]{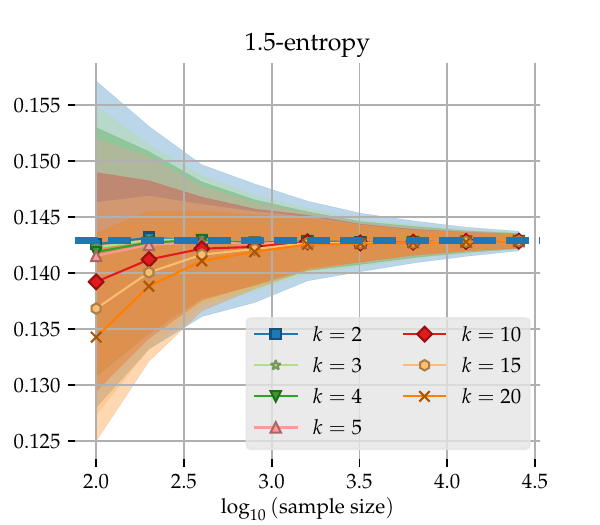}}
\end{minipage}
\qquad
\begin{minipage}[b]{0.28\linewidth}
\centerline{\includegraphics[width=\textwidth]{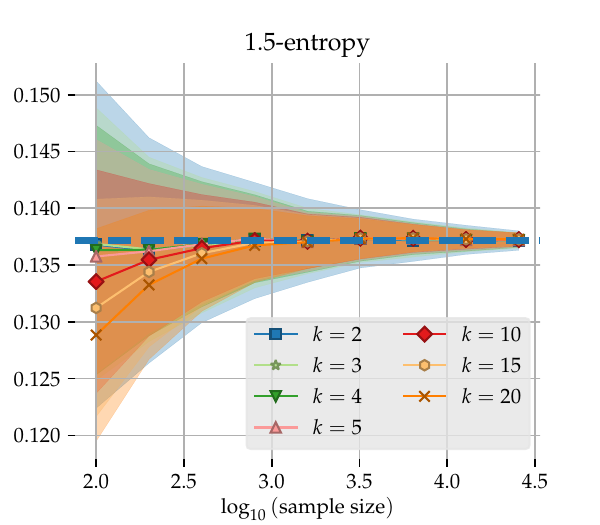}}
\end{minipage}
\begin{minipage}[b]{0.28\linewidth}
\centerline{\includegraphics[width=\textwidth]{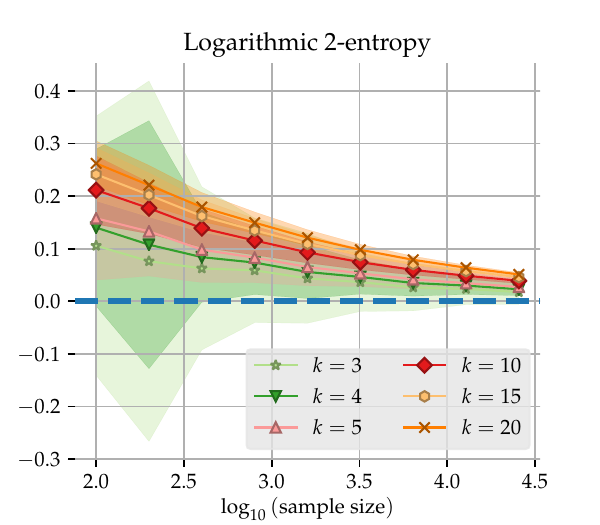}}
\end{minipage}
\qquad
\begin{minipage}[b]{0.28\linewidth}
\centerline{\includegraphics[width=\textwidth]{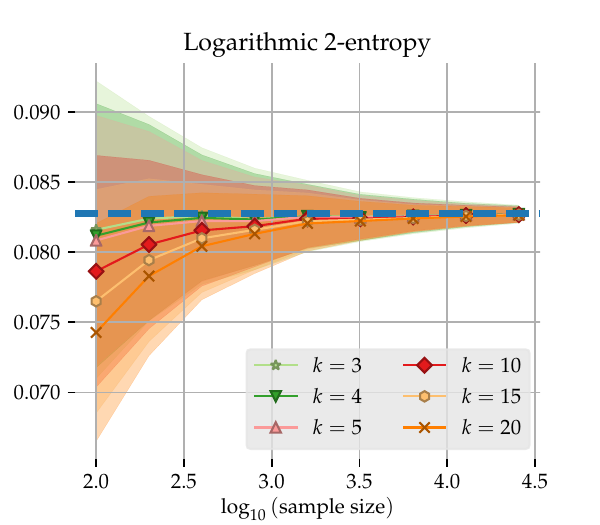}}
\end{minipage}
\qquad
\begin{minipage}[b]{0.28\linewidth}
\centerline{\includegraphics[width=\textwidth]{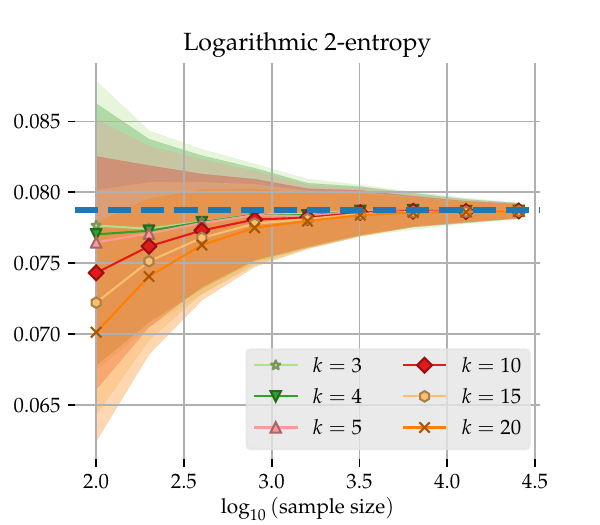}}
\end{minipage}
\begin{minipage}[b]{0.28\linewidth}
\centerline{\includegraphics[width=\textwidth]{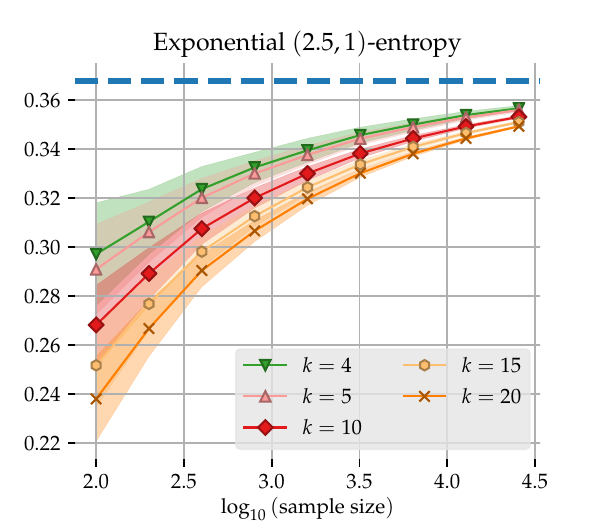}}
\end{minipage}
\qquad
\begin{minipage}[b]{0.28\linewidth}
\centerline{\includegraphics[width=\textwidth]{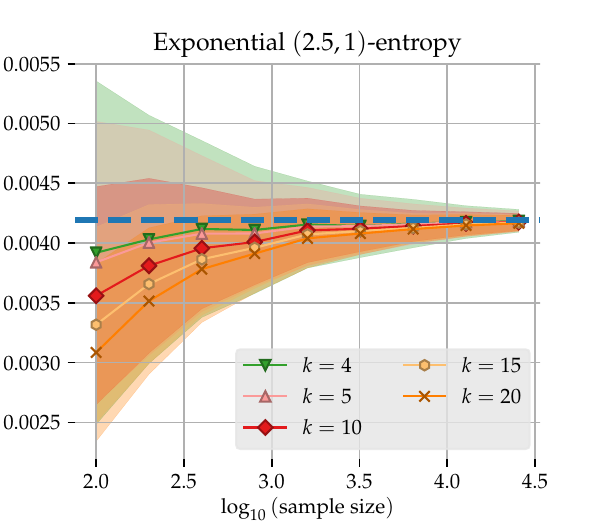}}
\end{minipage}
\qquad
\begin{minipage}[b]{0.28\linewidth}
\centerline{\includegraphics[width=\textwidth]{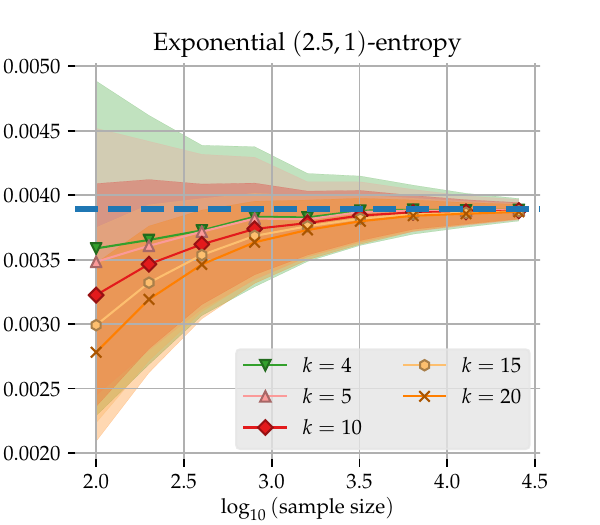}}
\end{minipage}
\vspace{-0.5em}
\caption{Convergence of the single-density functional estimator for differential entropy, $\a$-entropies $\a\in\{0.5,1.5\}$, logarithmic $2$-entropy, and exponential $(2.5,1)$-entropy for 3-dimensional densities.
The first, second, and third columns present simulation results with $\Unif([0,1]^3)$, $\Normal(0,I_3)$ restricted to $\|\xb\|\le 3$, and $\Normal(0,I_3)$, respectively.
The true functional values are indicated as dashed lines and one sample standard deviations of the estimates are indicated as shaded area.}
\label{fig:convergence:single}
\end{figure*}

\begin{figure*}[htp]
\centering
\begin{minipage}[b]{0.28\linewidth}
  \centerline{\includegraphics[width=\textwidth]{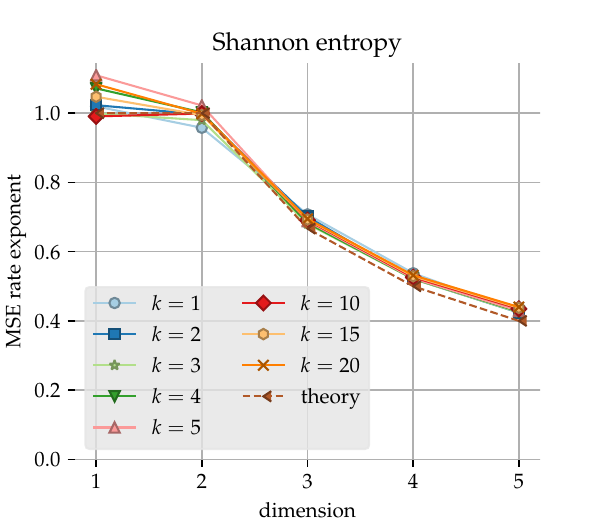}}
\end{minipage}
\qquad
\begin{minipage}[b]{0.28\linewidth}
  \centerline{\includegraphics[width=\textwidth]{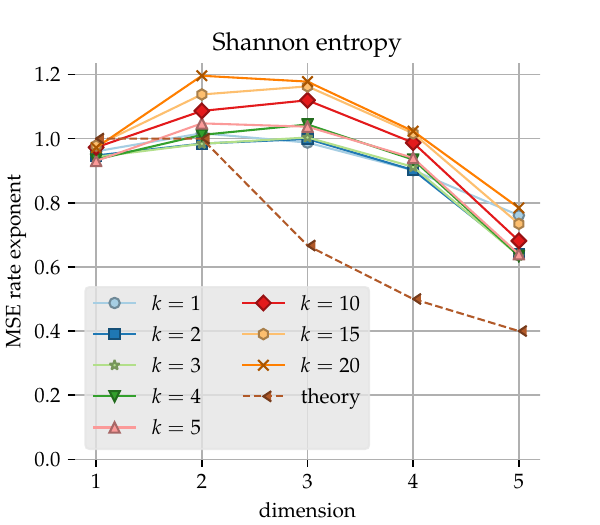}}
\end{minipage}
\qquad
\begin{minipage}[b]{0.28\linewidth}
  \centerline{\includegraphics[width=\textwidth]{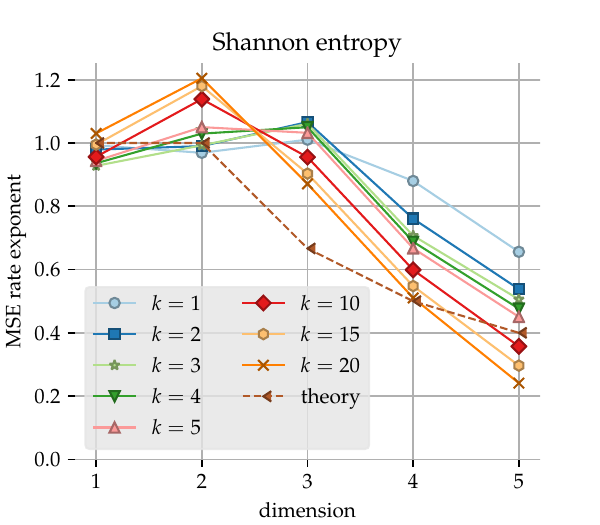}}
\end{minipage}\\
\begin{minipage}[b]{0.28\linewidth}
  \centerline{\includegraphics[width=\textwidth]{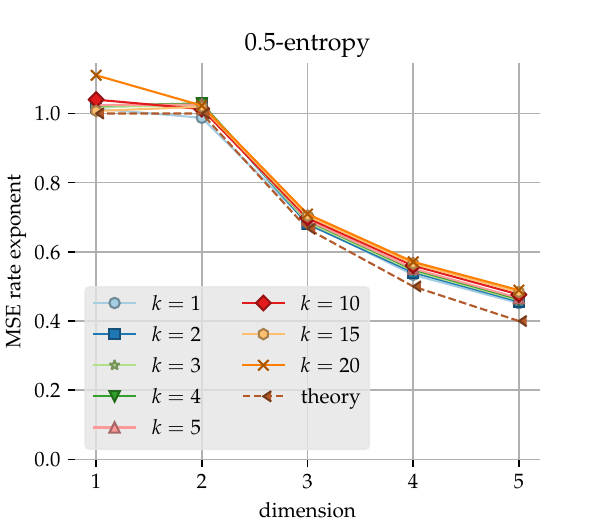}}
\end{minipage}
\qquad
\begin{minipage}[b]{0.28\linewidth}
  \centerline{\includegraphics[width=\textwidth]{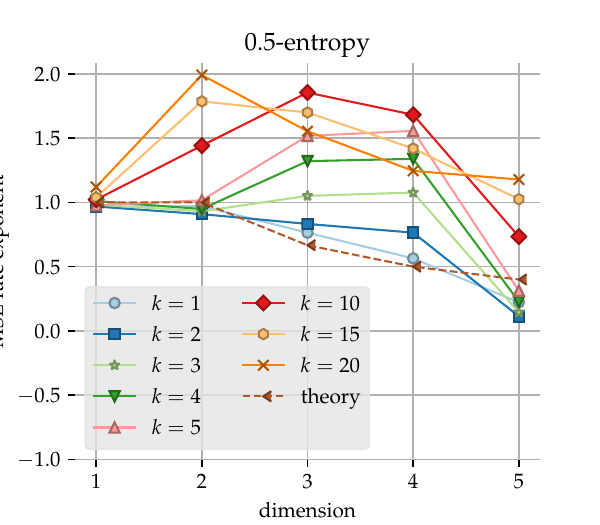}}
\end{minipage}
\qquad
\begin{minipage}[b]{0.28\linewidth}
  \centerline{\includegraphics[width=\textwidth]{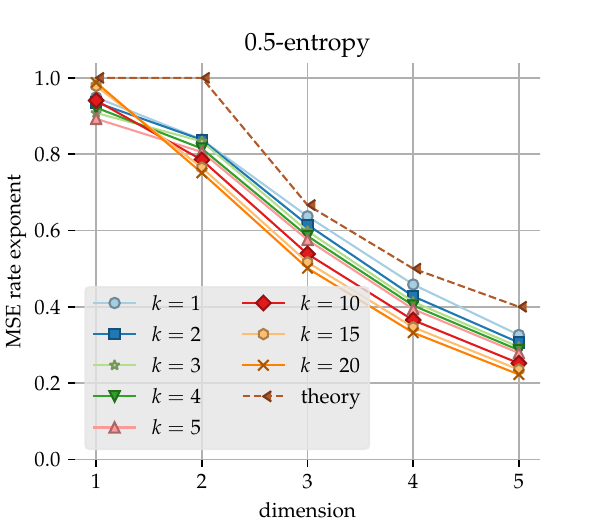}}
\end{minipage}
\begin{minipage}[b]{0.28\linewidth}
\centerline{\includegraphics[width=\textwidth]{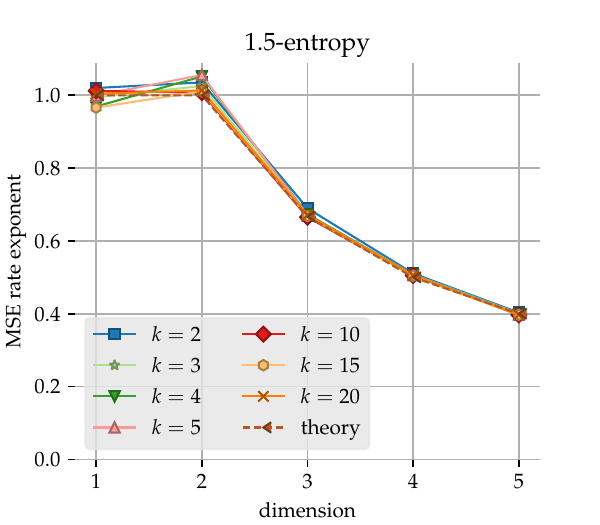}}
\end{minipage}
\qquad
\begin{minipage}[b]{0.28\linewidth}
\centerline{\includegraphics[width=\textwidth]{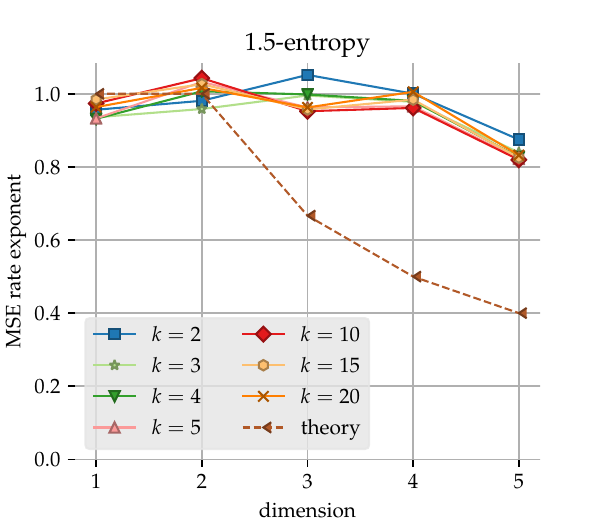}}
\end{minipage}
\qquad
\begin{minipage}[b]{0.28\linewidth}
\centerline{\includegraphics[width=\textwidth]{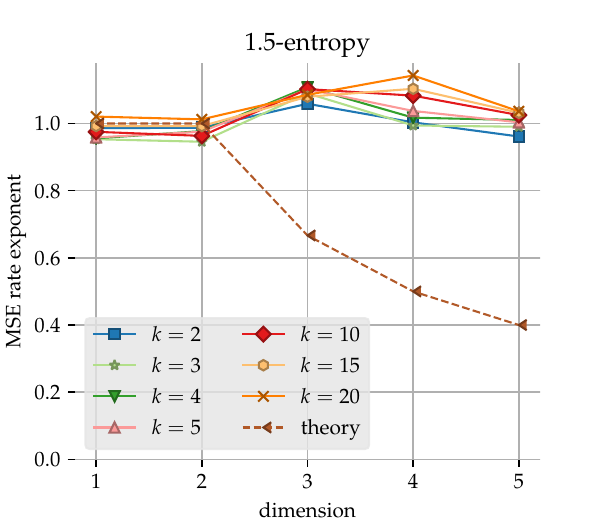}}
\end{minipage}
\begin{minipage}[b]{0.28\linewidth}
\centerline{\includegraphics[width=\textwidth]{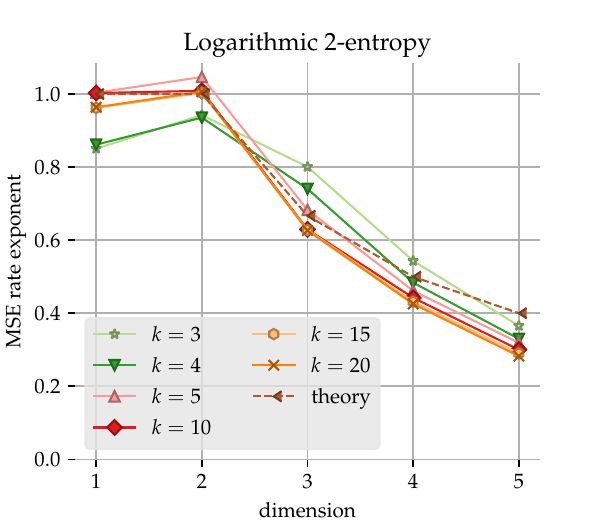}}
\end{minipage}
\qquad
\begin{minipage}[b]{0.28\linewidth}
\centerline{\includegraphics[width=\textwidth]{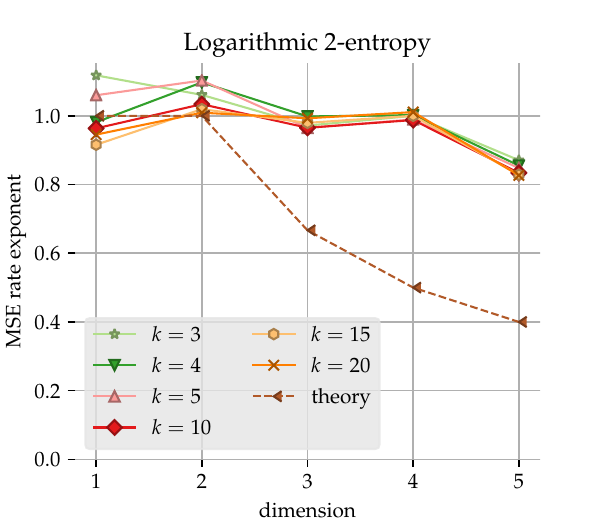}}
\end{minipage}
\qquad
\begin{minipage}[b]{0.28\linewidth}
\centerline{\includegraphics[width=\textwidth]{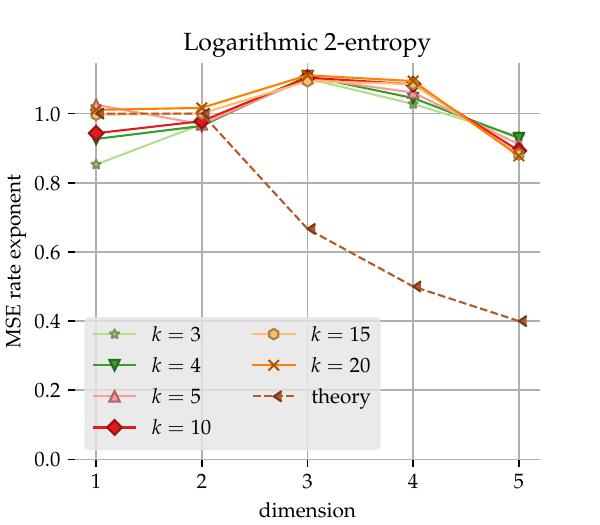}}
\end{minipage}
\begin{minipage}[b]{0.28\linewidth}
\centerline{\includegraphics[width=\textwidth]{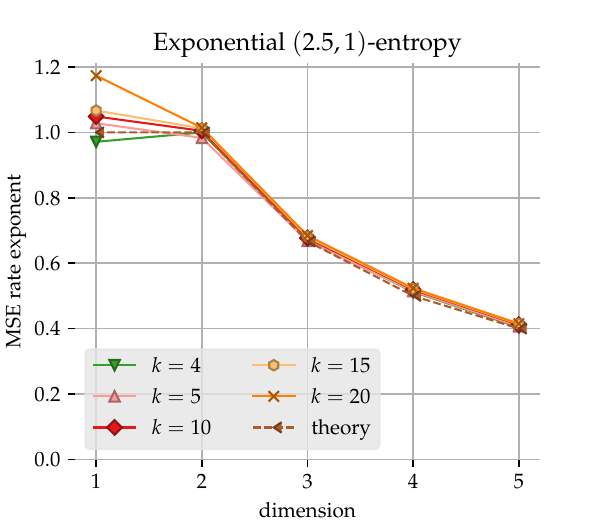}}
\end{minipage}
\qquad
\begin{minipage}[b]{0.28\linewidth}
\centerline{\includegraphics[width=\textwidth]{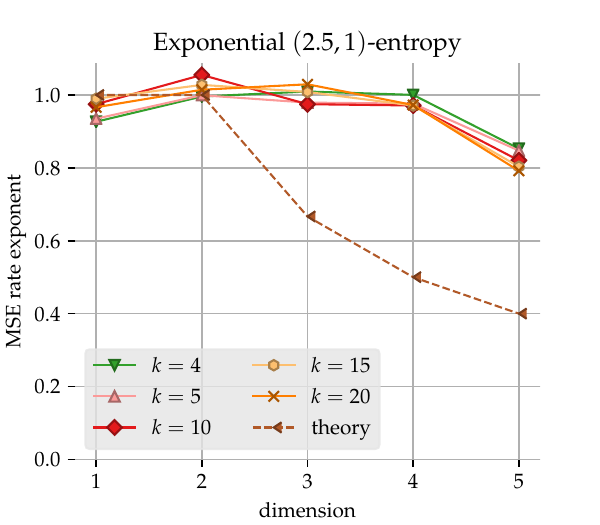}}
\end{minipage}
\qquad
\begin{minipage}[b]{0.28\linewidth}
\centerline{\includegraphics[width=\textwidth]{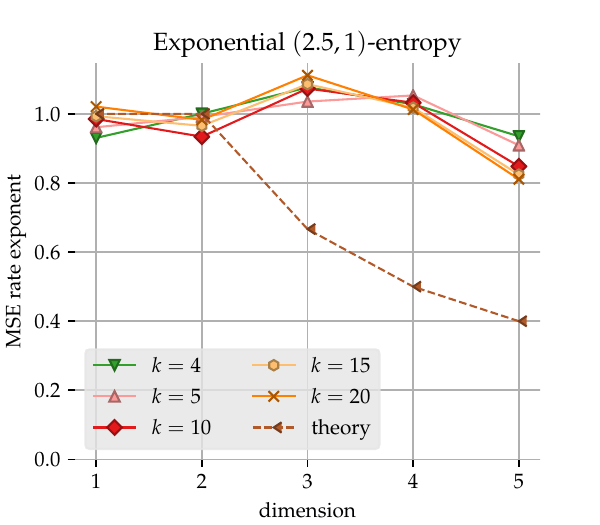}}
\end{minipage}
\vspace{-0.5em}
\caption{Simulated MSE rate exponents of the single-density functional estimator for differential entropy, $\a$-entropies for $\a\in\{0.5,1.5\}$, logarithmic $2$-entropy, and exponential $(2.5,1)$-entropy.
The first, second, and third columns present simulation results with $\Unif([0,1]^d)$, $\Normal(0,I_d)$ restricted to $\|\xb\|\le 3$, and $\Normal(0,I_d)$, respectively, for $d\in\{1,2,3,4,5\}$.}
\label{fig:exponent:single}
\end{figure*}

Figure~\ref{fig:convergence:single} presents the convergence of the estimator for differential entropy,  $\a$-entropies for $\a\in\{0.5,1.5\}$, logarithmic 2-entropy, and exponential $(2.5,1)$-entropy for 3-dimensional densities.
The simulation results show that smaller $k$ yields faster convergence while incurring larger variance, which suggests the use of a moderate size of $k$ in practice.
Figure~\ref{fig:exponent:single} summarizes the empirical exponents of the estimator for each functional and density.
A simple upper bound $(2/d)\wedge 1$ on the theoretical exponents established in Corollary~\ref{cor:single_mse} is also plotted for comparison; see also Examples~\ref{ex:estimating_differential_entropy} and \ref{ex:estimating_renyi_entropy}.
Empirical convergence rates are 
consistently better than theoretical bounds for the truncated densities.

\begin{figure*}[htp]
\centering
\begin{minipage}[b]{0.28\linewidth}
  \centerline{\includegraphics[width=\textwidth]{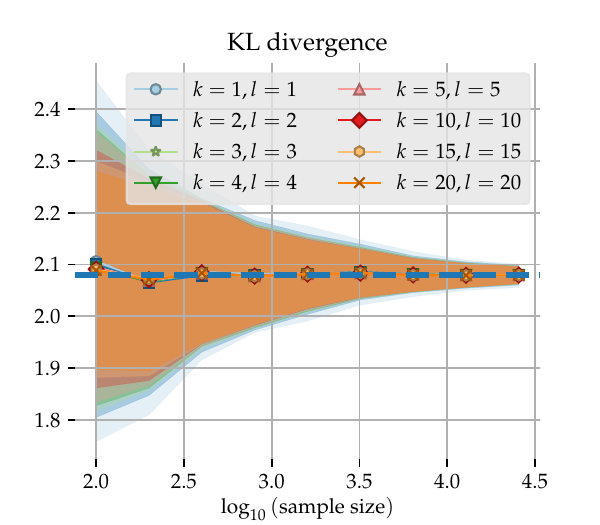}}
\end{minipage}
\qquad
\begin{minipage}[b]{0.28\linewidth}
  \centerline{\includegraphics[width=\textwidth]{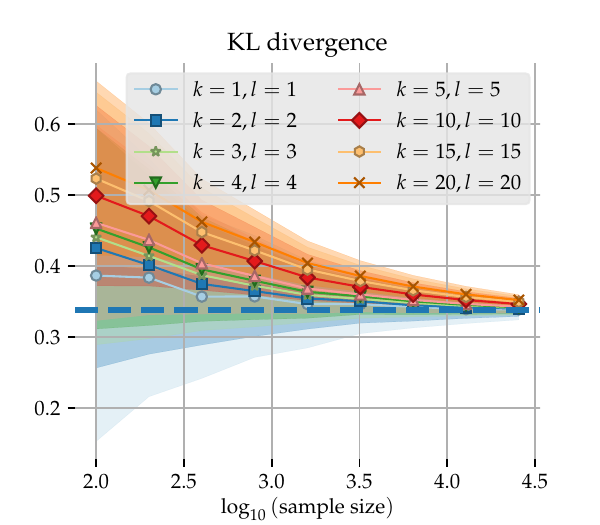}}
\end{minipage}
\qquad
\begin{minipage}[b]{0.28\linewidth}
  \centerline{\includegraphics[width=\textwidth]{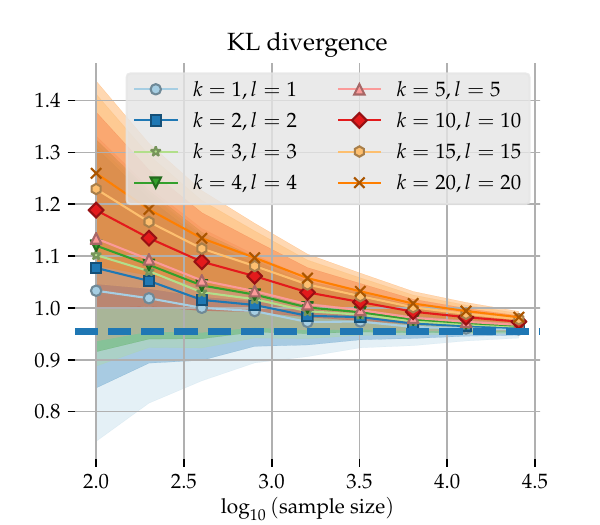}}
\end{minipage}\\
\begin{minipage}[b]{0.28\linewidth}
  \centerline{\includegraphics[width=\textwidth]{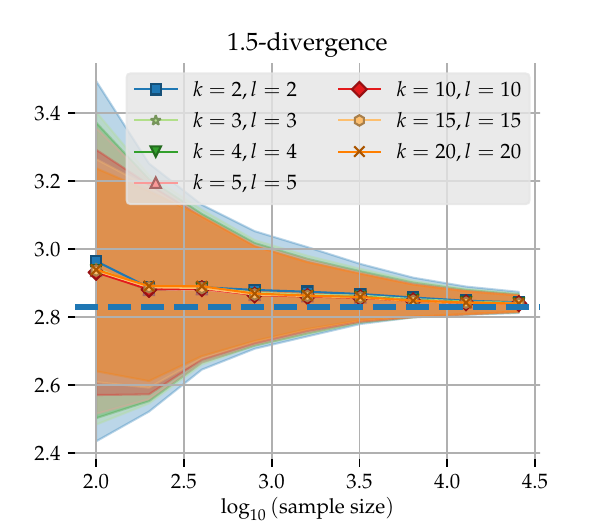}}
\end{minipage}
\qquad
\begin{minipage}[b]{0.28\linewidth}
  \centerline{\includegraphics[width=\textwidth]{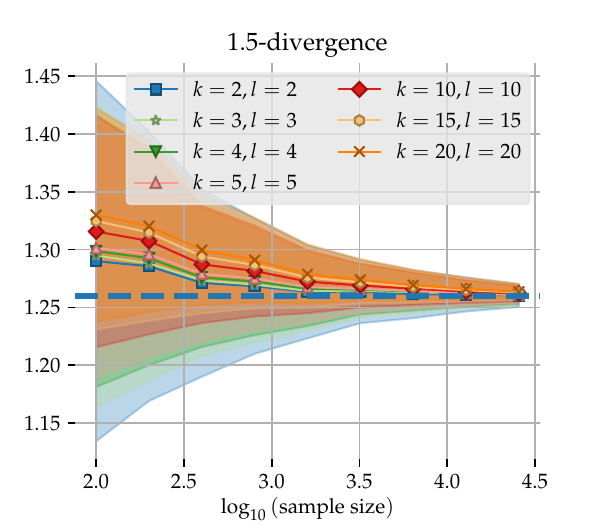}}
\end{minipage}
\qquad
\begin{minipage}[b]{0.28\linewidth}
  \centerline{\includegraphics[width=\textwidth]{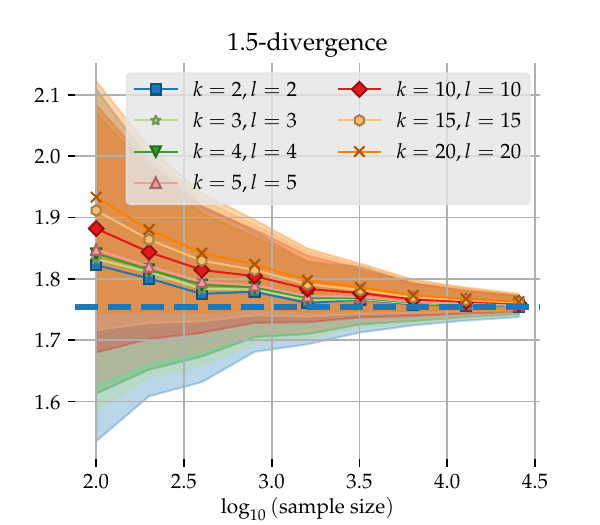}}
\end{minipage}
\begin{minipage}[b]{0.28\linewidth}
\centerline{\includegraphics[width=\textwidth]{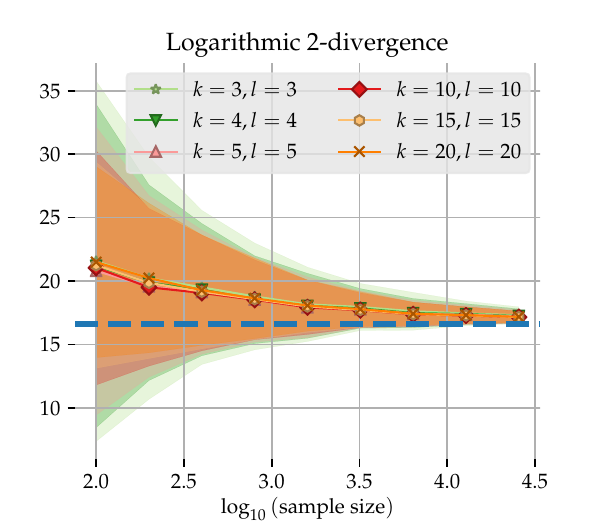}}
\end{minipage}
\qquad
\begin{minipage}[b]{0.28\linewidth}
\centerline{\includegraphics[width=\textwidth]{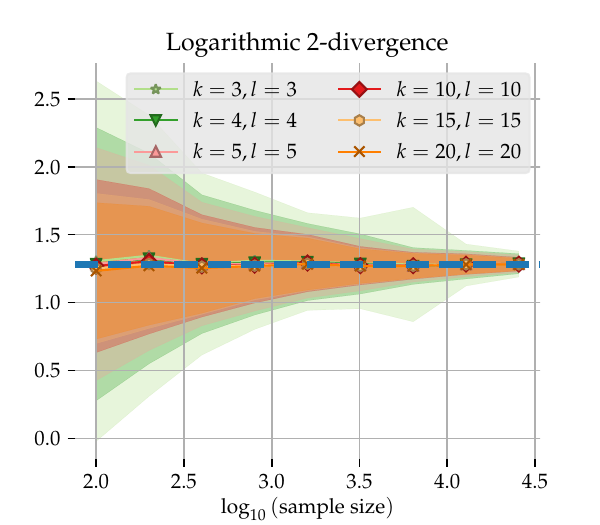}}\end{minipage}
\qquad
\begin{minipage}[b]{0.28\linewidth}
\centerline{\includegraphics[width=\textwidth]{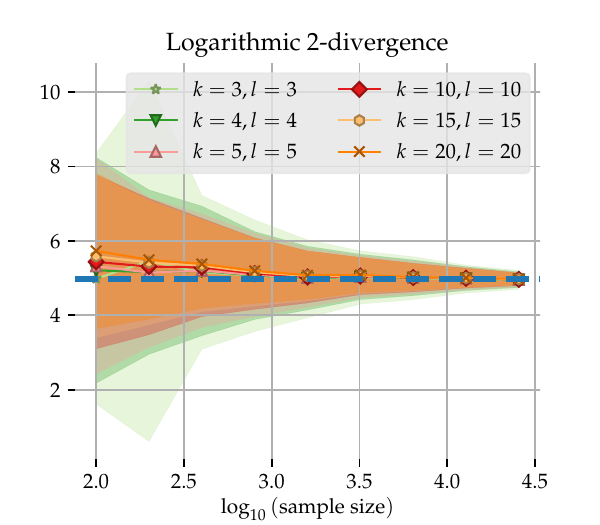}}
\end{minipage}
\begin{minipage}[b]{0.28\linewidth}
\centerline{\includegraphics[width=\textwidth]{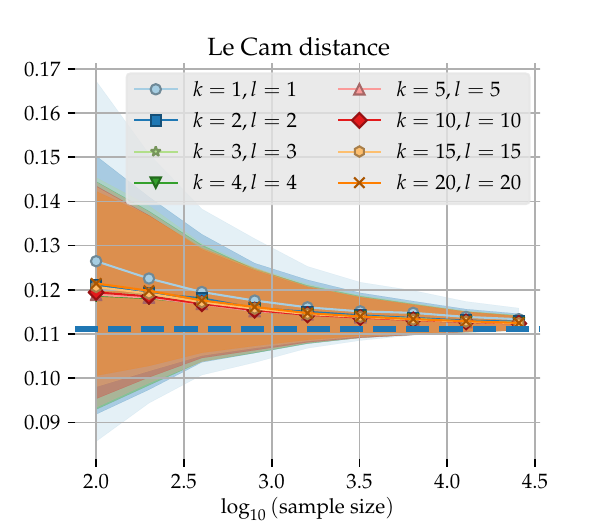}}
\end{minipage}
\qquad
\begin{minipage}[b]{0.28\linewidth}
\centerline{\includegraphics[width=\textwidth]{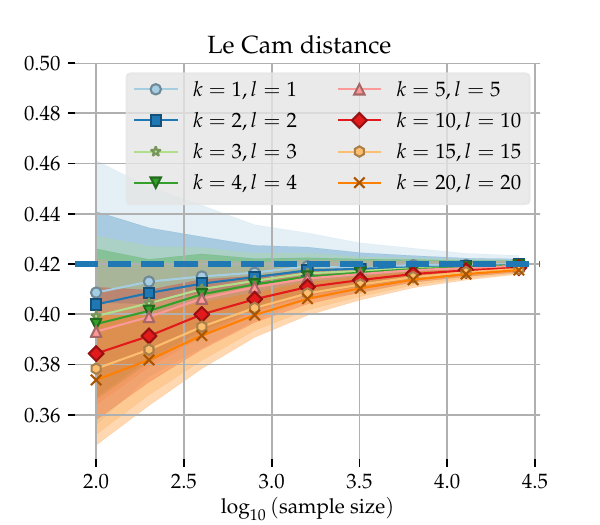}}\end{minipage}
\qquad
\begin{minipage}[b]{0.28\linewidth}
\centerline{\includegraphics[width=\textwidth]{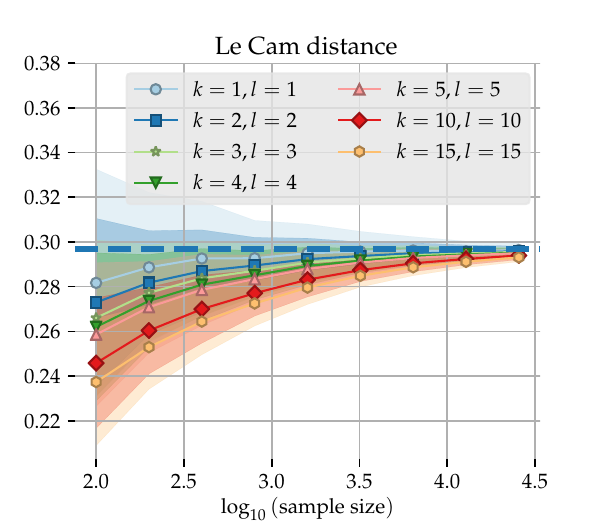}}
\end{minipage}
\begin{minipage}[b]{0.28\linewidth}
\centerline{\includegraphics[width=\textwidth]{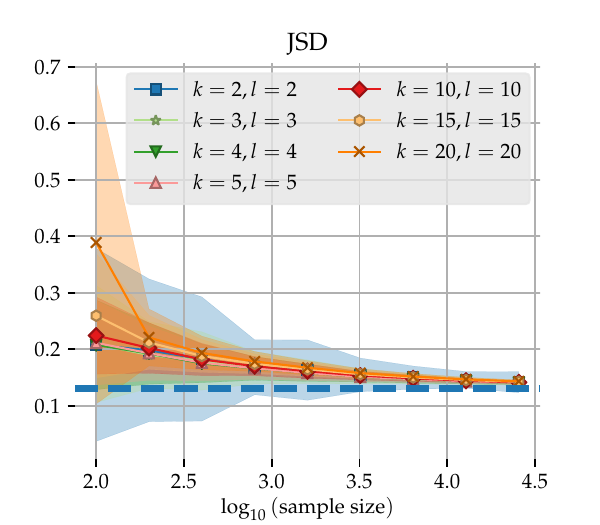}}
\end{minipage}
\qquad
\begin{minipage}[b]{0.28\linewidth}
\centerline{\includegraphics[width=\textwidth]{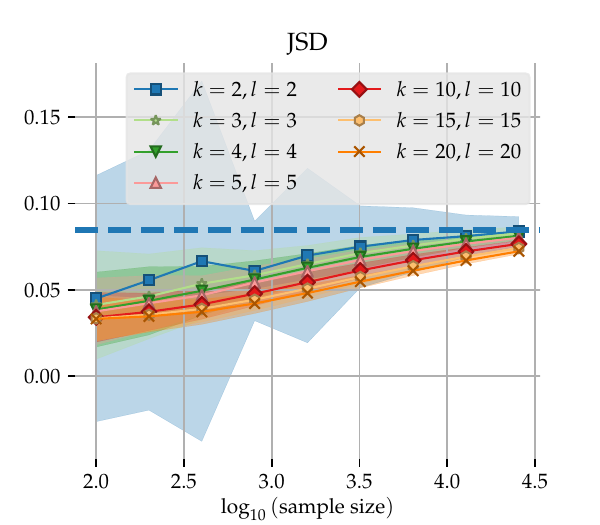}}\end{minipage}
\qquad
\begin{minipage}[b]{0.28\linewidth}
\centerline{\includegraphics[width=\textwidth]{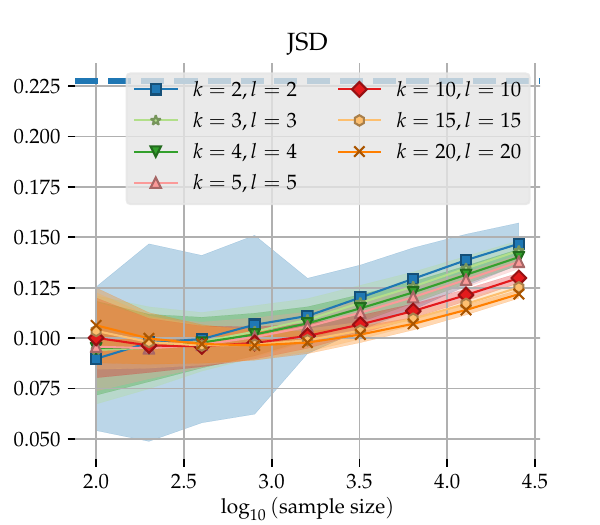}}
\end{minipage}
\vspace{-1.0em}
\caption{Convergence of the double-density functional estimator for KL divergence, $1.5$-divergence, and logarithmic $2$-divergence for 3-dimensional densities. 
The first, second, and third columns present simulation results for the densities $p$ and $q$ considered as $\Unif([0,1]^3)$ and $\Unif([0,2]^3)$, $\Normal(0,I_3)$ restricted to $\|\xb\|\le 3$ and $\Normal(0,4I_3)$ restricted to $\|\xb\|\le 3$, and $\Normal(0,I_3)$ and $\Normal(0,4I_3)$, respectively.
The true functional values are indicated as dashed lines and one sample standard deviations of the estimates are indicated as shaded area.
LCD and JSD are abbreviations for Le Cam distance and Jensen--Shannon divergence, respectively.
}
\label{fig:convergence:double}
\end{figure*}

\begin{figure*}[htp]
\centering
\begin{minipage}[b]{0.28\linewidth}
  \centerline{\includegraphics[width=\textwidth]{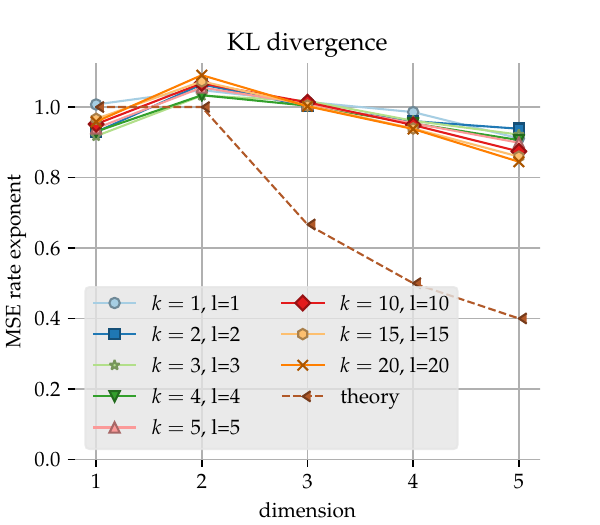}}
\end{minipage}
\qquad
\begin{minipage}[b]{0.28\linewidth}
  \centerline{\includegraphics[width=\textwidth]{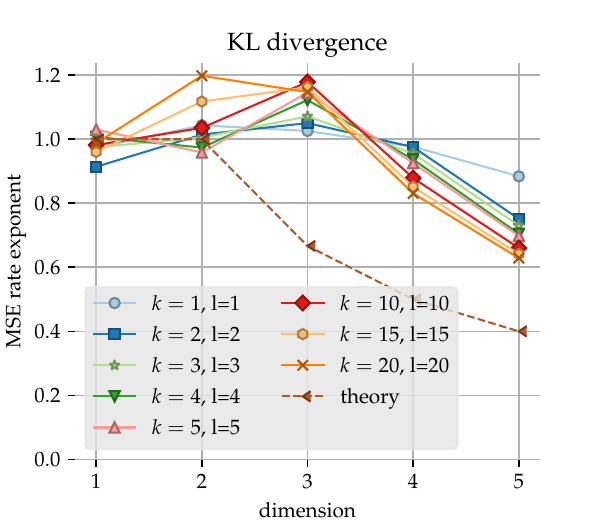}}
\end{minipage}
\qquad
\begin{minipage}[b]{0.28\linewidth}
  \centerline{\includegraphics[width=\textwidth]{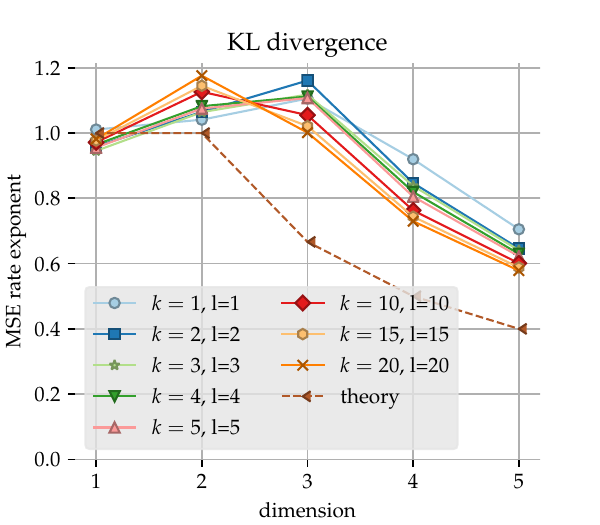}}
\end{minipage}\\
\begin{minipage}[b]{0.28\linewidth}
  \centerline{\includegraphics[width=\textwidth]{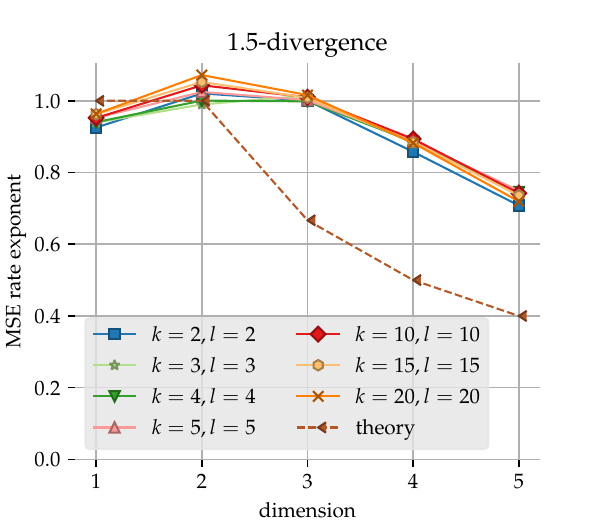}}
\end{minipage}
\qquad
\begin{minipage}[b]{0.28\linewidth}
  \centerline{\includegraphics[width=\textwidth]{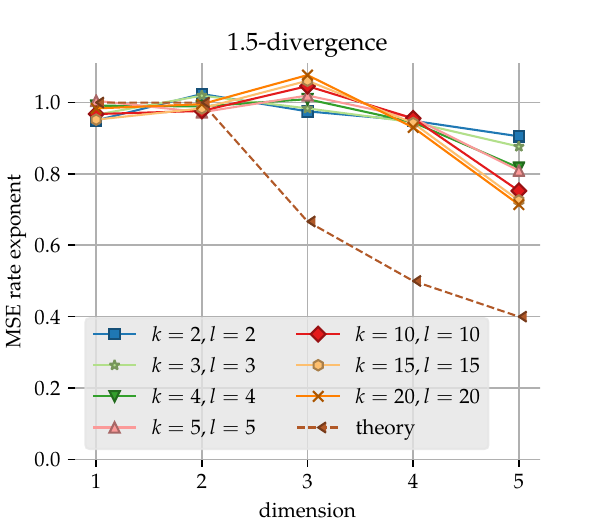}}
\end{minipage}
\qquad
\begin{minipage}[b]{0.28\linewidth}
  \centerline{\includegraphics[width=\textwidth]{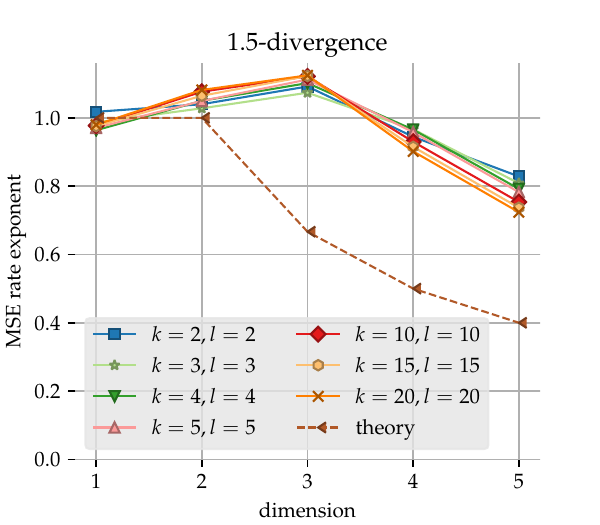}}
\end{minipage}
\begin{minipage}[b]{0.28\linewidth}
\centerline{\includegraphics[width=\textwidth]{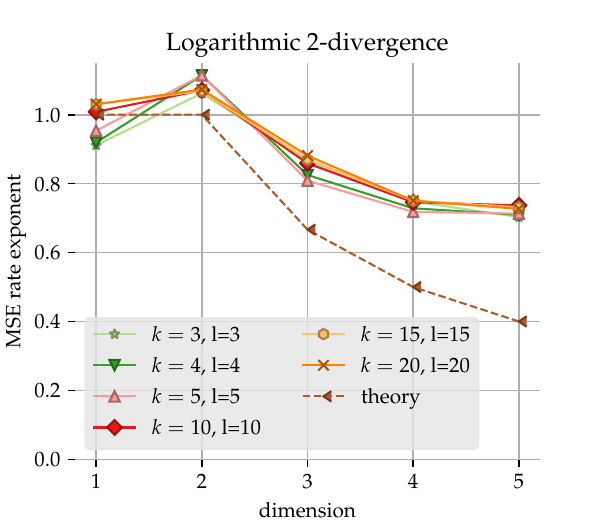}}
\end{minipage}
\qquad
\begin{minipage}[b]{0.28\linewidth}
\centerline{\includegraphics[width=\textwidth]{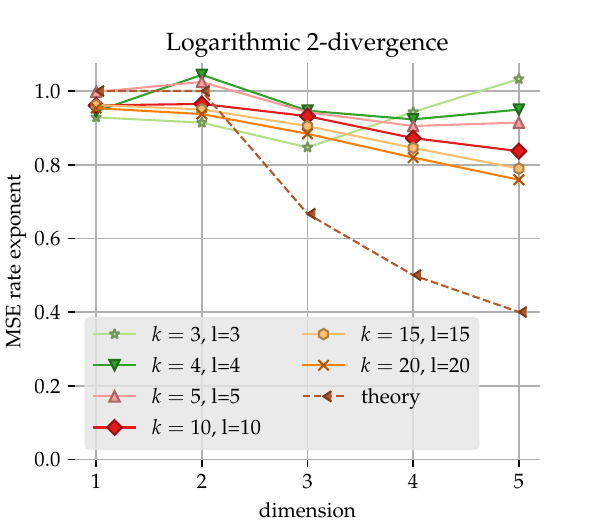}}
\end{minipage}
\qquad
\begin{minipage}[b]{0.28\linewidth}
\centerline{\includegraphics[width=\textwidth]{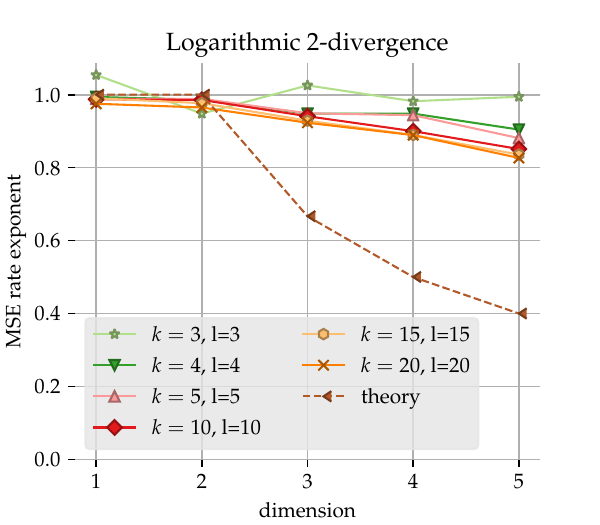}}
\end{minipage}
\begin{minipage}[b]{0.28\linewidth}
\centerline{\includegraphics[width=\textwidth]{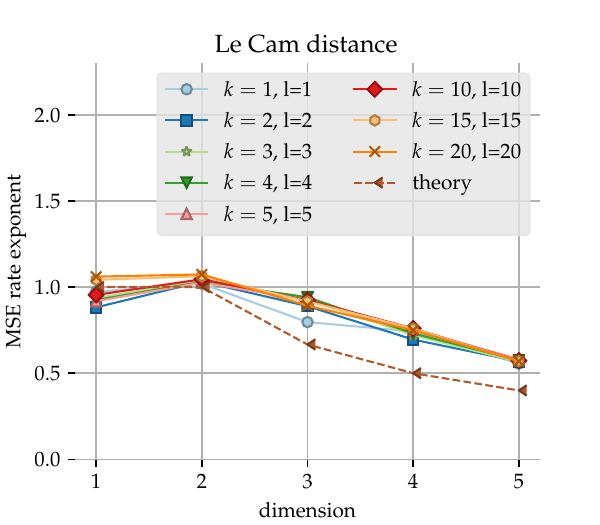}}
\end{minipage}
\qquad
\begin{minipage}[b]{0.28\linewidth}
\centerline{\includegraphics[width=\textwidth]{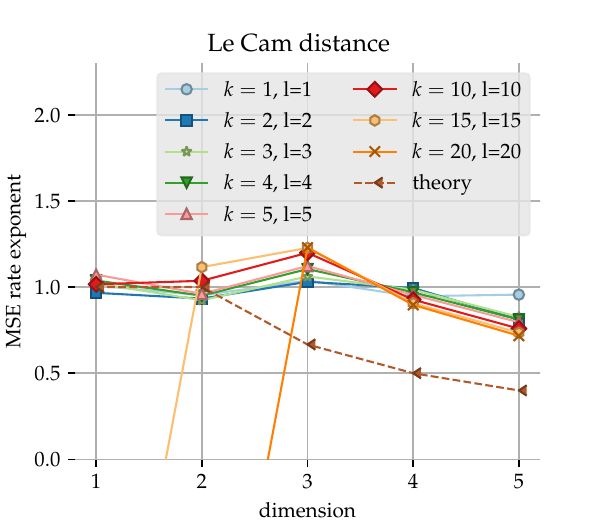}}
\end{minipage}
\qquad
\begin{minipage}[b]{0.28\linewidth}
\centerline{\includegraphics[width=\textwidth]{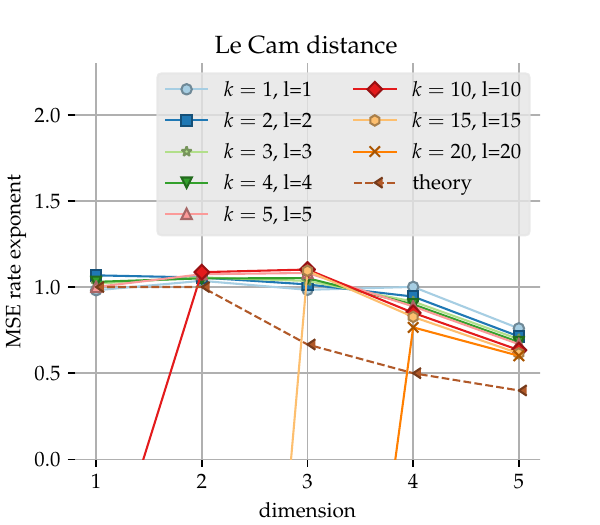}}
\end{minipage}
\begin{minipage}[b]{0.28\linewidth}
\centerline{\includegraphics[width=\textwidth]{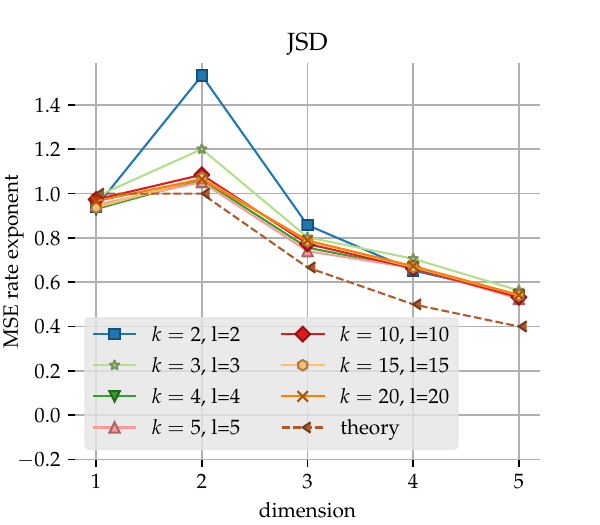}}
\end{minipage}
\qquad
\begin{minipage}[b]{0.28\linewidth}
\centerline{\includegraphics[width=\textwidth]{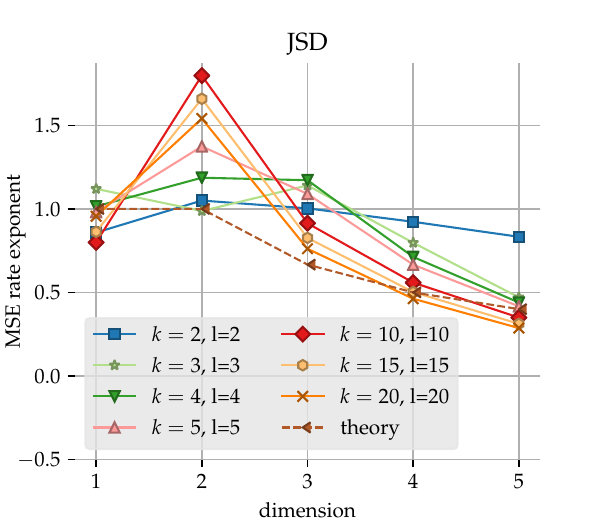}}\end{minipage}
\qquad
\begin{minipage}[b]{0.28\linewidth}
\centerline{\includegraphics[width=\textwidth]{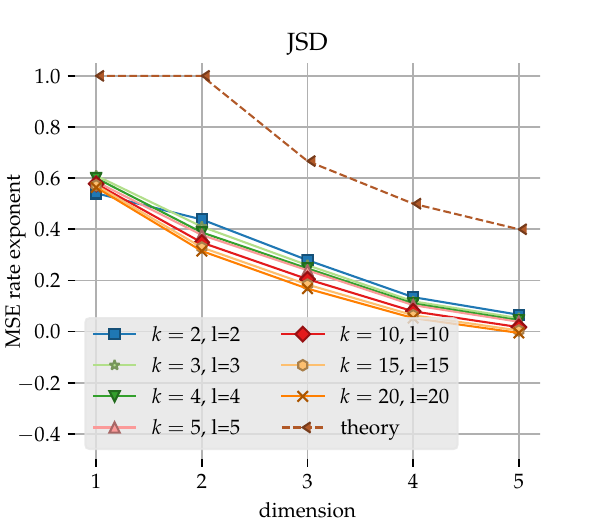}}
\end{minipage}
\vspace{-0.5em}
\caption{Simulated MSE rate exponents of the double-density functional estimator for KL divergence, $1.5$-divergence, and logarithmic $2$-divergence. 
The first, second, and third columns present simulation results for the densities $p$ and $q$ considered as $\Unif([0,1]^3)$ and $\Unif([0,2]^3)$, $\Normal(0,I_3)$ restricted to $\|\xb\|\le 3$ and $\Normal(0,4I_3)$ restricted to $\|\xb\|\le 3$, and $\Normal(0,I_3)$ and $\Normal(0,4I_3)$, respectively, for $d\in\{1,2,3,4,5\}$.
LCD and JSD are abbreviations for Le Cam distance and Jensen--Shannon divergence, respectively.}
\label{fig:exponent:double}
\end{figure*}

Corresponding simulation results for a few representative double-density functionals (KL divergence, $\a$-divergence, logarithmic $\a$-divergence, Le Cam distance, and Jensen--Shannon divergence) are presented in Figures~\ref{fig:convergence:double} and \ref{fig:exponent:double}.
These simulations indicate that the requirement $k>-4a$ and $l>-4\at$ in Theorem~\ref{thm:vanishing_variance:double} may be relaxed to the milder condition $k>-2a$ and $l>-2\at$.
For example, the estimator with $k=l=4$ for logarithmic $2$-divergence ($k=3\le -4(1-2)=4$ and $l=3\le -4(1-2)$) still exhibit consistency in Figure~\ref{fig:convergence:double}.
As presented in the last two rows in Figures~\ref{fig:convergence:double} and \ref{fig:exponent:double}, simulations also indicate that our estimator is consistent in practice for the exceptional examples of Le Cam distance and Jensen--Shannon divergence even without truncation.
For estimating Le Cam distance, we observed that using too large values for $k$ or $l$ lead to bad convergence behavior for small dimensions; see, e.g., the case of $k=l=15$ for $d=1$ at the second column of the fourth row in Figure~\ref{fig:exponent:double}.

\section{Concluding remarks}\label{sec:conclusion}
In this paper, we developed a systematic approach to designing $k$-NN based consistent estimators for a variety of functionals, starting from the fundamental requirement of asymptotic unbiasedness and utilizing the limiting behavior of the $k$-NN statistics (Proposition~\ref{prop:knndist}).
The proposed estimators rediscovered and unified several existing $k$-NN based estimators for Shannon entropy, KL divergence, $\a$-entropies and $\a$-divergences, and polynomial functionals, which have been sporadically studied and individually analyzed in the literature. It demystified the need of the known, but rather ad-hoc ``bias corrections'' for some functionals, providing an alternative, principled recipe to identify  $L_2$-consistent estimators.
Our list of examples is not exhaustive; other density functionals in the same form may exist or may be discovered in future, and our recipe will furnish consistent $k$-NN estimators for the same, with nonasymptotic performance predicted by our current analysis.

\revision{We remark that the established convergence rates are not minimax optimal; see Remark~\ref{rem:minimax}.
As further noted in Remark~\ref{rem:smoothness}, the proposed estimators cannot adapt to a higher order of smoothness $\sigma>2$, due to the inherent limitation of positive-valued kernels.
One possible solution to both problems is the ensemble approach~\cite{Sricharan--Wei--Hero2013,Moon--Hero2014a} that takes a weighted average of multiple estimators based on the asymptotic bias expansion of each density functional estimator. 
Studying the ensemble version of the estimators is beyond our scope and left as a future direction; see \cite{Berrett--Samworth2019} for a weighted version of the proposed divergence functional estimator with local minimax optimality.}

Throughout the paper, we assumed the Euclidean distance $\rho(\xb,\yb)=\|\xb-\yb\|$.
We conclude the paper with specifying technical issues one needs to address in order to extend the results of this paper to a general metric measure space $(\Xc,\rho,\mu)$, where $(\Xc,\rho)$ is a complete separable metric space and $\mu$ is a locally finite measure on the Borel $\sigma$-algebra of $\Xc$ (see, \eg \citet{Sturm2006}).
Consider a $\mu$-absolutely continuous probability measure $\P$ with density $p$. 
In general, the weak convergence property in Proposition~\ref{prop:knndist} for asymptotic unbiasedness (Theorems~\ref{thm:vanishing_bias:single} and \ref{thm:vanishing_bias:double}) requires the Lebesgue differentiation theorem to hold in the metric measure space $(\Xc,\rho,\mu)$, \ie we need 
\begin{align*}
\lim_{r\to 0} \frac{\P(\Bb(x,r))}{\mu(\Bb(x,r))} = p(x)
\end{align*}
for $\mu$-a.e.\@ $x\in\Xc$. Further, for the bias rate analysis to work, we need to extend Lemma~\ref{supp:lem:GOVLemma4_smoothing}, which states that if $p$ is locally $\sigma$-H\"older smooth on $\Bb(x,R)$, then for $r<R,$
\begin{align*}
    \Bigl|\frac{\P(\Bb(x,r))}{\mu(\Bb(x,r))} - p(x)\Bigr| &\lesssim r^{\sigma}
    \text{ and }
    \Bigl|\frac{\diff\P(\Bb(x,r))}{\diff\mu(\Bb(x,r))} - p(x)\Bigr| \lesssim r^{\sigma}.
\end{align*}
If there exists a nonsmooth boundary, we then further need Lemma~\ref{supp:lem:hausdorff} to hold in the metric measure space.
For the variance analysis to hold under $p$-norm and other norms, we can apply and extend the analysis in \cite{Gao--Oh--Viswanath2018tit} as pointed out earlier in Remark~\ref{rem:variance_analysis}.

\appendices
\section{Notation}
In what follows, let $\Rho_U(u)=\Pr\{U\le u\}$ and $\rho_U(u)=\diff \Rho_U(u) / \diff u$ denote the cumulative distribution function (cdf) and the density of a random variable $U$, respectively.
We use $B_{n,P}$ to denote a binomial random variable with parameters $n$ and $P$.
We use $P_q$ to denote a Poisson random variable with rate $q>0$.
We use $X_{\a,\b}$ to denote a beta random variable with parameters $\a,\b>0$ for $\a,\b>0$, whose density is
\begin{align*}
\frac{t^{\a-1}(1-t)^{\b-1}}{\Beta(\a,\b)},\quad 0\le t\le 1.
\end{align*}
Here $\Beta(\a,\b)\defeq \int_0^1 t^{\a-1}(1-t)^{\b-1}\diff t$ denotes the beta function.
Finally, we use 
$H^{d-1}$ to denote the $(d-1)$-dimensional Hausdorff measure.

\section{Technical lemmas}\label{supp:sec:technical_lemmas}

\subsection{Auxiliary lemmas}
\begin{lemma}\label{supp:lem:consistency_measure_zero}
Assume that $\P$ and $\Pt$ have densities $p$ and $\pt$, respectively, with respect to the Lebesgue measure $\Leb$. 
If $\P\ll \Pt$, then
$\P(\{\xb\suchthat\minimal_r{\pt}(\xb)>0\})=1$ for any $r>0$.
\end{lemma}
\begin{proof}
Let $r>0$ be fixed.
We first observe that $\P(\supp(\pt))=1$, since
\begin{align*}
1-\P(\supp(\pt))
&= \int p(\xb)(1-\ones_{\supp(\pt)}(\xb))\diff\xb\\
&= \int p(\xb)\ones_{\{\exists \d>0 \text{ s.t.} \Pt(\Bb(\xb,\d))=0\}}\diff\xb\\
&\stackrel{(a)}{\le} \int p(\xb)\ones_{\{\exists \d>0 \text{ s.t.} \P(\Bb(\xb,\d))=0\}}\diff\xb,\\
&\stackrel{(b)}{=}0.
\end{align*}
Here, (a) follows from the absolute continuity $\P\ll\Pt$, and (b) follows since $p(\xb)=0$ for a.e.\@ $\xb$ over the set $\{\xb\suchthat \exists\d>0\text{ s.t.} \P(\Bb(\xb,\d))=0\}$, by the Lebesgue differentiation theorem.

Now, define $A_{\d}\pt(\xb)=\Pt(\Bb(\xb,\d))/\Leb(\Bb(\xb,\d))$ for each $\d>0$ and $\xb\in\Real^d$.
On the one hand, we have
\[
\lim_{\d\to 0} A_{\d}\pt(\xb)
=\pt(\xb)
\]
for $\Leb$-a.e. $\xb$ by the Lebesgue differentiation theorem.
On the other hand, for $\xb\in T\cap\supp(\pt)$ where $T\defeq\{\xb\suchthat\minimal_r{\pt}(\xb)=0\}$, we have
$A_{\d}\pt(\xb)>0$ for every $\d > 0$ and
\[
0 = \minimal_r{\pt}(\xb) = \inf_{0<\d\le r}A_{\d}\pt(\xb)
\]
for any $r>0$.
Hence, we must have
\[
\pt(\xb) = \lim_{\d\to 0} A_{\d}\pt(\xb) =0
\]  
for $\Leb$-a.e. $\xb\in T\cap \supp(\pt)$, which implies that $\Pt(T\cap\supp(\pt)) = 0$, and thus  $\P(T\cap\supp(\pt)) = 0$ since $\P\ll\Pt$. 
This, together with $\P(\supp(\pt))=1$, establishes that $\P(T)=0$.
\end{proof}

\begin{lemma}\label{supp:lem:incomplete_gamma}
For the \emph{lower incomplete gamma function}
$\gamma(s,x)\defeq \int_0^x t^{s-1}e^{-t}\diff t$
and the \emph{upper incomplete gamma function} $\Gamma(s,x)\defeq \int_x^\infty t^{s-1}e^{-t}\diff t$,
we have
\begin{align}
\gamma(s,x)
&\le \Gamma(s)\wedge\frac{x^s}{s}, &&\forall s>0, x> 0,
\label{supp:eq:incomplete_gamma_lower}\\
\Gamma(s,x)
&\le \Gamma(s)x^{s-1}e^{-x+1}, &&\forall s\ge 1, x\ge 1.
\label{supp:eq:incomplete_gamma_upper}
\end{align}
\end{lemma}
\begin{proof}
As $\Gamma(s,x)/\Gamma(s)$ is decreasing in $s$ for fixed $x\ge 1$, we have that for $s\ge 1$,
\begin{align*}
    \frac{\Gamma(s,x)}{\Gamma(s)}
    \le \frac{\Gamma(\floor{s},x)}{\Gamma(\floor{s})}
    &= e^{-x}\sum_{k=0}^{\floor{s}-1} \frac{x^k}{k!}\\
    &\le e^{-x}x^{\floor{s}-1}\sum_{k=0}^\infty \frac{1}{k!}
    \le e^{-x+1}x^{s-1}.
\end{align*}

The second inequality follows since, for any $x>0$, letting $t=xe^{-u}$, we have
\begin{align*}
\gamma(s,x)
=\int_0^x t^{s-1}e^{-t}\diff t
&=x^s\int_0^{\infty} e^{-(su+xe^{-u})}\diff u\\
&\le x^s \int_0^{\infty} e^{-su}\diff u 
= \frac{x^s}{s}.\qedhere
\end{align*}
\end{proof}

\subsection{Convergence of distribution of \texorpdfstring{$k$}{k}-NN statistics}
\label{supp:sec:technical_lemmas:knn_stats}
We first state a basic statistical property of $k$-NN statistics. 

\begin{lemma}[Distribution of $k$-NN distance]
\label{supp:lem:knn_ball_cdf}
The cdf of $r_{km}(\xb)$ is 
\begin{align*}
\Rho_{r_{km}(\xb)}(r)
&= \Pr\{B_{m,\P(\Bb(\xb,r))}\ge k\}
=\Rho_{X_{k,m-k+1}}(\P(\Bb(\xb,r))).
\end{align*}
\end{lemma}
\begin{proof}
Consider
\begin{align*}
\Rho_{r_{km}(\xb)}(r)
&=\Pr\{r_{km}(\xb)\le r\}\\
&= \Pr\{\rho(\xb,\Xb_{(k)}(\xb))\le r\}\\
&=\Pr\bigl\{\abs{\{i\in[m]\suchthat \Xb_i\in \Bb(\xb,r)\}} \ge k\bigr\}\\
&= \Pr\{B_{m,\P(\Bb(\xb,r))}\ge k\}\\
&= \Rho_{X_{k,m-k+1}}(\P(\Bb(\xb,r))).
\end{align*}
The last equality follows from the identity
\begin{equation*}
\Pr\{B_{m,P}\ge k\}=\Rho_{X_{k,m-k+1}}(P).\qedhere
\end{equation*}
\end{proof}

Using this fact, Proposition~\ref{prop:knndist}, which claims the weak convergence of the $k$-NN statistics $U_{km}(\xb)$ to a Gamma random variable, readily follows. 
\begin{proof}[Proof of Proposition~\ref{prop:knndist}]
Fix $\xb\in\Real^d$ and $u>0$, and let $P_m\defeq \P(\Bb(\xb,\rvol(\frac{u}{m})))$.
Since $\Rho_{U_{km}(\xb)}(u)=\Rho_{r_{km}(\xb)}(\rvol(\frac{u}{m}))$, we have
\begin{align*}
\Rho_{U_{km}(\xb)}(u)
&= \Pr\{B_{m,P_m}\ge k\}
\end{align*}
from Lemma~\ref{supp:lem:knn_ball_cdf}.
By the Lebesgue differentiation theorem (see, \eg \cite{Rudin1987}), for $\Leb$-a.e.\@ $\xb$,
\begin{align*}
\lim_{m\to\infty}mP_m
=\lim_{m\to\infty}
u\frac{\P(\Bb(\xb,\rvol(\frac{u}{m})))}{\Leb(\Bb(\xb,\rvol(\frac{u}{m})))}
=up(\xb).
\end{align*}
Therefore, for each $i=0,\ldots,k-1$, we have
\begin{align*}
&\lim_{m\to\infty}\binom{m}{i}P_m^i(1-P_m)^{m-i}\\
&=\lim_{m\to\infty}\frac{i!}{m^i}\binom{m}{i}\bigl(1-P_m\bigr)^{m-i}\frac{(mP_m)^i}{i!}\\
&= e^{-u p(\xb)}\frac{(u p(\xb))^i}{i!},
\end{align*}
since 
\begin{align*}
\lim_{m\to\infty}\frac{i!}{m^i}\binom{m}{i}&=1
\text{ and }
\lim_{m\to\infty}(1-P_m)^{m-i}=e^{-up(\xb)}.
\end{align*}
This leads us to concludes that 
\begin{align*}
\pushQED{\qed} 
\lim_{m\to\infty}\Pr\{U_{km}(\xb)>u\}
&=\sum_{i=0}^{k-1}e^{-up(\xb)}\frac{up(\xb)^i}{i!}\\
&=\Pr\{U_{k\infty}(\xb)>u\},
\end{align*}
where $U_{k\infty}(\xb)$ is a $\GammaDist(k,p(\xb))$ random variable.
\end{proof}

Moreover, if the density $p$ is locally smooth, then one can establish a polynomial convergence rate of the density of $U_{km}(\xb)$ to $U_{k\infty}(\xb)$ as follows.

\begin{lemma}[Generalization of {\cite[Lemma~2]{Gao--Oh--Viswanath2018tit}}]\label{supp:lem:GOVLemma2_pdf_gap_bound}
Suppose that $\nu_m=o(\sqrt{m})$ and $k=k_m=o(\sqrt{m})$ as $m\to\infty$.
For $\xb\in\supp(p)$, if $p(\xb)\le \Cab_p<\infty$ and $p$ is $\sigma_p$-H\"older continuous ($\sigma_p\in[0,2]$) over $\Bb(\xb,\rvol(\frac{u}{m}))$ with H\"older constant $L$,
we have
\begin{align*}
&\bigl|\rho_{U_{km}(\xb)}(u)-\rho_{U_{k\infty}(\xb)}(u)\bigr|\\
&\lesssim_{\sigma_p,L,\Cab_p,d} 
(1+u)\Bigl(\frac{u}{m}\Bigr)^{\frac{\sigma_p}{d}}
+ k^{-k}\frac{(k^2+u^{2})u^{k-1}e^{-up(\xb)}}{m}
\end{align*}
for $u\in [0,\nu_m]$ and $m$ sufficiently large.
\end{lemma}

We first state two technical lemmas required to prove Lemma~\ref{supp:lem:GOVLemma2_pdf_gap_bound}, whose proofs are omitted here; we refer the interested readers to \cite{Gao--Oh--Viswanath2018tit}.
The first lemma in the following establishes a rate of convergence of a Poisson binomial random variable $B_{m,Q/m}\sim\Binom(m,Q/m)$ to a Poisson random variable $P_Q\sim\Poisson(Q)$ in distribution.

\begin{lemma}[Generalization of {\cite[Lemma~5]{Gao--Oh--Viswanath2018tit}}]\label{supp:lem:GOVLemma5_poisson}
For any $Q,k=o(\sqrt{m})$ as $m\to\infty$, there exists a constant $C_0>0$ such that for $m$ sufficiently large
\begin{align}
\bigl|\Pr\{B_{m,\frac{Q}{m}}=k\}-\Pr\{P_Q=k\}\bigr|
\le C_0\frac{Q^ke^{-Q}}{k!}\frac{(k^2+Q^2)}{m}.\nonumber
\end{align}
\end{lemma}
%

The second lemma establishes the speed of convergence of $\P(\Bb(\xb,r))/\Leb(\Bb(\xb,r))$ and $\diff\P(\Bb(\xb,r))/\diff\Leb(\Bb(\xb,r))$ to $p(\xb)$ as $r\to 0$, when $p$ is locally smooth at $\xb$.

\begin{lemma}[Generalization of {\cite[Lemma~4]{Gao--Oh--Viswanath2018tit}}]\label{supp:lem:GOVLemma4_smoothing}
If a density $p$ is $\sigma_p$-H\"older continuous with constant $L>0$ over $\Bb(\xb,R)$ for $\xb\in\Real^d$ and some $\sigma_p\in[0,2]$,
we have for any $0<r< R$,
\begin{align*}
\Bigl|\frac{\P(\Bb(\xb,r))}{\lambda(\Bb(\xb,r))}-p(\xb)\Bigr|
&\le \frac{d}{\sigma_p+d}Lr^{\sigma_p},\\
\Bigl|\frac{\diff\P(\Bb(\xb,r))}{\diff\lambda(\Bb(\xb,r))}-p(\xb)\Bigr| 
&\le Lr^{\sigma_p}.
\end{align*}
\end{lemma}

The proof of the first inequality can be found in \cite{Jiao--Gao--Han2018} and the second inequality can be proved by a similar argument.

\begin{remark}
If $g$ is bounded above over $\Bb(\xb,R)$, then $g$ is $\sigma_p$-H\"older continuous over $\Bb(\xb,R)$ with $\sigma_p=0$.
The convergence of $U_{km}(\xb)$ to a $\GammaDist(k,p(\xb))$ random variable as $m\to\infty$ can be quantified in terms of a gap between the densities using this lemma and the order of smoothness $\sigma_p$ of the underlying density $p$; however, the bounds in Lemma~\ref{supp:lem:GOVLemma4_smoothing} cannot be improved further beyond $\bigO(r^2)$.
It is consistent with the observation that the higher-order smoothness beyond 2 cannot be exploited with $k$-NN methods \cite{Tsybakov--van-der-Meulen1996,Han--Jiao--Weissman--Wu2017}.
\end{remark}

Now we are ready to present the proof of Lemma~\ref{supp:lem:GOVLemma2_pdf_gap_bound}.
\begin{proof}[Proof of Lemma~\ref{supp:lem:GOVLemma2_pdf_gap_bound}]
First note that the density of the $k$-th NN statistics $r_{km}(\xb)$ is
\begin{align*}
\rho_{r_{km}(\xb)}(r)
&= m\Pr\{B_{m-1,\P(\Bb(\xb,r))}=k-1\} \frac{\diff \P(\Bb(\xb,r))}{\diff r}\\
&= g_{km}(\P(\Bb(\xb,r))) \frac{\diff \P(\Bb(\xb,r))}{\diff r}
\end{align*}
from Lemma~\ref{supp:lem:knn_ball_cdf} in Appendix~\ref{supp:sec:technical_lemmas:knn_stats}.
Here we define
\begin{align*}
g_{km}(P)\defeq m\Pr\{B_{m-1,P}=k-1\}
\end{align*}
for $p\in [0,1]$, which is the density of the $k$-th order statistic from among $m$ random samples drawn from the uniform distribution over $[0,1]$.
It is easy to check that
$g_{km}(P)\le m$ and $g_{km}'(P)\le 2m(m-1)\le 2m^2$ for any $P\in [0,1]$.
Recall that $P_m(u|\xb)\defeq \P(\Bb(\xb,\rvol(\frac{u}{m})))$.
The density of $U_{km}(\xb)$ can then be written as
\begin{align*}
\rho_{U_{km}(\xb)}(u) 
&= \rho_{r_{km}(\xb)}\bigl(\rvol\bigl(\frac{u}{m}\bigr)\bigr) \frac{\diff \rvol(\frac{u}{m})}{\diff u}\\
&= g_{km}(P_m(u|\xb)) \frac{\diff P_m(u|\xb)}{\diff u}.
\end{align*}
We define an intermediate density approximation
\begin{align*}
\rho_{km}(u)
&\defeq g_{km}\Bigl(\frac{up(\xb)}{m}\Bigr)\frac{p(\xb)}{m}
\end{align*}
for $u\le m/\Cab_p$, and bound the density gap by
\begin{align*}
&\abs*{\rho_{U_{km}(\xb)}(u)-\rho_{U_{k\infty}(\xb)}(u)}\\
&\le \abs*{\rho_{U_{km}(\xb)}(u)-\rho_{km}(u)}
+ \abs*{\rho_{km}(u)-\rho_{U_{k\infty}(\xb)}(u)}.
\end{align*}
We bound each term on the right hand side.

For the first term, consider
\begin{align*}
&|\rho_{U_{km}(\xb)}(u)-\rho_{km}(u)|\\
&\le g_{km}(P_m(u|\xb))
\Bigl|\frac{\diff P_m(u|\xb)}{\diff u}-\frac{p(\xb)}{m}\Bigr|
    \\&\quad 
    +\Bigl|g_{km}(P_m(u|\xb))-g_{km}\Bigl(\frac{up(\xb)}{m}\Bigr)\Bigr|\frac{p(\xb)}{m}\\
&\le g_{km}(P_m(u|\xb))\Bigl|\frac{\diff P_m(u|\xb)}{\diff u}-\frac{p(\xb)}{m}\Bigr|
	\\&\quad 
	+\max_{p\in(0,1)}|g_{km}'(p)| \Bigl|P_m(u|\xb)-\frac{up(\xb)}{m}\Bigr| \frac{p(\xb)}{m}\\
&\le m\Bigl|\frac{dP_m(u|\xb)}{\diff u}-\frac{p(\xb)}{m}\Bigr|
	+2m^2\Bigl|P_m(u|\xb)-\frac{up(\xb)}{m}\Bigr|\frac{p(\xb)}{m}\\
&= 
\Bigl|
\frac{\diff \P(\Bb(\xb,\rvol(\frac{u}{m})))}{\diff \lambda(\Bb(\xb,\rvol(\frac{u}{m}))}
-p(\xb)
\Bigr|
	+2up(\xb)\Bigl|\frac{\P(\Bb(\xb,\rvol(\frac{u}{m})))}{\Leb(\Bb(\xb,\rvol(\frac{u}{m}))}-p(\xb)\Bigr|\\
&\le \Bigl(1+2\Cab_p\frac{d}{\sigma_p+d} u\Bigr)L\rvol^{\sigma_p}\bigl(\frac{u}{m}\bigr)\\
&\lesssim_{\sigma_p,L,\Cab_p,d} (1+u)\Bigl(\frac{u}{m}\Bigr)^{\frac{\sigma_p}{d}}.
\end{align*}
The second last inequality follows from Lemma~\ref{supp:lem:GOVLemma4_smoothing}.
Note that this term is independent of $k$.

The second term can be bounded using Lemma~\ref{supp:lem:GOVLemma5_poisson}.
For $m$ sufficiently large, we have
\begin{align*}
&\abs*{\rho_{km}(u)-\rho_{U_{k\infty}(\xb)}(u)}\\
&=\frac{k}{u}\abs*{\Pr\bigl\{B_{m,up(\xb)/m}=k\bigr\}-\Pr\{P_{up(\xb)}=k\}}\\
&\le \frac{k}{u}C_0\frac{(up(\xb))^ke^{-up(\xb)}}{k!}\frac{k^2+u^2p^2(\xb)}{m}\\
&= \frac{C_0}{\Gamma(k)}\frac{(k^2+(up(\xb))^2) (up(\xb))^{k}e^{-up(\xb)}}{mu}\\
&\lesssim_{C_0,\Cab_p} k^{-k}\frac{(k^2+u^2)u^{k-1}e^{-up(\xb)}}{m},
\end{align*}
which holds uniformly for all $u,k=o(\sqrt{m})$ as $m\to\infty$. Here we use the Stirling approximation $\Cab_p^k/k!\sim (e\Cab_p)^k/k^{k+\half}$.
\end{proof}

\begin{remark}
This proof closely follows the one in \citep{Gao--Oh--Viswanath2018tit}, while keeping track of the explicit dependence on the constants $C_0, \Cab_p$ and $k$.
\end{remark}

The following lemma quantifies the convergence of the cdf of $U_{km}(\xb)$ to the cdf of $U_{k\infty}(\xb)$ when the underlying density $p$ is smooth.

\begin{lemma}[Generalization of {\cite[Lemma~3]{Gao--Oh--Viswanath2018tit}}]\label{supp:lem:GOVLemma3_cdf_gap_bound}
Suppose that $\nu_m=o(\sqrt{m})$ and $k=k_m=o(\sqrt{m})$ as $m\to\infty$.
For $\xb\in\supp(p)$, if $p(\xb)\le \Cab_p<\infty$ and $p$ is $\sigma_p$-H\"older continuous ($\sigma_p\in[0,2]$) over $\Bb(\xb,\rvol(u/m))$ with H\"older constant $L$,
we have
\begin{align}
&\abs*{\Rho_{U_{km}(\xb)}(u)-\Rho_{U_{k\infty}(\xb)}(u)}\nonumber\\
&\lesssim_{\sigma_p,L,\Cab_p,d}
  ku\Bigl(\frac{u}{m}\Bigr)^{\frac{\sigma_p}{d}}
  +\frac{(k^2+u^2)u^{k-1}e^{-up(\xb)}}{m}
  \label{supp:eq:cdf_bound},
\end{align}
for $u\in [1/(p(\xb)),\nu_m]$ for $m$ sufficiently large.
\end{lemma}

\begin{proof}
First, note that
\begin{align*}
\Rho_{U_{k\infty}(\xb)}(u)
&=1-\sum_{j=0}^{k-1} \Pr\{P_{up(\xb)}=j\}
\end{align*}
and
\begin{align*}
\Rho_{U_{km}(\xb)}(u)
&= 1-\sum_{j=0}^{k-1} \Pr\{B_{m,P_m(u|\xb)}=j\},
\end{align*}
from Lemma~\ref{supp:lem:knn_ball_cdf} in Appendix~\ref{supp:sec:technical_lemmas:knn_stats}.
By triangle inequality, we have
\begin{align*}
&\bigl|\Rho_{U_{km}(\xb)}(u)-\Rho_{U_{k\infty}(\xb)}(u)\bigr|\\
&\le \sum_{j=0}^{k-1} \bigl|\Pr\{P_{up(\xb)}=j\}-\Pr\{B_{m,P_m(u|\xb)}=j\}\bigr|\\
&\le \sum_{j=0}^{k-1} \Bigl\{\bigl|\Pr\{P_{up(\xb)}=j\}-\Pr\{B_{m,\frac{up(\xb)}{m}}=j\}\bigr|
\\&\qquad\quad
+\bigl|\Pr\{B_{m,\frac{up(\xb)}{m}}=j\}-\Pr\{B_{m,P_m(u|\xb)}=j\}\bigr|\Bigr\}.
\end{align*}

For the first term, using Lemma~\ref{supp:lem:GOVLemma5_poisson}, we obtain
\begin{align*}
&\bigl|\Pr\{P_{up(\xb)}=j\}-\Pr\bigl\{B_{m,\frac{up(\xb)}{m}}=j\bigr\}\bigr|\\
&\le C_0\frac{(up(\xb))^je^{-up(\xb)}}{j!}\frac{j^2+(up(\xb))^2}{m},
\end{align*}
for each $j=0,\ldots,k-1$,
which implies that
\begin{align*}
&\sum_{j=0}^{k-1}\bigl|\Pr\{P_{up(\xb)}=j\}-\Pr\bigl\{B_{m,\frac{up(\xb)}{m}}=j\bigr\}\bigr|\\
&\le C_0\frac{k^2+(up(\xb))^2}{m}e^{-up(\xb)}\sum_{j=0}^{k-1}\frac{(up(\xb))^j}{j!}\\
&=C_0\frac{k^2+(up(\xb))^2}{m}\frac{\Gamma(k,up(\xb))}{\Gamma(k)}\\
&\le C_0\frac{k^2+(up(\xb))^2}{m} (up(\xb))^{k-1}e^{-up(\xb)+1},
\end{align*}
where the last inequality follows from Lemma~\ref{supp:lem:incomplete_gamma}.

For the second term, we have
\begin{align*}
&\bigl|\Pr\bigl\{B_{m,\frac{up(\xb)}{m}}=j\bigr\}-\Pr\{B_{m,P_m(u|\xb)}=j\}\bigr|\\
&\le 2m\Bigl|P_m(u|\xb)-\frac{up(\xb)}{m}\Bigr|\\
&= 2u\Bigl|\frac{\P(\Bb(\xb,\rvol(\frac{u}{m})))}{\Leb(\Bb(\xb,\rvol(\frac{u}{m})))}-p(\xb)\Bigr|\\
&\le 2u\frac{d}{\sigma_p+d}L\rvol^{\sigma_p}\bigl(\frac{u}{m}\bigr),
\end{align*}
for each $j=0,\ldots,k-1$, from Lemma~\ref{supp:lem:GOVLemma4_smoothing}.

Putting the bounds together and using the triangle inequality, we have that for $k,u=o(\sqrt{m})$
\begin{align*}
&\bigl|\Rho_{U_{km}(\xb)}(u)-\Rho_{U_{k\infty}(\xb)}(u)\bigr|\\
&\le 
\frac{2kud}{\sigma_p+d}L\rvol^{\sigma_p}\bigl(\frac{u}{m}\bigr)
+C_0\frac{k^2+(up(\xb))^2}{m} (up(\xb))^{k-1}e^{-up(\xb)+1}\\
&\lesssim_{\sigma_p,d,L,C_0,\Cab_p} 
    ku\Bigl(\frac{u}{m}\Bigr)^{\frac{\sigma_p}{d}}
    +\frac{(k^2+u^2)u^{k-1}e^{-up(\xb)}}{m},
\end{align*}
which concludes the proof.
\end{proof}

\subsection{Bounds on distribution of \texorpdfstring{$k$}{k}-NN statistics}
We now present several bounds on 
\begin{align*}
F_{km}(u|\xb)
&\defeq \Pr\{U_{km}(\xb)\le u\}\\
&=\Pr\bigl\{r_{km}(\xb)\le \rvol\bigl(\frac{u}{m}\bigr)\bigr\}\\
&=\Pr\{B_{m,P_m(u|\xb)}\ge k\},
\end{align*}
which is the cdf of $U_{km}(\xb)$.
Here and henceforth, for $\xb\in\Real^d$ and $u\ge 0$, we define
\[
P_m(u|\xb)
\defeq \P\bigl(\Bb\bigl(\xb,\rvol\bigl(\frac{u}{m}\bigr)\bigr)\bigr)
= \frac{u}{m}\frac{\P(\Bb(\xb,\rvol(\frac{u}{m})))}{\lambda(\Bb(\xb,\rvol(\frac{u}{m})))}.
\]
Note that by the definitions of $\minimal_r{p}(\xb)$ and $\maximal_r{p}(\xb)$, we have
\begin{align*}
u'\minimal_{r}{p}(\xb) 
\le mP_m(u'|\xb) 
\le m\wedge (u'\maximal_{r}{p}(\xb))
\end{align*}
for $r=\rvol(\frac{u}{m})$ and for any $0< u'\le u$.

The following lemma presents an upper bound on the cdf $F_{km}(u|\xb)$.
\begin{lemma}[Generalization of {\cite[Eq.~(3.19)]{Bulinski--Dimitrov2019b}}]\label{supp:lem:bounds_cdf_of_u_km}
For any $\xb\in\Real^d$ and $u>0$, we have
\begin{align}
F_{km}(u|\xb)&\le \frac{(mP_m(u|\xb))^k}{k!}.
\label{supp:eq:upper_bound_cdf_of_u_km}
\end{align}
\end{lemma}
\begin{proof}
Since $F_{km}(u|\xb)=\Rho_{T_{k,m-k+1}}(P_m(u|\xb))$ from Lemma~\ref{supp:lem:knn_ball_cdf}, we have
\begin{align*}
F_{km}(u|\xb)
&=\int_0^{P_m(u|\xb)} \frac{t^{k-1}(1-t)^{m-k}}{\Beta(k,m-k+1)}\diff t\\
&\le \frac{P_m^k(u|\xb)}{k\Beta(k,m-k+1)}\\
&= \binom{m}{k}P_m^k(u|\xb)\\
&\le \frac{(mP_m(u|\xb))^k}{k!},
\end{align*}
which concludes the proof.
\end{proof}

We present two upper bounds on the complementary cdf $1-F_{km}(u|\xb)$.

\begin{lemma}[{\cite[Eq.~(3.23)]{Bulinski--Dimitrov2019b}}]
\label{supp:lem:bounds_ccdf_of_u_km1}
For any $\xb\in\Real^d$, $0<D<1$, and $u\ge 0$, 
we have
\begin{align}
1-F_{km}(u|\xb)&\le (1-D)^{-k+1}e^{-DmP_m(u|\xb)}.
\label{supp:eq:lower_bound_cdf_of_u_km_11}
\end{align}
In particular, if $mP_m(u|\xb)>k$, we have
\begin{align}
1-F_{km}(u|\xb)\le \Bigl(\frac{emP_m(u|\xb)}{k}\Bigr)^k e^{-mP_m(u|\xb)}.
\label{supp:eq:lower_bound_cdf_of_u_km_12}
\end{align}
\end{lemma}
\begin{proof}
Since we can write $1-F_{km}(u|\xb)=\Pr\{B_{m,P_m(u|\xb)}<k\}$ from Lemma~\ref{supp:lem:knn_ball_cdf}, 
the bound follows immediately from a Chernoff bound on a binomial random variable.
For any $\lambda>0$,
\begin{align*}
\Pr\{B_{m,P}<k\} 
&\le e^{\lambda k}\E[e^{-\lambda B_{m,P}}]\\
&= e^{\lambda k}(1-P+Pe^{-\lambda})^m\\
&\le e^{\lambda k}e^{-mP(1-e^{-\lambda})},
\end{align*}
and this proves \eqref{supp:eq:lower_bound_cdf_of_u_km_11} if we set $D\defeq 1-e^{-\lambda}\in(0,1)$.
If $mP> k$, we then can minimize the right hand side by plugging in $\lambda=\ln \frac{mp}{k}$, which obtains
\[
\Pr\{B_{m,P}<k\} \le
\Bigl(\frac{emP}{k}\Bigr)^k e^{-mP}.\qedhere
\]
\end{proof}

\begin{lemma}[{\cite[Eq.~(3.32)]{Bulinski--Dimitrov2019b}}]\label{supp:lem:bounds_ccdf_of_u_km2}
For any $\xb\in\Real^d$, $\d>0$, $m\ge (1+1/\d)(k-1)$, and $u\ge 0$, we have
\begin{align}
1-F_{km}(u|\xb)&\le (1+\d)(1-P_m(u|\xb)).
\label{supp:eq:lower_bound_cdf_of_u_km_2}
\end{align}
\end{lemma}
\begin{proof}
Consider
\begin{align*}
&1-F_{km}(u|\xb)\\
&= \sum_{j=0}^{k-1}\binom{m}{j}P_m^j(u|\xb)(1-P_m(u|\xb))^{m-j}\\
&= (1-P_m(u|\xb))\\
&\qquad\times\sum_{j=0}^{k-1}\frac{m}{m-j}\binom{m-1}{j}P_m^j(u|\xb)(1-P_m(u|\xb))^{m-j-1}.
\end{align*}
For any fixed $\d>0$, if $m\ge (1+\d^{-1})(k-1)$, then
\[
\frac{m}{m-j} \le \frac{m}{m-k+1} \le 1+\d
\]
for $j=0,\ldots,k-1$.
Therefore, we have
\[
1-F_{km}(u|\xb)
\le (1+\d)(1-P_m(u|\xb)).\qedhere
\]
\end{proof}

\begin{lemma}\label{supp:lem:bound_kNN_ball_density}
If $p(\zb)\le \Cab_p$ for $\zb\in\widebar{\Bb}(\xb,r)$, we have
\begin{equation*}
\rho_{U_{km}(\xb)}(u)\le \frac{\Cab_p^k u^{k-1}}{\Gamma(k)}.
\end{equation*}
\end{lemma}

We first prove the following lemma.
Let us denote the sphere centered at $\xb\in\Real^d$ of radius $r>0$ by $\Sb(\xb,r)\defeq \{\yb\suchthat \rho(\xb,\yb)=r\}$.
Note that the the Hausdorff measure $H^{d-1}(\Sb(\xb,r))$ of the sphere is $d\ups_d r^{d-1}$.
\begin{lemma}\label{supp:lem:bound_p_and_dpdu}
If $p(\zb)\le \Cab_p$ for $\zb\in\Sb(\xb,r)$, we have
\begin{align*}
\frac{\diff \P(\Bb(\xb,r))}{\diff \Leb(\Bb(\xb,r))}
\le \Cab_p.
\end{align*}
\end{lemma}
\begin{proof}[Proof of Lemma~\ref{supp:lem:bound_p_and_dpdu}]
It is easy to see that $p(\xb)\le \maximal_r{p}(\xb)$ for any $r>0$ by contradiction.
From the coarea formula \cite{Evans--Gariepy2015}, we have
\begin{align*}
\frac{\diff \P(\Bb(\xb,r))}{\diff r}
&=\frac{\diff}{\diff r}
	\int_{\Bb(\xb,r)}p(\yb)\diff \yb\\
&=\int_{\Sb(\xb,r)}p(\yb)H^{d-1}(\diff \yb)\\
&\le \Cab_p(d\ups_d r^{d-1})
\end{align*}
since $p(\xb)\le \Cab_p$ for $\xb\in\Sb(\xb,r)$.
Therefore, we have
\begin{equation*}
\frac{\diff \P(\Bb(\xb,r))}{\diff \Leb(\Bb(\xb,r))}
= \frac{ \frac{\diff }{\diff r} \P(\Bb(\xb,r))}{ \frac{\diff }{\diff r} \Leb(\Bb(\xb,r))}
\le \Cab_p.\qedhere
\end{equation*}
\end{proof}

\begin{proof}[Proof of Lemma~\ref{supp:lem:bound_kNN_ball_density}]
Now, from Lemma~\ref{supp:lem:knn_ball_cdf} and Lemma~\ref{supp:lem:bound_p_and_dpdu}, if $p(\yb)\le\Cab_p$ for $\yb\in\Bb(\xb,r)$, then
\begin{align*}
\rho_{r_{km}(\xb)}(r)
&=\rho_{X_{k,m-k+1}}(\P(\Bb(\xb,r)))\frac{\diff\P(\Bb(\xb,r))}{\diff r}\\
&\le \frac{m^k}{\Gamma(k)}\P^{k-1}(\Bb(\xb,r))\frac{\diff\P(\Bb(\xb,r))}{\diff r}\\
&\le \frac{(\Cab_pm)^k}{\Gamma(k)}\frac{d}{r}\Leb^k(\Bb(\xb,r)).
\end{align*}
We then bound the density of $U_{km}(\xb)$ as
\begin{align*}
\rho_{U_{km}(\xb)}(u)
&= \rho_{r_{km}(\xb)}\bigl(\rvol\bigl(\frac{u}{m}\bigr)\bigr) \frac{\diff \rvol(\frac{u}{m})}{\diff u}\\
&\le \frac{(\Cab_pm)^k}{\Gamma(k)}\frac{d}{\rvol(\frac{u}{m})}\bigl(\frac{u}{m}\bigr)^k \frac{\rvol(\frac{u}{m})}{du}
= \frac{\Cab_p^k}{\Gamma(k)}u^{k-1},
\end{align*}
which concludes the proof.
\end{proof}

\subsection{Bounds on expected values of \texorpdfstring{$k$}{k}-NN statistics}
\label{supp:sec:technical_lemmas:generic}

Let $\tilde{f}_{km}(u|\xb)\defeq \rho_{\Ut_{km}(\xb)}(v)$ denote the density of the normalized volume $\Ut_{km}(\xb)=\Leb(\Bb(\xb,r_k(\xb|\tilde{\Xb}_{1:m})))$, where $\tilde{\Xb}_{1:m}$ is drawn \iid from density $\pt$. 
Later, the density $\pt$ may be identified as the density $p$ for $\Xb_{1:m}$ or the density $q$ for $\Yb_{1:n}$.
Pick any numbers $0\le \tau_m\le 1\le \nu_m \le \kappa_m <\infty$.
Suppose that we are given a nondecreasing function $\xi\in\Xi$.
For $(a,b)\in\Real^2$ and $k\in\Natural$, we define, for each $\xb\in\Real^d$
\begin{align*}
A_{km}(\xb;\pt;\xi)
&\defeq \int_0^{\tau_m} \xi(u^a)\ft_{km}(u|\xb)\diff u,
\numberthis\label{supp:eq:akm}
\\
B_{km}^{(1)}(\xb;\pt;\xi)
&\defeq \int_1^{\nu_m} \xi(u^b)\ft_{km}(u|\xb)\diff u,
\numberthis
\label{supp:eq:bkm1}
\\
B_{km}^{(2)}(\xb;\pt;\xi)
&\defeq \int_{\nu_m}^{\kappa_m} \xi(u^b)\ft_{km}(u|\xb)\diff u,
\numberthis
\label{supp:eq:bkm2}
\end{align*}
and
\begin{align*}
B_{km}^{(3)}(\xb;\pt;\xi)
&\defeq \int_{\kappa_m}^{\infty} \xi(u^b)\ft_{km}(u|\xb)\diff u.
\numberthis
\label{supp:eq:bkm3}
\end{align*}
\begin{lemma}\label{supp:lem:akm}
For $r=\varrho(\frac{\tau_m}{m})$, we have
\begin{align*}
&A_{km}(\xb;\pt;\xi)\\
&\le \frac{(\maximal_{r}{\pt}(\xb))^k}{k!}
\Bigl(\tau_m^k\xi(\tau_m^a) - \ones_{(-\infty,0)}(a)\int_0^{\tau_m} u^k\diff\xi(u^a)\Bigr).
\end{align*}
In particular, if $\tau_m=1$ and $-\int_0^1 u^k\diff\xi(u^a)<\infty$, we have for $r=\varrho(\frac{1}{m})$,
\begin{align*}
A_{km}(\xb;\pt;\xi)
&\lesssim \frac{(\maximal_{r}{\pt}(\xb))^k}{k!}.
\end{align*}
\end{lemma}
\begin{proof}
Integrating by parts and applying Lemma~\ref{supp:lem:bounds_cdf_of_u_km}, we have
\begin{align*}
&A_{km}(\xb;\pt;\xi)\\
&= \int_0^{\tau_m}\xi(u^a)\diff \Ft_{km}(u|\xb)\\
&\le \xi(\tau_m^{a})\Ft_{km}(\tau_m|\xb) - \int_0^{\tau_m} \Ft_{km}(u|\xb)\diff\xi(u^a)\\
&\le \frac{(\maximal_{\varrho(\frac{\tau_m}{m})}{\pt}(\xb))^k}{k!}\tau_m^k\xi(\tau_m^a) - \int_0^{\tau_m} \Ft_{km}(u|\xb)\diff\xi(u^a).
\end{align*}
If $a<0$, we again apply Lemma~\ref{supp:lem:bounds_cdf_of_u_km} again to the remaining integral and obtain 
\begin{align*}
&A_{km}(\xb;\pt;\xi)\\
&\le \frac{(\maximal_{\varrho(\frac{\tau_m}{m})}{\pt}(\xb))^k}{k!}\Bigl(\tau_m^k\xi(\tau_m^a) - \int_0^{\tau_m} u^k\diff\xi(u^a)\Bigr).\qedhere
\end{align*}
\end{proof}
\begin{lemma}\label{supp:lem:bkm1:consistency}
If $b\le 0$, we have
\begin{align*} B_{km}^{(1)}(\xb;\pt;\xi) \lesssim 1.
\end{align*}
If $b>0$ and $\int_0^\infty e^{-t}\xi(t^b)\diff t<\infty$, then for any $0<D<1$ and $r=\rvol(\frac{\nu_m}{m})$, we have
\begin{align*} B_{km}^{(1)}(\xb;\pt;\xi)
&\lesssim_{k,D}
\xi(\nu_m^b) e^{-D\nu_m(\minimal_{r}p(\xb))}+\xi((D\minimal_r{\pt}(\xb))^{-b}).
\end{align*}
\end{lemma}
\begin{proof}
By definition, if $b\le 0$, we have
\begin{align*}
B_{km}^{(1)}(\xb;\pt;\xi)
&=\int_1^{\nu_m}\xi(u^b)\ft_{km}(u|\xb)\diff u\\
&\le \xi(1)\int_1^{\nu_m} \ft_{km}(u|\xb)\diff u\\
&\le \xi(1).
\end{align*}
We now assume $b>0$.
Integrating by parts, we have
\begin{align*}
&B_{km}^{(1)}(\xb;\pt;\xi)\\
&=-\int_1^{\nu_m}\xi(u^b)\diff(1-\Ft_{km}(u|\xb))\\
&\le \xi(1)(1-\Ft_{km}(u|\xb))
+\int_1^{\nu_m} (1-\Ft_{km}(u|\xb))\diff\xi(u^b).
\end{align*}
Applying Lemma~\ref{supp:lem:bounds_ccdf_of_u_km1} yields, for any $0<D<1$, that
\begin{align*}
&B_{km}^{(1)}(\xb;\pt;\xi)\\
&\le \xi(1) + (1-D)^{-k+1}\int_1^{\nu_m} e^{-Dm\Pt_m(u|\xb)}\diff\xi(u^b)\\
&\le \xi(1) + (1-D)^{-k+1}\int_1^{\nu_m} e^{-Du(\minimal_r{\pt}(\xb))}\diff\xi(u^b)\numberthis\label{supp:eq:lem:bkm1:consistency:temp1}
\end{align*}
for $r=\rvol(\frac{\nu_m}{m})$.
Integrating by parts again, we thus obtain
\begin{align*}
&\int_1^{\nu_m} e^{-Du(\minimal_r{\pt}(\xb))}\diff\xi(u^b)\\
&\le \xi(\nu_m^b)e^{-D\nu_m(\minimal_r{\pt}(\xb))} 
\\&\qquad
+ D(\minimal_r{\pt}(\xb))\int_1^{\nu_m} e^{-Du(\minimal_r{\pt}(\xb))}\xi(u^b)\diff u\\
&\le \xi(\nu_m^b)e^{-D\nu_m(\minimal_r{\pt}(\xb))} 
\\&\qquad
+ \int_{D(\minimal_r{\pt}(\xb))}^{D\nu_m(\minimal_r{\pt}(\xb))} e^{-t}\xi(t^b(D\minimal_r{\pt}(\xb))^{-b})\diff t.
\numberthis\label{supp:eq:lem:bkm1:consistency:temp2}
\end{align*}
Here, using the property that $\xi(xy)\le \xi(x)\xi(y)$ for any $x,y>t_0$ for some $t_0\ge 0$, it is easy to show that 
\begin{align*}
&\int_{D(\minimal_r{\pt}(\xb))}^{D\nu_m(\minimal_r{\pt}(\xb))} e^{-t}\xi(t^b(D\minimal_r{\pt}(\xb))^{-b})\diff t\\
&\le \Bigl(t_0\xi(t_0^b)+\int_0^{\infty}e^{-t}\xi(t^b)\diff t\Bigr)\xi((D\minimal_r{\pt}(\xb))^{-b}) 
\\&\qquad
+ \xi(t_0)\int_0^{\infty} e^{-t}\xi(t^b)\diff\xb\\
&\lesssim 1+\xi((D\minimal_r{\pt}(\xb))^{-b}).
\numberthis\label{supp:eq:lem:bkm1:consistency:temp3}
\end{align*}
Putting \eqref{supp:eq:lem:bkm1:consistency:temp1}, \eqref{supp:eq:lem:bkm1:consistency:temp2}, and \eqref{supp:eq:lem:bkm1:consistency:temp3} together, we obtain the desired bound.
\end{proof}
\begin{lemma}\label{supp:lem:bkm2:consistency}
For any $0<D<1$ and $r=\rvol(\frac{\nu_m}{m})$, we have
\begin{align*}
B_{km}^{(2)}(\xb;\pt;\xi)
&\lesssim_{k,D} \xi(\nu_m^b\vee \kappa_m^b) e^{-D\nu_m(\minimal_{r}\pt(\xb))}.
\end{align*}
\end{lemma}
\begin{proof}
Integrating by parts, we have
\begin{align*}
&B_{km}^{(2)}(\xb;\pt;\xi)\\
&= -\int_{\nu_m}^{\kappa_m}\xi(u^b)\diff(1-\Ft_{km}(u|\xb))\\
&\le \xi(\nu_m^b)(1-\Ft_{km}(\nu_m|\xb)) + \int_{\nu_m}^{\kappa_m} (1-F_{km}(u|\xb))\diff\xi(u^b)\\
&\le 2\xi(\nu_m^b\vee \kappa_m^b)(1-\Ft_{km}(\nu_m|\xb)).\numberthis\label{supp:eq:lem:bkm2:consistency:intermediate}
\end{align*}
Applying Lemma~\ref{supp:lem:bounds_ccdf_of_u_km1}, we have that for any $0<D<1$ and $r=\rvol(\frac{\nu_m}{m})$
\begin{align*}
&B_{km}^{(2)}(\xb;\pt;\xi)\\
&\le 2(1-D)^{-k+1}\xi(\nu_m^b\vee\kappa_m^b)
e^{-Dm\Pt_m(\nu_m|\xb)}\\
&\le 2(1-D)^{-k+1}\xi((\nu_m^b\vee\kappa_m^b)
e^{-D\nu_m(\minimal_r{p}(\xb))}.
\qedhere
\end{align*}
\end{proof}
\begin{lemma}\label{supp:lem:bkm3}
For any $\d>0$ and $m$ sufficiently large, we have
\begin{align*}
&B_{km}^{(3)}(\xb;\pt;\xi)\\
&\lesssim_{\d} \xi(m^b)
\int p(\yb)\xi(\ups^{b}(\rho(\xb,\yb)))\ones_{\{\rho(\xb,\yb)>\rvol(\frac{\kappa_m}{m})\}}\diff\yb.
\end{align*}
\end{lemma}
\begin{proof}
We recall the following bound~\eqref{supp:eq:lower_bound_cdf_of_u_km_2} on the complementary cdf $1-\Ft_{km}(u|\xb)$ from Lemma~\ref{supp:lem:bounds_cdf_of_u_km}: for any $\d>0$ and $m\ge (1+1/\d)(k-1)$, we have
\begin{align*}
1-\Ft_{km}(u|\xb)
&\le (1+\d)(1-\Pt_m(u|\xb))\\
&=(1+\d)\int \pt(\yb)\ones_{\{\rho(\xb,\yb)>\rvol(\frac{u}{m})\}}\diff\yb.
\end{align*}
Integrating by parts, we first obtain 
\begin{align*}
&B_{km}^{(3)}(\xb;\pt;\xi)\\
&=-\int_{\kappa_m}^{\infty} \xi(u^b)\diff(1-\Ft_{km}(u|\xb))\\
&\le \xi(\kappa_m^b)(1-\Ft_{km}(\kappa_m|\xb)) 
\\&\qquad
+ \int_{\kappa_m}^{\infty}(1-\Ft_m(u|\xb)) \diff\xi(u^b)\\
&\le \xi(\kappa_m^b)(1-\Ft_{km}(\kappa_m|\xb)) 
\\&\qquad
+ (1+\d)\int_{\kappa_m}^{\infty}(1-\Pt_m(u|\xb)) \diff\xi(u^b).\numberthis\label{supp:eq:proof:lem:bkm3:temp1}
\end{align*}
Integrating the second term by parts leads to
\begin{align*}
&\int_{\kappa_m}^{\infty}(1-\Pt_m(u|\xb)) \diff\xi(u^b)
\numberthis\label{supp:eq:proof:lem:bkm3:temp2}\\
&\le \lim_{u\to\infty} \xi(u^b) (1-\Pt_m(u|\xb)) 
+ \int_{\kappa_m}^{\infty} \xi(u^b)\diff\Pt_m(u|\xb).
\end{align*}
For the first term in \eqref{supp:eq:proof:lem:bkm3:temp2}, since for $m$ sufficiently large with $m^b> t_0$ and $(\kappa_m/m)^b>t_0$, we have $\xi(u^b)\le \xi(m^b)\xi((u/m)^b)$ for $u\ge \kappa_m$, it follows that
\begin{align*}
&\xi(u^b)(1-\Pt_m(u|\xb))\\
&=\xi(u^b)\int \pt(\yb)\ones_{\{\rho(\xb,\yb)>\rvol(\frac{u}{m})\}}\diff\yb\\
&\le\xi(m^b) \int \pt(\yb)\xi\bigl(\bigl(\frac{u}{m}\bigr)^b\bigr)\ones_{\{\rho(\xb,\yb)>\rvol(\frac{u}{m})\}}\diff\yb\\
&\le\xi(m^b) \int \pt(\yb)\xi(\ups^b(\rho(\xb,\yb)))\ones_{\{\rho(\xb,\yb)>\rvol(\frac{u}{m})\}}\diff\yb\\
&\le\xi(m^b) \int \pt(\yb)\xi(\ups^b(\rho(\xb,\yb)))\ones_{\{\rho(\xb,\yb)>\rvol(\frac{\kappa_m}{m})\}}\diff\yb.
\end{align*}
Therefore, 
\begin{align*}
&\lim_{u\to\infty} \xi(u^b)(1-\Pt_m(u|\xb)) 
\numberthis\label{supp:eq:proof:lem:bkm3:temp3}
\\
&\le \xi(m^b) \int \pt(\yb)\xi(\ups^b(\rho(\xb,\yb)))\ones_{\{\rho(\xb,\yb)>\rvol(\frac{\kappa_m}{m})\}}\diff\yb.
\end{align*}
The second term in \eqref{supp:eq:proof:lem:bkm3:temp2} can be bounded similarly as
\begin{align*}
&\int_{\kappa_m}^{\infty} \xi(u^b)\diff \Pt_m(u|\xb)
\numberthis\label{supp:eq:proof:lem:bkm3:temp4}\\
&= \int \pt(\yb)\xi((m\ups(\rho(\xb,\yb)))^b)\ones_{\{\rho(\xb,\yb)>\rvol(\frac{\kappa_m}{m})\}}\diff\yb\\
&\le \xi(m^b)
\int \pt(\yb)\xi(\ups^{b}(\rho(\xb,\yb)))\ones_{\{\rho(\xb,\yb)>\rvol(\frac{\kappa_m}{m})\}}\diff\yb.
\end{align*}
Plugging \eqref{supp:eq:proof:lem:bkm3:temp2}, \eqref{supp:eq:proof:lem:bkm3:temp3}, and \eqref{supp:eq:proof:lem:bkm3:temp4} into \eqref{supp:eq:proof:lem:bkm3:temp1} establishes the desired bound.
\end{proof}

The following is the key lemma in establishing vanishing bias and vanishing variance for single- and double-density cases.

\begin{lemma}\label{supp:lem:generic:boundedness}
Assume that $-\int_0^1 u^k\diff\xi(u^{a\wedge 0})<\infty$ and $\int_0^{\infty} e^{-t}\xi(t^{b\vee 0})\diff t<\infty$.
If the densities $p$ and $\pt$ satisfy $\Pt\ll\P$, \ref{cond:low}, and \ref{cond:high}, 
we have
\begin{align*}
\limsup_{m\to\infty} \int p(\xb) \int_{0}^{\infty} \xi(\psi_{a,b}(u))\diff \Ft_{km}(u|\xb)\diff\xb
<\infty.
\end{align*}
\end{lemma}
\begin{proof}
Let $\tau_m=1$ and $\kappa_m=e^{o(m)}$.
Then, there exists $\nu_m$ such that $\nu_m\to\infty$, $\nu_m/m \to 0$, and for any $c>0$, $e^{-c\nu_m}\xi(\kappa_m^b)\to 0$, as $m\to\infty$.
Consider
\begin{align*}
&\int_{0}^{\infty} \xi(\psi_{a,b}(u))\diff F_{km}(u|\xb)\\
&=\int_0^1\xi(u^a)\diff\Ft_{km}(u|\xb)+\int_1^{\infty}\xi(u^b)\diff\Ft_{km}(u|\xb)\\
&=A_{km}(\xb;\pt;\xi) + B_{km}(\xb;\pt;\xi),
\end{align*}
where 
\[
B_{km}(\xb;\pt;\xi)\defeq
B_{km}^{(1)}(\xb;\pt;\xi)
+B_{km}^{(2)}(\xb;\pt;\xi)
+B_{km}^{(3)}(\xb;\pt;\xi).
\]
Recall the definitions of $A_{km}(\xb;\pt;\xi)$, $B_{km}^{(1)}(\xb;\pt;\xi)$, $B_{km}^{(2)}(\xb;\pt;\xi)$, and $B_{km}^{(3)}(\xb;\pt;\xi)$ in \eqref{supp:eq:akm}, \eqref{supp:eq:bkm1}, \eqref{supp:eq:bkm2}, and \eqref{supp:eq:bkm3}, respectively.
Letting
\[
A_{km}(p,\pt;\xi)\defeq \int p(\xb)A_{km}(\xb;\pt;\xi)\diff\xb
\]
and
\[
B_{km}(p,\pt;\xi)\defeq \int p(\xb)B_{km}(\xb;\pt;\xi)\diff\xb,
\]
we show separately that $\limsup_{m\to\infty}A_{km}(p,\pt;\xi)<\infty$ and $\limsup_{m\to\infty}B_{km}(p,\pt;\xi)<\infty$.

\textbf{Step 1. Bounding $A_{km}(p,\pt;\xi)$.} 
If $a\ge 0$, we trivially have 
$A_{km}(p,\pt;\xi)
\le \xi(1)$.
If $a<0$, by Lemma~\ref{supp:lem:akm}, we have
\begin{align*}
A_{km}(p,\pt;\xi)
&\le \frac{\maxf(p,\pt;k,\rvol(\frac{1}{m}))}{k!}
\Bigl(\xi(1) - \int_0^{1} u^k\diff\xi(u^a)\Bigr)\\
&\lesssim_k \maxf\bigl(p,\pt;k,\rvol\bigl(\frac{1}{m}\bigr)\bigr).
\end{align*}
Hence, since there exists $r'>0$ such that $\maxf(p,\pt;k,r')<\infty$ by the the condition~\ref{cond:low} and $\maxf(p,\pt;k,r)$ is nonincreasing as $r\to 0$, we conclude that
$A_{km}(p,\pt;\xi)< \infty$ for $m$ sufficiently large such that $\rvol(1/m)<r'$.

\textbf{Step 2. Bounding $B_{km}(p,\pt;\xi)$.}
If $b\le 0$, then we trivially have $B_{km}(p,\pt;\xi)\le \xi(1)$.
If $b>0$, by applying Lemmas~\ref{supp:lem:bkm1:consistency}, \ref{supp:lem:bkm2:consistency}, and \ref{supp:lem:bkm3}, we have that for any $0<D<1$ and $m$ sufficiently large
\begin{align*}
B_{km}(p,\pt;\xi)
&\lesssim \xi(\kappa_m^b)\int e^{-D\nu_m(\minimal_{r_1}{\pt}(\xb))} p(\xb)\diff\xb 
\\&\qquad
+ \minf(p,\pt;\xi,b,r_1)
\\&\qquad
+ \xi(m^b)\expvolfunc(p,\pt;\xi,b,r_2),
\end{align*}
where $r_1=\rvol(\nu_m/m)$ and $r_2=\rvol(\kappa_m/m)$.
\begin{itemize}
    \item For the first term, since $\P\ll\Pt$ implies that $\P(\{\xb\suchthat \minimal_{r_1}{\pt}(\xb)>0\})=1$ (Lemma~\ref{supp:lem:consistency_measure_zero}), we have $\xi(\kappa_m^b)e^{-\nu_m(\minimal_{r_1}{\pt}(\xb))}\to 0$ as $m\to\infty$ for $\P$-a.e. $\xb$ by definition of $\nu_m$ and $\kappa_m$. 
    Therefore, by the dominated convergence theorem,
    \[
    \lim_{m\to\infty} \int \xi(\kappa_m^b) e^{-\nu_m(\minimal_{r_1}{\pt}(\xb))} p(\xb)\diff\xb = 0.
    \]
    \item Since there exists $r''>0$ such that $\minf(p,\pt;\xi,b,r'')<\infty$ by the condition~\ref{cond:high} and $\minf(p,\pt;\xi,b,r)$ is nonincreasing as $r\to 0$, the second term is bounded for $m$ sufficiently large such that $\rvol(\frac{\kappa_m}{m})<r''$.
    \item The limit superior of the last term $\xi(m^b)R(p,\pt;\xi,b,r_2)$ as $m\to\infty$ is bounded by the condition~\ref{cond:high}.
\end{itemize}
Overall, we conclude that
\[
\limsup_{m\to\infty} B_{km}(p,\pt;\xi) <\infty.\qedhere
\]
\end{proof}

Following the proof of Lemma~\ref{supp:lem:generic:boundedness} with the stronger assumptions establishes the following bound.
\begin{lemma}\label{supp:lem:bound_ukm_integral}
Assume that $-\int_0^1 u^k\diff\xi(u^{a\wedge 0})<\infty$ and $\int_0^{\infty} e^{-t}\xi(t^{b\vee 0})\diff t<\infty$.
If $\pt$ satisfies the conditions \ref{cond:bounded_above}, \ref{cond:bounded_below}, \ref{cond:bounded_support}, and \ref{cond:regular_support},
we have
\begin{align*}
\int_{0}^{\infty} \xi(\psi_{a,b}(u))\diff \Ft_{km}(u|\xb)
\lesssim 1
\end{align*}
for $\P$-a.e.\@ $\xb$.
\end{lemma}

Continuing from \eqref{supp:eq:lem:bkm2:consistency:intermediate} and applying \eqref{supp:eq:lower_bound_cdf_of_u_km_12} in Lemma~\ref{supp:lem:bounds_ccdf_of_u_km1} yield the following bound, which is required for establishing performance guarantees with adaptive choices of $k$ and $l$.
\begin{lemma}\label{supp:lem:bkm2:rate:bounded}
For $r=\rvol(\nu_m/m)$, we have
\begin{align*}
B_{km}^{(2)}(\xb;\pt;\xi)
&\le 2\xi(\nu_m^b\vee \kappa_m^b) \Bigl(\frac{e\nu_m\maximal_{r}{\pt}(\xb)}{k}\Bigr)^ke^{-\nu_m \minimal_{r}{\pt}(\xb)}.
\end{align*}
\end{lemma}

\begin{lemma}\label{supp:lem:bound_gamma_integral}
If $k+a>0$, for $\xb\in\supp(p)$, we have
\begin{align*}
&\int_0^\infty
	\psi_{a,b}(u)\rho_{U_{k\infty}(\xb)}(u)\diff u\\
&\le \frac{p^k(\xb)}{(k+a)\Gamma(k)} +
\frac{\Gamma((k+b)\vee 1)}{\Gamma(k)}(p(\xb))^{(k-1)\wedge (-b)}.
\end{align*}
In particular, if $\Cbe_p\le p(\xb)\le\Cab_p$, then 
\begin{align*}
\int_0^\infty
	\psi_{a,b}(u)\rho_{U_{k\infty}(\xb)}(u)\diff u
\lesssim 1.
\end{align*}
\end{lemma}
\begin{proof}
First, consider
\begin{align*}
\int_0^1 u^a \rho_{U_{k\infty}(\xb)}(u)\diff u
&= \frac{p^k(\xb)}{\Gamma(k)}\int_0^1 u^{k+a-1}e^{-up(\xb)}\diff u\\
&= \frac{(p(\xb))^{-a}}{\Gamma(k)}\int_0^{p(\xb)} t^{k+a-1}e^{-t}\diff t\\
&\le \frac{p^k(\xb)}{(k+a)\Gamma(k)},
\end{align*}
where the last inequality follows from the bound on the lower incomplete gamma in Lemma~\ref{supp:lem:incomplete_gamma}.
Similarly, we consider
\begin{align*}
\int_0^1 u^b \rho_{U_{k\infty}(\xb)}(u)\diff u
&= \frac{(p(\xb))^{-b}}{\Gamma(k)}\int_0^{p(\xb)} t^{k+b-1}e^{-t}\diff t.
\end{align*}
On the one hand, if $k+b>1$, by bounding the integral by $\Gamma(k+b)$, we have
\[
\int_0^1 u^b \rho_{U_{k\infty}(\xb)}(u)\diff u
\le \frac{\Gamma(k+b)}{\Gamma(k)}(p(\xb))^{-b}.
\]
On the other hand, if $k+b\le 1$, we have
\begin{align*}
\int_0^1 u^b \rho_{U_{k\infty}(\xb)}(u)\diff u
&\le \frac{(p(\xb))^{k-1}}{\Gamma(k)} \int_{p(\xb)}^{\infty}e^{-t}\diff t \\
&\le \frac{(p(\xb))^{k-1}}{\Gamma(k)}.
\end{align*}
Therefore, we obtain
\[
\int_0^1 u^b \rho_{U_{k\infty}(\xb)}(u)\diff u
\le \frac{\Gamma((k+b)\vee 1)}{\Gamma(k)}(p(\xb))^{(k-1)\wedge (-b)},
\]
which completes the proof.
\end{proof}

\subsection{Generic bias bounds}
\label{supp:sec:generic_bias_bounds}

\begin{lemma}[Generic inner bias bound]\label{supp:lem:generic_inner_bias:single_case1}
Suppose that the density $p$ satisfies the conditions \ref{cond:bounded_above}, \ref{cond:smoothness_int}, and \ref{cond:boundary}, and let $k=o(\sqrt{m})$ as $m\to\infty$.
\begin{enumerate}
\item We have
\begin{align}\label{supp:eq:inner_first}
I_{\textin,1}
&= \bigO\Bigl(
\frac{\tau_m^{(a+\frac{\sigma_p}{d}+1)\wedge 0}}{m^{\frac{\sigma_p}{d}}}
+\frac{k^{-k}}{m}
\nonumber\\
&\quad\qquad
+\Bigl(\frac{1}{m}\Bigr)^{\frac{1}{d}}\tau_m^{(a+1)\wedge 0}
\Bigr).
\end{align}
\item If $\nu_m=o(\sqrt{m})$ as $m\to\infty$, we have
\begin{align}
\label{supp:eq:inner_second}
I_{\textin,2}
&= \bigO\Bigl(
\frac{\nu_m^{(b+\frac{\sigma_p}{d}+2)\vee 0}}{m^{\frac{\sigma_p}{d}}}
+\frac{k^{-k}\nu_m^{(b+k+2)\vee 0}}{m} \nonumber\\
&\quad\qquad +\Bigl(\frac{\nu_m}{m}\Bigr)^{\frac{1}{d}}\nu_m^{(b+2)\vee 0}
\Bigr).
\end{align}
\end{enumerate}
\end{lemma}
\begin{proof}
We establish each bound separately. 

\textbf{Bounding the lower inner bias $I_{\mathrm{in},1}$.}
For each $r>0$, define a set 
\begin{align*}
S_p(r)
&\defeq \{\xb\in\supp(p)\suchthat \text{$p$ is $\sigma_p$-H\"older continuous}\\
&\qquad\qquad\qquad\qquad
\text{over $\Bb(\xb,r)$}\}.
\end{align*}
By the smoothness assumption~\ref{cond:smoothness_int}, we can bound the inner bias incurred at the ``smooth region'', \ie
\[
I_{\textin,1,\text{smooth}}=\int_{S_p(\rvol(\frac{1}{m}))}I_{\textin,1}(\xb)p(\xb)\diff\xb,
\]
by applying Lemma~\ref{supp:lem:GOVLemma2_pdf_gap_bound}.
Since $p(\xb)\le C_p<\infty$ for $\P$-a.e.\@ $\xb$, this lemma holds for $m$ sufficiently large uniformly over $\P$-a.e.\@ $\xb$.
Applying Lemma~\ref{supp:lem:GOVLemma2_pdf_gap_bound} for $\xb\in S_p(\rvol(\frac{1}{m}))$, we have
\begin{align*}\numberthis\label{supp:eq:inner_first_smooth_pointwise}
&I_{\textin,1}(\xb)\\
&\lesssim_{\sigma_p,L,\Cab_p,C_0,d}
\int_{\tau_m}^1
u^a\Bigl\{(1+u)\Bigl(\frac{u}{m}\Bigr)^{\frac{\sigma_p}{d}}
\\&\qquad\qquad\qquad\qquad\quad
+k^{-k}\frac{(k^2+u^{2})u^{k-1}e^{-up(\xb)}}{m}\Bigr\}\diff u.
\end{align*}
It is easy to see that the first term is bounded by $\bigO\bigl(\tau_m^{(a+\frac{\sigma_p}{d}+1)\wedge 0}m^{-\frac{\sigma_p}{d}}\bigr).$\footnote{Here $a+\frac{\sigma_p}{d}+1\neq 0$ is implicitly assumed. If $a+\frac{\sigma_p}{d}+1=0$, then the first term behaves as $O((\ln \tau_m)m^{-\frac{\sigma_p}{d}})$} 
To bound the second term, we use the upper bound on the lower incomplete gamma function (Lemma~\ref{supp:lem:incomplete_gamma}). Since we always assume that $k+a>0$, we have
\begin{align*}
&\int_{\tau_m}^1\frac{k^{-k}}{m}(k^2 + u^2)u^{k+a-1}e^{-up(\xb)}\diff u\\
&\le \frac{k^{-k}}{m}\bigl\{k^2p(\xb)^{-(k+a)}\gamma(k+a, p(\xb))
\\&\qquad
+ p(\xb)^{-(k+a+2)}\gamma(k+a+2, p(\xb))\bigr\}
=\bigO\Bigl(\frac{k^{-k}}{m}\Bigr).
\end{align*}
Hence, we conclude that
\begin{align*}
&I_{\textin,1,\text{smooth}}
\numberthis\label{supp:eq:inner_first_smooth}
\\
&=\bigO(\tau_m^{(a+\frac{\sigma_p}{d}+1)\wedge 0}m^{-\frac{\sigma_p}{d}} + k^{-k}m^{-1}).
\end{align*}
To control the inner bias incurred at $\xb \in\supp(p)\backslash S_p(\rvol(m^{-1}))$, \ie
\begin{align*}
I_{\textin,1,\text{nonsmooth}}
&=\int_{\supp(p)\backslash S_p(\rvol(\frac{1}{m}))}I_{\textin,1}(\xb)p(\xb)\diff\xb,
\end{align*}
we first note that the bound~\eqref{supp:eq:inner_first_smooth_pointwise} on $I_{\textin,1}(\xb)$ holds with $\sigma_p=0$ from the upper boundedness assumption~\ref{cond:bounded_above}, which implies that
\begin{align*}
&I_{\textin,1,\text{nonsmooth}}
\numberthis\label{supp:eq:inner_first_nonsmooth}
\\
&=\bigO(
\Leb(\supp(p)\backslash S_p(\rvol(m^{-1})))(\tau_m^{(a+1)\wedge 0}+k^{-k}m^{-1})
).
\end{align*}
We now only need to bound the Lebesgue measure of the set where $\supp(p)\backslash S_p(\rvol(m^{-1}))$.
Observe that for any $r>0$
\[
\supp(p)\backslash S_p(r)
\subseteq \{\xb\in\Real^d\suchthat \Bb(\xb,r)\cap \partial(\supp(p))\neq \emptyset\},
\]
where $\partial A$ denotes the boundary of a set $A$. 
Using the following lemma with the condition~\ref{cond:boundary} on the finiteness of the Hausdorff measure of the boundary of the support, we can bound the Lebesgue measure of 
$\Real^d\backslash S_p(\rvol(m^{-1}))$ by $\bigO(\rvol(1/m))=\bigO(m^{-\frac{1}{d}})$.

\begin{lemma}[{\cite[Section~A]{Gao--Oh--Viswanath2018tit}}]\label{supp:lem:hausdorff}
For $S\subset \Real^d$, suppose that $0<H^{d-1}(S)<\infty$.
Let $T(r)\defeq \{\xb\in\Real^d\suchthat \Bb(\xb,r)\cap S\neq \varnothing\}$ for $r>0$.
Then $\Leb(T(r))=2rH^{d-1}(S)+o(r)$ for $r$ sufficiently small.
\end{lemma}
Combining \eqref{supp:eq:inner_first_smooth} and \eqref{supp:eq:inner_first_nonsmooth} establishes the desired bound~\eqref{supp:eq:inner_first}. \medskip

\textbf{Bounding the upper inner bias $I_{\mathrm{in},2}$.}
The proof follows a similar line of argument as that of \eqref{supp:eq:inner_first}. 
We first apply Lemma~\ref{supp:lem:GOVLemma2_pdf_gap_bound} for $\xb\in S_p(\rvol(\frac{\nu_m}{m}))$ and obtain
\begin{align*}
&I_{\textin,2}(\xb)\\
&\lesssim_{\sigma_p,L,\Cab_p,C_0,d}
\int_{1}^{\nu_m}
u^b\Bigl\{(1+u)\Bigl(\frac{u}{m}\Bigr)^{\frac{\sigma_p}{d}}
\\&\qquad\qquad\qquad\qquad\qquad
+ k^{-k}\frac{(k^2+u^{2})u^{k-1}e^{-up(\xb)}}{m}\Bigr\}\diff u.
\end{align*}
The first term is bounded by $\bigO(m^{-\frac{\sigma_p}{d}}\nu_m^{(b+\frac{\sigma_p}{d}+2)\vee 0}).$
The second term is again bounded by the upper bound on the lower incomplete gamma function. If $b+k>0$, we have
\begin{align*}
&\int_1^{\nu_m}\frac{k^{-k}}{m}(k^2 + u^2)u^{b+k-1}e^{-up(\xb)}\diff u\\
&\le \frac{k^{-k}}{m}(k^2p^{-(b+k)}(\xb)\gamma(b+k, \nu_m p(\xb))
\\&\qquad\qquad
+ p^{-(b+k+2)}(\xb)\gamma(b+k+2, \nu_m p(\xb)))\\
&= \bigO\Bigl(
k^{-k}\frac{(k^2\nu_m^{(b+k)\vee 0} + \nu_m^{(b+k+2)\vee 0})}{m}
\Bigr)\\
&= \bigO\Bigl(
k^{-k}\frac{\nu_m^{(b+k+2)\vee 0}}{m}
\Bigr).
\end{align*} 
One can easily show that the bound also holds when $b+k\le 0$. 
Hence, we conclude that
\begin{align*}
&I_{\textin,2,\text{smooth}}
\numberthis\label{supp:eq:inner_second_smooth}
\\
&=\int_{S_p(\rvol(\frac{\nu_m}{m}))}I_{\textin,2}(\xb)p(\xb)\diff\xb\\
&=\bigO(m^{-\frac{\sigma_p}{d}}\nu_m^{(b+\frac{\sigma_p}{d}+2)\vee 0} + m^{-1}\nu_m^{(b+k+2)\vee 0}).
\end{align*}
Similar to \eqref{supp:eq:inner_first_nonsmooth}, we have 
\begin{align*}
&I_{\textin,2,\text{nonsmooth}}
\numberthis\label{supp:eq:inner_second_nonsmooth}\\
&=\int_{\supp(p)\backslash S_p(\rvol(\frac{\nu_m}{m}))}I_{\textin,2}(\xb)p(\xb)\diff\xb,
\\&
=\bigO((\nu_m/m)^{\frac{1}{d}} (\nu_m^{(b+2)\vee 0}+m^{-1}\nu_m^{(b+k+2)\vee 0})
),
\end{align*}
since $\Leb(\supp(p)\backslash S_p(\rvol(\nu_m/m)))=O(\rvol(\nu_m/m))=O((\nu_m/m)^{\frac{1}{d}})$ by Lemma~\ref{supp:lem:hausdorff}.
Putting \eqref{supp:eq:inner_second_smooth} and \eqref{supp:eq:inner_second_nonsmooth} together establishes the desired bound~\eqref{supp:eq:inner_second}.
\end{proof}

\begin{lemma}[Generic outer bias bound]
\label{supp:lem:generic_outer_bias_case1}
Suppose that the density $p$ satisfies \ref{cond:bounded_above}.
\begin{enumerate}
\item If $k>-a$, we have
\begin{align}\label{supp:eq:outer_first}
I_{\textout,1}
&= \bigO\bigl(
k^{-k}\tau_m^{k+a}
\bigr).
\end{align}
\item If $p$ satisfies \ref{cond:bounded_below}, \ref{cond:bounded_support}, and \ref{cond:regular_support}, 
then, for $m$ sufficiently large, we have
\begin{align*}
I_{\textout,2}
&=\bigO\bigl(
k^{b}\nu_m^{b+k-1}e^{-\Cbe_p\nu_m}
\\&\qquad~
+(\nu_m^b\vee \kappa_m^b)\bigl(\frac{\nu_m}{k}\bigr)^ke^{-\eta_p\Cbe_p\nu_m}\bigr).
\numberthis
\label{supp:eq:outer_second}
\end{align*}
\end{enumerate}
\end{lemma}
\begin{proof}
Recall that
\[
\rho_{U_{k\infty}(\xb)}(u) = \frac{p^k(\xb)}{\Gamma(k)}u^{k-1}e^{-up(\xb)}.
\]
Define
\begin{align*}
A_{k\infty}(\xb;p)&\defeq \int_0^{\tau_m} u^a \rho_{U_{k\infty}(\xb)}(u)\diff u
\end{align*}
and
\begin{align*}
B_{k\infty}(\xb;p)&\defeq \int_{\nu_m}^{\infty} u^b\rho_{U_{k\infty}(\xb)}(u)\diff u.
\end{align*}
For some $\kappa_m=\omega(m)$ such that $\kappa_m\ge \nu_m$, we also let $A_{km}(\xb;p)\defeq A_{km}(\xb;p;\xi)$, $B_{km}^{(2)}(\xb;p)\defeq B_{km}^{(2)}(\xb;p;\xi)$, and $B_{km}^{(3)}(\xb;p)\defeq B_{km}^{(3)}(\xb;p;\xi)$ for $\xi(t)=t$; recall the definitions in Appendix~\ref{supp:sec:technical_lemmas:generic}.
Now we can write the lower outer bias as
\[
I_{\textout,1}=\int p(\xb) (A_{km}(\xb;p)+A_{k\infty}(\xb;p))\diff\xb
\]
and the upper outer bias as
\[
I_{\textout,2}=\int p(\xb) (B_{km}^{(2)}(\xb;p)+B_{km}^{(3)}(\xb;p)+B_{k\infty}(\xb;p))\diff\xb
\]


\textbf{Bounding the lower outer bias $I_{\textout,1}$.}
On the one hand, by invoking the lower incomplete gamma function in Lemma~\ref{supp:lem:incomplete_gamma}, we obtain
\begin{align*}
A_{k\infty}(\xb;p)
&= \frac{p^k(\xb)}{\Gamma(k)}\int_0^{\tau_m} u^{k+a-1}e^{-up(\xb)}\diff u\\
&= \frac{p^{-a}(\xb)}{\Gamma(k)}\gamma(k+a, \tau_m p(\xb)) \\
&\le \frac{p^k(\xb)\tau_m^{k+a}}{\Gamma(k)(k+a)}\\
&\le \frac{\Cab_p^k\tau_m^{k+a}}{\Gamma(k)(k+a)}
= \bigO(k^{-k}\tau_m^{k+a}).
\end{align*}
On the other hand, by applying Lemma~\ref{supp:lem:akm} with the upper boundedness condition~\ref{cond:bounded_above}, we obtain
\begin{align*}
\int p(\xb) A_{km}(\xb;p)\diff\xb
&\le 
\frac{\Cab_p^k\tau_m^{k+a}}{k!}\Bigl(1\vee \frac{k}{k+a}\Bigr)\\
&=\bigO(k^{-k}\tau_m^{k+a}).
\end{align*}
Combining the two bounds, we conclude that $I_{\textout,1}=\bigO(k^{-k}\tau_m^{k+a})$.

\textbf{Bounding the upper outer bias $I_{\textout,2}$.}
For the $B_{k\infty}(\xb;p)$ term in the upper outer bias $I_{\textout,2}$, we apply the bound~\eqref{supp:eq:incomplete_gamma_upper} on the upper incomplete gamma function in Lemma~\ref{supp:lem:incomplete_gamma}. 
Consider
\begin{align*}
B_{k\infty}(\xb;p)
&= \frac{p^k(\xb)}{\Gamma(k)}
\int_{\nu_m}^\infty u^{k+b-1}e^{-up(\xb)}\diff u\\
&= \frac{p^{-b}(\xb)}{\Gamma(k)} \int_{\nu_m p(\xb)}^\infty t^{k+b-1}e^{-t}\diff t.
\end{align*}
If $\nu_m p(\xb)<1$, we have
\begin{align*}
B_{k\infty}(\xb;p)
&\le \frac{p^{-b}(\xb)}{\Gamma(k)} \int_0^{\infty} t^{k+b-1} e^{-t}\diff t\\
&\le \frac{\Gamma((k+b)\vee 1)}{\Gamma(k)} p^{-b}(\xb).
\end{align*}
We now assume that $\nu_m p(\xb)\ge 1$.
If $k+b\ge 1$, we have
\begin{align*}
B_{k\infty}(\xb;p)&= \frac{p^{-b}(\xb)}{\Gamma(k)}\Gamma(k+b,\nu_m p(\xb))\\
&\le \frac{p^{-b}(\xb)}{\Gamma(k)}\Gamma(k+b)(\nu_m p(\xb))^{k+b-1} e^{-\nu_m p(\xb)+1}\\
&= \frac{\Gamma(k+b)}{\Gamma(k)}\nu_m^{k+b-1}p^{k-1}(\xb) e^{-\nu_m p(\xb)+1},
\end{align*}
where the inequality follows from Lemma~\ref{supp:lem:incomplete_gamma}.
For $k+b<1$, a similar bound can be derived:
\begin{align*}
B_{k\infty}(\xb;p)
&= \frac{p^{-b}(\xb)}{\Gamma(k)} (\nu_m p(\xb))^{k+b-1}\int_{\nu_m p(\xb)}^{\infty} e^{-t}\diff t\\
&=\frac{1}{\Gamma(k)}\nu_m^{k+b-1}p^{k-1}(\xb)e^{-\nu_m p(\xb)}.
\end{align*}
To sum up, we can bound $B_{k\infty}(\xb;p)$ as
\begin{align*}
&B_{k\infty}(\xb;p)\\
&\le \frac{\Gamma((k+b)\vee 1)}{\Gamma(k)} (p^{-b}(\xb)\ones_{\{\nu_m p(\xb)<1\}}
\\&\qquad\qquad\qquad\qquad
+ \nu_m^{k+b-1}p^{k-1}(\xb)e^{-\nu_m p(\xb)+1})\\
&\stackrel{(a)}{\le} \frac{\Gamma((k+b)\vee 1)}{\Gamma(k)} (p^{-b}(\xb)+ \nu_m^{k+b-1}p^{k-1}(\xb)) e^{-\nu_m p(\xb)+1}\\
&\stackrel{(b)}{\le} \frac{\Gamma((k+b)\vee 1)}{\Gamma(k)} ((\Cab_p^{-b}\vee \Cbe_p^{-b})+ \nu_m^{k+b-1}\Cab_p^{k-1}) e^{-\nu_m \Cbe_p+1}\\
&= \bigO(k^{b}\nu_m^{k+b-1}e^{-\Cbe_p\nu_m}).
\end{align*}
Here, (a) follows from the inequality $\ones_{\{t\le 1\}}\le e^{-t+1}$, and (b) follows from the boundedness conditions~\ref{cond:bounded_above} and \ref{cond:bounded_below}.
Therefore, we conclude that
\begin{align*}
\int p(\xb) B_{k\infty}(\xb;p)\diff\xb 
&= \bigO(k^{b}\nu_m^{k+b-1}e^{-\Cbe_p\nu_m}).
\numberthis\label{supp:eq:bkinfty}
\end{align*}

Next, we bound $\int p(\xb)(B_{km}^{(2)}(\xb;p)+B_{km}^{(3)}(\xb;p))\diff\xb$. 
On the one hand,
applying Lemma~\ref{supp:lem:bkm2:rate:bounded} with the upper boundedness condition~\ref{cond:bounded_above}, we first have
\begin{align*}
&\int p(\xb) B_{km}^{(2)}(\xb;p)\diff\xb\\
&\le 2(\nu_m^b\vee \kappa_m^b)\Bigl(\frac{e\Cab_p\nu_m}{k}\Bigr)^k \int p(\xb) e^{-\nu_m \minimal_{r}{p}(\xb)}\diff\xb
\end{align*}
for $r=\rvol(\frac{\nu_m}{m})$.
Further, since we have 
\[
\eta_p= \inf_{\xb\in\supp(p)}\inf_{r'\in(0,r]} \frac{\Leb(\Bb(\xb,r)\cap \supp(p))}{\Leb(\Bb(\xb,r))}>0
\]
from condition~\ref{cond:regular_support}, it follows that
$\minimal_r{p}(\xb)\ge \Cbe_p\eta_p$ for $\xb\in\supp(p)$, leading to
\begin{align*}
\int p(\xb) B_{km}^{(2)}(\xb;p)\diff\xb
&\le 2(\nu_m^b\vee \kappa_m^b)\Bigl(\frac{e\nu_m\Cab_p}{k}\Bigr)^k e^{-\eta_p\Cbe_p\nu_m}.
\end{align*}
On the other hand, since the support of the density $p$ is bounded by the condition~\ref{cond:bounded_support}, $\expvolfunc(p,p;\xi,b,\rvol(\kappa_m/m))$ becomes 0 for $m$ sufficiently large, since $\kappa_m/m\to\infty$ as $m\to\infty$.
Hence, by applying Lemma~\ref{supp:lem:bkm3} for a fixed $\d>0$, we have
\begin{align*}
&\int p(\xb) B_{km}^{(3)}(\xb;p)\diff\xb\\
&\le 3(1+\d)m^b\expvolfunc\bigl(p,p;\xi,b,\rvol\bigl(\frac{\kappa_m}{m}\bigr)\bigr)
=0
\end{align*}
for $m$ sufficiently large.
Therefore, we conclude that
\begin{align*}
&\int p(\xb)(B_{km}^{(2)}(\xb;p)+B_{km}^{(3)}(\xb;p))\diff\xb\\
&=\bigO(\nu_m^b\vee \kappa_m^b)\bigl(\frac{\nu_m}{k}\bigr)^ke^{-\eta_p\Cbe_p\nu_m}\bigr).
\numberthis\label{supp:eq:bkm23}
\end{align*}

Combining the bounds \eqref{supp:eq:bkinfty} and \eqref{supp:eq:bkm23} establishes the desired bound~\eqref{supp:eq:outer_second}.
\end{proof}

\begin{remark}
\label{supp:rem:tail_condition}
A more general condition, namely, that
\begin{enumerate}[label=\textbf{(B\arabic*$_p^\prime$)},leftmargin=*]
\setcounter{enumi}{0}
\item\label{supp:cond:exp_integral_exp}
there exists $E_0,E_1>0$ such that $\int p(\xb)e^{-\b p(\xb)}\diff\xb \le E_0e^{-E_1\b}$ for all $\b>1$,
\end{enumerate}
was originally assumed in \cite{Gao--Oh--Viswanath2018tit}. 
Known examples of densities that satisfy the condition~\ref{supp:cond:exp_integral_exp} satisfy the more intuitive condition~\ref{cond:bounded_below}.
We remark, however, that it is nontrivial to adapt the proofs in this paper to work with \ref{supp:cond:exp_integral_exp} in place of \ref{cond:bounded_below}, as the lower boundedness condition~\ref{cond:bounded_below} is explicitly utilized to remove the upper truncation of the estimator in the analysis of \cite{Gao--Oh--Viswanath2018tit}.
\end{remark}

\subsection{Generic variance bounds}\label{supp:sec:tech_lem_var}
\begin{lemma}\label{supp:lem:generic_variance:single}
For a given function $\phi\suchthat\Real_+\to\Real$, let $\zeta_k(\xb|\xb_{1:m})\defeq \phi(r_k(\xb|\xb_{1:m}))$ for any points $\xb,\xb_{1:m}$ in the $d$-dimensional Euclidean space $(\Real^d, \|\cdot\|)$.
Let
\begin{align}
\Phi(\xb_{1:m}) 
=\frac{1}{m}\sum_{i=1}^m \zeta_k(\xb_i|\xb_{1:m}^{\sim i}).\label{supp:eq:generic_est_form_single}
\end{align}
If the samples $\Xb_{1:m}$ are \iid, then 
\begin{align*}
&\Var(\Phi(\Xb_{1:m}))\\
&\le\frac{2(1+k\gamma_d)}{m}
\{(2k+1)\E[\zeta_{k}^2(\Xb_m|\Xb_{1:m-1})]
\\&\qquad\qquad\qquad\quad
+2k\E[\zeta_{k+1}^2(\Xb_m|\Xb_{1:m-1})]\},
\end{align*}
where $\gamma_d\in \Natural$ is a constant which depends only  on $d$.
\end{lemma}
Before we prove Lemma~\ref{supp:lem:generic_variance:single}, we introduce two technical lemmas.

\begin{lemma}[Efron--Stein inequality \cite{Efron--Stein1981,Steele1986}]\label{supp:lem:efron_stein}
Let $X_1,\ldots,X_n$ be independent random variables, and let $g(X_{1:n})=g(X_1,\ldots,X_n)$ be a square-integrable function of $X_1,\ldots,X_n$.
Then if $X_1',\ldots,X_n'$ are independent copies of $X_1,\ldots,X_n$, we have
\begin{align*}
&\Var(g(X_{1:n}))\\
&\le \half \sum_{i=1}^n
\E\bigl[|g(X_{1:n})-g(X_{1:i-1}X_i'X_{i+1:n})|^2\bigr].
\end{align*}
\end{lemma}
The proof of this lemma can be found in \cite{Steele1986}.

We need another fact on $k$-nearest neighbors in the Euclidean space, stated below in Lemma~\ref{supp:lem:generic_variance:single}. 
Informally speaking, given a finite collection $S$ of points in $\Real^d,$ each fixed point in $\Real^d$ can be one of the $k$ nearest neighbors of at most $\gamma_d$ points in $S$, where $\gamma_d$ depends only on $d.$ 
Henceforth, for a set of points $A$ such that $\xb\notin A$, we use $N_k(\xb|A)$ to denote the $k$-nearest neighbors of $\xb$ in $A$.
\begin{lemma}[{\cite[Lemma~20.6]{Biau--Devroye2015}, \cite[Ch.~5.3]{Devroye--Gyorfi--Lugosi2013}}]\label{supp:lem:knn_cone_covering}
In the $d$-dimensional Euclidean space $(\Real^d, \|\cdot\|)$
there exists a constant $\gamma_d>0$ which depends only on $d$ such that
for any $m\in\Natural$ and for any distinct points $\xb,\xb_1,\ldots,\xb_m\in\Real^d$,
\begin{align*}
\sum_{i=1}^m \ones_{\{\xb\in N_k(\xb_i|\xb_{1:m}^{\sim i},\xb)\}}\le k\gamma_d.
\end{align*}
\end{lemma}
\iftrue
\begin{proof}
We follow the proof of Stone's lemma in \citet[Ch.~5.3]{Devroye--Gyorfi--Lugosi2013}.
\newcommand{\Cc}{\mathcal{C}}
For $\zb\in\Real^d\backslash\{\mathbf{0}\}$ and $\theta\in(0,\pi/2]$, we define a cone $\Cc(\zb,\theta)\defeq \{\yb\in\Real^d\suchthat \yb=\mathbf{0} \text{ or }\angle(\zb,\yb)\le \theta\}$.
It is well known~\cite[Theorem~20.16]{Biau--Devroye2015} that there exists a constant $\gamma_d>0$, which depends only on the dimension $d$, such that there exist $\gamma_d$ cones $\Cc(\zb_1,\pi/6),\ldots,\Cc(\zb_{\gamma_d},\pi/6)$ which  cover the entire space $\Real^d$.
Furthermore, it is easy to see that ($\star$) if $\yb_1,\yb_2\in\Cc(\xb,\pi/6)$ and $\norm{\yb_1}<\norm{\yb_2}$, then $\norm{\yb_1-\yb_2} <\norm{\yb_2}$; see, e.g., \cite[Lemma~20.5]{Biau--Devroye2015}.

Now, for each $j\in[\gamma_d]$, \emph{mark} all $\xb_i$'s (if any) among the $k$-nearest neighbors of $\xb$ in $\xb+\Cc(\zb_j,\pi/6)$.
If $\xb_i\in\xb+\Cc(\zb_j,\pi/6)$ for some $j\in [\gamma_d]$ and $\xb_i$ is not marked, then $\xb$ is not among the $k$-nearest neighbors of $\xb_i$ in $\xb_{1:i-1},\xb_{i+1:m},\xb$, \ie $\xb\notin N_k(\xb_i|\xb_{1:m}^{\sim i},\xb)$, by the property ($\star$).
Therefore, we have
\begin{align*}
\sum_{i=1}^n \ones_{\{\xb\in N_k(\xb_i|\xb_{1:m}^{\sim i},\xb)\}}
\le \sum_{i=1}^n \ones_{\{\xb_i \text{ is marked}\}}
\le k\gamma_d,
\end{align*}
since there exist at most $k\gamma_d$ marked points.
\end{proof}
\fi

We are now ready to prove Lemma~\ref{supp:lem:generic_variance:single}.
\begin{proof}[Proof of Lemma~\ref{supp:lem:generic_variance:single}]
Let $\Xb_1'$ be an independent copy of $\Xb_1$.
Then, by applying the Efron--Stein inequality (Lemma~\ref{supp:lem:efron_stein}), we have
\begin{align}
&\Var\bigl(\Phi(\Xb_{1:m})\bigr)\nonumber\\
&\le \frac{m}{2}\E\bigl[
\bigl(\Phi(\Xb_{1:m})
-\Phi(\Xb_1'\Xb_{2:m})\bigr)^2
\bigr]
\nonumber\\
&\stackrel{(a)}{\le} m\E\bigl[
\bigl(
\Phi(\Xb_{1:m})
-\frac{m-1}{m}\Phi(\Xb_{2:m})
\bigr)^2
\\&\qquad\qquad
+\bigl(
\Phi(\Xb_1'\Xb_{2:m})
-\frac{m-1}{m}\Phi(\Xb_{2:m})
\bigr)^2
\bigr]
\nonumber\\
&=2m\E\bigl[
\bigl(
\Phi(\Xb_{1:m})
-\frac{m-1}{m}\Phi(\Xb_{2:m})
\bigr)^2
\bigr],\label{supp:eq:intermediate_label}
\end{align}
where (a) follows from the elementary inequality $(a-b)^2\le 2((a-x)^2+(b-x)^2)$.

Define
\[
E_i\defeq \{\Xb_1\text{~is one of the $k$-NNs of }\Xb_i\text{ in }\Xb_{1:m}^{\sim i}\}
\]
for $2\le i\le m$. 
Applying Lemma~\ref{supp:lem:knn_cone_covering}, we obtain 
\[
\sum_{i=2}^m \ones_{E_i} \le k\gamma_d.
\]
Further, note that if $E_i^c$ occurs, \ie $\Xb_1$ is not among the $k$ nearest neighbors of $\Xb_i$ in $\Xb_{1:m}^{\sim i}$, then $\zeta_k(\Xb_i|\Xb_{1:m}^{\sim i})=\zeta_k(\Xb_i|\Xb_{2:m}^{\sim i})$.
We thus obtain \eqref{eq:eq1_proof_III_5}, where (b) follows from Cauchy--Schwarz inequality.
\begin{table*}
\centering
\begin{minipage}{\textwidth}
\begin{align*}
m^2\bigl(
  \Phi(\Xb_{1:m})
  -\frac{m-1}{m}
      \Phi(\Xb_{2:m})
    \bigr)^2
&=\Bigl(\zeta_k(\Xb_1|\Xb_{2:m})
    + \sum_{i=2}^m \ones_{E_i}
    \bigl(
      \zeta_k(\Xb_i|\Xb_{1:m}^{\sim i})
      -\zeta_k(\Xb_i|\Xb_{2:m}^{\sim i})
    \bigr)
    \Bigr)^2\\
&\stackrel{(b)}{\le} \Bigl(1+\sum_{i=2}^m \ones_{E_i}\Bigr)
        \Bigl(
        \zeta_k^2(\Xb_1|\Xb_{2:m})
        + \sum_{i=2}^m \ones_{E_i}
            \bigl(
            \zeta_k(\Xb_i|\Xb_{1:m}^{\sim i})-\zeta_k(\Xb_i|\Xb_{2:m}^{\sim i})
        \bigr)^2
    \Bigr)\\
&\le (1+k\gamma_d)
        \Bigl(
        \zeta_k^2(\Xb_1|\Xb_{2:m})
        + 2\sum_{i=2}^m \ones_{E_i}
            \bigl(
            \zeta_k^2(\Xb_i|\Xb_{1:m}^{\sim i})
            +\zeta_k^2(\Xb_i|\Xb_{2:m}^{\sim i})
            \bigr)
        \Bigr).
    \numberthis\label{eq:eq1_proof_III_5}
\end{align*}
\end{minipage}
\medskip
\hrule
\end{table*}
By taking expectations with respect to $\Xb_{1:m}$ on both sides and multiplying by $2/m$, we can continue from \eqref{supp:eq:intermediate_label} to obtain
\begin{align*}
\numberthis\label{supp:eq:intermediate_label1}
&\Var\bigl(\Phi(\Xb_{1:m})\bigr)\\
&\le \frac{2(1+k\gamma_d)}{m}
\\&\quad\times
\Bigl\{
\E\bigl[\zeta_k^2(\Xb_1|\Xb_{2:m})\bigr]
\\&\quad\qquad
+2\E\Bigl[\sum_{i=2}^m \ones_{E_i}(\zeta_k^2(\Xb_i|\Xb_{1:m}^{\sim i})
+\zeta_k^2(\Xb_i|\Xb_{2:m}^{\sim i}))\Bigr]    \Bigr\}.
\end{align*}
Note that if $E_i$ occurs, \ie $\Xb_1$ is among the $k$ nearest neighbors of $\Xb_i$ in $\Xb_{1:m}^{\sim i}$, we have 
$\zeta_k(\Xb_i|\Xb_{2:m}^{\sim i})=\zeta_{k+1}(\Xb_i|\Xb_{1:m}^{\sim i})$.
Therefore, it follows that
\begin{align*}
&\E\Bigl[\sum_{i=2}^m \ones_{E_i}(\zeta_k^2(\Xb_i|\Xb_{1:m}^{\sim i})
+\zeta_k^2(\Xb_i|\Xb_{2:m}^{\sim i}))\Bigr]\\
&= \E\Bigl[\sum_{i=2}^m \ones_{E_i}(\zeta_k^2(\Xb_i|\Xb_{1:m}^{\sim i})
            +\zeta_{k+1}^2(\Xb_i|\Xb_{1:m}^{\sim i}))\Bigr]\\
&\stackrel{(c)}{=} \E\Bigl[\sum_{i=2}^m \ones_{\{\Xb_i\text{ is among the $k$-NNs of $\Xb_1$ in $\Xb_{2:m}$}\}}
    \\&\qquad\qquad\times
    (\zeta_k^2(\Xb_1|\Xb_{2:m})
            +\zeta_{k+1}^2(\Xb_1|\Xb_{2:m}))\Bigr]\\
&= k\E[\zeta_k^2(\Xb_1|\Xb_{2:m})+\zeta_{k+1}^2(\Xb_1|\Xb_{2:m})],\numberthis\label{supp:eq:intermediate_label2}
\end{align*}
where (c) follows by exchanging $\Xb_1$ and $\Xb_i$ in each summand $2\le i\le m$.
Therefore, plugging the equation in \eqref{supp:eq:intermediate_label2}  into \eqref{supp:eq:intermediate_label1} proves the desired bound.
\end{proof}

For the double-density case, we can establish a similar variance bound.

\begin{lemma}\label{supp:lem:generic_variance:double}
For a given function $\phi\suchthat\Real_+\times\Real_+\to\Real$, let $\zeta_{kl}(\xb|\xb_{1:m},\yb_{1:n})\defeq \phi(r_k(\xb|\xb_{1:m}),r_l(\xb|\yb_{1:n}))$ for any points $\xb,\xb_{1:m},\yb_{1:n}$ in the $d$-dimensional Euclidean space $(\Real^d, \|\cdot\|)$.
Let
\begin{align}\label{supp:eq:generic_est_form_double}
    \Phi(\xb_{1:m},\yb_{1:n})
    \defeq \frac{1}{m}\sum_{i=1}^m \zeta_{kl}(\xb_i|\xb_{1:m}^{\sim i},\yb_{1:n}).
\end{align}
If $\Xb_{1:m}$ and $\Yb_{1:n}$ are independent \iid samples, we have
\begin{align*}
&\Var(\Phi(\Xb_{1:m},\Yb_{1:n}))\\
&\le
\frac{2(1+k\gamma_d)}{m}
    \{(2k+1)\E[\zeta_{kl}^2(\Xb_m|\Xb_{1:m-1},\Yb_{1:n})]
        \\&\qquad\qquad\qquad\quad
        +2k\E[\zeta_{k+1,l}^2(\Xb_m|\Xb_{1:m-1},\Yb_{1:n})]
    \}.
\end{align*}
\end{lemma}
\begin{proof}
Given $\Yb_{1:n}=\yb_{1:n}$, we can show that
\begin{align*}
&\Var\bigl(\Phi\bigl(\Xb_{1:m},\yb_{1:n}\bigr)\bigr)\\
&\le 2m\E\bigl[
            \bigl(
            	\Phi(\Xb_{1:m},\yb_{1:n})
            	-\frac{m-1}{m}\Phi(\Xb_{2:m},\yb_{1:n})
            \bigr)^2
        \bigr]\\
&\le\frac{2(1+k\gamma_d)}{m}
    \{
        (2k+1)\E[\zeta_{kl}^2(\Xb_m|\Xb_{1:m-1},\yb_{1:n})]
        \\&\qquad\qquad\qquad\quad
        +2k\E[\zeta_{k+1,l}^2(\Xb_m|\Xb_{1:m-1},\yb_{1:n})]
    \}
\end{align*}
by following the same line of reasoning as in the proof of Lemma~\ref{supp:lem:generic_variance:single}.
Since $\Yb_{1:n}$ is independent of $\Xb_{1:m}$, taking expectation on both sides with respect to $\Yb_{1:n}$ establishes the desired bound.
\end{proof}

\section{Deferred proofs of main results}
\label{supp:sec:proof:main_results}

\subsection{Detailed proof of Theorem~\ref{thm:bias_class1_fixed_k_single}}\label{supp:sec:proof:thm:bias_class1_fixed_k_single}
We continue the proof from \eqref{eq:bias_decomposition}.
\begin{align*}
\bigl|\E[\That_f^{(k)}]-T_f(p)\bigr|
&\lesssim 
I_{\textout,1}
+I_{\textin,1}
+I_{\textin,2}
+I_{\textout,2}.
\tag{\ref{eq:bias_decomposition}}
\end{align*}
Applying the bounds in Lemmas~\ref{supp:lem:generic_inner_bias:single_case1} and \ref{supp:lem:generic_outer_bias_case1},
we obtain the following bias bound
for an underlying density $p$ satisfying the conditions \ref{cond:bounded_above}, \ref{cond:bounded_below}, \ref{cond:smoothness_int}, and \ref{cond:boundary}, provided that $\nu_m=o(\sqrt{m})$ as $m\to\infty$ and $k\in\Natural$ is fixed:
\begin{align*}
&|\E[\That_f^{(k)}]-T_f(p)|
\lesssim_{\sigma_p,L,\Cab_p,C_0,d,k}\\
&\quad
m^{-\frac{\sigma_p}{d}} \tau_m^{(a+\frac{\sigma_p}{d}+1)\wedge 0}+m^{-1}
+m^{-\frac{1}{d}} \tau_m^{(a+1)\wedge 0}\\
&\quad +m^{-\frac{\sigma_p}{d}} \nu_m^{(b+\frac{\sigma_p}{d}+2)\vee 0}
+m^{-1}\nu_m^{(b+k+2)\vee 0}
+m^{-\frac{1}{d}}\nu_m^{(b+2)\vee 0+\frac{1}{d}}\\
&\quad
+\tau_m^{k+a}
+\nu_m^{b+k-1}e^{-\Cbe_p\nu_m}.
\end{align*}
First, by choosing $\nu_m=\Theta((\ln m)^{1+\d})$ for some $\d>0$, we make the last term $\nu_m^{b+k-1}e^{-\Cbe_p\nu_m}$ decay faster than any polynomial rate. 
With this choice, the bound can be simplified as
\begin{align*}
&|\E[\That_f^{(k)}]-T_f(p)|\\
&=\bigOtilde_{\sigma_p,L,\Cab_p,C_0,d,k}(
\tau_m^{(a+\frac{\sigma_p}{d}+1)\wedge 0}m^{-\frac{\sigma_p}{d}}+\tau_m^{(a+1)\wedge 0}m^{-\frac{1}{d}}
\\&\qquad\qquad\qquad\qquad
+m^{-\frac{\sigma_p\wedge 1}{d}}+\tau_m^{k+a}).
\end{align*}
We consider three different ranges of the lower tail exponent $a$. 
\begin{enumerate}
\item If $a\le-\sigma_p/d-1$, we have
\begin{align*}
|\E[\That_f^{(k)}]-T_f(p)|
=\bigOtilde(
	\tau_m^{a+1}m^{-\frac{\sigma_p\wedge 1}{d}}
    +\tau_m^{k+a})
\end{align*}
as a suboptimal bound.
By equating the two terms, we establish a rate $\bigOtilde(m^{-\frac{(\sigma_p\wedge 1)}{d}\frac{k+a}{k-1}})$ with $\tau_m=\Theta(m^{-\frac{(\sigma_p\wedge 1)}{d}\frac{1}{k-1}})$.

    \item If $-\sigma_p/d-1< a \le-1$, the rate becomes
\begin{align*}
|\E[\That_f^{(k)}]-T_f(p)|
=\bigOtilde(
    &\tau_m^{a+1} m^{-\frac{1}{d}}
	+m^{-\frac{\sigma_p\wedge 1}{d}}
    +\tau_m^{k+a}).
\end{align*}
Equating $\tau_m^{a+1}m^{-\frac{1}{d}}$ and $\tau_m^{k+a}$ as a suboptimal choice, we obtain $\tau_m = \Theta(m^{-\frac{1}{d}\frac{1}{k-1}})$, which results in the final rate
\begin{align*}
|\E[\That_f^{(k)}]-T_f(p)|
&=\bigOtilde(
m^{-\frac{1}{d}\frac{k+a}{k-1}}
+m^{-\frac{\sigma_p\wedge 1}{d}}
)\\
&=\bigOtilde(
m^{-\frac{1}{d}(\sigma_p\wedge \frac{k+a}{k-1})})
\end{align*}

    \item If $a> -1$, we can attain the bias rate $\bigOtilde(m^{-\frac{\sigma_p\wedge 1}{d}})$ by using $\tau_m=O(m^{-\frac{1}{d(a+1)}})$.
\end{enumerate}
To sum up, by choosing
\begin{align}
\tau_m&=\tau(m,d,\sigma_p,a,k)
\label{supp:eq:alpha_classP_w_sigma}
\\
&= 
\begin{cases}
\Theta\bigl(m^{-\frac{\sigma_p\wedge 1}{d(k-1)}}\bigr)
	& \text{if }a\le -\frac{\sigma_p}{d}-1,\\
\Theta\bigl(m^{-\frac{1}{d(k-1)}}\bigr)
	&\text{if }-\frac{\sigma_p}{d}-1 < a \le -1,\\
O\bigl(m^{-\frac{1}{d(a+1)}}\bigr)
	&\text{if }a > -1,
\end{cases}
\nonumber
\end{align}
we establish the bias bound in Theorem~\ref{thm:bias_class1_fixed_k_single}.
\qed

\subsection{\secondrevision{Proof of Theorem~\ref{thm:vanishing_bias:double}}}
\label{supp:sec:proof:thm:vanishing_variance:double}
Following a similar line of reasoning as in the proof of Proposition~\ref{prop:knndist} and using the continuous mapping theorem, 
it is easy to show that
$\phi_k(U_{k,m-1}(\Xb_m),V_{ln}(\Xb_m))$
converges to $\phi_{kl}(U_{k\infty}(\Xb),V_{l\infty}(\Xb))$ in distribution as $m,n\to\infty$, where $U_{k\infty}(\xb)$ and $V_{l\infty}(\xb)$ are a $\GammaDist(k,p(\xb))$ random variable and a $\GammaDist(l,q(\xb))$ random variable, respectively, which are independent of each other and of $\Xb\sim p$, for $\P$-a.e.\@ $\xb$.
Hence, if we can only show that the collection of random variables $(\phi_{kl}(U_{k,m-1}(\Xb_m),V_{ln}(\Xb_m)))_{m,n\ge 1}$ is uniformly integrable, we can readily establish the asymptotic unbiasedness as follows:
\begin{align*}
&\lim_{m,n\to\infty} \E[\That_f^{(kl)}(\Xb_{1:m},\Yb_{1:n})]\\
&=\lim_{m,n\to\infty} \E[\phi_{kl}(U_{k,m-1}(\Xb_m),V_{ln}(\Xb_m))]\\
&=\E[\phi_{kl}(U_{k\infty}(\Xb),V_{l\infty}(\Xb))] \\
&=T_f(p,q).
\end{align*}
Consider
\begin{align*}
&\E[\xi(\abs{\phi_{kl}(U_{k,m-1}(\Xb_m),V_{ln}(\Xb_m))})]\\
&=\int p(\xb)\E[\xi(\abs{\phi_{kl}(U_{k,m-1}(\xb),V_{ln}(\xb))})]\diff\xb.
\end{align*}
By invoking the polynomial bound $|\phi_{kl}(u,v)|\lesssim \psi_{a,b}(u)\psi_{\at,\bt}(v)$ and using the independence of $U_{k,m-1}(\xb)$ and $V_{ln}(\xb)$, we have
\begin{align*}
&\E[\xi(\abs{\phi_{kl}(U_{k,m-1}(\Xb_m),V_{ln}(\Xb_m))})]\numberthis\label{supp:eq:proof:thm:vanishing_variance:double:temp1}\\
&\lesssim_{\xi(t_0)}
1+\E[\xi(\psi_{a,b}(U_{k,m-1}(\Xb_m)))]
\\&\qquad\qquad
+\E[\xi(\psi_{\at,\bt}(V_{ln}(\Xb_m)))]
\\&\qquad\qquad
+\{\E[(\E[\xi(\psi_{a,b}(U_{km}(\Xb_m)))|\Xb_m]
\\&\qquad\qquad\qquad\quad\times
\E[\xi(\psi_{\at,\bt}(V_{ln}(\Xb_m)))|\Xb_m])]\}^2,
\end{align*}
since $\xi(xy)\le \xi(x)\xi(y)$ for any $x,y>t_0$. 
We can bound the last term as
\begin{align*}
&\{\E[(\E[\xi(\psi_{a,b}(U_{km}(\Xb_m)))|\Xb_m]
\\&\qquad\times
\E[\xi(\psi_{\at,\bt}(V_{ln}(\Xb_m)))|\Xb_m])]\}^2\\
&\stackrel{(a)}{\le}
\E[(\E[\xi^2(\psi_{a,b}(U_{km}(\Xb_m)))|\Xb_m])^2]
\\&\qquad\times
\E[(\E[\xi^2(\psi_{\at,\bt}(V_{ln}(\Xb_m)))|\Xb_m])^2]\\
&\stackrel{(b)}{\le} \E[\xi^2(\psi_{a,b}(U_{k,m-1}(\Xb_m)))]
\E[\xi^2(\psi_{\at,\bt}(V_{ln}(\Xb_m)))],
\end{align*}
where (a) and (b) follow from Cauchy--Schwarz inequality and Jensen's inequality.
We thus only need to show that
\begin{align*}
\limsup_{m\to\infty} 
\E[\xi^2(\psi_{a,b}(U_{k,m-1}(\Xb_m)))]
<\infty
\end{align*}
and
\begin{align*}
\limsup_{n\to\infty}
\E[\xi^2(\psi_{\at,\bt}(V_{ln}(\Xb_m)))]
<\infty,
\end{align*}
since they would imply that all the terms in \eqref{supp:eq:proof:thm:vanishing_variance:double:temp1} are bounded.
By applying Lemma~\ref{supp:lem:generic:boundedness} to both integrals for $k>-2a\omega(\xi)$ and $l>-2\at\omega(\xi)$, we conclude the proof by the de la Vall\'{e}e Poussin theorem (Lemma~\ref{lem:de_la_vallee}).
\qed

\subsection{Proof of Theorem~\ref{thm:vanishing_variance:double}}
Recall from the generic variance bound (Lemma~\ref{supp:lem:generic_variance:double}) that we have
\begin{align*}
&\Var(T_f^{(kl)})\\
&\le
\frac{2(1+k\gamma_d)}{m}
    \{(2k+1)\E[\phi_{kl}^2(U_{k,m-1}(\Xb_m),V_{ln}(\Xb_m))]
    \\&\qquad\qquad\qquad
    +2k\E[\phi_{kl}^2(U_{k+1,m-1}(\Xb_m),V_{ln}(\Xb_m))]
    \}.
\end{align*}
Hence, following the same logic as in Section~\ref{supp:sec:proof:thm:vanishing_variance:double}, in order to ensure that $\Var(\That_f^{(kl)})=\bigO(m^{-1})$ for $m$ and $n$ sufficiently large,
it is enough to show that
\begin{align*}
\limsup_{m\to\infty} 
\E[\xi^2(\psi_{a,b}(U_{k',m-1}(\Xb_m)))]
<\infty
\end{align*}
and
\begin{align*}
\limsup_{n\to\infty}
\E[\xi^2(\psi_{\at,\bt}(V_{ln}(\Xb_m)))]
<\infty
\end{align*}
for $\xi(t)=t^2$ and for $k'\in\{k,k+1\}$.
By applying Lemma~\ref{supp:lem:generic:boundedness} to both integrals for $k>-4a$ and $l>-4\at$ with $\xi(t)=t^2$, we conclude the proof.
\qed

\subsection{Proof of Theorem~\ref{thm:bias_class1_fixed_kl_double}}
\label{sec:proof:thm:two}
Let $k>-a$ and $l>-\at$ be fixed.
First, following similar steps as in \eqref{eq:expectation_estimator_single}, we can write the expected value of $\That_{f}^{(kl)}(\Xb_{1:m},\Yb_{1:n})$ as
\begin{align*}
&\E\bigl[\That_{f}^{(kl)}(\Xb_{1:m},\Yb_{1:n})\bigr]\\
&=\int p(\xb)\E\bigl[\phi_{kl}(U_{k,m-1}(\xb),V_{ln}(\xb))\bigr]\diff\xb,
\end{align*}
since $U_{k,m-1}(\xb)$ and $V_{ln}(\xb)$ are independent of $\Xb_m=\xb$ for $\P$-a.e.\@ $\xb$.
Moreover, similar to \eqref{eq:target_functional_single}, we can write the target density functional as
\begin{align*}
T_f(p,q)
&= \int p(\xb) \E[\phi_{kl}(U_{k\infty}(\xb),V_{l\infty}(\xb))] \diff\xb,
\end{align*}
where $U_{k\infty}(\xb)\sim\GammaDist(k,p(\xb))$ and
$V_{l\infty}(\xb)\sim\GammaDist(l,q(\xb))$ are independent each other, and of
$\Xb\sim p$ for $\P$-a.e.\@ $\xb$.
Consider real numbers $\tau_m$, $\nu_m, \taut_n$, and $\nut_n$, to be determined later, such that $0\le \tau_m \le 1 \le \nu_m <\infty$ and $0\le \taut_n \le 1 \le \nut_n <\infty$.
Using the polynomial bound $\abs{\phi_{kl}(u,v)}\lesssim \psi_{a,b}(u)\psi_{\at,\bt}(v)$ and the triangle inequality, we then have
\begin{align}
|\E[\That_{f}^{(kl)}] - T_f(p,q)|
&\lesssim \int
  (I_{\textin}(\xb)+I_{\textout}(\xb))p(\xb)
  \diff\xb\nonumber\\
&=I_{\textin}+I_{\textout},
\label{supp:eq:abs_bias_two}
\end{align}
where $I_{\textin}(\xb)$ and $I_{\textout}(\xb)$ are defined in \eqref{eq:def_textin} and \eqref{eq:def_textout}, where $\square_{m,n}\defeq (\tau_m,\nu_m)\times(\taut_n,\nut_n)$.
\begin{table*}
\centering
\begin{minipage}{\textwidth}
\begin{align*}
I_{\textin}(\xb)
&\defeq \int_{\square_{m,n}}
\psi_{a,b}(u)\psi_{\at,\bt}(v)
\bigl|\rho_{U_{k\infty}(\xb)}(u)\rho_{V_{l\infty}(\xb)}(v)
-\rho_{U_{k,m-1}(\xb)}(u)\rho_{V_{ln}(\xb)}(v)\bigr|
\diff u\diff v,
\numberthis\label{eq:def_textin}\\
I_{\textout}(\xb)
&\defeq \int_{\Real^2_+\backslash \square_{m,n}}
\psi_{a,b}(u)\psi_{\at,\bt}(v)
(\rho_{U_{k\infty}(\xb)}(u)
\rho_{V_{l\infty}(\xb)}(v)
+\rho_{U_{k,m-1}(\xb)}(u)
\rho_{V_{ln}(\xb)}(v))
\diff u\diff v.
\numberthis\label{eq:def_textout}
\end{align*}
\end{minipage}
\medskip
\hrule
\end{table*}
We bound the inner bias $I_{\textin}=\int I_{\textin}(\xb) p(\xb)\diff\xb$ and the outer bias $I_{\textout}=\int I_{\textout}(\xb) p(\xb)\diff\xb$ separately.
Henceforth, we use the following shorthand notation:
\begin{align*}
\psit_{a,b}(u;\tau,\nu)=\psi_{a,b}(u)\ones_{(\tau,\nu)}(u)
\end{align*}
and
\begin{align*}
\psitt_{a,b}(u;\tau,\nu)=\psi_{a,b}(u)(1-\ones_{(\tau,\nu)}(u)).
\end{align*}

\textbf{Step 1: Bounding the inner bias.}
For $\xb\in\Real^d$, let $\d_{km}^{(p)}(u|\xb)\defeq 
\abs{\rho_{U_{k,m-1}(\xb)}(u)-\rho_{U_{k\infty}(\xb)}(u)}$ and $\d^{(q)}_{ln}(v|\xb)\defeq
\abs{\rho_{V_{ln}(\xb)}(v)-\rho_{V_{l\infty}(\xb)}(v)}$.
By the triangle inequality, we have
\begin{align*}
&\abs{\rho_{U_{k,m-1}(\xb)}(u)\rho_{V_{ln}(\xb)}(v)-\rho_{U_{k\infty}(\xb)}(u)\rho_{V_{l\infty}(\xb)}(v)}\\
&\le \d^{(p)}_{km}(u|\xb)\rho_{V_{ln}(\xb)}(v)
    +\d^{(q)}_{ln}(v|\xb)\rho_{U_{k\infty}(\xb)}(v)\\
&\le
\d^{(p)}_{km}(u|\xb)\d^{(q)}_{ln}(v|\xb)+\d^{(p)}_{km}(u|\xb)\rho_{V_{l\infty}}(\xb)
\\&\qquad\qquad\qquad\qquad\quad
+\d^{(q)}_{ln}(v|\xb)\rho_{U_{k\infty}}(\xb).
\end{align*}
Therefore, for each $\xb\in\supp(p)$, we can bound $I_{\textin}(\xb)$ as
\begin{align*}
&I_{\textin}(\xb)\\
&\le
    \int_{\tau_m}^{\nu_m}\psi_{a,b}(u)\d^{(p)}_{km}(u|\xb)\diff u
    \int_{\taut_n}^{\nut_n}\psi_{\at,\bt}(v)\d^{(q)}_{ln}(v|\xb)\diff v
\\&\quad
    +\E[\psit_{\at,\bt}(V_{l\infty}(\xb);\taut_n,\nut_n)]
    \int_{\tau_m}^{\nu_m}\psi_{a,b}(u)\d^{(p)}_{km}(u|\xb)\diff u
\\&\quad
    +\E[\psit_{a,b}(U_{k\infty}(\xb);\tau_m,\nu_m)]
    \int_{\taut_n}^{\nut_n}\psi_{\at,\bt}(v)\d^{(q)}_{ln}(v|\xb)\diff v\\
&\stackrel{(a)}{\lesssim} \int_{\tau_m}^{\nu_m}\psi_{a,b}(u)\d^{(p)}_{km}(u|\xb)\diff u
    +\int_{\taut_n}^{\nut_n}\psi_{\at,\bt}(v)\d^{(q)}_{ln}(v|\xb)\diff v,
\end{align*}
where (a) follows by applying Lemma~\ref{supp:lem:bound_gamma_integral} with the assumptions \ref{cond:bounded_above} and \ref{cond:bounded_below}.
Therefore, we have
\begin{align*}
I_{\textin}
&\lesssim 
\int p(\xb)
    \Bigl(\int_{\tau_m}^{\nu_m}\psi_{a,b}(u)\d^{(p)}_{km}(u|\xb)\diff u
    \\&\qquad\qquad\quad
    +\int_{\taut_n}^{\nut_n}\psi_{\at,\bt}(v)\d^{(q)}_{ln}(v|\xb)\diff v
    \Bigr)\diff\xb,
\end{align*}
and we can now apply the generic inner bias bounds in Lemma~\ref{supp:lem:generic_inner_bias:single_case1} to bound the inner bias.

\textbf{Step 2: Bounding the outer bias.}
We first consider the upper bound of $I_{\textout}(\xb)$ in \eqref{eq:bound_outer_bias}.
\begin{table*}
\centering
\begin{minipage}{\textwidth}
\begin{align*}
\numberthis\label{eq:bound_outer_bias}
I_{\textout}(\xb)
&\le
\int_{\Real\backslash(\tau_m,\nu_m)}
\int_{\taut_n}^{\nut_n} 
(\rho_{U_{k\infty}(\xb)}(u)\rho_{V_{l\infty}(\xb)}(v)
    +\rho_{U_{k,m-1}(\xb)}(u)\rho_{V_{ln}(\xb)}(v))
    \psi_{a,b}(u)\psi_{\at,\bt}(v) \diff u\diff v
\\&\quad
+ \int_{\tau_m}^{\nu_m}
\int_{\Real\backslash(\taut_n,\nut_n)}
(\rho_{U_{k\infty}(\xb)}(u)\rho_{V_{l\infty}(\xb)}(v)
+\rho_{U_{k,m-1}(\xb)}(u)\rho_{V_{ln}(\xb)}(v))
\psi_{a,b}(u)\psi_{\at,\bt}(v) \diff u\diff v.
\end{align*}
\end{minipage}
\medskip
\hrule
\end{table*}
For the first integral, we have
\begin{align*}
&\int_{\Real\backslash(\tau_m,\nu_m)}
\int_{\taut_n}^{\nut_n} 
\{\rho_{U_{k\infty}(\xb)}(u)\rho_{V_{l\infty}(\xb)}(v)
\\&\qquad\qquad\qquad\quad
+\rho_{U_{k,m-1}(\xb)}(u)\rho_{V_{ln}(\xb)}(v)\}
\\&\qquad\qquad\qquad\quad\times
\psi_{a,b}(u)\psi_{\at,\bt}(v) \diff u\diff v\\
&= \E[\psitt(U_{k\infty}(\xb);\tau_m,\nu_m)+\psitt(U_{k,m-1}(\xb);\tau_m,\nu_m)]
\\&\quad\times
\E[\psit(V_{l\infty}(\xb);\taut_n,\nut_n)+\psit(V_{ln}(\xb);\taut_n,\nut_n)]\\
&\stackrel{(b)}{\lesssim} \E[\psitt(U_{k\infty}(\xb);\tau_m,\nu_m)+\psitt(U_{k,m-1}(\xb);\tau_m,\nu_m)],
\end{align*}
where (b) follows from Lemmas~\ref{supp:lem:bound_gamma_integral} and \ref{supp:lem:bound_ukm_integral}.
The second integral can be bounded similarly. 
Overall, we have
\begin{align*}
I_{\textout}
&\lesssim \int p(\xb) \E[\psitt(U_{k\infty}(\xb);\tau_m,\nu_m)
\\&\qquad \qquad \qquad 
+\psitt(U_{k,m-1}(\xb);\tau_m,\nu_m)]\diff\xb
\\&\quad 
+\int p(\xb)\E[ \psitt(V_{l\infty}(\xb);\taut_n,\nut_n)
\\&\qquad \qquad \qquad 
+\psitt(V_{ln}(\xb);\taut_n,\nut_n)]\diff\xb,
\end{align*}
and we can now apply the generic outer bias bounds in Lemma~\ref{supp:lem:generic_outer_bias_case1}.

\textbf{Step 3: Choosing break points.}
Putting the bounds on the inner and outer bias together and choosing the break points $(\tau_m,\nu_m,\taut_n,\nut_n)$ as in the proof of Theorem~\ref{thm:bias_class1_fixed_kl_double}, we obtain the desired bias rates.
\qed

\subsection{Proof of Theorem~\ref{thm:variance_rate:double}}
\label{sec:proof:lem:variance:two}
By Lemma~\ref{supp:lem:generic_variance:double}, we have
\begin{align*}
&\Var(T_f^{(kl)}(\Xb_{1:m},\Yb_{1:n}))\\
&\le
\frac{2(1+k\gamma_d)}{m}
    \{
        (2k-2)\E[\phi_{kl}^2(U_{k-1,m-1}(\Xb_m),V_{ln}(\Xb_m))]\\
        &\qquad\qquad\qquad+(2k+1)\E[\phi_{kl}^2(U_{k,m-1}(\Xb_m),V_{ln}(\Xb_m))]\\
        &\qquad\qquad\qquad+\E[\phi_{kl}^2(U_{k+1,m-1}(\Xb_m),V_{ln}(\Xb_m))]
    \}.
\end{align*}
Using Lemma~\ref{supp:lem:bound_ukm_integral}, we have
\begin{align*}
&\E[\phi_{kl}^2(U_{k',m-1}(\Xb_m),V_{ln}(\Xb_m))]\\
&= \int p(\xb)\E[\phi_{kl}^2(U_{k',m-1}(\xb),V_{ln}(\xb))]\diff\xb\\
&\lesssim \int p(\xb)\E[\psi_{a,b}^2(U_{k',m-1}(\xb))]
\E[\psi_{\at,\bt}^2(V_{ln}(\xb))]\diff\xb\\
&\lesssim 1
\end{align*}
for
\[
k\in\begin{cases}
\{1,2\} & \text{if }k=1,\\
\{k-1,k,k+1\} & \text{if }k\ge2,
\end{cases}
\]
and for $m$ and $n$ sufficiently large, which concludes the proof.
\qed

\section{Deferred proofs of auxiliary results}

\subsection{Proof of Proposition~\ref{prop:bias_class2_fixed_k_single}}

Similar to Lemmas~\ref{supp:lem:generic_inner_bias:single_case1} and \ref{supp:lem:generic_outer_bias_case1}, we establish the following bounds.

\begin{lemma}[Generic inner bias bound under \ref{cond:smoothness_Rd}]\label{supp:lem:generic_inner_bias:single_case2}
Suppose that the density $p$ satisfies the conditions \ref{cond:bounded_above} 
and \ref{cond:smoothness_Rd} ,
and let $k=o(\sqrt{m})$ as $m\to\infty$.
\begin{enumerate}

\item We have
\begin{align*}
I_{\textin,1} &= \bigO\Bigl(
\frac{\tau_m^{(a+\frac{\sigma_p}{d}+1)\wedge 0}}{m^{\frac{\sigma_p}{d}}}
+\frac{k^{-k}}{m}
\Bigr).
\end{align*}

\item 
Suppose that $\phi_{k}(u)$ is differentiable at every $u>0$ and $\abs{\phi_k'(u)}\lesssim \psi_{a-1,b-1}(u)$.
If $\nu_m=o(\sqrt{m})$ as $m\to\infty$, then we have
\begin{align*}
I_{\textin,2}
&= \bigO\Bigl(
k\frac{\nu_m^{(b+\frac{\sigma_p}{d}+1)\vee 0}}{m^{\frac{\sigma_p}{d}}}
+\frac{\nu_m^{(b+k+1)\vee 0}}{m}
\Bigr).
\end{align*}
\end{enumerate}
\end{lemma}
\begin{proof}
To establish the second bound, we invoke Lemma~\ref{supp:lem:GOVLemma3_cdf_gap_bound} instead of Lemma~\ref{supp:lem:GOVLemma2_pdf_gap_bound}; this helps us obtain a tighter bias bound by reducing the exponent of $\nu_m$ by at most 1, which comes at the cost of additional factors in $k$.
Let 
\begin{equation*}
\Delta_{km}(u)\defeq \abs{\Rho_{U_{km}(\xb)}(u)-\Rho_{U_{k\infty}(\xb)}(u)}.
\end{equation*}
Since we assume that $\phi_{k}(u)$ is differentiable at any $u>0$ and $\abs{\phi_k'(u)}\lesssim \psi_{a-1,b-1}(u)$,
integration by parts leads to
\begin{align*}
I_{\textin,2}(\xb)
&= \Bigl|[\phi_{k}(u)\Delta_{km}(u)]_1^{\nu_m}
	+\int_1^{\nu_m}\phi_{k}'(u)
    		\Delta_{km}(u)
    	\diff u
    	\Bigr| \\
&\le \abs{\phi_{k}(\nu_m)}\cdot
	\Delta_{km}(\nu_m)
	+\abs{\phi_{k}(1)}\cdot
    \Delta_{km}(1)
    \\&\qquad
    +\int_1^{\nu_m}\abs{\phi_k'(u)}
    \cdot\Delta_{km}(u)\diff u\\
&=\bigOtilde_{\sigma_p,L,d}\Bigl(
k\frac{\nu_m^{(b+\frac{\sigma_p}{d}+1)\vee 0}}{m^{\frac{\sigma_p}{d}}}
    +\frac{\nu_m^{(k+b+1)\vee 0}}{m}
\Bigr)
\end{align*}
for $\xb\in\supp(p)$, 
establishing the second bound.
\end{proof}

Assuming \ref{cond:exp_integral_poly} in place of \ref{cond:bounded_below}, we obtain a different generic bound on the upper outer bias $I_{\textout,2}$ than that of Lemma~\ref{supp:lem:generic_outer_bias_case1}; see also Remark~\ref{supp:rem:tail_condition}.

\begin{lemma}[Generic outer bias bound under \ref{cond:exp_integral_poly} and \ref{cond:additional_regularity}]
\label{supp:lem:generic_outer_bias_case2}
Suppose that the density $p$ satisfies the conditions \ref{cond:bounded_above},
\ref{cond:exp_integral_poly},
and \ref{cond:additional_regularity},
we have
\begin{align*}
I_{\textout,2}
&=\bigO(\nu_m^{b+k-1-\th}).
\end{align*}
\end{lemma}

For any density $p$ satisfying the conditions \ref{cond:bounded_above}, 
\ref{cond:exp_integral_poly},
\ref{cond:additional_regularity}, and
\ref{cond:smoothness_Rd}, if $\nu_m=o(\sqrt{m})$ and $k$ is fixed, we have the bias bound from Lemmas~\ref{supp:lem:generic_inner_bias:single_case2} and \ref{supp:lem:generic_outer_bias_case2}:
\begin{align*}
&\bigl|\E\bigl[\That_f^{(k)}\bigr]-T_f(p)\bigr|\\
&\lesssim_{\sigma_p,L,\Cab_p,C_0,d,k}
\frac{\tau_m^{(a+\frac{\sigma_p}{d}+1)\wedge 0}}{m^{\frac{\sigma_p}{d}}}
+\frac{\nu_m^{(b+\frac{\sigma_p}{d}+1)\vee 0}}{m^{\frac{\sigma_p}{d}}}
\\&\qquad\qquad\qquad
+\frac{\nu_m^{(b+k+1)\vee 0}}{m}
+\tau_m^{k+a}+\nu_m^{b+k-1-\th}.
\end{align*}
Since $\nu_m\to\infty$ as $m\to\infty$, we require $b+k-1-\th<0$ to guarantee that the bias vanishes in our analysis, which forces us to choose a fixed $k$.

We first choose $\tau_m$. 
If $a+\frac{\sigma_p}{d}+1>0$, we can take $\tau_m=\bigO(m^{-\frac{\sigma_p}{d}\frac{1}{k+a}})$.
Otherwise, we take $\tau_m=\Th(m^{-\frac{\sigma_p}{d}\frac{1}{k-1-\frac{\sigma_p}{d}}})$ to make the first and the fourth terms decay with the same speed.
To summarize, we choose
\begin{align}
\tau_m&=\begin{cases}
\Th(m^{-\frac{\sigma_p}{d}\frac{1}{k-\frac{\sigma_p}{d}-1}}) & \text{if }a\le -\frac{\sigma_p}{d}-1,\\
O(m^{-\frac{\sigma_p}{d}\frac{1}{k+a}}) & \text{o.w.}
\label{supp:eq:tau_classQ}
\end{cases}
\end{align}
to bound the first and the fourth terms as
\begin{align*}
&\frac{\tau_m^{(a+\frac{\sigma_p}{d}+1)\wedge 0}}{m^{\frac{\sigma_p}{d}}}+\tau_m^{k+a}\\
&=\begin{cases}
\bigO(m^{-\frac{\sigma_p}{d}\frac{k+a}{k-\frac{\sigma_p}{d}-1}}) & \text{if }a\le -\frac{\sigma_p}{d}-1,\\
\bigO(m^{-\frac{\sigma_p}{d}}) & \text{o.w.}
\end{cases}
\end{align*}

\begin{table*}
\begin{minipage}{\textwidth}
\begin{align}\label{supp:eq:beta_classQ}
\nu_m
&=\begin{cases}
\Th(m^{(\frac{\sigma_p}{d}\wedge 1)\frac{1}{\th-k-b+1}})& 
    \text{if }k\le -b-1, b\le-\frac{\sigma_p}{d}-1,\\
\Th(m^{\frac{\sigma_p}{d}\frac{1}{\th-k+\frac{\sigma_p}{d}+2}}\bigr)
    & \text{if }k\le -b-1, b>-\frac{\sigma_p}{d}-1,\\
\Th(m^{\frac{1}{\th+2}})
    & \text{if } k>-b-1, b\le-\frac{\sigma_p}{d}-1,\\
\Th(m^{(\frac{\sigma_p}{d}\wedge 1)\frac{1}{\th+2}})
    & \text{if } k>-b-1, b>-\frac{\sigma_p}{d}-1
\end{cases}
\end{align}
\end{minipage}
\medskip
\hrule
\end{table*}
Similarly, by choosing $\nu_m$ as defined in \eqref{supp:eq:beta_classQ}
with $\nu_m=o(\sqrt{m})$ as $m\to\infty$,
we bound
the second, third, and last terms as
\[
\frac{1}{m^{\frac{\sigma_p}{d}}}+\frac{\nu_m^{(b+k+2)\vee 0}}{m}+\nu_m^{b+k-\th-1}
=\bigO(m^{-\lambda_\nu}),
\]
where $\lambda_\nu$ is as defined in \eqref{eq:lambda_beta_classQ}.
\qed

\subsection{Proof of Proposition~\ref{prop:bias_class1_increasing_k_single}}
\label{supp:sec:proofs:adaptive:single}

For any density $p$ satisfying the conditions \ref{cond:bounded_above}, \ref{cond:bounded_below}, \ref{cond:smoothness_int}, and \ref{cond:boundary}, if $\nu_m=o(\sqrt{m})$ and $k\to\infty$ with $k=o(\sqrt{m})$ as $m\to\infty$, we have the bias bound from Lemma~\ref{supp:lem:generic_inner_bias:single_case1}:
\begin{align*}
&\bigl|\E[\That_f^{(k)}]-T_f(p)\bigr|\\
&\lesssim_{\sigma_p,L,\Cab_p,C_0,d}
\frac{\tau_m^{(a+\frac{\sigma_p}{d}+1)\wedge 0}}{m^{\frac{\sigma_p}{d}}}+\frac{k^{-k}}{m}+\frac{\tau_m^{(a+1)\wedge 0}}{m^{\frac{1}{d}}}\\
&\qquad\qquad+\frac{\nu_m^{(b+\frac{\sigma_p}{d}+2)\vee 0}}{m^{\frac{\sigma_p}{d}}}
+k^{-k}\frac{\nu_m^{(b+k+2)\vee 0}}{m}
+\frac{\nu_m^{(b+2)\vee 0+\frac{1}{d}}}{m^{\frac{1}{d}}}\\
&\qquad\qquad+k^{-k}\tau_m^{k+a}
+k^{(b\vee 0)}\nu_m^{b+k-1}e^{-\Cbe_p\nu_m}.
\end{align*}
Setting $\nu_m=\Theta((\ln m)^{1+\d})$ and $k=\Theta((\ln m)^{1+\d'})$ for some $0<\d'<\d$, the last term $k^{(b\vee 0)}\nu_m^{b+k-1}e^{-\Cbe_p\nu_m}$ decays faster than any polynomial rate, that is, for any $C>0$,
\[
(b\vee 0)\ln k+ (b+k-1)\ln \nu_m -\Cbe_p\nu_m < -C\ln m
\]
for $m$ sufficiently large. 
With these choices of $\nu_m$ and $k$, the bias bound then can be simplified as
\begin{align*}
&\bigl|\E[\That_f^{(k)}]-T_f(p)\bigr|\\
&=\bigOtilde_{\sigma_p,L,\Cab_p,C_0,d}
\Bigl(
\frac{\tau_m^{(a+\frac{\sigma_p}{d}+1)\wedge 0}}{m^{\frac{\sigma_p}{d}}}+\frac{\tau_m^{(a+1)\wedge 0}}{m^{\frac{1}{d}}}+\frac{1}{m^{\frac{\sigma_p\wedge 1}{d}}}\Bigr).
\end{align*}
By choosing
\begin{align}
\tau_m
&=\tau'(m,a_k)
\label{eq:alpha_classP_varying_k}\\
&=\begin{cases}
\bigO((\polyln m)^{-1}) &\text{if }a_k\le -1\\
0 & \text{if }a_k>-1,
\end{cases}\nonumber
\end{align} 
we obtain
\begin{align*}
\bigl|\E[\That_f^{(k)}]-T_f(p)\bigr|
=\bigOtilde_{\sigma_p,L,\Cab_p,C_0,d}
\bigl(m^{-\frac{\sigma_p\wedge 1}{d}}\bigr).
\end{align*}

Now, we show that $\Var(\That_f^{(k)})=\bigOtilde(m^{-1})$ if $k=\Th((\ln m)^{1+\d})$ as $m\to\infty$ for some $\d>0$.
Using Lemmas~\ref{supp:lem:akm}, \ref{supp:lem:bound_kNN_ball_density}, \ref{supp:lem:bkm2:rate:bounded}, and \ref{supp:lem:bkm3}, if we choose $\nu_m$ and $\kappa_m$ such that $\nu_m/m\to 0$ and $\kappa_m/m\to\infty$ as $m\to\infty$, we have
\begin{align*}
&\Var(\That_f^{(k)})\\
&=\bigO\Bigl(\frac{k^2}{m}\Bigl\{\frac{\Cab_p^k}{k!}+\nu_m^{2b\vee 0} 
\\&\qquad\qquad\quad
+ (\nu_m^{2b}\vee \kappa_m^{2b})e^{-\nu_m\eta_p\Cbe_p}\Bigl(\frac{e\Cab_p\nu_m}{k}\Bigr)^k\Bigr\}\Bigr)
\end{align*}
for $m$ sufficiently large.
Letting $\nu_m=(2b/(\eta_p\Cbe_p))(\ln m)^{1+\d/2}$ and $\kappa_m=e^{(\ln m)^{1+\d/4}}$ ensures that the bound is $\bigOtilde(m^{-1})$.
\qed

\section{Derivation of estimator functions}\label{supp:sec:exmps}
In this section, we present derivations of some selected examples of estimator functions $\phi_{kl}(u,v)$ for some functions $f(p,q)$ in Table~\ref{table:estimator_functions_two}. 
Estimator functions $\phi_k(u)$ for the single-density case can be computed in a similar manner.
In particular, we present the examples of KL divergence (Example~\ref{supp:ex:kl}), logarithmic $\a$-divergences (Example~\ref{supp:ex:log_alpha_div}), entropy difference (Example~\ref{supp:ex:ent_diff}), reverse KL divergence (Example~\ref{supp:ex:rev_kl}), polynomial functionals (Example~\ref{supp:ex:poly}), Le Cam distance (Example~\ref{supp:ex:NNclassification}), and Jensen--Shannon divergence (Example~\ref{supp:ex:JSD}).

We remark that as alluded to in the main text, the estimator function $\phi_{kl}(u,v)$ is a function of $u/v$ if $f(p,q)$ is a function of $q/p$.

\begin{proposition}\label{supp:prop:gandh}
If $f(p,q)$ is a function of $q/p$, then there exists a function $\varphi_{kl}\suchthat \Real_+\to \Real$ such that $\phi_{kl}(u,v)=\varphi_{kl}(u/v)$.
\end{proposition}
\begin{proof}
Suppose that we can write $f(p,q)=g(q/p)$ for some function $g\suchthat\Real_+\to\Real$.
Recall that we have
\begin{align*}
&\Lc\{u^{k-1}v^{l-1}\phi_{kl}(u,v)\}(p,q)\\
&=\iint_{\mathbbm{R}^2_+}u^{k-1}v^{l-1}e^{-pu}e^{-qv}\phi_{kl}(u,v)\diff u \diff v \\
&= \frac{\Gamma(k)\Gamma(l)}{p^k q^l}g\Bigl(\frac{q}{p}\Bigr).  
\end{align*}
Now, for any $c>0$, we consider
\begin{align*}
&\Lc\{u^{k-1}v^{l-1}\phi_{kl}(cu,cv)\}(p,q)\\
&=\iint_{\mathbbm{R}^2_+}u^{k-1}v^{l-1}e^{-pu}e^{-qv}\phi_{kl}(cu,cv)\diff u \diff v\\
    &= \frac{1}{c^{k+l}}\iint_{\mathbbm{R}^2_+}\ut^{k-1}\vt^{l-1}e^{-p\ut/c}e^{-q\vt/c}\phi_{kl}(\ut,\vt)\diff\ut \diff\vt\\
    &= \frac{1}{c^{k+l}}\cdot\frac{\Gamma(k)\Gamma(l)}{(p/c)^k(q/c)^l}g\Bigl(\frac{q/c}{p/c}\Bigr)\\
    &= \frac{\Gamma(k)\Gamma(l)}{p^kq^l}g\Bigl(\frac{q}{p}\Bigr).
\end{align*}
Thus, by the (a.e.) uniqueness of Laplace transform, we have $\phi_{kl}(cu,cv) = \phi_{kl}(u,v),$ whence $\phi_{kl}(u,v)$ can be written as 
$\phi_{kl}(u,v)=\varphi_{kl}(u/v)$ for some function $\varphi\suchthat\Real_+\to\Real$.
\end{proof}

In what follows, for the one-dimensional inverse Laplace transform of two-variable functions, we will specify the transformed variable by a subscript of the inverse Laplace operator.
For example, $\Lc_p^{-1}\{G(p,q)\}(u)$ denotes the inverse Laplace transform of $G(p,q)$ along the $p$-axis with a corresponding time-domain variable $u$.

\begin{example}[KL divergence; Example~\ref{ex:kl_divergence}]\label{supp:ex:kl}
For $f(p,q)=\ln (p/q)$, the corresponding functional $T_f(p,q)=\D{p}{q}$ is the KL divergence.
This is one of the simplest cases, as we only need to deal with one-dimensional inverse Laplace transforms by linearity:
\begin{align*}
&\Lc^{-1}\Bigl\{\frac{1}{p^kq^l}\ln \frac{p}{q}\Bigr\}\\
&= \Lc^{-1}\Bigl\{\frac{\ln p}{p^k}\Bigr\}\Lc^{-1}\Bigl\{\frac{1}{q^l}\Bigr\}
- \Lc^{-1}\Bigl\{\frac{1}{p^k}\Bigr\}\Lc^{-1}\Bigl\{\frac{\ln q}{q^l}\Bigr\}.
\end{align*}
Note that for any $\kappa>0$,
\begin{align}\label{supp:eq:inv_laplace_lnp_over_pk}
\Lc^{-1}\Bigl\{\frac{\ln p}{p^\kappa}\Bigr\}
= \frac{u^{\kappa-1}}{\Gamma(\kappa)}\bigl(\digamma(\kappa)-\ln u\bigr).
\end{align}
This can be verified by taking Laplace transform of the right-hand expression.
From the definition of the estimator function $\phi_{kl}(u,v)$ in \eqref{eq:estimator_double}, we obtain
\begin{align}
\phi_{kl}(u,v)
&=\ln\frac{v}{u}+\digamma(k)-\digamma(l).
\label{supp:eq:est_function_KL}
\end{align}
\end{example}

\begin{example}[Polynomial functionals; Example~\ref{ex:polynomial_functional}]
\label{supp:ex:poly}
Consider $f(p,q)=p^{\a-1}q^{\b}$ for some $\a,\b\in\Real$, which corresponds to the functional
\begin{align*}
T_f(p,q)=\E\bigl[p^{\a-1}(\Xb)q^{\b}(\Xb)\bigr]
=\int p^{\a}(\xb)q^{\b}(\xb)\diff \xb.
\end{align*}
This includes many special cases such as \Renyi{} entropies, \Renyi{} divergences, Hellinger distance, and $\chi^2$-divergence.
The estimator function is
\begin{align*}
\phi_{kl}(u,v)=\frac{\Gamma(k)\Gamma(l)}{\Gamma(k-\a+1)\Gamma(l-\b)}u^{1-\a}v^{-\b}
\end{align*}
for $k>\a-1$ and $l>\b$.
We remark that our estimator recovers the bias-corrected estimator presented in \cite{Poczos--Xiong--Sutherland--Schneider2012}.
\end{example}

\begin{example}[Logarithmic $\a$-divergence; Example~\ref{ex:logarithmic_alpha_divergence}]
\label{supp:ex:log_alpha_div}
For $\a\in\Real$, consider a function
$f(p,q)=(p/q)^{\a-1}\ln\frac{p}{q}$, which corresponds to the functional
\begin{align*}
T_f(p,q)
&=\E\Bigl[\Bigl(\frac{p(\Xb)}{q(\Xb)}\Bigr)^{\a-1}\ln\frac{p(\Xb)}{q(\Xb)}\Bigr]\\
&=\int p^{\a}(\xb)q^{1-\a}(\xb)\ln\frac{p(\xb)}{q(\xb)}\diff \xb.
\end{align*}
Similar to KL divergence, the estimator function can be found immediately from \eqref{supp:eq:inv_laplace_lnp_over_pk}, \ie
\begin{align*}
\phi_{kl}(u,v)
&=\frac{\Gamma(k)\Gamma(l)}{\Gamma(k-\a+1)\Gamma(l+\a-1)}\Bigl(\frac{v}{u}\Bigr)^{\a-1}
\\&\quad\times
\Bigl(
\ln\frac{v}{u}+\digamma(k-\a+1)-\digamma(l+\a-1)
\Bigr),
\end{align*}
for $k>\a-1$ and $l>-\a+1$.
Note that $\a=1$ recovers the estimator function for the KL divergence \eqref{supp:eq:est_function_KL}.
\end{example}

\begin{example}[Le Cam distance; Example~\ref{ex:asymp_nn_classfication}]\label{supp:ex:NNclassification}
For $f(p,q)=1-2q/(p+q)$, we wish to compute the estimator function $\phi_{kl}(u,v)$, that is,
\begin{align*}
\phi_{kl}(u,v)=2\frac{\Gamma(k)\Gamma(l)}{u^{k-1}v^{l-1}}\Lc^{-1}\Bigl\{\frac{1}{p^kq^l}\frac{1}{1+\frac{q}{p}}\Bigr\} - 1.
\end{align*}
The two-dimensional inverse Laplace transform can be peeled off dimension by dimension as follows:
\begin{align*}
&\Lc^{-1}_{p,q}\Bigl\{\frac{1}{p^kq^l}\frac{1}{1+\frac{q}{p}}\Bigr\}(u,v)\\
&=\Lc^{-1}_p\Bigl\{\frac{1}{p^{k+l}}
    \Lc^{-1}_q\Bigl\{\frac{1}{(\frac{q}{p})^{l}(1+\frac{q}{p})}\Bigr\}(v)
    \Bigr\}(u).
    \numberthis\label{supp:eq:calc1}
\end{align*}
Letting $\qtil=q/p$, we first find the inverse Laplace transform of
\begin{align*}\numberthis\label{supp:eq:useful_identity}
\frac{1}{\qtil^l(1+\qtil)}
=(-1)^l\Bigl(\sum_{i=1}^l\frac{(-1)^i}{\qtil^i}+\frac{1}{1+\qtil}\Bigr),
\end{align*}
which is
\begin{align*}
\Lc_{\qtil}^{-1}\Bigl\{\frac{1}{\qtil^l(1+\qtil)}\Bigr\}(v)
=(-1)^l\Bigl(e^{-v}-\sum_{i=0}^{l-1}\frac{(-v)^i}{i!}\Bigr),
\end{align*}
since we have
\begin{align*}
\Lc_p^{-1}\Bigl\{\frac{1}{p^{n+1}}\Bigr\}(u) = \frac{u^{n}}{n!}\ones_{[0,\infty)}(u)
\end{align*}
for $n\in\Natural\cup\{0\}$ and
\begin{align*}
\Lc_p^{-1}\Bigl\{\frac{1}{s+a}\Bigr\}(u) = e^{-au}\ones_{[0,\infty)}(u).
\end{align*}
Moreover, by the time-scaling property, we have
\begin{align*}
&\Lc^{-1}_q\Bigl\{\frac{1}{(\frac{q}{p})^{l}(1+\frac{q}{p})}\Bigr\}(v)\\
&= (-1)^l\Bigl(pe^{-pv}-\sum_{i=0}^{l-1}\frac{(-v)^i}{i!}p^{i+1}\Bigr).
\end{align*}
Now, continuing from \eqref{supp:eq:calc1}, we have \eqref{eq:interim_le_cam}, which leads to the estimator function~\eqref{eq:le_cam_estimator_function}.
As a bound on the estimator function $\phi_{kl}(u,v)$, we observe that
\begin{align*}
|\phi_{kl}(u,v)|
&\lesssim
\Bigl(\frac{u}{v}\Bigr)^{l-1}\Bigl(\sum_{i=0}^{l-1}\Bigl(\frac{v}{u}\Bigr)^i + \sum_{j=0}^{k+l-2}\Bigl(\frac{v}{u}\Bigr)^{j}\Bigr)\\
&\lesssim \psi_{-k+1,l-1}(u)\psi_{-l+1,k-1}(v).
\end{align*}



\begin{table*}
\begin{minipage}{\textwidth}

\begin{align*}
\Lc^{-1}_{p,q}\Bigl\{\frac{1}{p^kq^l}\frac{1}{1+\frac{q}{p}}\Bigr\}(u,v)
&= \Lc^{-1}_{p}\Bigl\{(-1)^l\Bigl(\frac{e^{-pv}}{p^{k+l-1}}-\sum_{i=0}^{l-1}\frac{(-v)^i}{i!}\frac{1}{p^{k+l-i-1}}\Bigr)\Bigr\}(u)\\
&= (-1)^l\Bigl(\frac{(u-v)^{k+l-2}}{(k+l-2)!}\ones_{[v,\infty)}(u)-\sum_{i=0}^{l-1}\frac{(-v)^i}{i!}\frac{u^{k+l-i-2}}{(k+l-i-2)!}\Bigr)\\
&= (-1)^l\frac{u^{k+l-2}}{(k+l-2)!}\Bigl(\Bigl(1-\frac{v}{u}\Bigr)^{k+l-2}\ones_{[v,\infty)}(u)-\sum_{i=0}^{l-1}\binom{k+l-2}{i}\Bigl(\frac{-v}{u}\Bigr)^i\Bigr).
\numberthis\label{eq:interim_le_cam}
\end{align*}
\end{minipage}
\medskip
\hrule
\end{table*}

\begin{table*}
\begin{minipage}{\textwidth}
\begin{align*}
\phi_{kl}(u,v)&=
2\binom{k+l-2}{k-1}^{-1}\Bigl(-\frac{u}{v}\Bigr)^{l-1}
\Bigl(\sum_{i=0}^{l-1}\binom{k+l-2}{i}\Bigl(-\frac{v}{u}\Bigr)^i-\Bigl(1-\frac{v}{u}\Bigr)^{k+l-2}\ones_{[v,\infty)}(u)\Bigr)-1.
\numberthis\label{eq:le_cam_estimator_function}
\end{align*}
\end{minipage}
\medskip
\hrule
\end{table*}

\end{example}

For the remaining examples, we assume that $\Q\ll\P$.

\begin{example}[Entropy difference]
\label{supp:ex:ent_diff}
For $f(p,q)=\ln(1/p)-(q/p)\ln(1/q)$, the corresponding functional $T_f(p,q)=h(p)-h(q)$ becomes the difference of the differential entropies $h(p)$ and $h(q)$.
It is easy to show that
\begin{align*}
\phi_{kl}(u,v) 
    &= \frac{(l-1)}{k}\frac{u}{v}(\digamma(l-1)-\ln v)
        -(\digamma(k)-\ln u).
\end{align*}
As a bound on the estimator function $\phi_{kl}(u,v)$, we have
\begin{align*}
|\phi_{kl}(u,v)|
&\lesssim \frac{u}{v}(1+|\ln v|)+(1+|\ln u|)\\
&\lesssim \psi_{1,1}(u)\psi_{-1-\eps,-1+\eps}(v)+\psi_{-\eps,\eps}(u)\\
&\lesssim \psi_{-\eps,1}(u)\psi_{-1-\eps,-1+\eps}(v).
\end{align*}
\end{example}

\begin{example}[Reverse KL divergence]
\label{supp:ex:rev_kl}
When $\Q\ll\P$, we can write the reverse KL divergence as
\begin{align*}
\D{q}{p}
&=\int q(\xb) \ln\frac{q(\xb)}{p(\xb)}\diff\xb\\
&=\int p(\xb) \frac{q(\xb)}{p(\xb)} \ln\frac{q(\xb)}{p(\xb)}\diff\xb=T_f(p,q)
\end{align*}
for $f(p,q)=(q/p)\ln (q/p)$.
Then, for $k\ge 1$ and $l\ge 2$, we have
\begin{align*}
\Lc^{-1}\Bigl\{\frac{f(p,q)}{p^kq^l}\Bigr\}
&=\Lc^{-1}\Bigl\{\frac{1}{p^{k+1}}\Bigr\} \Lc_q^{-1}\Bigl\{\frac{\ln q}{q^{l-1}}\Bigr\}
\\&\quad
-\Lc^{-1}\Bigl\{\frac{\ln p}{p^{k+1}}\Bigr\} \Lc_q^{-1}\Bigl\{\frac{1}{q^{l-1}}\Bigr\}\\
&=\frac{u^k}{\Gamma(k+1)}\frac{v^{l-2}}{\Gamma(l-1)}\bigl(\digamma(l-1)-\ln v\bigr)
\\&\quad
-\frac{u^k}{\Gamma(k+1)}\bigl(\digamma(k+1)-\ln u\bigr)\frac{v^{l-2}}{\Gamma(l-1)}.
\end{align*}
Here, the case $l=1$ is excluded, since $\Lc^{-1}\{\ln s\}$ is ill-defined.
Finally, we have
\begin{align*}
\phi_{kl}(u,v)
&= \frac{\Gamma(k)\Gamma(l)}{u^{k-1}v^{l-1}}
    \frac{u^k}{\Gamma(k+1)}\frac{v^{l-2}}{\Gamma(l-1)}
\\&\quad\times
    \bigl\{
    \bigl(\Psi(l-1)-\ln v\bigr)
    -\bigl(\Psi(k+1)-\ln u\bigr)
    \bigr\}\\
&= \frac{l-1}{k}\frac{u}{v}\bigl(\ln\frac{u}{v}+\Psi(l-1)-\Psi(k+1)\bigr).
\end{align*}
As a bound on the estimator function $\phi_{kl}(u,v)$, we have
\begin{align*}
|\phi_{kl}(u,v)|
&\lesssim \frac{u}{v}(1+|\ln u|+|\ln v|)\\
&\lesssim \frac{u}{v}(1+|\ln u|)(1+|\ln v|)\\
&\lesssim \psi_{1-\eps,1+\eps}(u)\psi_{-1-\eps,-1+\eps}(v).
\end{align*}
\end{example}

\begin{table}
	\caption{Inverse Laplace transforms of few elementary functions and basic operations.}
    \centering
    \begin{adjustbox}{center}
    \begin{tabular}{c c}
     \toprule
     \makecell{Frequency domain\\$F(p)=\Lc\{f(u)\}$}
      & \makecell{Time domain\\$f(u)=\Lc^{-1}\{F(p)\}$}\\
      \midrule
     $p^{-k}~(k>0)$
      & $ u^{k-1}/\Gamma(k)$\\
     $ \ln p/p$
      & $ -(\ln u+\gamma)$\\
     $ 1/(p+\a)$
      & $ e^{-\a u}$\\
     \midrule
     $F(ap)$
      & $ f\bigl(u/a\bigr)/a$\\
     $e^{-ap}F(p)$
      & $f(u-a)\ones_{[a,\infty)}(u)$\\
     $F^{(n)}(p)$
      & $(-1)^n u^n f(u)$\\
     $ F(p)/p$
      & $ \int_0^u f(t)\diff t$\\
     $F(p)G(p)$
        & $ (f*g)(u) =\int_0^u f(\ut)g(u-\ut)\diff\ut$\\
     $pF(p)$
        & $ f'(u)-f(0)$\\
     \bottomrule
    \end{tabular}
    \end{adjustbox}
    \label{supp:table:inverse_laplace}
\end{table}
\begin{example}[Jensen--Shannon divergence; Example~\ref{ex:JSD}]
\label{supp:ex:JSD}
We wish to compute the estimator function $\phi_{kl}(u,v)$ for
\begin{equation*}
f(p,q)=\half\Bigl(\frac{q}{p}+1\Bigr)\ln\frac{2}{(q/p)+1}+\frac{q}{2p}\ln \frac{q}{p}.
\end{equation*} 
For $l\ge 2,$ we have
\begin{align*}
\frac{2f(p,q)}{p^kq^l}
&=\Bigl(\frac{1}{p^{k+1}q^{l-1}}+\frac{1}{p^kq^l}\Bigr)\ln 2+\frac{1}{p^{k+1}q^{l-1}}\ln\frac{q}{p}
\\&\quad
    +\frac{G_{l-1}(\frac{q}{p})+G_l(\frac{q}{p})}{p^{k+l}},
\end{align*}
where we define $G_l(q)\defeq -\ln(q+1)/q^l$.
Using the identity~\eqref{supp:eq:useful_identity}, we can show that for $l\in\Natural$
\begin{align*}
g_l(v)
&= \Lc_q^{-1}\{G_l(q)\}(v)\\
&=(-1)^{l}\Bigl(
    \int_1^\infty \frac{e^{-vx}}{x^l}\diff x
    -\sum_{j=0}^{l-2}\frac{(-v)^{j}}{(l-1-j)j!}
    \Bigr).
\end{align*}
Now the desired estimator function can be written as
\begin{align*}
2\phi_{kl}(u,v)
&= \frac{\Gamma(k)\Gamma(l)}{u^{k-1}v^{l-1}}\Lc^{-1}\Bigl\{\frac{f(p,q)}{p^kq^l}\Bigr\}(u,v)\nonumber \\
&= \frac{l-1}{k}\frac{u}{v} \Bigl(\Psi(l-1)-\Psi(k+1)+\ln\frac{u}{v}\Bigr)
\\&\quad
+\Bigl(\frac{l-1}{k}\frac{u}{v}+1\Bigr)\ln 2
+A_{kl}(u,v),
\numberthis
\label{supp:eq:JS_finalexp}
\end{align*}
where we define
\begin{align}
&A_{kl}(u,v)\nonumber\\
&=\frac{\Gamma(k)\Gamma(l)}{u^{k-1}v^{l-1}}\Lc_p^{-1}\Bigl\{\frac{\Lc_q^{-1}\{G_{l-1}(\frac{q}{p})+G_l(\frac{q}{p})\}(v)}{p^{k+l}}\Bigr\}(u)\nonumber \\
&\stackrel{(a)}{=} \frac{\Gamma(k)\Gamma(l)}{u^{k-1}v^{l-1}}\Lc_p^{-1}\Bigl\{\frac{g_{l-1}(pv)+g_{l}(pv)}{p^{k+l-1}}\Bigr\}(u)\nonumber \\
&= B_{kl}(u,v) + \frac{l-1}{k}\frac{u}{v} B_{k+1,l-1}(u,v),\label{supp:eq:JS_Akl}
\end{align}
where
\begin{align*}
B_{kl}(u,v)
&= \frac{\Gamma(k)\Gamma(l)}{u^{k-1}v^{l-1}}\Lc_p^{-1}\Bigl\{\frac{g_l(pv)}{p^{k+l-1}}\Bigr\}(u).
\end{align*}
Here, (a) follows by the time scaling property, that is, $\Lc_q^{-1}\{G_l(q/p)\}(v)=pg_l(pv)$.
\begin{table*}
\begin{minipage}{\textwidth}
\begin{align*}
\Lc_p^{-1}\Bigl\{\frac{g_l(pv)}{p^{k+l-1}}\Bigr\}
&= \int_1^\infty \frac{1}{x^l}\Lc_p^{-1}\Bigl\{\frac{e^{-pvx}}{p^{k+l-1}}\Bigr\}\diff x 
    - \sum_{j=0}^{l-2}\frac{(-v)^{j}}{(l-1-j)j!} \Lc_p^{-1}\Bigl\{\frac{1}{p^{k+l-1-j}}\Bigr\}\\
&= \int_1^\infty \frac{1}{x^l}
    \ones_{[vx,\infty)}(u)\frac{(u-vx)^{k+l-2}}{(k+l-2)!}
    \diff x 
    - \sum_{j=0}^{l-2}\frac{(-v)^{j}}{(l-1-j)j!} \frac{u^{k+l-2-j}}{(k+l-2-j)!},
\numberthis
\label{eq:interim_jsd}
\end{align*}
\end{minipage}
\medskip
\hrule
\end{table*}
Now, since we have \eqref{eq:interim_jsd},
it follows that
\begin{align*}
&\binom{k+l-2}{k-1}B_{kl}(u,v)\\
&=-\ones_{[1,\infty)}(w) (-w)^{-k+1}\int_1^w \frac{(x-w)^{k+l-2}}{x^l}\diff x
\\&\quad
    + \sum_{j=0}^{l-2} \binom{k+l-2}{j}\frac{(-w)^{l-1-j}}{l-1-j},
\end{align*}
where $w\defeq u/v$.

Rearranging the integral in the parenthesis as
\begin{align*}
&(-w)^{k+1}\int_1^w \frac{(x-w)^{k+l-2}}{x^l}\diff x\\
&= \sum_{\substack{i=0\\i\neq k-1}}^{k+l-2}\binom{k+l-2}{i} \frac{(-1)^{k-1-i}-(-w^{-1})^{k-1-i}}{k-1-i}
\\&\qquad
+ \binom{k+l-2}{k-1}\ln w,
\end{align*}
we finally obtain
\begin{align*}
\numberthis\label{supp:eq:JS_Bkl_1}
&B_{kl}(u,v)\\
&= \binom{k+l-2}{k-1}^{-1}\sum_{j=0}^{l-2}\binom{k+l-2}{j}\frac{(-u/v)^{l-1-j}}{l-1-j} 
\end{align*}
if $\frac{u}{v}< 1$,
and 
\begin{align*}
\numberthis\label{supp:eq:JS_Bkl_2}
B_{kl}(u,v)
&= -\ln \frac{u}{v}
+ \binom{k+l-2}{k-1}^{-1}
\\&\qquad\times
\Bigl\{\sum_{i=0}^{k-2}\binom{k+l-2}{i}\frac{(-v/u)^{k-1-i}}{k-1-i}
\\&\qquad\qquad
-\sum_{\substack{i=0\\i\neq k-1}}^{k+l-2}\binom{k+l-2}{i} \frac{(-1)^{k-1-i}}{k-1-i} \Bigr\}
\end{align*}
if $\frac{u}{v}\ge 1$.
Substituting the expressions for $B_{kl}(u, v)$ from~\eqref{supp:eq:JS_Bkl_1} and~\eqref{supp:eq:JS_Bkl_2} into~\eqref{supp:eq:JS_Akl} and then into~\eqref{supp:eq:JS_finalexp} yields the final expression for the estimator function as
\begin{align*}
&\phi_{kl}(u,v)\\
&= \half\Bigl\{\ln 2
+\frac{l-1}{k}\frac{u}{v} \Bigl(\ln 2 + \Psi(l-1)-\Psi(k+1)+\ln\frac{u}{v}\Bigr)
\\&\qquad\quad
+ B_{kl}(u,v) + \frac{l-1}{k}\frac{u}{v} B_{k+1,l-1}(u,v)\Bigr\}.
\end{align*}
As a bound on the estimator function $\phi_{kl}(u,v)$, we have
\[
|\phi_{kl}(u,v)|
\lesssim \psi_{-k+1,l-1}(u)\psi_{-l+1,k-1}(v).
\]
\end{example}

\section{Examples of smooth densities}
\label{supp:sec:examples_densities}
In this section, we show that the $d$-dimensional truncated Gaussian, Cauchy, and exponential distributions, as well as the uniform distribution and the $d$-dimensional product of identical beta distributions with parameters $\a\ge 3$ and $\b\ge 3$ satisfy the conditions \ref{cond:bounded_above}, \ref{cond:bounded_below}, \ref{cond:smoothness_int}, and \ref{cond:boundary} with $\sigma_p=2$, and the $d$-dimensional truncated Laplace distribution satisfies the conditions with $\sigma_p=1$.
We remark that the boundedness of the Hessian of the density $p$ over a compact set implies 2-H\"older continuity, if the Hessian is integrable.
Since we have considered that the Hessian is integrable, we only need to prove the boundedness of the Hessian in order to demonstrate the 2-H\"older continuity.

\begin{example}[Truncated Gaussian]
Consider the {\it truncated $d$-dimensional Gaussian distribution} defined by the density
\begin{equation*} p(\xb)\defeq \frac{\Gamma(d/2+1)}{\pi^{d/2}K_d(R)}e^{-\|\xb\|_2^2/2}\ones_{(-\infty,R]}(\|\xb\|_2),
\end{equation*}
where $K_d(R)\defeq\int_0^Rdr^{d-1}e^{-r^2/2}\diff r.$ Then,
$\supp(p) = \{\xb\in\Real^d: \|\xb\|\le R\}$
and
\begin{equation*} \frac{\Gamma(d/2+1)}{\pi^{d/2}K_d(R)}e^{-R^2/2}\le p(\xb)\le \frac{\Gamma(d/2+1)}{\pi^{d/2}K_d(R)}\end{equation*}
for $\xb\in\supp(p).$ Moreover, on $\supp(p)^\mathrm{o},$
\begin{equation*} \nabla^2p(\xb)_{ij} = \frac{\Gamma(d/2+1)}{\pi^{d/2}K_d(R)}(x_ix_j - \d_{ij})e^{-\|\xb\|_2^2/2},\end{equation*}
whence,
\begin{equation*} \|\nabla^2p(\xb)\|\le\|\nabla^2p(\xb)\|_F\le\frac{\Gamma(d/2+1)}{\pi^{d/2}K_d(R)}\sqrt{R^4+d}.
\end{equation*}
Finally, $\partial\supp(p) = \Sb(\mathbf{0}, R)$ satisfies
\begin{equation*} H^{d-1}(\Sb(\mathbf{0}, R)) = d\ups_dR^{d-1}.\end{equation*}
Therefore, this density satisfies the conditions \ref{cond:bounded_above}, \ref{cond:bounded_below}, \ref{cond:smoothness_int}, and \ref{cond:boundary} with $\sigma_p=2$ and
\begin{align*}
\sup_{\xb}p(\xb) &= \frac{\Gamma(d/2+1)}{\pi^{d/2}K_d(R)},\\
L(p;\supp(p)^{\mathrm{o}}) &= \frac{\Gamma(d/2+1)}{\pi^{d/2}K_d(R)}\sqrt{R^4+d}, \\
H^{d-1}(\partial\supp(p)) &= d\ups_dR^{d-1}.
\end{align*}
\end{example}

\begin{example}[Truncated exponential]
Let $S_R\defeq \{\xb\in\Real^d\suchthat x_1,\ldots,x_d\ge 0,x_1+\ldots+x_d\le R\}$.
The {\it truncated $d$-dimensional exponential distribution} defined by the density
\begin{equation*} p(\xb)\defeq \frac{e^{-(x_1+\cdots+x_d)}}{1-\bigl(\sum_{i=0}^{d-1}\frac{R^i}{i!}\bigr)e^{-R}}\ones_{S_R}(\xb)
\end{equation*}
is $2$-H\"older continuous over $\supp(p)$ and satisfies
\begin{align*}
\sup_{\xb}p(\xb) &= \Bigl(1-\Bigl(\sum_{i=0}^{d-1}\frac{R^i}{i!}\Bigr)e^{-R}\Bigr)^{-1},\\
L(p;\supp(p)^{\mathrm{o}}) 
	&= d\sup_{\xb}p(\xb),
\end{align*}
and
\begin{align*}
H^{d-1}(\partial\supp(p))
    &=\Bigl(\frac{\sqrt{d}}{(d-1)!}+d\Bigr)R^{d-1},
\end{align*}
as can be seen by an analysis similar to that in the previous example.
\end{example}
\begin{example}[Truncated Laplace]
Consider the {\it truncated $d$-dimensional Laplace distribution} defined by the density
\begin{equation*} 
p(\xb)\defeq \frac{e^{-(|x_1|+\cdots+|x_d|)}}{2^d\bigl(1-\bigl(\sum_{i=0}^{d-1}\frac{R^i}{i!}\bigr)e^{-R}\bigr)}
\ones_{(-\infty,R]}(\|\xb\|_1).
\end{equation*}
Then, \ref{cond:bounded_above}, \ref{cond:bounded_below}, and \ref{cond:boundary} can be demonstrated similarly to the previous examples. 
For \ref{cond:smoothness_int}, note that for $x,y\in\Real,$
\begin{equation*} \bigl|e^{-|x|}-e^{-|y|}\bigr|\le|x-y|.\end{equation*}
Generalizing this to $d$ dimensions, we have
\begin{align*} \bigl|e^{-(|x_1|+\cdots+|x_d|)} - e^{-(|y_1|+\cdots+|y_d|)}\bigr|
&\le\|\xb - \yb\|_1\\
&\le\sqrt{d}\|\xb - \yb\|_2.\end{align*}
Therefore, the truncated $d$-dimensional Laplace distribution is $1$-H\"older continuous over $\supp(p)$ and satisfies
\begin{align*}
\sup_{\xb}p(\xb)  
    &= \Bigl(2^d\Bigl(1-\Bigl(\sum_{i=0}^{d-1}\frac{R^i}{i!}\Bigr)e^{-R}\Bigr)\Bigr)^{-1},\\
L(p;\supp(p)^{\mathrm{o}}) 
	&= \sqrt{d}\sup_{\xb}p(\xb),
\end{align*}
and
\begin{align*}
H^{d-1}(\partial\supp(p)) 
	&= \frac{2^d\sqrt{d}}{(d-1)!}R^{d-1}.
\end{align*}
\end{example}

\begin{example}[Truncated Cauchy]
Consider the {\it truncated $d$-dimensional Cauchy distribution} defined by the density
\begin{equation*} p(\xb)\defeq \frac{\Gamma\bigl((d+1)/2\bigr)}{\pi^{(d+1)/2}L_d(R)\bigl(1+\|\xb\|_2^2\bigr)^{(d+1)/2}}\ones_{(-\infty,R]}(\|\xb\|_2),
\end{equation*}
where
\begin{equation*} L_d(R)\defeq \frac{\int_0^{\arctan R}\sin^{d-1}\theta \diff\theta}{\int_0^{\pi/2}\sin^{d-1}\theta \diff\theta}\in[0,1].\end{equation*}
Then, we have
\begin{align*} \nabla^2p(\xb)_{ij} &= \frac{(d+1)\Gamma((d+1)/2)}{\pi^{(d+1)/2}L_d(R)\bigl(1+\|\xb\|_2^2\bigr)^{(d+5)/2}}
\\&\quad\times
\bigl((d+3)x_ix_j-\bigl(1+\|\xb\|_2^2\bigr)\d_{ij}\bigr),\end{align*}
which leads to the bound
\begin{equation*} \norm{\nabla^2p(\xb)}\le\frac{(d+1)\Gamma((d+1)/2)}{\pi^{(d+1)/2}L_d(R)}\sqrt{R^4(d+1)(d+3)+d}\end{equation*}
on $\supp(p)^\mathrm{o}.$ 
Therefore, the truncated $d$-dimensional Cauchy distribution is $2$-H\"older continuous over $\supp(p)$ and satisfies
\begin{align*}
\sup_{\xb}p(\xb)  
	&= \frac{\Gamma\bigl((d+1)/2\bigr)}{\pi^{(d+1)/2}L_d(R)},\\
L(p;\supp(p)^{\mathrm{o}}) 
	&= \frac{(d+1)\Gamma((d+1)/2)}{\pi^{(d+1)/2}L_d(R)}
	\\&\quad\times
	\sqrt{R^4(d+1)(d+3)+d},
\end{align*}
and
\begin{align*}
H^{d-1}(\partial\supp(p)) 
	&= d\ups_dR^{d-1}.
\end{align*}
\end{example}


\ifCLASSOPTIONcaptionsoff
  \newpage
\fi



\bibliographystyle{IEEEtranN}
%
\bibliography{ref}

%

\begin{IEEEbiographynophoto}{J. Jon Ryu}
J. Jon Ryu (S'18) received the B.S. (Hons.) degrees in electrical and computer engineering and mathematical science (double major) from Seoul National University, Seoul, South Korea, in 2015. He is pursuing the Ph.D. degree in the Department of Electrical and Computer Engineering from the University of California San Diego (UCSD), La Jolla, CA, USA. 
He was a recipient of Kwanjeong Scholarship for graduate study from 2015 to 2020.
His research interests include information theory, data science, and statistical machine learning.
\end{IEEEbiographynophoto}

\begin{IEEEbiographynophoto}{Shouvik Ganguly}
Shouvik Ganguly (S'17–M'21) received the B.Tech. degree in electrical
engineering from Indian Institute of Technology, Kanpur in 2013, and the
Ph.D. degree in electrical engineering from the University of California San
Diego (UCSD) in 2020. In 2020, he joined XCOM Labs, San Diego, CA, USA, where he is currently a Member, Technical Staff. His research interests
include network information theory and communication theory.
\end{IEEEbiographynophoto}

\vfill

\begin{IEEEbiographynophoto}{Young-Han Kim}
Young-Han Kim (S'99–M'06–SM'12–F'15) received the B.S. degree (Hons.) in electrical engineering from Seoul National University, Seoul, South Korea, in 1996, and the M.S. degrees in electrical engineering and in statistics and the Ph.D. degree in electrical engineering from Stanford University, Stanford, CA, USA, in 2001, 2006, and 2006, respectively. In 2006, he joined the University of California San Diego, La Jolla, CA USA, where he is currently a Professor in the Department of Electrical and Computer Engineering. Since 2020, he has also been a founding CEO of Gauss Labs Inc., an industrial AI startup company in Silicon Valley and Seoul, South Korea. He has co-authored the book Network Information Theory (Cambridge University Press, 2011) and the monograph Fundamentals of Index Coding (Now Publishers, 2018). His current research interests include data science, machine learning, information theory, and their applications in manufacturing, microelectronics, communications, networking, cryptography, and bioinformatics. Prof. Kim was a recipient of the 2008 NSF Faculty Early Career Development Award, the 2009 US–Israel Binational Science Foundation Bergmann Memorial Award, the 2012 IEEE Information Theory Paper Award, and the 2015 IEEE Information Theory Society James L. Massey Research and Teaching Award for Young Scholars. He served as an Associate Editor of the IEEE Transactions on Information Theory and a Distinguished Lecturer for the IEEE Information Theory Society. He is a foreign member of the National Academy of Engineering of Korea.
\end{IEEEbiographynophoto}

\begin{IEEEbiographynophoto}{Yung-Kyun Noh}
Yung-Kyun Noh (M'19) is an Associate Professor in the Department of Computer Science at Hanyang University and an Affiliate Professor in the School of Computational Sciences at the Korea Institute for Advanced Study. 
He received the BS degree in physics from POSTECH, and the PhD degree in computer science from Seoul National University. His research interests include metric learning and dimensionality reduction in machine learning, and he is especially interested in applying statistical theory of nearest neighbors to real, large datasets. He worked in the GRASP Robotics Laboratory, University of Pennsylvania in Philadelphia as a visiting researcher. He is currently a visiting scientist at the RIKEN Center for Advanced Intelligence Project in Tokyo and a visiting scholar at the Mayo Clinic Gastroenterology and Hepatology in Rochester.
\end{IEEEbiographynophoto}

\begin{IEEEbiographynophoto}{Daniel D. Lee}
Dr. Daniel Dongyuel Lee (F'14) is the Tisch University Professor in Electrical and Computer Engineering at Cornell Tech and Executive Vice President and Head of the Global AI Center for Samsung Research. He received his B.A. summa cum laude in Physics from Harvard University and his Ph.D. in Condensed Matter Physics from the Massachusetts Institute of Technology. He was also a researcher at Bell Labs in the Theoretical Physics and Biological Computation departments. He is a Fellow of the IEEE and AAAI and has received the NSF CAREER award and the Lindback award for distinguished teaching. He was also a fellow of the Hebrew University Institute of Advanced Studies in Jerusalem, an affiliate of the Korea Advanced Institute of Science and Technology, and organized the US-Japan National Academy of Engineering Frontiers of Engineering symposium and Neural Information Processing Systems (NeurIPS) conference. His group focuses on understanding general computational principles in biological systems and on applying that knowledge to build autonomous systems.
\end{IEEEbiographynophoto}





\vfill


\end{document}